\documentclass[onefignum,onetabnum]{siamart220329}

\usepackage{lipsum}
\usepackage{amsfonts}
\usepackage{graphicx}
\usepackage{epstopdf}
\usepackage{algorithmic}
\usepackage{multirow}
\ifpdf
  \DeclareGraphicsExtensions{.eps,.pdf,.png,.jpg}
\else
  \DeclareGraphicsExtensions{.eps}
\fi

\newsiamremark{remark}{Remark}
\newsiamremark{hypothesis}{Hypothesis}
\crefname{hypothesis}{Hypothesis}{Hypotheses}
\newsiamthm{claim}{Claim}

\headers{Potential Singularity of the Axisymmetric Euler Equations}{Thomas Y. Hou and Shumao Zhang}

\title{Potential Singularity of the Axisymmetric Euler Equations with $C^\alpha$ Initial Vorticity for A Large Range of $\alpha$}

\author{Thomas Y. Hou\thanks{Department of Computing and Mathematical Sciences, California Institute of Technology, Pasadena, CA (\email{hou@cms.caltech.edu}).}
\and Shumao Zhang\thanks{Department of Computing and Mathematical Sciences, California Institute of Technology, Pasadena, CA (\email{shumaoz@caltech.edu}).}}

\usepackage{amsopn}

\makeatletter
\newcommand*{\addFileDependency}[1]{
  \typeout{(#1)}
  \@addtofilelist{#1}
  \IfFileExists{#1}{}{\typeout{No file #1.}}
}
\makeatother

\ifpdf
\hypersetup{
  pdftitle={Potential Singularity of the axisymmetric Euler equations},
  pdfauthor={Thomas Y. Hou, and Shumao Zhang}
}
\fi

\begin{document}

\maketitle

\begin{abstract}
    We provide numerical evidence for a potential finite-time self-similar singularity of the 3D axisymmetric Euler equations with no swirl and with $C^\alpha$ initial vorticity for a large range of $\alpha$. We employ a highly effective adaptive mesh method to resolve the potential singularity sufficiently close to the potential blow-up time. Resolution study shows that our numerical method is at least second-order accurate. Scaling analysis and the dynamic rescaling method are presented to quantitatively study the scaling properties of the potential singularity. We demonstrate that this potential blow-up is stable with respect to the perturbation of initial data. Our numerical study shows that the 3D axisymmetric Euler equations with our initial data develop finite-time blow-up when the H\"older exponent $\alpha$ is smaller than some critical value $\alpha^*$, which has the potential to be $1/3$. We also study the $n$-dimensional axisymmetric Euler equations with no swirl, and observe that the critical H\"older exponent $\alpha^*$ is close to $1-\frac{2}{n}$. Compared with Elgindi's blow-up result in a similar setting \cite{elgindi2021finite}, our potential blow-up scenario has a different H\"older continuity property in the initial data and the scaling properties of the two initial data are also quite different. We also propose a relatively simple one-dimensional model and numerically verify its approximation to the $n$-dimensional axisymmetric Euler equations. This one-dimensional model sheds useful light to our understanding of the blow-up mechanism for the $n$-dimensional Euler equations. 
\end{abstract}



\section{Introduction}

The three-dimensional (3D) incompressible Euler equations in fluid dynamics describe the motion of inviscid incompressible flows and are one of the most fundamental equations in fluid dynamics. Despite their wide range of applications, the question regarding the global regularity of the Euler equations has been widely recognized as a major open problem in partial differential equations (PDEs) and is closely related to the Millennium Prize Problem on the Navier-Stokes equations listed by the Clay Mathematics Institute \cite{fefferman2000existence}. In 2014, Luo and Hou \cite{luo2014potentially, luo2014toward} considered the 3D axisymmetric Euler equations with smooth initial data and boundary, and presented strong numerical evidences that they can develop potential finite-time singularity. The presence of the boundary, the symmetry properties and the direction of the flow in the initial data collaborate with each other in the formation of a sustainable finite-time singularity. Recently, Chen and Hou \cite{chen2022stable,chen2023stable} provided a rigorous justification of the Luo-Hou blow-up scenario.

In 2021, Elgindi \cite{elgindi2021finite} showed that given appropriate $C^{\alpha}$ initial vorticity with $\alpha>0$ sufficiently small, the 3D axisymmetric Euler equations with no swirl can develop finite-time singularity. In Elgindi's work, the initial data for the vorticity $\omega$ have $C^\alpha$ H\"{o}lder continuity near $r=0$ and $z=0$. When $\alpha$ is small enough, Elgindi approximated the 3D axisymmetric Euler equations by a fundamental model that develops a self-similar finite-time singularity. The blow-up result obtained in \cite{elgindi2021finite} has infinite energy. In a subsequent paper \cite{drivas2022singularity}, the authors improved the result obtained in \cite{elgindi2021finite} to have finite energy blow-up.

In this work we study potential finite-time singularity of the 3D axisymmetric Euler equations with no swirl and $C^\alpha$ initial vorticity for a large range of $\alpha$. Define $\omega=\nabla\times u$ as the vorticity vector, and then the 3D incompressible Euler equations can be written in the vorticity stream function formulation:
\begin{equation}
\label{eq: vort_stream_euler}
    \begin{aligned}
    \omega_t+u\cdot\nabla\omega&=\omega\cdot\nabla u,\\
    -\nabla\psi&=\omega,\\
    u&=\nabla\times\psi,
    \end{aligned}
\end{equation}
where $\psi$ is the vector-valued stream function. Let us use $x=(x_1, x_2, x_3)$ to denote a point in $\mathbb{R}^3$, and let $e_r$, $e_\theta$, $e_z$ be the unit vectors of the cylindrical coordinate system
$$e_r=\left(x_1/r, x_2/r, 0\right),\quad e_\theta=\left(x_2/r, -x_1/r, 0\right),\quad e_z=\left(0, 0, 1\right),$$
where $r=\sqrt{x_1^2+x_2^2}$ and $z=x_3$. We say a vector field $v$ is axisymmetric if it admits the decomposition
$$v=v^r(r, z)e_r + v^\theta(r, z)e_\theta + v^z(r, z)e_z,$$
namely, $v^r$, $v^\theta$ and $v^z$ are independent of the angular variable $\theta$. Denote by $u^\theta$, $\omega^\theta$, and $\psi^\theta$ the angular velocity, vorticity and stream function, respectively. The axisymmetric condition can then simplify the 3D Euler equations \eqref{eq: vort_stream_euler} to \cite{majda2002vorticity}:
\begin{subequations}
\label{eq: vort_stream_3d}
\begin{align}
    u^\theta_t + u^ru^\theta_r + u^zu^\theta_z &= -\frac{1}{r}u^ru^\theta, \label{eq: velo_theta_3d}\\
    \omega^\theta_t + u^r\omega^\theta_r + u^z\omega^\theta_z &= \frac{2}{r}u^\theta u^\theta_z + \frac{1}{r}u^r\omega^\theta, \label{eq: vort_theta_3d}\\
    -\psi^\theta_{rr}-\psi^\theta_{zz}-\frac{1}{r}\psi^\theta_{r}+\frac{1}{r^2}\psi^\theta&=\omega^\theta, \label{eq: stream_theta_3d}\\
    u^r=-\psi^\theta_z,\quad u^z&=\frac{1}{r}\psi^\theta+\psi^\theta_r. \label{eq: velo_rz_3d}
\end{align}
\end{subequations}
In the case of no swirl, i.e. $u^\theta \equiv 0$, the axisymmetric Euler equations are further simplified into:
\begin{subequations}
\label{eq: vort_stream_3d_noswirl}
\begin{align}
    \omega^\theta_t + u^r\omega^\theta_r + u^z\omega^\theta_z &= \frac{1}{r}u^r\omega^\theta, \label{eq: vort_theta_3d_noswirl}\\
    -\psi^\theta_{rr}-\psi^\theta_{zz}-\frac{1}{r}\psi^\theta_{r}+\frac{1}{r^2}\psi^\theta&=\omega^\theta, 
    \label{eq: stream_theta_3d_noswirl}\\
    u^r=-\psi^\theta_z,\quad u^z&=\frac{1}{r}\psi^\theta+\psi^\theta_r. \label{eq:velo_rz_3d_noswirl}
\end{align}
\end{subequations}
When the initial condition for the angular vorticity $\omega^\theta$ is smooth, it is well known that the 3D axisymmetric Euler equations with no swirl \eqref{eq: vort_stream_3d_noswirl} will not develop finite-time blow-up \cite{ukhovskii1968axially}. Therefore, we consider \eqref{eq: vort_stream_3d_noswirl} when the initial condition for the angular vorticity $\omega^\theta$ is $C^\alpha$ H\"{o}lder continuous for a large range of $\alpha$. By using an effective adaptive mesh method, we will provide convincing numerical evidence that the 3D axisymmetric Euler equations with no swirl and $C^\alpha$ initial voriticity with $\alpha$ greater or equal to $0$, and smaller than some critical value $\alpha^*$ can develop potential finite-time self-similar blow-up. The critical H\"{o}lder exponent $\alpha^*$ is observed to be larger than $0.3$ and close to $1/3$. Our result serves as an example to support Conjecture 8 of \cite{drivas2022singularity} that the critical value $\alpha^*$ is equal to $1/3$.

We perform scaling analysis and use the dynamic rescaling formulation \cite{hou2018potential, chen2021finite, chen2021finite2} to study the behavior of the potential self-similar blow-up. An operator splitting method is proposed to solve the dynamic rescaling formulation and the late time solution from the adaptive mesh method is used as our initial condition for the dynamic rescaling formulation. We observe rapid convergence to a steady state, which implies that this potential singularity is self-similar. We notice that the size of the finite computational domain needs to be large enough to approximate steady state, due to a scale-invariant property of the dynamic rescaling equation. So we conduct domain size study to verify the accuracy of our results. We will also demonstrate that this potential blow-up is stable with respect to the perturbation of initial data, suggesting that the underlying blow-up mechanism is generic and insensitive to the initial data. 

We choose the following $C^\alpha$ initial data for the angular vorticity $\omega^\theta$:
\[
\omega^\theta_0 =\frac{-12000 \;r^\alpha \left(1-r^2\right)^{18}\sin(2\pi z)}{1+12.5\cos^2(\pi z)}\,.
\]
The initial condition is a smooth and periodic function in $z$ and is $C^\alpha$ in $r$.
The velocity field $u$ becomes $C^{1,\alpha}$ continuous. We further introduce the new variables:
\begin{align}
\label{eq: new_variable_3d_noswirl}
    \omega_1(r,z)=\frac{1}{r^\alpha}\omega^\theta(r,z),\quad\psi_1(r,z)=\frac{1}{r}\psi^\theta(r,z),
\end{align}
to remove the formal singularity in \eqref{eq: vort_stream_3d_noswirl} near $r=0$. 
In terms of the new variables $\left(\omega_1, \psi_1\right)$, the 3D axisymmetric Euler equations with no swirl have the following equivalent form
\begin{subequations}
\label{eq: vort_stream_1_3d_noswirl}
\begin{align}
    \omega_{1,t} + u^r\omega_{1,r} + u^z\omega_{1,z} &= -(1-\alpha)\psi_{1,z}\omega_1, \label{eq: vort_1_3d_noswirl}\\
    -\psi_{1,rr}-\psi_{1,zz}-\frac{3}{r}\psi_{1,r}&=\omega_1r^{\alpha-1}, \label{eq: stream_1_3d_noswirl}\\
    u^r=-r\psi_{1,z},\quad u^z&=2\psi_1+r\psi_{1,r}. \label{eq: velo_rz_1_3d_noswirl}
\end{align}
\end{subequations}
The above reformulation is crucial for us to perform accurate numerical computation of the potential singular solution and allow us to push the computation sufficiently close to the singularity time.

It is important to note that the initial condition for the rescaled vorticity field $\omega_1$ is a smooth function of $r$ and $z$. Using the above reformulation enables us to resolve the potential singular solution sufficiently close to the potential singularity time. If we solve the original 3D Euler equations \eqref{eq: vort_theta_3d_noswirl}--\eqref{eq:velo_rz_3d_noswirl}, it is extremely difficult to resolve the H\"older continuous vorticity even with an adaptive mesh, especially for small $\alpha$. For this reason, we have not been able to compute the finite-time singularity in Elgindi's work \cite{elgindi2021finite} since such a reformulation is not available for his initial data.

Compared with Elgindi's blow-up result \cite{elgindi2021finite}, our potential blow-up scenario has very different scaling properties. In our scenario, the scaling factor $c_l$, defined in \eqref{eq: self_similar}, increases with $\alpha$ and tends to infinity as $\alpha$ approaches $\alpha^*$, see Table \ref{tab: cl_n3}. In contrast, the scaling factor $c_l$ in Elgindi's scenario is $1/\alpha$, which decreases with $\alpha$ and tends to infinity as $\alpha$ approaches $0$. Another difference is that Elgindi's initial vorticity is $C^{2\alpha}$ in $r$ and $C^\alpha$ in $z$ near the origin, while our initial vorticity is $C^{\alpha}$ in $r$, but smooth in $z$. We discuss the comparison in details in Section \ref{sec: compare_elgindi}.

We also consider the $n$-dimensional axisymmetric Euler equations (see definitions in \eqref{eq: vort_stream_nd}). We observe the same potential self-similar finite-time blow-up and find that the critical value $\alpha^*$ is close to $1-\frac{2}{n}$. The self-similar profiles for high dimensional Euler equations are qualitatively similar to those of the 3D Euler equations. We observe that the blow-up is more robust for higher space dimensions and the scaling factor $c_l$ decreases as $n$ increases, which may be partially due to the stronger nonlinearity in the vortex stretching term with larger $n$. Similar to the 3D Euler equations, we see $c_l$ quickly increases with $\alpha$, and has the trend to go to infinity as $\alpha \rightarrow \alpha^*$. We observe that the stream function $\psi_1$ becomes almost linear in the near field along the rescaled $z$ variable (denoted as $\zeta$) as $\alpha \rightarrow \alpha^*$. Another interesting observation is that $\omega_1$ becomes increasingly flat as a function of the rescaled $r$ variable $\xi$. Based on this observation, we propose a simplified one-dimensional model along the $z$-direction by extending $\omega_1$ as a constant in the $r$-direction. This is equivalent to approximating $\omega^\theta (r,z,t) \approx r^\alpha \omega_1 (0,z,t)$, which still captures the effect of the H\"older continuity of $\omega^\theta$ through $r^\alpha$. Although this 1D model is supposed to give a good approximation of the 3D Euler equation as $\alpha \rightarrow \alpha^*$, we also observe that it can approximate the self-similar profile of the 3D Euler equation along the $z$-axis and its scaling factor $c_l$ very well even for small $\alpha$. The analysis of the 1D model should shed useful light on the blow-up mechanism of the 3D and the $n$-D Euler equations.
 
There have been some exciting recent progress on the global regularity of the high dimensional axisymmetric Euler equations with no swirl with smooth initial data under some assumptions, see e.g. \cite{Choi2022,Miller2022,Lim2023}. In particular, when $n = 4$, Choi-Jeong-Lim \cite{Choi2022} proved global regularity of the 4D axisymmetric Euler equation with no swirl under the assumption that the initial vorticity satisfies some decay condition at infinity and is vanishing at the symmetry axis. Further, if the initial vorticity is of one sign, they proved global regularity for $n \leq 7$. This result is further improved in a subsequent paper by Lim \cite{Lim2023} to any dimension $ n \geq 4$ under a similar decay assumption and for one-signed vorticity. In \cite{Miller2022}, Miller showed that the four and higher dimensional axisymmetric Euler equations with no swirl have properties which could lead to finite-time blow-up that is excluded for the 3D Euler equation. The author also considered a model for the infinite-dimensional vorticity equation, which exhibits finite-time blow-up of a Burgers shock type. The blow-up result of this model equation seems to suggest that the Euler equation in sufficiently high dimension is likely to develop a finite-time blow-up with smooth initial data.

In an excellent survey paper by Drivas and Elgindi \cite{drivas2022singularity}, the authors discussed singularity formation in the high dimensional incompressible Euler equation in some details. In particular, the authors asked in their Open Question 7 that if singularities can form from smooth data for the axisymmetric no-swirl Euler equations on $\mathbb{R}^n$ when $n\geq 4$. To the best of our knowledge, there has been no strong numerical evidence for potential finite-time blow-up for the $n$-dimensional axisymmetric Euler equation with no swirl and smooth initial data. We have performed such computation by ourselves and did not find any evidence for finite-time blow-up. One of the reasons for the non-blow-up is that the quantity $\omega^\theta/r^{n-2}$ satisfies a transport equation with no vortex stretching, thus $\omega^\theta/r^{n-2}$ satisfies a maximum principle. As long as $\omega^\theta/r^{n-2}$ is well defined at $t=0$ and $\omega^\theta$ decays rapidly at infinity, we will have control of the maximum growth of vorticity, thus there is no finite-time blow-up \cite{Choi2022,Miller2022,Lim2023}.

Theoretical analysis of the 3D Euler equations have been studied for long. The Beale-Kato-Majda (BKM) blow-up criterion \cite{beale1984remarks, ferrari1993blow} gives a necessary and sufficient condition for the finite-time singularity for the smooth solutions of the 3D Euler equations at time $T$ if and only if $  \int_0^T\|\omega(\cdot,t)\|_{L^\infty}\mathrm{d}t=+\infty$. This result also holds true for H\"{o}lder continuous initial data, see Theorem 4.3 of \cite{majda2002vorticity} and can be easily generalized to high dimensional Euler equations. In \cite{constantin1996geometric}, Constantin, Fefferman and Majda asserted that there will be no finite-time blow-up if the velocity $u$ is uniformly bounded and the direction of vorticity $\xi=\omega/|\omega|$ is sufficiently regular (Lipschitz continuous) in an $O(1)$ domain containing the location of the maximum vorticity. Inspired by the work of \cite{constantin1996geometric}, Deng-Hou-Yu developed a more localized non-blow-up criterion using a Lagrangian approach in \cite{deng2005geometric}.

There have been a number of numerical attempts in search of the potential finite-time blow-up. 
The finite-time blow-up in the numerical study was first reported by Grauer and Sideris \cite{grauer1991numerical} and Pumir and Siggia \cite{pumir1992development} for the 3D axisymmetric Euler equations. However, the later work of E and Shu \cite{weinan1992numerical} suggested that the finite-time blow-up in \cite{grauer1991numerical, pumir1992development} could be caused by numerical artifact. 
Kerr and his collaborators \cite{kerr1993evidence, bustamante20083d} presented finite-time singularity formation in the Euler flows generated by a pair of perturbed anti-parallel vortex tubes. In \cite{hou2006dynamic}, Hou and Li reproduced Kerr's computation using a similar  initial condition with much higher resolutions and did not observe finite-time blow-up. The maximum vorticity grows slightly slower than double exponential in time. Later on, Kerr  confirmed in \cite{kerr2013bounds} that the solutions from \cite{kerr1993evidence} eventually converge to a super-exponential growth and are unlikely to lead to a finite-time singularity.

In \cite{caflisch1993singularity, siegel2009calculation}, Caflisch and his collaborators studied axisymmetric Euler flows with complex initial data and reported singularity formation in the complex plane. The review paper \cite{gibbon2008three} lists a more comprehensive collection of interesting numerical results with more detailed discussions.

Due to the lack of stable structure in the potentially singular solutions, the previously mentioned numerical results remain inconclusive. In \cite{luo2014potentially, luo2014toward}, Luo and Hou reported that the 3D axisymmetric Euler equations with a smooth initial condition developed a self-similar finite-time blow-up in the meridian plane on the boundary of $r=1$, see also \cite{luo2019formation}. The Hou-Luo blow-up scenario has generated a great deal of interests in both the mathematics and fluid dynamics communities, and inspired a number of subsequent developments \cite{kiselev2014small, kiselev2016finite, kiselev2018small, choi2017finite, chen2021asymptotically, chen2021finite, chen2021finite2, chen2022stable,chen2023stable,cordoba2023finite,cordoba2023blow,chen2023remarks}.

We remark that there has been some recent exciting progress on the theoretical study of  singularity formation of the 2D Boussinesq system and the 3D Euler equations. In \cite{Elgindi2023}, Elgindi and Pasqualotto established finite-time singularity formation for $C^{1,\alpha}$ solutions to the Boussinesq system that are compactly supported on $\mathbb{R}^2$ and infinitely smooth except in the radial direction at the origin. 
In \cite{cordoba2023finite}, Cordoba, Martinez-Zorao and Zheng constructed non-self-similar blow-up solutions of the 3D axisymmetric Euler equation with no swirl and $C^\alpha$ initial vorticity, which is smooth except at the origin. Inspired by the work of \cite{cordoba2023finite}, Chen in \cite{chen2023remarks} showed that Elgindi's $C^{1,\alpha}$ self-similar blow-up can be improved to be smooth except at the origin. Further, Cordoba and Martinez-Zorao constructed finite-time singularity of 3D incompressible Euler equations with velocity in $C^{3,1/2}\cap L^2$ and uniform $C^{1,1/2-\epsilon} \cap L^2$ force.

The rest of this paper is organized as follows.  In Section \ref{sec: setting_numeric}, we briefly introduce the numerical method. We present the evidence of the potential self-similar blow-up in Section \ref{sec: evidence}, and provide the resolution study and scaling analysis for the case of $\alpha=0.1$. In Section \ref{sec: dynamic_rescaling} we use the dynamic rescaling method to provide further evidence of the potential blow-up. In Section \ref{sec: holder_dimension}, we consider the potential finite-time blow-up in the general case of the H\"{o}lder exponent $\alpha$ and the dimension $n$. The sensitivity of the potential blow-up to the initial data is considered in Section \ref{sec: sensitivity}, and the comparison of our potential blow-up scenario with Elgindi's scenario in \cite{elgindi2021finite} is discussed in Section \ref{sec: compare_elgindi}. Section \ref{sec: 1d_model} is devoted to a one-dimensional model to study the potential self-similar blow-up of the $n$-dimensional axisymmetric Euler equations. Some concluding remarks are made in Section \ref{sec: conclusions}. 

\section{Problem set up and numerical method}
\label{sec: setting_numeric}

In this section, we give details about the setup of the problem, the initial data, the boundary conditions, and some basic properties of the equations, and our numerical method.

\subsection{Boundary conditions and symmetries}
\label{sec: boundary_symmetry}

We consider \eqref{eq: vort_stream_1_3d_noswirl} in a cylinder region
$$\mathcal{D}_{\text{cyl}}=\left\{(r,z): 0\leq r\leq 1\right\},$$
We impose a periodic boundary condition in $z$ with period $1$:
\begin{align}
\label{eq: z_periodic}
    \omega_1(r,z)=\omega_1(r,z+1),\quad \psi_1(r,z)=\psi_1(r, z+1).
\end{align}
In addition, we enforce that $(\omega_1, \psi_1)$ are odd in $z$ at $z=0$:
\begin{align}
\label{eq: z_symmetric}
    \omega_1(r,z)=-\omega_1(r,-z),\quad \psi_1(r,z)=-\psi_1(r, -z).
\end{align}
And this symmetry will be preserved dynamically by the 3D Euler equations.

At $r=0$, it is easy to see that $u^r(0, z)=0$, so there is no need for the boundary condition for $\omega_1$ at $r=0$. Since $\psi^\theta=r\psi_1$ will at least be $C^2$-continuous, according to \cite{liu2006convergence, liu2009characterization}, $\psi^\theta$ must be an odd function of $r$. Therefore, we impose the following pole condition for $\psi_1$
\begin{align}
\label{eq: r_pole}
    \psi_{1,r}(0,z)=0.
\end{align}
We impose the no-flow boundary condition at the boundary $r=1$:
\begin{align}
\label{eq: r_noflow}
    \psi_{1}(1,z)=0.
\end{align}
This implies that $u^r(1,z)=0$. So there is no need to introduce a boundary condition for $\omega_1$ at $r=1$.

Due to the periodicity and the odd symmetry along the $z$-direction, the equations \eqref{eq: vort_stream_3d_noswirl} only need to be solved on the half-periodic cylinder
$$\mathcal{D}=\left\{(r,z): 0\leq r\leq 1, 0\leq z\leq1/2\right\}.$$
The above boundary conditions of $\mathcal{D}$ show that there is no transport of the flow across its boundaries. Indeed, we have
$$u^r=0\quad\text{ on }\quad r=0\text{ or } 1,\quad\text{ and }\quad u^z=0\quad\text{ on }\quad z=0\text{ or }1/2.$$
Thus, the boundaries of $\mathcal{D}$ behave like ``impermeable walls''.

\subsection{Initial data}
\label{sec: initial_data}

\begin{figure}[hbt!]
\centering
\includegraphics[width=.4\textwidth]{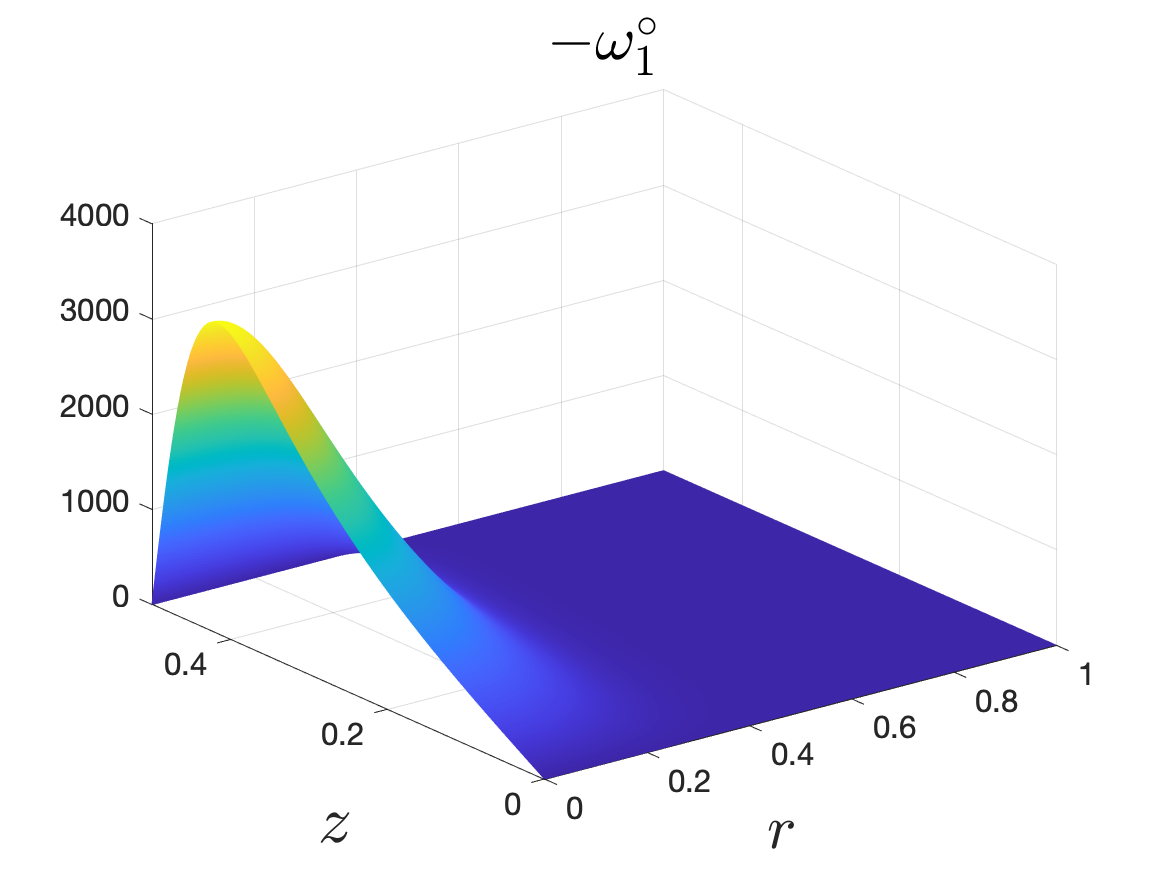}\includegraphics[width=.4\textwidth]{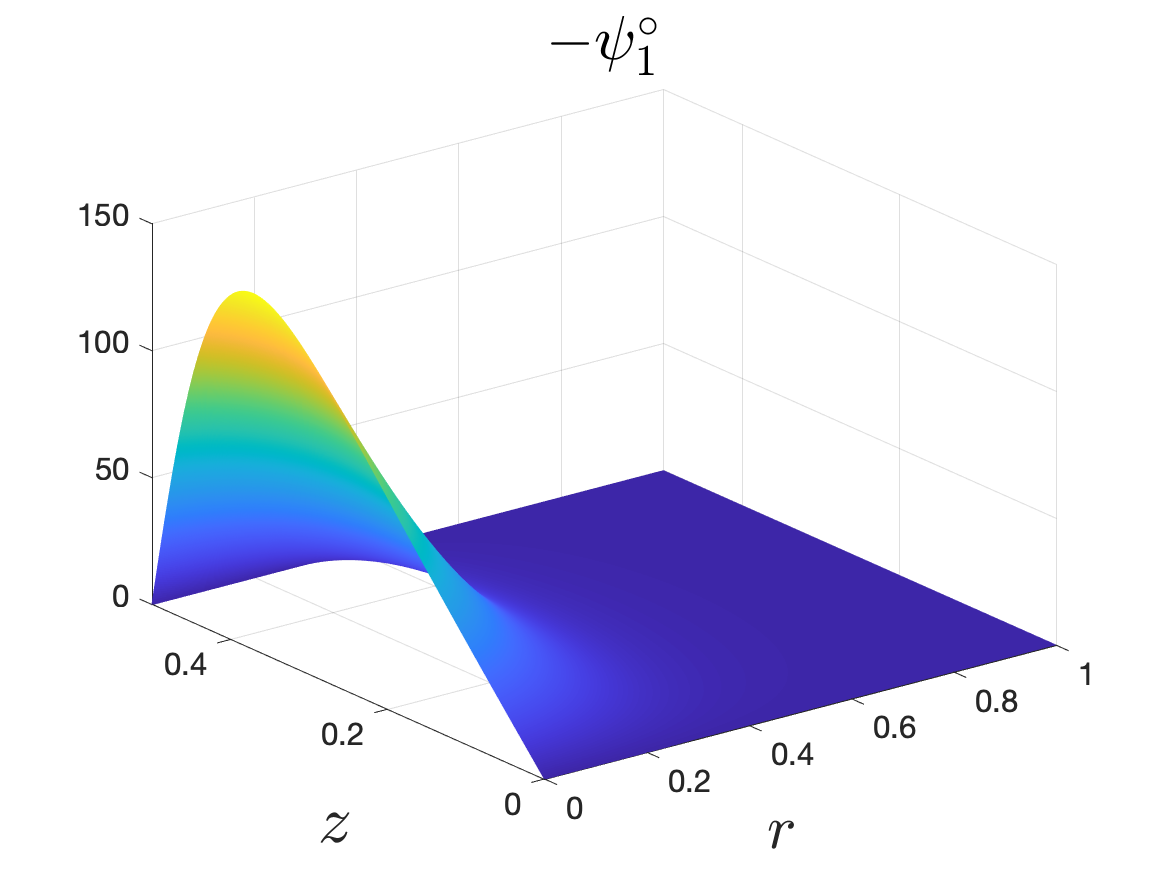}\\
\caption{3D profiles of the initial value $-\omega^\circ_1$ and $-\psi^\circ_1$.}\label{fig: initial_data}
\end{figure}

\begin{figure}[hbt!]
\centering
\includegraphics[width=.4\textwidth]{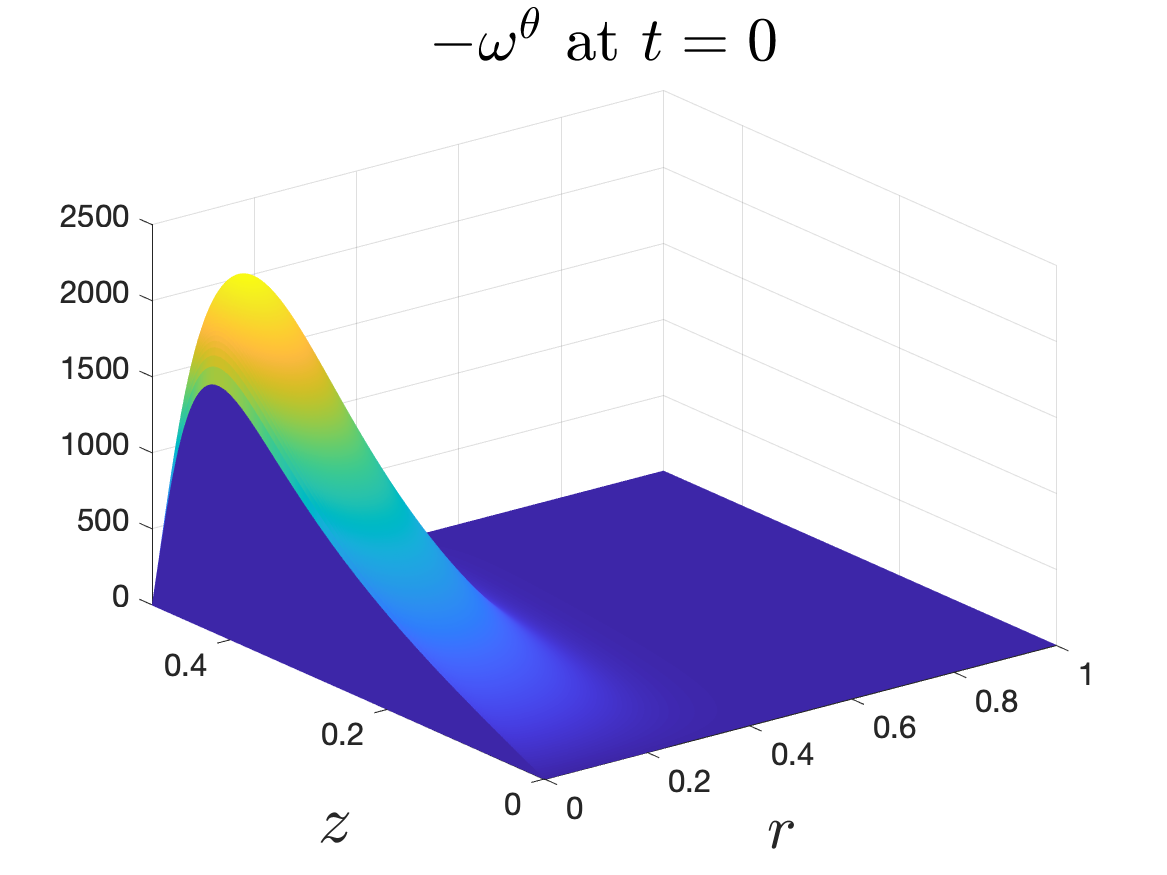}\includegraphics[width=.4\textwidth]{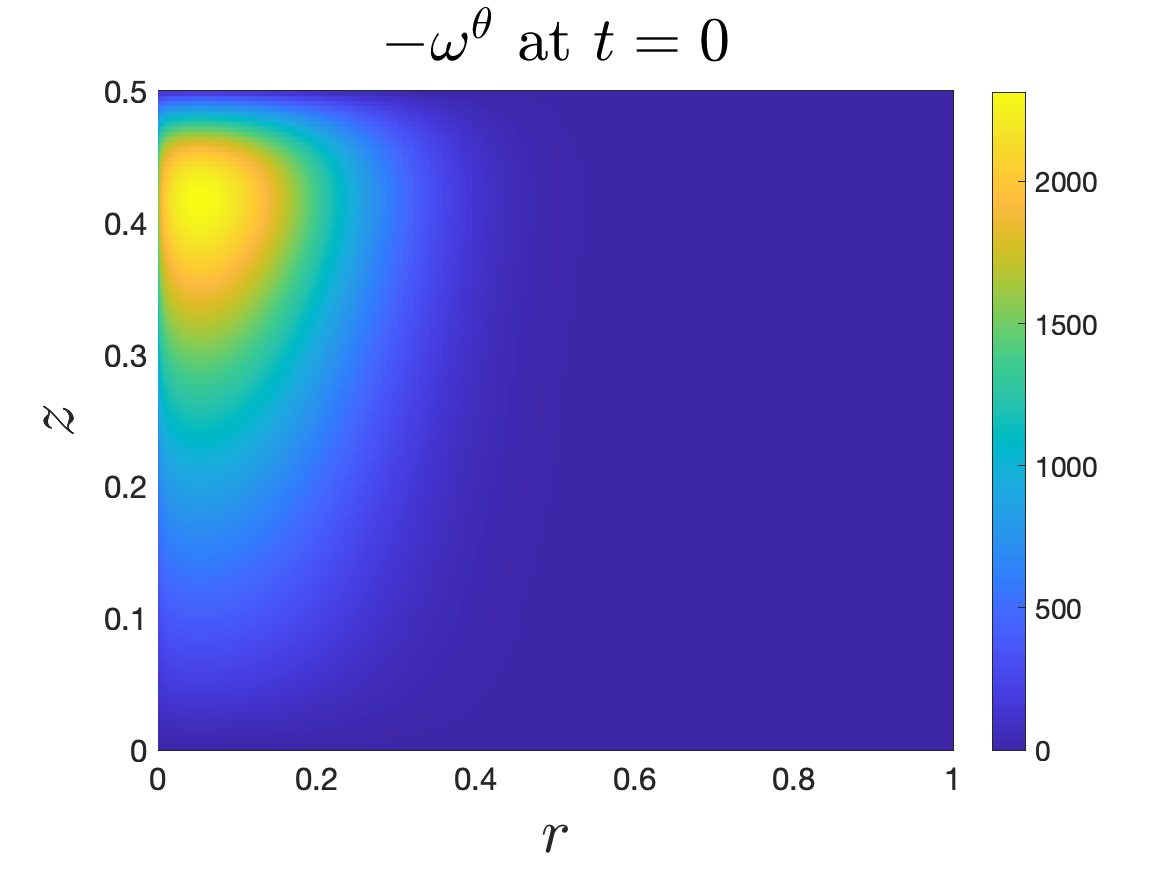}\\
\caption{The initial data for the angular vorticity $\omega^\theta$.}\label{fig: initial_omega}
\end{figure}

Inspired by the potential blow-up scenario in \cite{hou2021nearly}, we propose the following initial data for $\omega_1$ in $\mathcal{D}$,
\begin{align}
\label{eq: inital_data}
    \omega^\circ_1=\frac{-12000\left(1-r^2\right)^{18}\sin(2\pi z)}{1+12.5\cos^2(\pi z)}.
\end{align}
Later we will see in Section \ref{sec: sensitivity} that the self-similar singularity formation has some robustness to the choice of initial data. We solve the Poisson equation \eqref{eq: stream_1_3d_noswirl} to get the initial value $\psi^\circ_1$ of $\psi_1$.

\begin{figure}[hbt!]
\centering
\includegraphics[width=.4\textwidth]{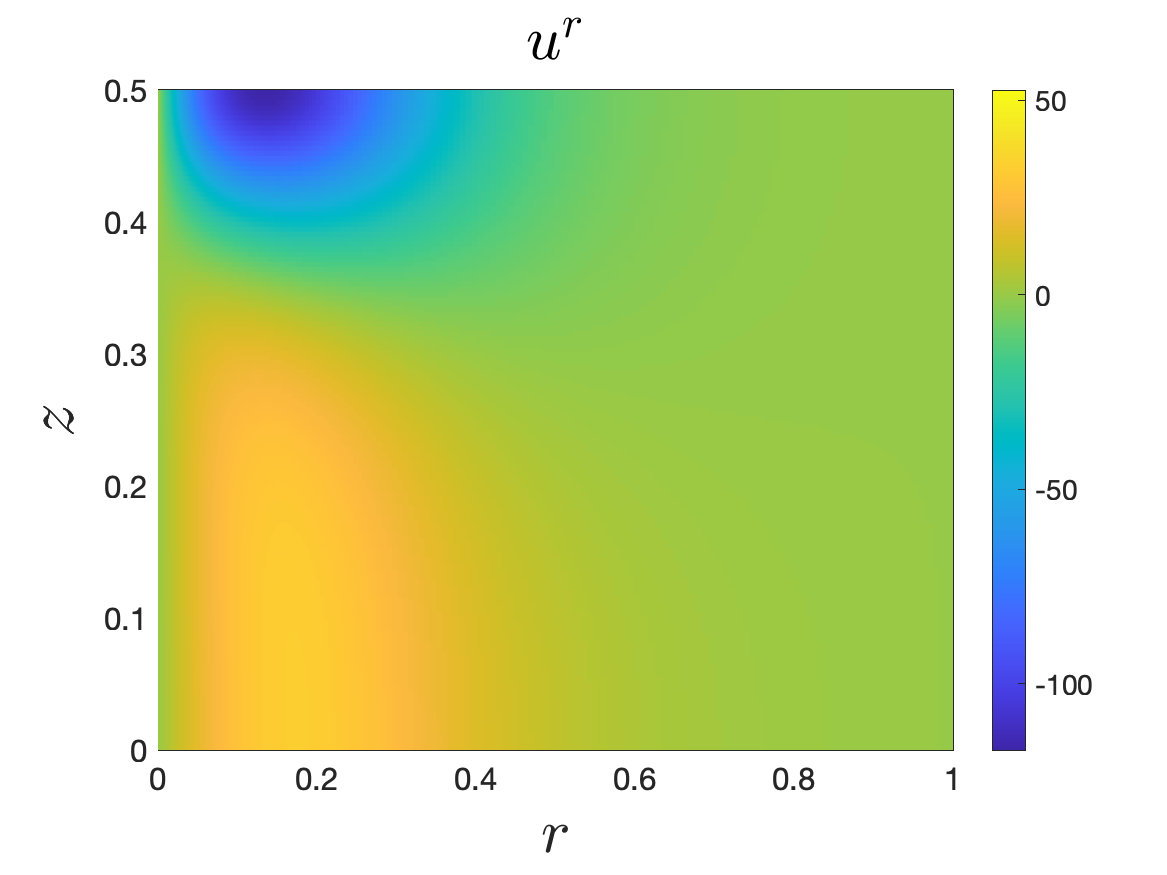}\includegraphics[width=.4\textwidth]{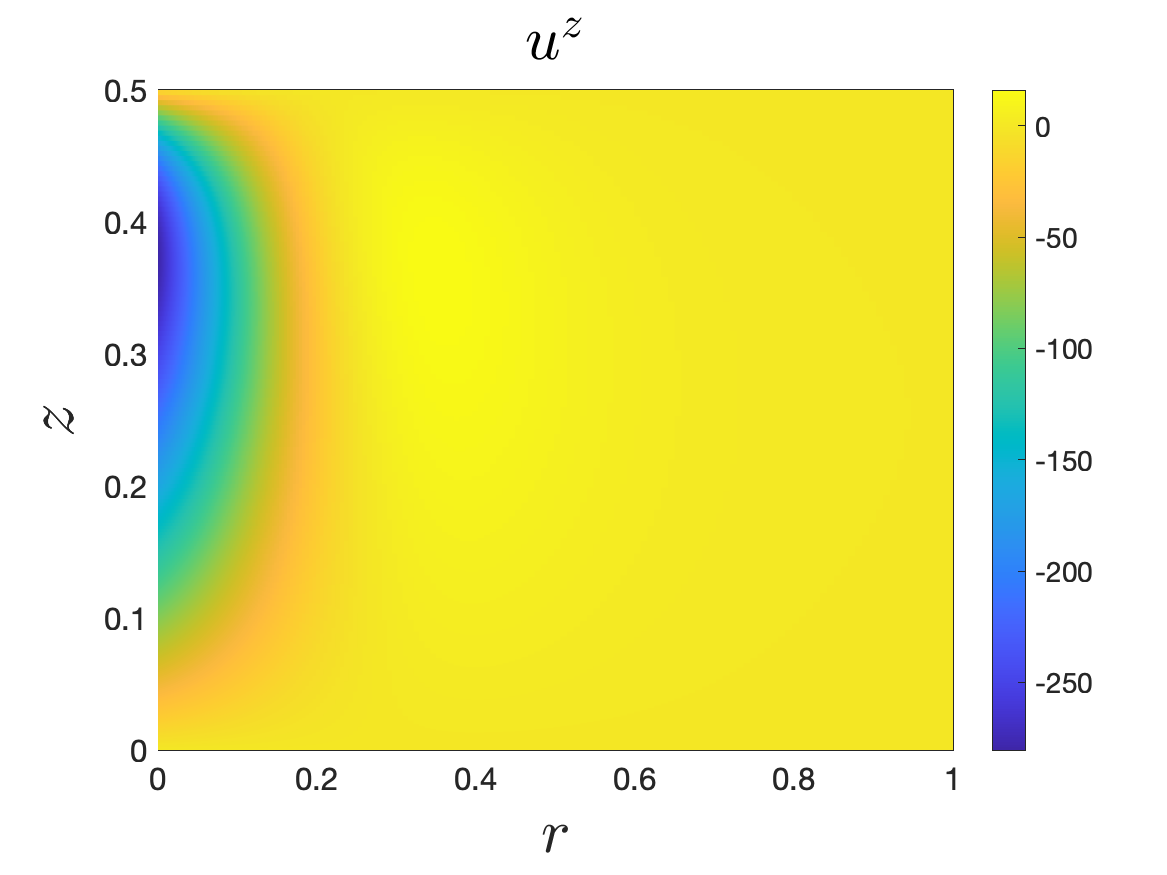}\\
\caption{Initial velocity fields $u^r$ and $u^z$.}\label{fig: initial_velo}
\end{figure}

\begin{figure}[hbt!]
\centering
\includegraphics[width=.5\textwidth]{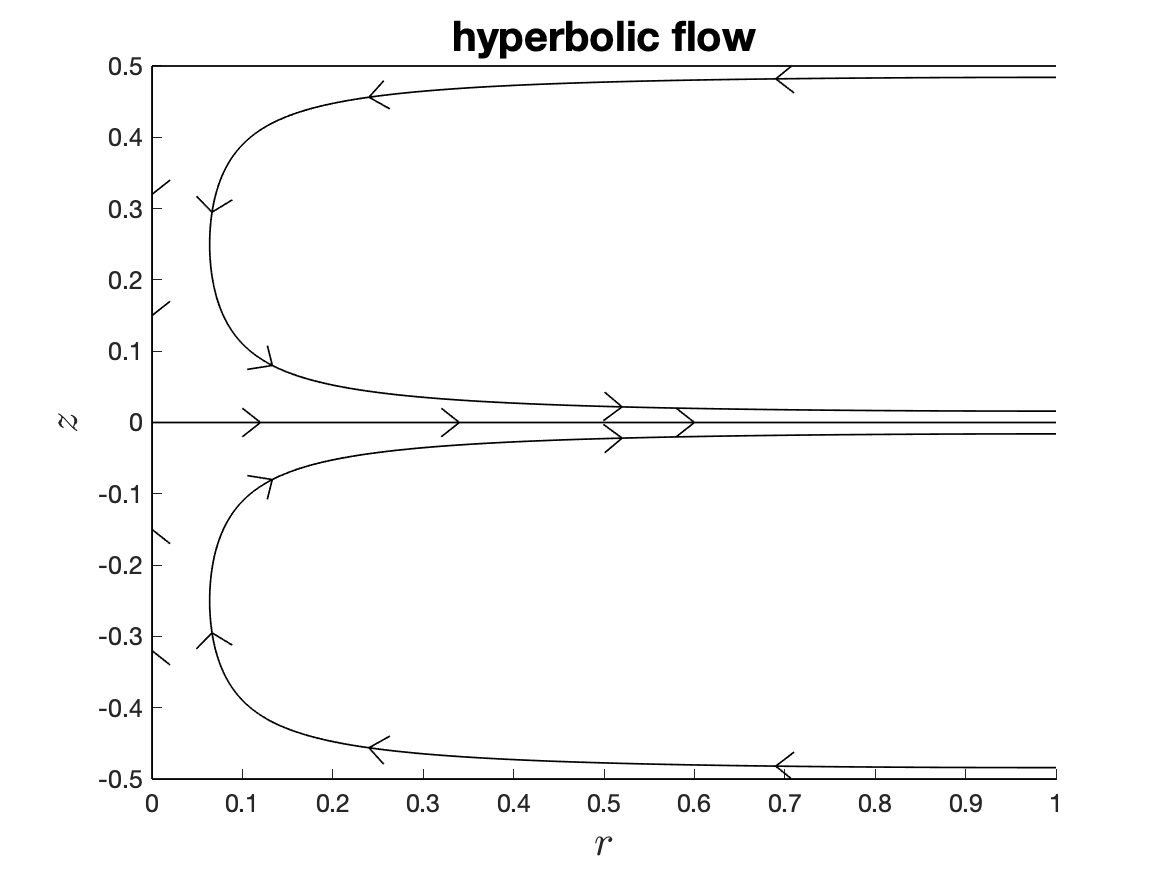}
\caption{A heuristic diagram of the hyperbolic flow.}\label{fig: hyper_flow}
\end{figure}

The 3D profiles of $(\omega^\circ_1, \psi^\circ_1)$ can be found in Figure \ref{fig: initial_data}. Since most parts of $\omega^\circ_1$ and $\psi^\circ_1$ are negative, we negate them for better visual effect when generating figures. In Figure \ref{fig: initial_omega}, we show the 3D profile and pseudocolor plot of the angular vorticity $\omega^\theta$ at $t=0$. We can see that there is a sharp drop to zero of $-\omega^\theta$ near $r=0$, which is due to the H\"{o}lder continuous of $\omega^\theta$ at $r=0$.

We plot the initial velocity field $u^r$ and $u^z$ in Figure \ref{fig: initial_velo}. We can see that $u^r$ is primarily positive near $z=0$ and negative near $z=1/2$ when $r$ is small, and $u^z$ is mainly negative when $r$ is small. Such a pattern suggests a hyperbolic flow near $(r,z)=(0,0)$ as depicted in the heuristic diagram Figure \ref{fig: hyper_flow}, which will extend periodically in $z$.

\subsection{Self-similar solution}
\label{sec: self_similar_solution}

Self-similar solutions are a common and important class of solutions to nonlinear PDEs that demonstrates their intrinsic structure and properties. A self-similar solution is when the local profile of the solution remains nearly unchanged in time after rescaling the spatial and the temporal variables of the physical solution. For example, for \eqref{eq: vort_stream_1_3d_noswirl}, the self-similar profile is the ansatz
\begin{equation}
\label{eq: self_similar}
    \begin{aligned}
    \omega_1(x,t)&\approx\frac{1}{(T-t)^{c_\omega}}\Omega\left(\frac{x-x_0}{(T-t)^{c_l}}\right),\\
    \psi_1(x,t)&\approx\frac{1}{(T-t)^{c_\psi}}\Psi\left(\frac{x-x_0}{(T-t)^{c_l}}\right),
    \end{aligned}
\end{equation}
for some parameters $c_\omega$, $c_\psi$, $c_l$, $x_0$ and $T$. Here $T$ is considered as the blow-up time, and $x_0$ is the location of the self-similar blow-up. The parameters $c_\omega$, $c_\psi$, $c_l$ are called scaling factors.

It is also important to notice that the 3D Euler equations \eqref{eq: vort_stream_euler} enjoy the following scaling invariant property: if $(u, \omega, \psi)$ is a solution to \eqref{eq: vort_stream_euler}, then $(u_{\lambda,\mu}, \omega_{\lambda,\mu}, \psi_{\lambda,\mu})$ is also a solution, where

\begin{align*}
    u_{\lambda,\mu}(x, t)= \frac{\lambda}{\mu}u\left(\frac{x}{\lambda}, \frac{t}{\mu}\right),~\omega_{\lambda,\mu}(x, t)= \frac{1}{\mu}\omega\left(\frac{x}{\lambda}, \frac{t}{\mu}\right),~\psi_{\lambda,\mu}(x, t)= \frac{\lambda^2}{\mu}\psi\left(\frac{x}{\lambda}, \frac{t}{\mu}\right),
\end{align*}
and $\lambda>0$, $\mu>0$ are two constant scaling factors. In the case of the 3D axisymmetric Euler equations with no swirl \eqref{eq: vort_stream_1_3d_noswirl}, the scaling invariant property can be equivalently translated to: if $(\omega_1, \psi_1)$ is a solution of \eqref{eq: vort_stream_1_3d_noswirl}, then 
\begin{align}
\label{eq: vort_stream_invariant}
\left\{\frac{1}{\lambda^\alpha\mu}\omega_1\left(\frac{x}{\lambda},\frac{t}{\mu}\right),~ \frac{\lambda}{\mu}\psi_1\left(\frac{x}{\lambda}, \frac{t}{\mu}\right)\right\}
\end{align}
is also a solution.

If we assume the existence of the self-similar solution \eqref{eq: self_similar}, then the new solutions in \eqref{eq: vort_stream_invariant} should also admit the same ansatz, regardless of the values of $\lambda$ and $\mu$. As a result, we must have
\begin{align}
\label{eq: scaling_relation}
    c_\omega=1+\alpha c_l,\quad c_\psi=1-c_l.
\end{align}
Therefore, the self-similar profile \eqref{eq: vort_stream_invariant} of \eqref{eq: vort_stream_1_3d_noswirl} only has one degree of freedom, for example $c_l$, in the scaling factors. In fact, $c_l$ cannot be determined by straightforward dimensional analysis.

As a consequence of the ansatz \eqref{eq: self_similar} and the scaling relation \eqref{eq: scaling_relation}, we have
\begin{align}
\label{eq: scaling_conclusion}
    \|\omega^\theta(x,t)\|_{L^\infty}\sim\frac{1}{T-t},\quad\|\psi_{1,z}(x,t)\|_{L^\infty}\sim\frac{1}{T-t},
\end{align}
which always holds true regardless of the value of $c_l$.

\subsection{Numerical method}
\label{sec: numeric}

Although the initial data are very smooth, the solutions of Euler equations quickly become very singular and concentrate in a rapidly shrinking region. Therefore, we use the adaptive mesh method to resolve the singular profile of the solutions. A detailed description of the adaptive mesh method can be found in \cite{hou2018potential,luo2019formation,zhang2023singularity}. Here we briefly introduce the idea behind the adaptive mesh method. The specific parameter setting used for the experiments in this work can be found in the appendix of \cite{zhang2023singularity}.

The Euler equations \eqref{eq: vort_stream_1_3d_noswirl} are originally posted as an initial-boundary value problem on the computational domain $(r, z)\in\left[0, 1\right]\times\left[0, 1/2\right]$. To capture the singular part of the solution, we introduce two variables $(\kappa, \eta)\in\left[0, 1\right]\times\left[0, 1\right]$, and the maps
$$r=r(\kappa),\quad z=z(\eta),$$
where we assume these two maps and their derivatives are all analytically known. We update these two maps from time to time according to some criteria and construct these two maps as  monotonically increasing functions. We will use these two maps to map the physical domain in $(r, z)$ to a computational domain in $(\kappa,\eta)$, so that $\omega_1(r(\kappa), z(\eta))$ and $\psi_1(r(\kappa), z(\eta))$ as functions of $(\kappa, \eta)$ are relatively smooth. Let $n_\kappa$, $n_\eta$ be the number of grid points along the $r$- and $z$- directions, respectively. And let $h_\kappa=1/n_\kappa$, $h_\eta=1/n_\eta$ be the mesh sizes along the $r$- and $z$- directions respectively. We place a uniform mesh on the computation domain of $(\kappa,\eta)$:
$$\mathcal{M}_{(\kappa,\eta)}=\left\{(ih_\kappa, jh_\eta): 0\leq i\leq n_\kappa, 0\leq j\leq n_\eta\right\}.$$
This is equivalent to covering the physical domain of $(r, z)$ with the tensor-product mesh:
$$\mathcal{M}_{(r,z)}=\left\{(r(ih_\kappa), z(jh_\eta)): 0\leq i\leq n_\kappa, 0\leq j\leq n_\eta\right\}.$$
With properly chosen maps of $r=r(\kappa)$ and $z=z(\eta)$, the mesh $\mathcal{M}_{(r,z)}$ can focus on the singular part of the solution, so that the accuracy of the numerical solution can be greatly improved.

As we will see in the following sections, the singular part of the solutions will gradually move towards the origin. Thus we dynamically update the maps to accommodate the movement of the focused region. The update of the maps is based on an adaptive strategy that quantitatively locates the singular part of the solution and then decides the necessity to change the maps, as well as the parameters for the new maps. Once we update the maps, we interpolate the solutions from the old mesh to the new mesh and use the new computational domain.
In our algorithm, we adopt a second-order implementation for our adaptive mesh method. In Section \ref{sec: resolution_study}, we will perform resolution study to confirm the order of accuracy of our numerical method.

\section{Numerical evidence for a potential self-similar singularity}
\label{sec: evidence}

In this section, we will focus on the case with H\"{o}lder exponent $\alpha=0.1$, and provide numerical evidences for the potential self-similar singularity observed from the 3D axisymmetric Euler equations with no swirl and with H\"{o}lder continuous initial data. For the cases with different values of H\"{o}lder exponent $\alpha$, we will present the results in Section \ref{sec: holder_dimension}.

\begin{figure}[hbt!]
\centering
\includegraphics[width=.33\textwidth]{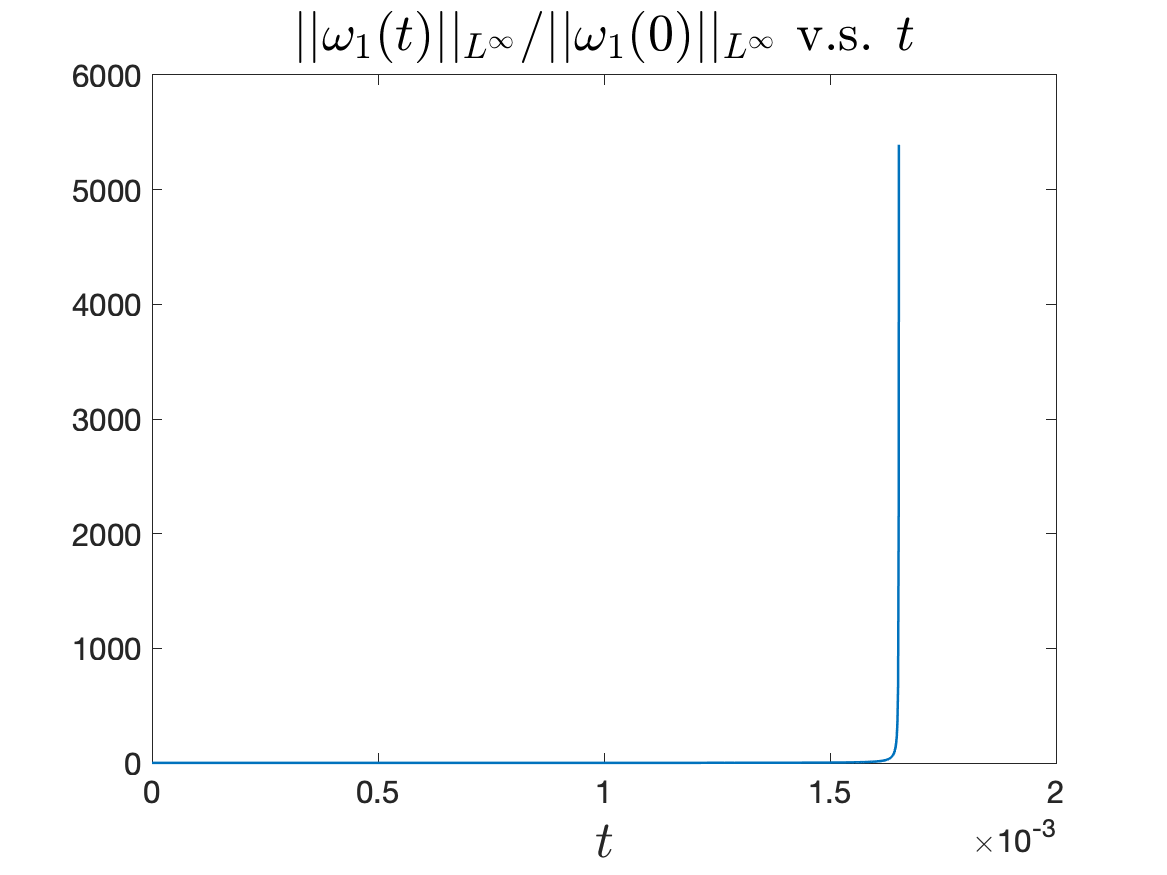}\includegraphics[width=.33\textwidth]{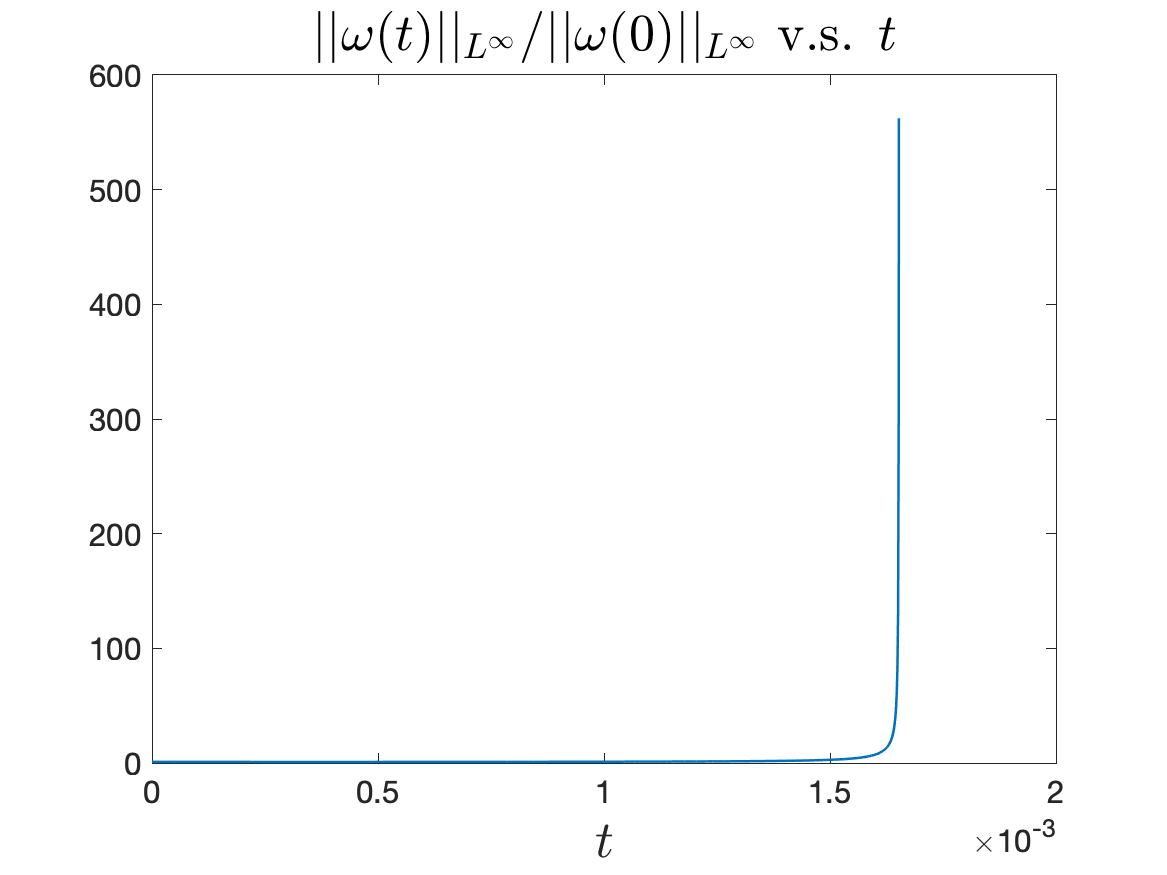}\includegraphics[width=.33\textwidth]{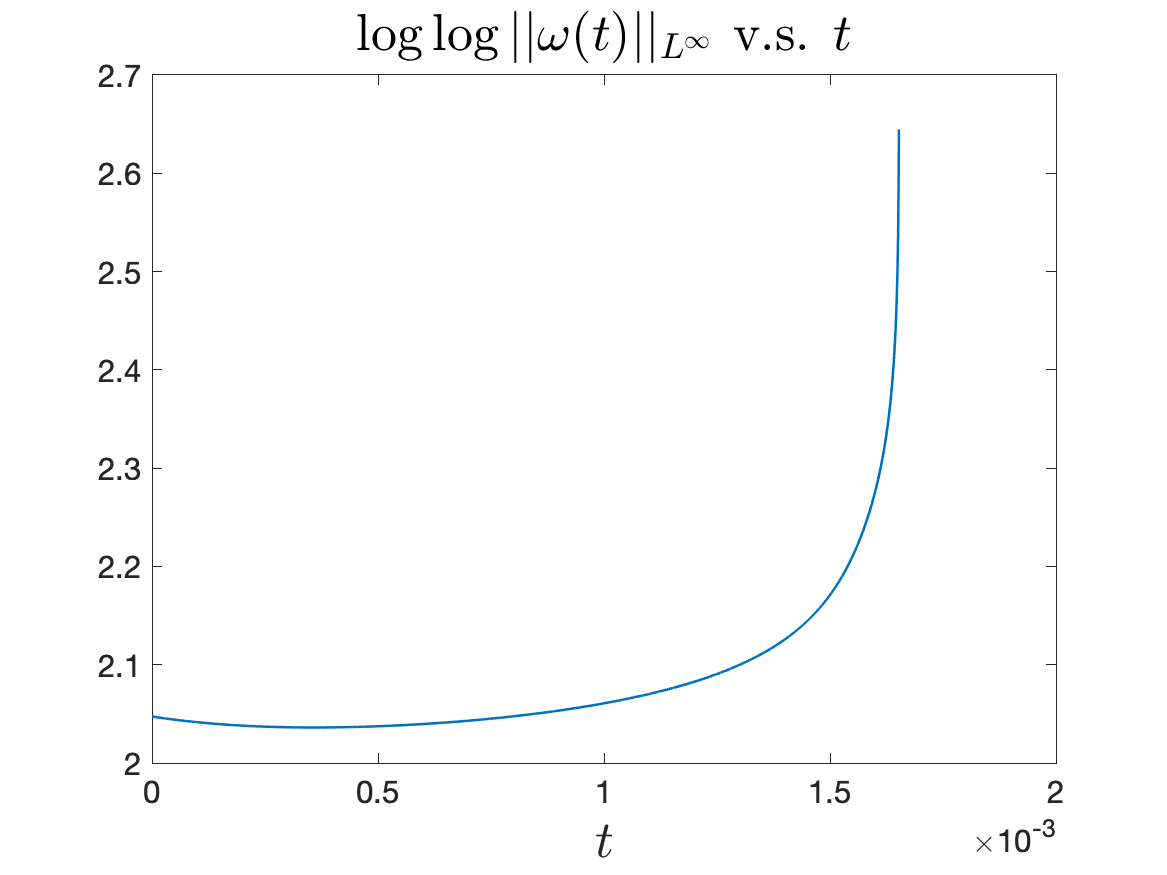}\\
\includegraphics[width=.33\textwidth]{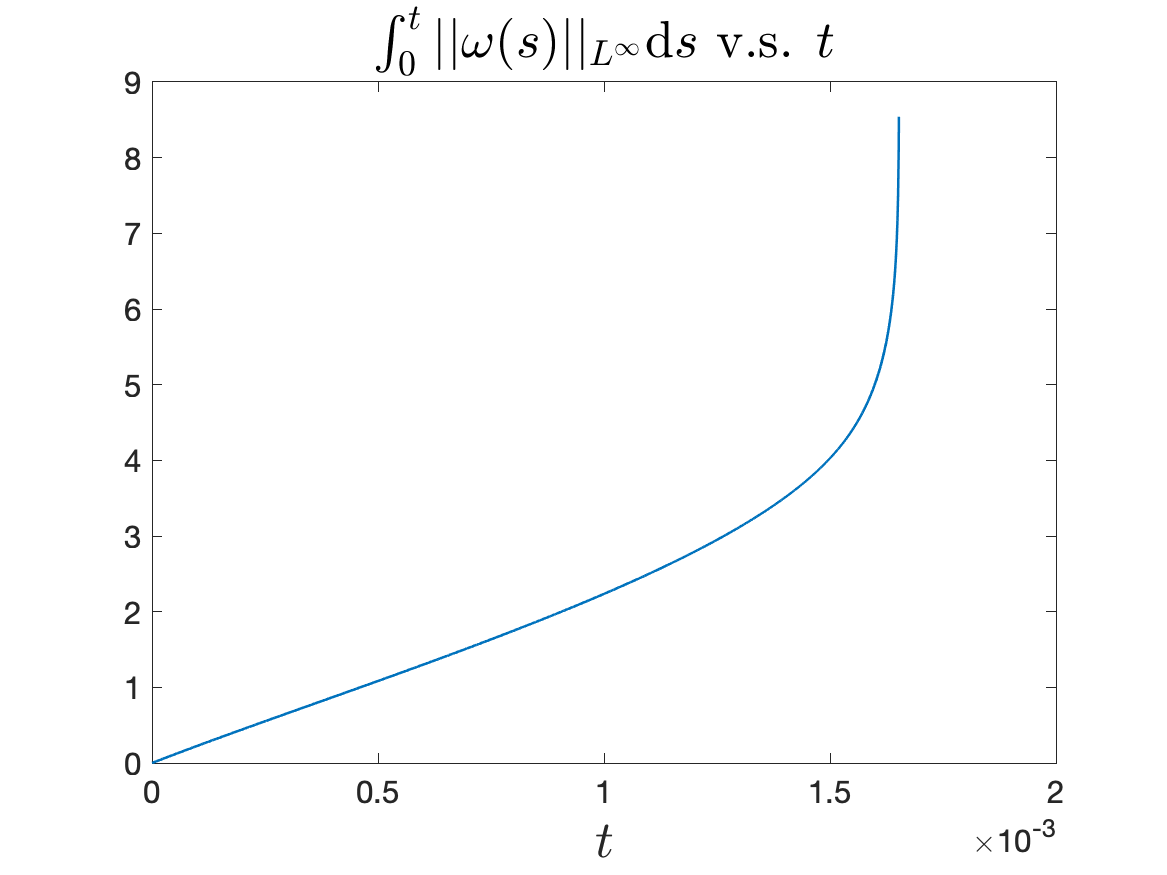}
\includegraphics[width=.33\textwidth]{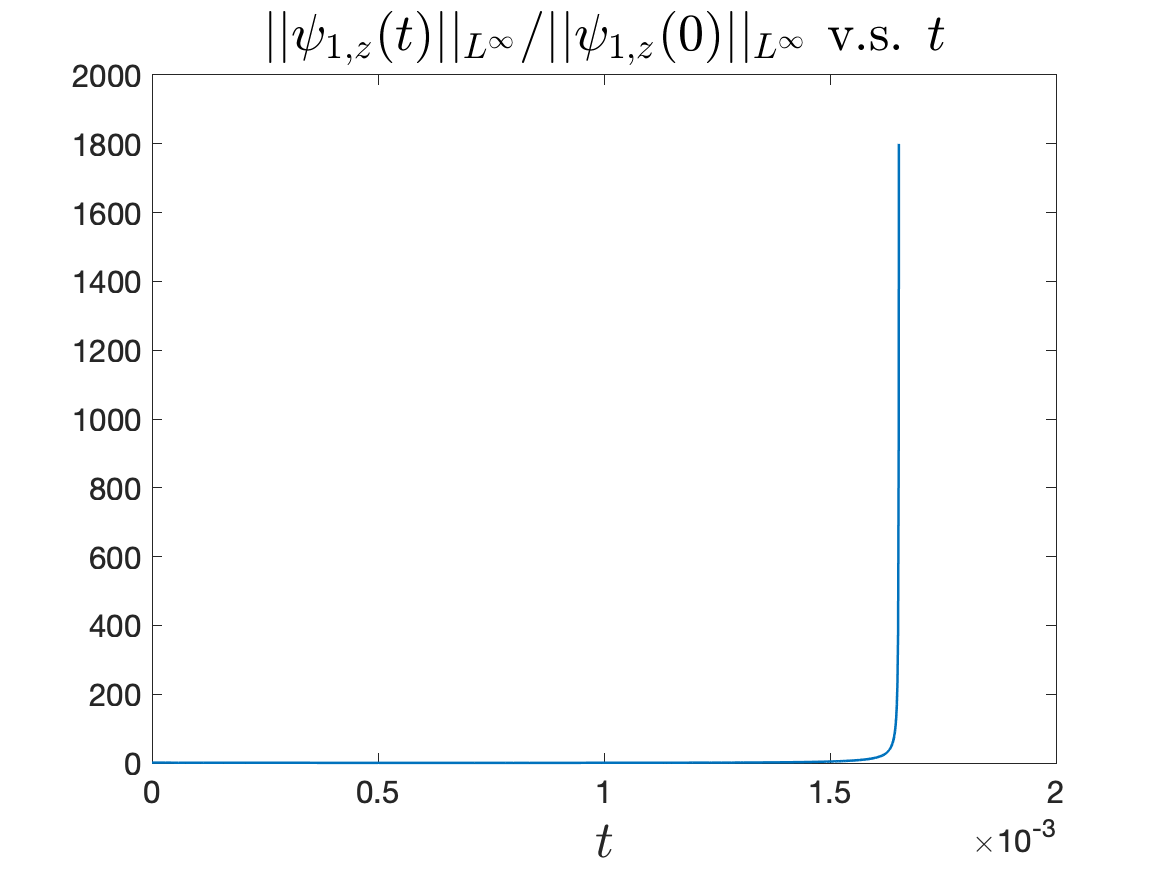}
\caption{Curves of $\|\omega_1\|_{L^\infty}$, $\|\omega\|_{L^\infty}$, $\log\log\|\omega\|_{L^\infty}$, $\int_0^t\|\omega(s)\|_{L^\infty}\mathrm{d}s$, $\|\psi_{1,z}\|_{L^\infty}$ as functions of time $t$.}\label{fig: stats_curve_1}
\end{figure}

\subsection{Evidence for a potential singularity}
\label{sec: evidence_singular}

\begin{figure}[hbt!]
\centering
\includegraphics[width=.33\textwidth]{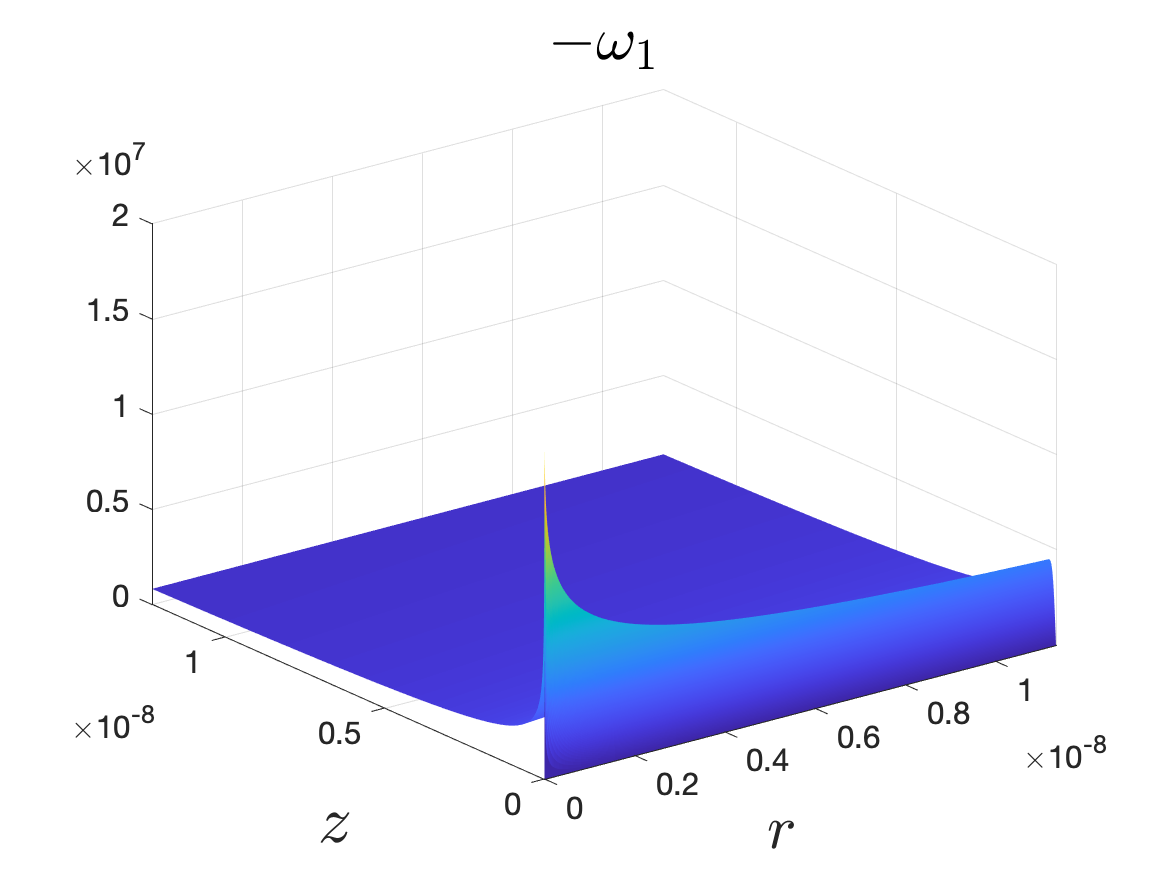}\includegraphics[width=.33\textwidth]{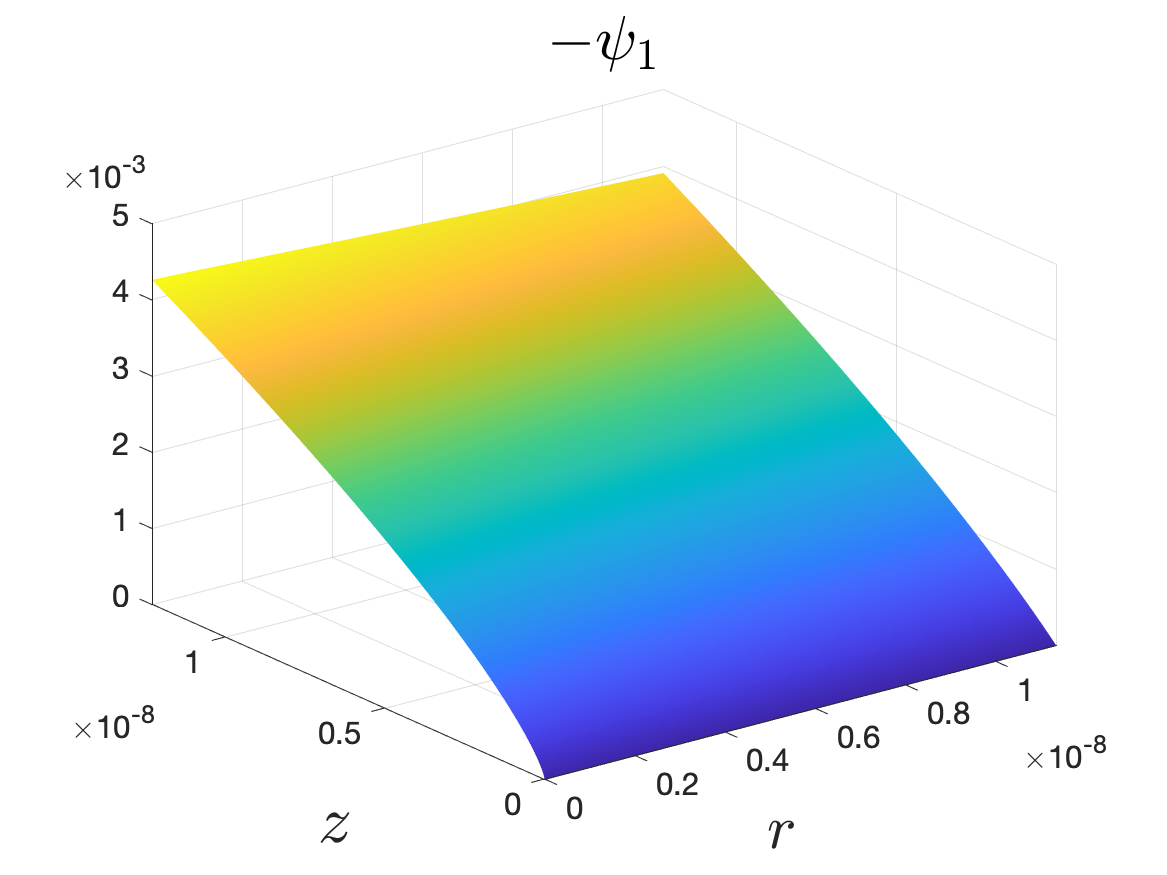}\includegraphics[width=.33\textwidth]{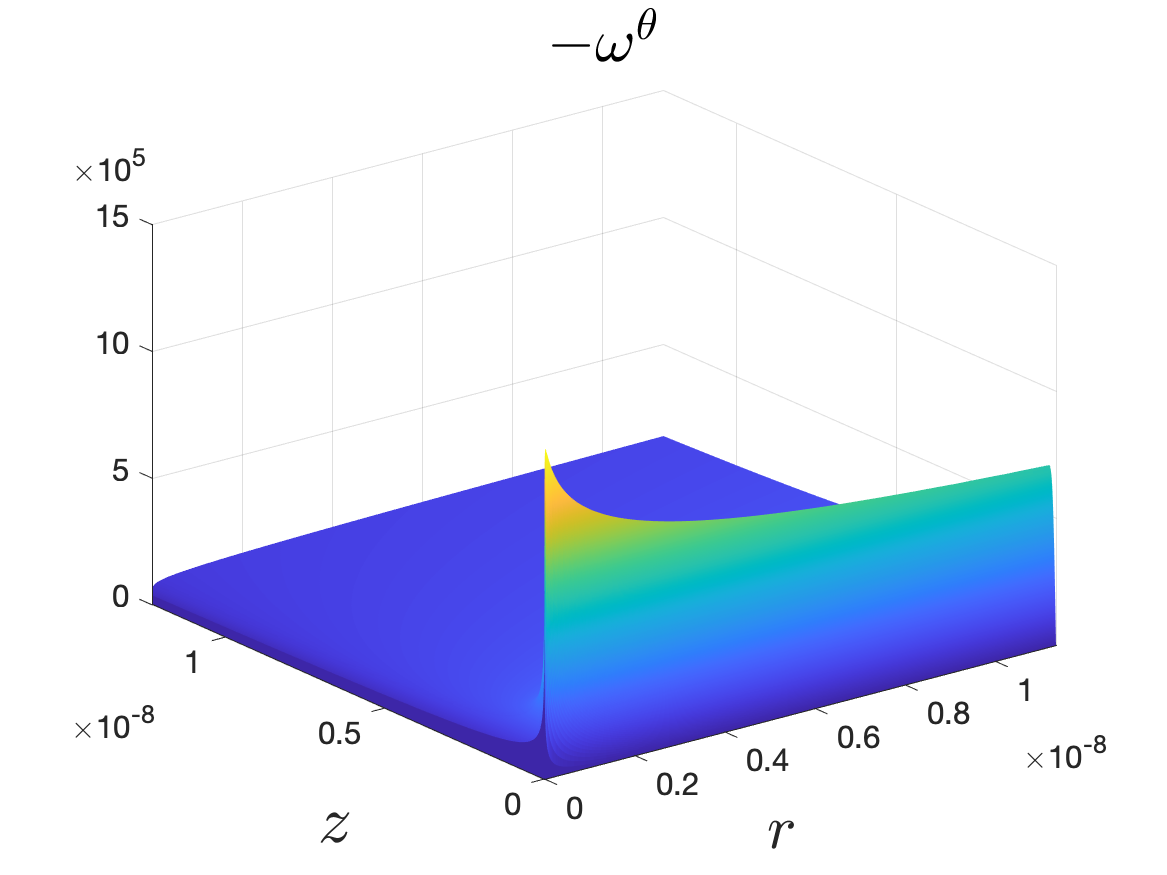}\\
\includegraphics[width=.33\textwidth]{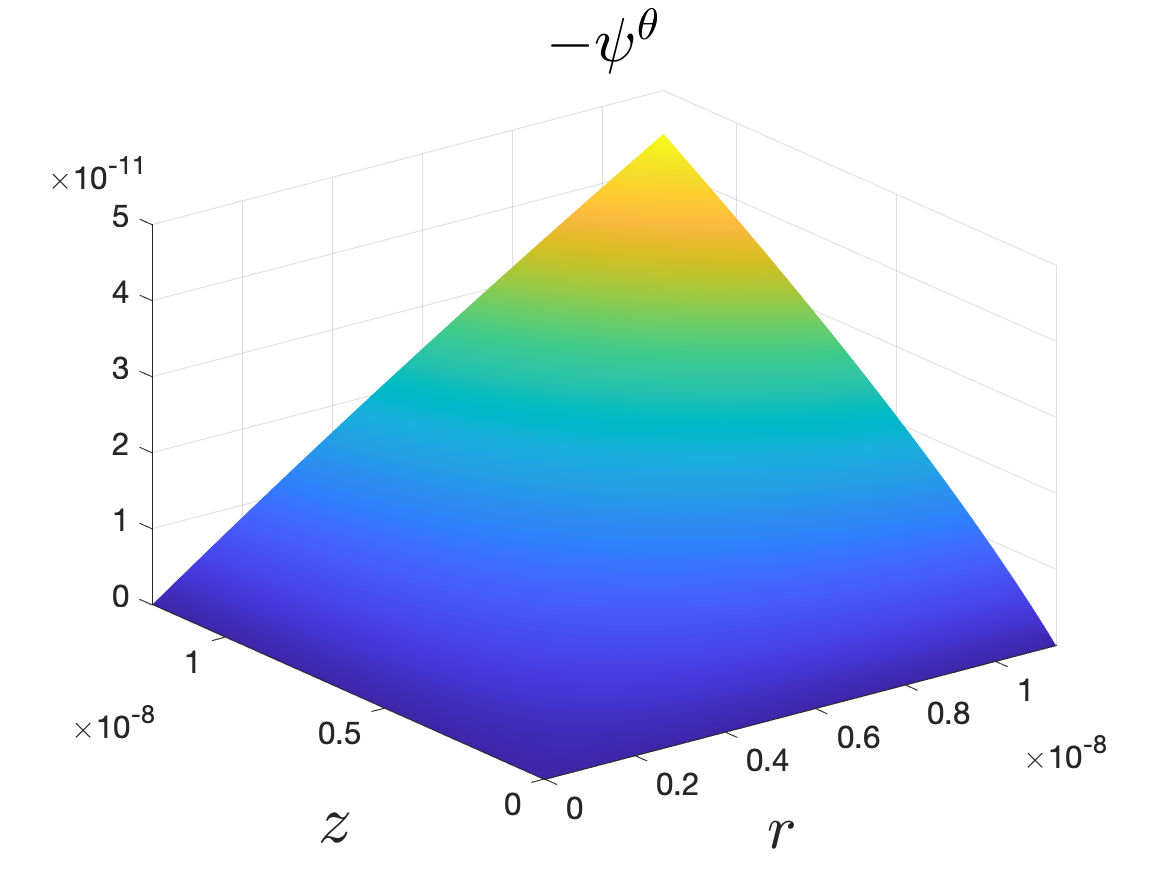}\includegraphics[width=.33\textwidth]{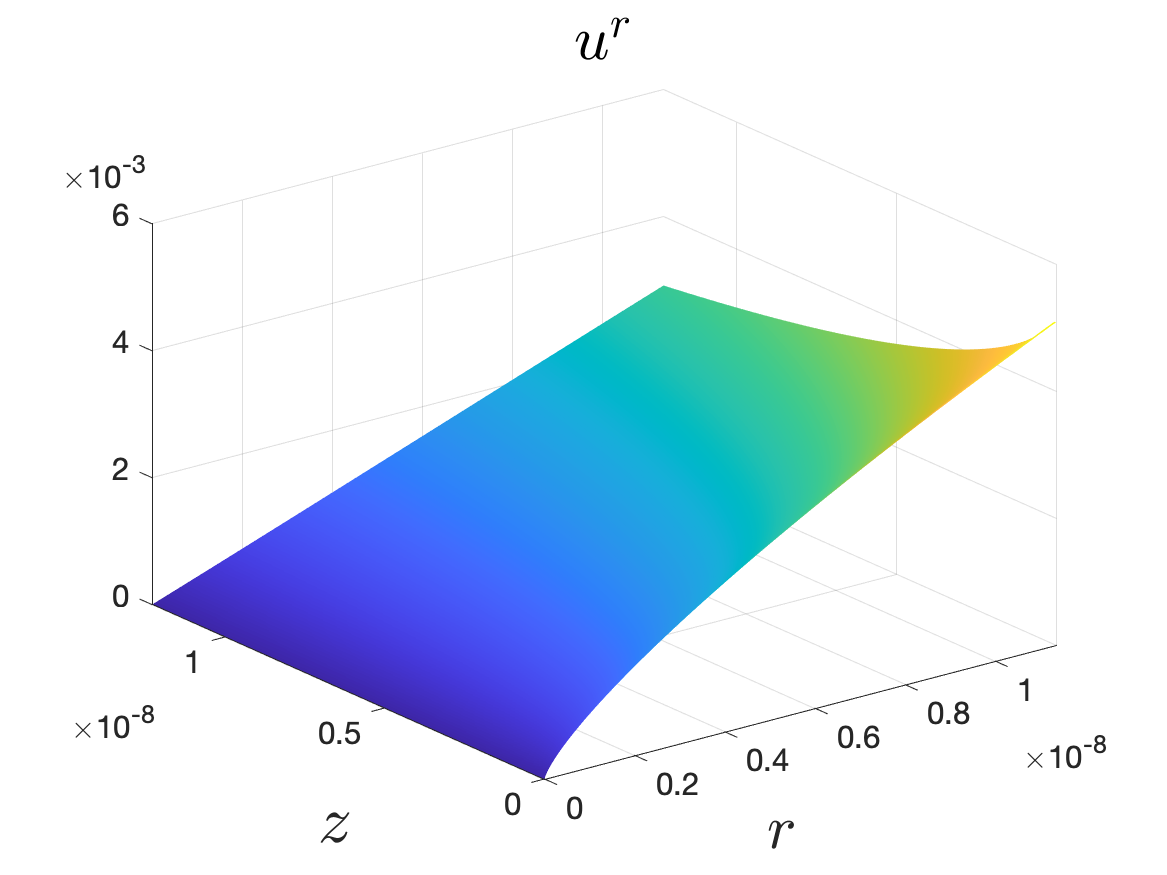}\includegraphics[width=.33\textwidth]{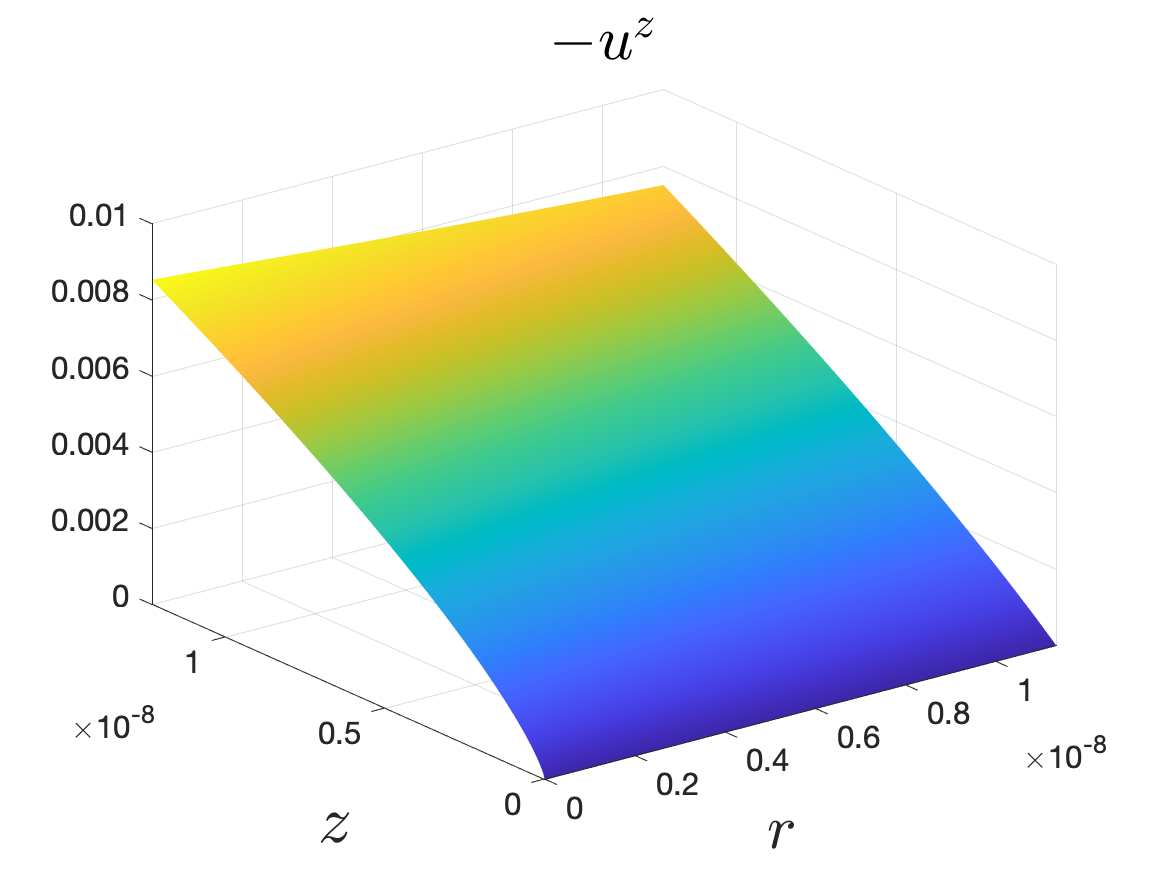}
\caption{Profiles of $-\omega_1$, $-\psi_1$, $-\omega^\theta$, $-\psi^\theta$, $u^r$ and $-u^z$ near the origin at $t=1.6524635\times10^{-3}$.}\label{fig: end_data}
\end{figure}

On $1024\times1024$ spatial resolution, we use the adaptive mesh method to solve \eqref{eq: vort_stream_1_3d_noswirl} with H\"{o}lder exponent $\alpha=0.1$, until the time when the smallest adaptive mesh size gets close to the machine precision. The final time of the computation is at $t=1.6524635\times10^{-3}$, after more than $6.5\times10^4$ iterations in time.

In Figure \ref{fig: stats_curve_1}, we plot the dynamic growth of several important quantities of the solution. The magnitude of $\omega_1$ has grown significantly, especially near the end of the computation. At the final time of the computation, $\|\omega_1\|_{L^\infty}$ has increased by a factor of around $5400$, and $\|\omega\|_{L^\infty}$ has increased by a factor of more than $560$. We also observe that the double logarithm curve of the maximum vorticity magnitude, $\log\log\|\omega\|_{L^\infty}$, maintains a super-linear growth, and the time integral $\int_0^t\|\omega(s)\|_{L^\infty}\mathrm{d}s$ has rapid growth with strong growth inertia close to the stopping time. This provides strong evidence for a potential finite-time blow-up of the 3D Euler equations by the Beale-Kato-Majda blow-up criterion. We see rapid increase in $\|\psi_{1,z}\|_{L^\infty}$ in time, which shows very strong vortex stretching term $-(n-2-\alpha)\psi_{1,z}\omega_1$ near the potential finite-time blow-up. We analyze their scaling time $t$ in Section \ref{sec: scaling-analysis}.

\begin{figure}[hbt!]
\centering
\includegraphics[width=.4\textwidth]{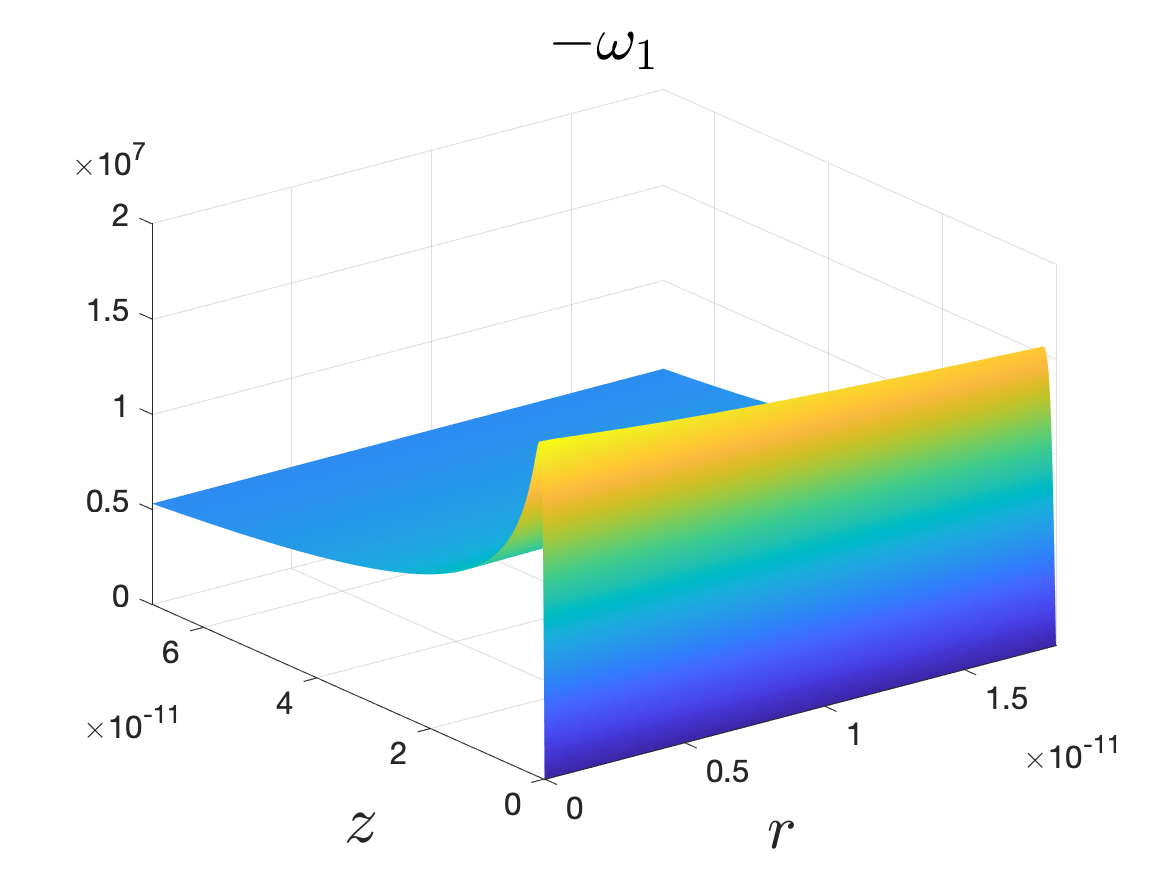}\includegraphics[width=.4\textwidth]{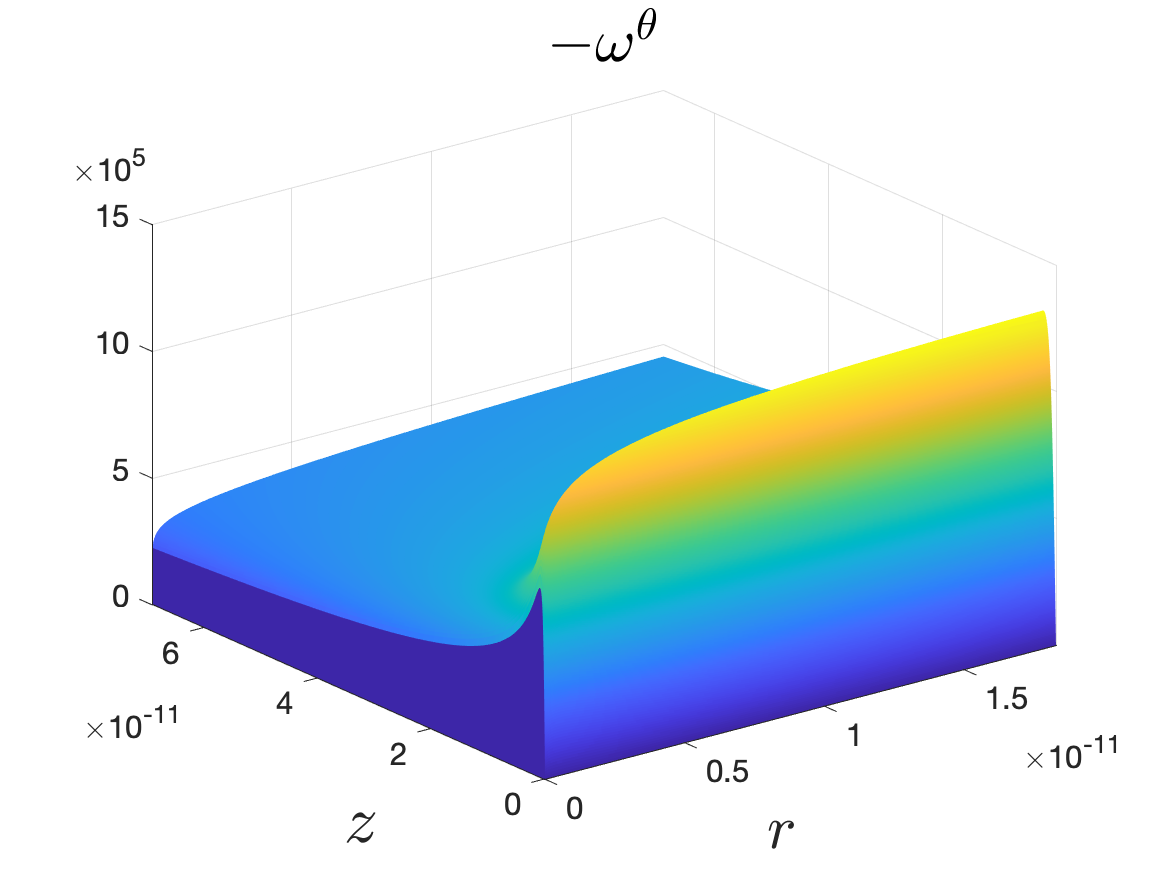}
\caption{Further zoomed-in profiles of $-\omega_1$ and $-\omega^\theta$ near the origin at $t=1.6524635\times10^{-3}$.}\label{fig: end_data_zoom}
\end{figure}

In Figure \ref{fig: end_data}, we plot the 3D profiles of $\omega_1$, $\psi_1$, $\omega^\theta$, $\psi^\theta$, $u^r$, and $u^z$ near the origin at end of our computation. We can see that $\omega_1$ is very concentrated near the origin, and so is $\omega^\theta$. Therefore, we further zoom-in around the origin and plot the local near field profiles of $\omega_1$ and $\omega^\theta$ in Figure \ref{fig: end_data_zoom}. We observe that the ``peak'' of $-\omega_1$ locates at the $z$-axis where $r=0$, and is being pushed toward the origin as implied by the velocity field $u^r$, $u^z$. We denote by $(R_1(t), Z_1(t))$ the position at which  $|\omega_1|$ achieves its maximum at time $t$. We have $R_1(t)=0$. At $(R_1(t), Z_1(t))$, the radial velocity $u^r$ is zero, and the axial velocity $u^z$ is negative, which pushes $(R_1(t), Z_1(t))$ toward the origin.

\begin{figure}[hbt!]
\centering
\includegraphics[width=.4\textwidth]{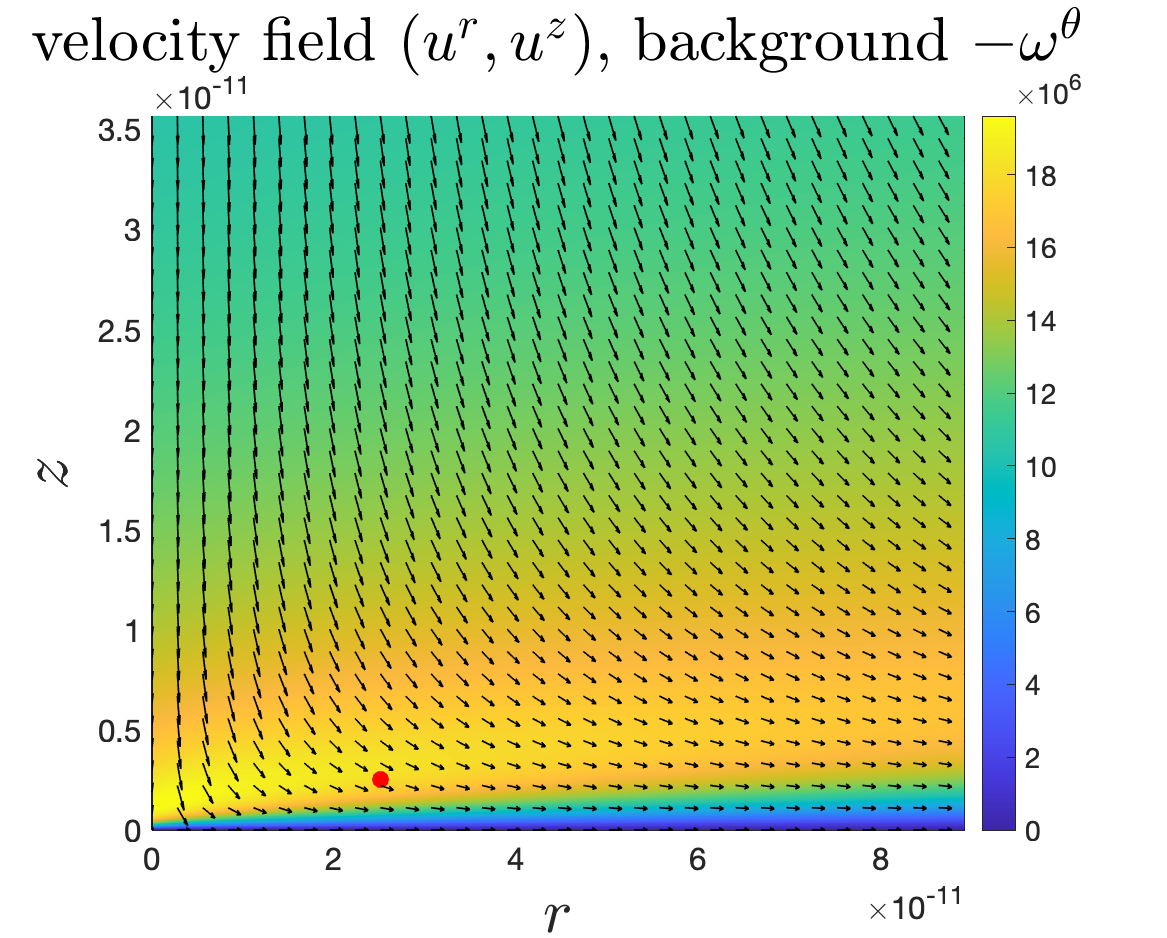}\includegraphics[width=.4\textwidth]{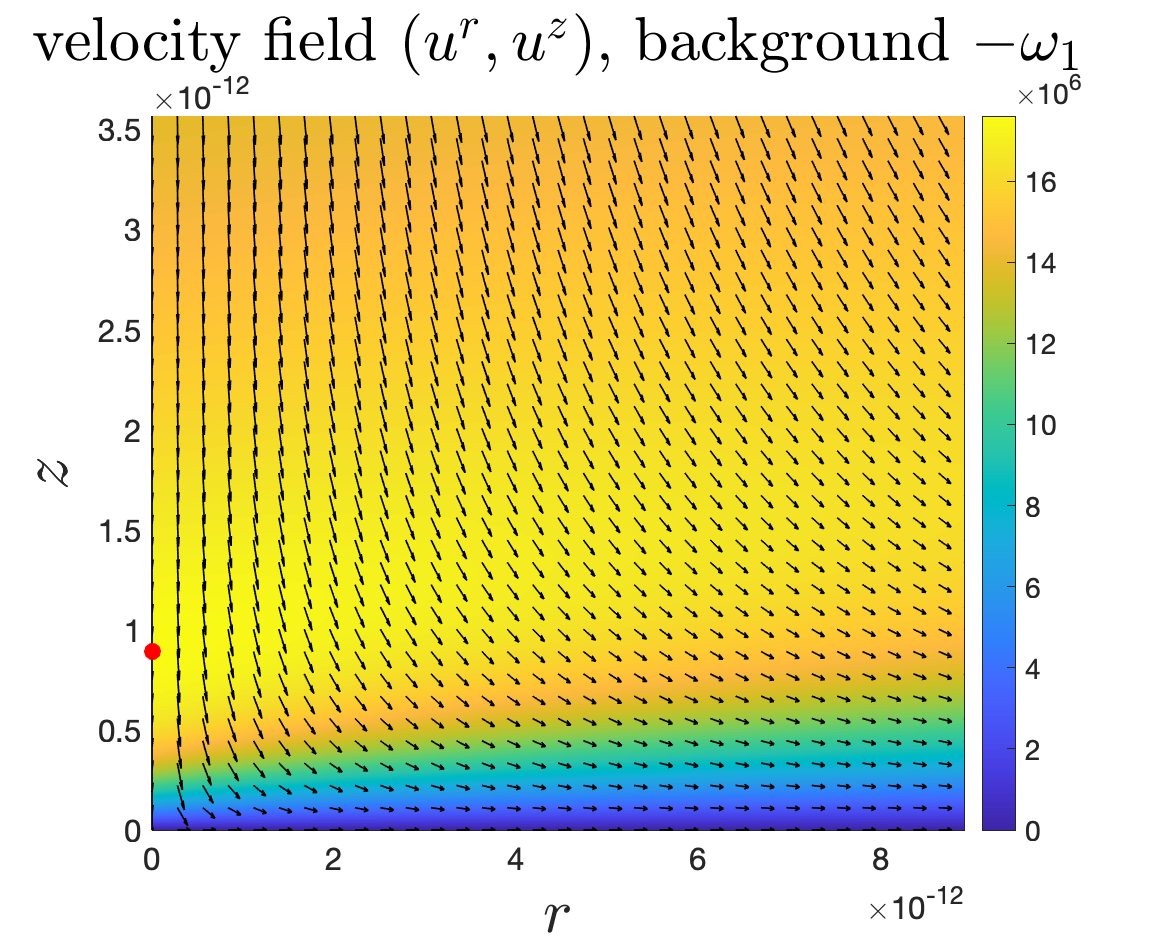}
\caption{The local velocity field near the maximum of $-\omega^\theta$ and $-\omega_1$. The pseudocolor plot of $-\omega^\theta$ or $-\omega_1$ is the background, and the red dot is its maximum.}\label{fig: velo_field}
\end{figure}

\begin{figure}[hbt!]
\centering
\includegraphics[width=.4\textwidth]{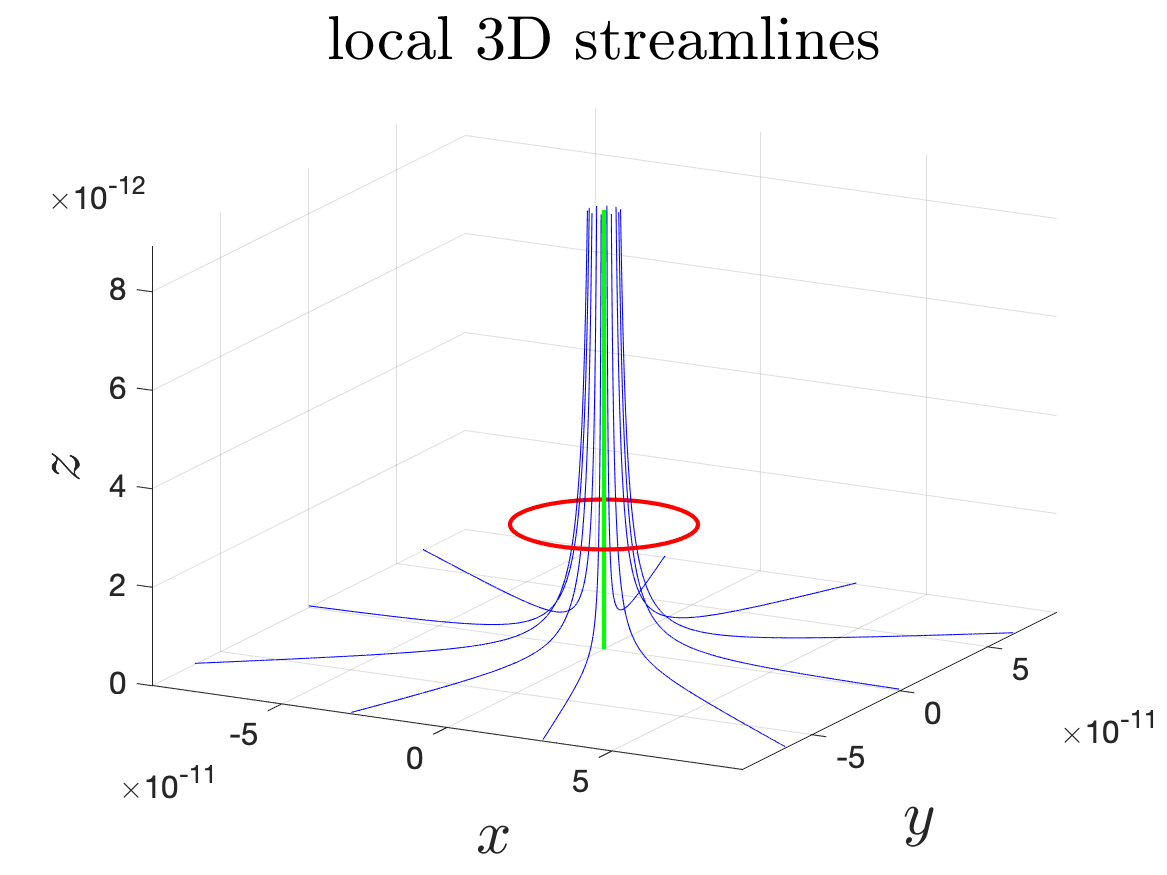}\includegraphics[width=.4\textwidth]{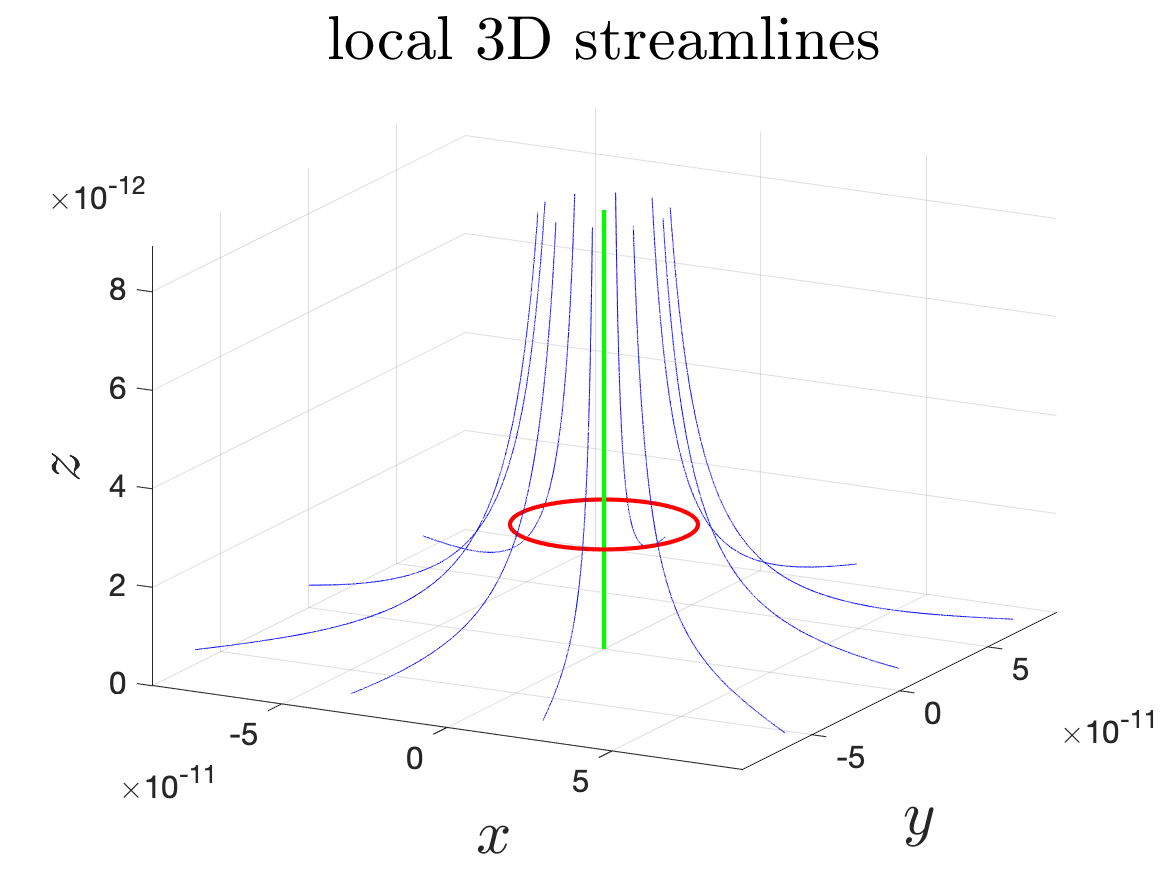}
\caption{The local streamlines near the origin. The green pole is the $z$-axis, and the red ring is where $-\omega^\theta$ achieves its maximum.}\label{fig: streamline}
\end{figure}

In Figure \ref{fig: velo_field}, we plot the local velocity field near the maximum of $-\omega^\theta$ and $-\omega_1$, respectively. We use the pseudocolor plots of $-\omega^\theta$ and $-\omega_1$ as the background respectively for the figure in the left and the right subplots, and mark the maximum of $-\omega^\theta$ or $-\omega_1$ with the red dot. The velocity field demonstrates a clear hyperbolic structure as depicted by Figure \ref{fig: hyper_flow}. And the velocity field clearly pushes the maximum $(R_1(t), Z_1(t))$ of $-\omega_1$ toward the origin.

In Figure \ref{fig: streamline}, we show the local streamlines near the maximum of $-\omega^\theta$ in $\mathbb{R}^3$. The maximum of $-\omega^\theta$ locates on the red ring centered at $(0,Z_1(t))$ along the $z$-axis. In the left figure, we plot a set of streamlines that travel through the maximum ring from top to bottom. And in the right figure, we plot a set of streamlines that travel around the maximum ring from top to bottom. From Figure \ref{fig: streamline}, we notice that the streamlines are axisymmetric, and do not form swirl around the $z$-axis.

\begin{figure}[hbt!]
\centering
\includegraphics[width=.4\textwidth]{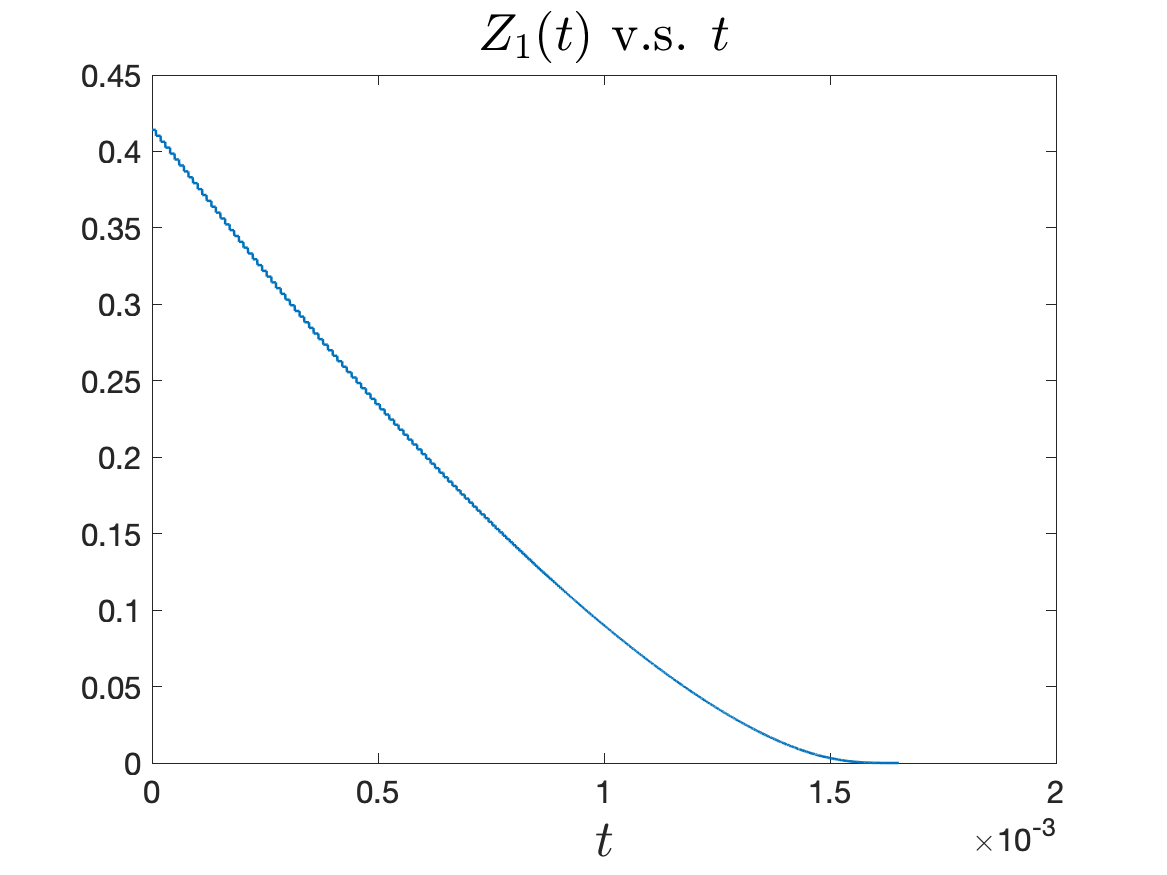}\includegraphics[width=.4\textwidth]{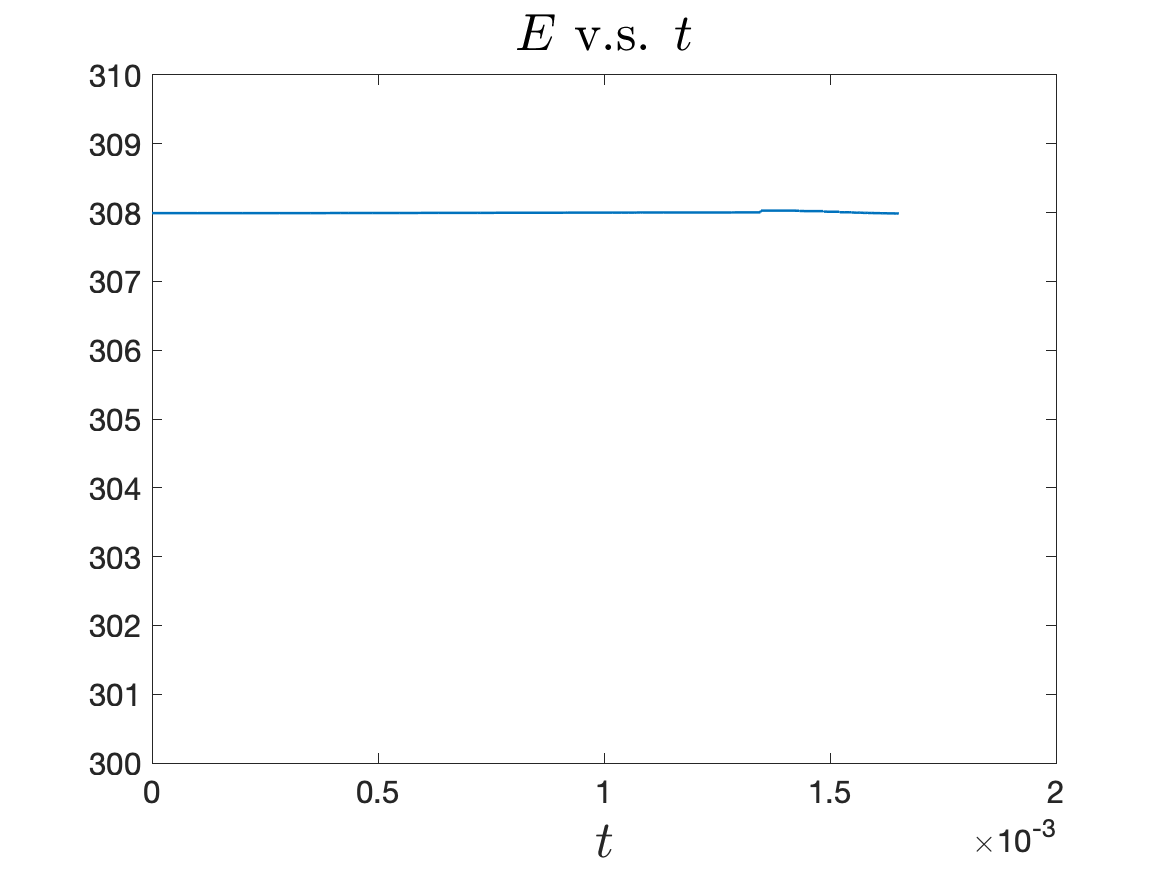}
\caption{Curves of $Z_1$ and $E$ as functions of time $t$.}\label{fig: stats_curve_2}
\end{figure}

In Figure \ref{fig: stats_curve_2}, we show the curve of the maximum location of $-\omega_1$, $Z_1$ and the kinetic energy, $E$ as functions of time. We can see that $Z_1(t)$ monotonically decreases to zero with $t$. The curve of $Z_1(t)$ seems to be convex, especially in time windows close to the stopping time. We refer to Section \ref{sec: scaling-analysis} for more study of the behavior of $Z_1(t)$. The kinetic energy $E$, which is defined as
$$E=\frac{1}{2}\int_{\mathcal{D}}\left|u\right|^2\mathrm{d}x=\pi\int_0^1\int_0^{1/2}\left(\left|u^r\right|^2+\left|u^z\right|^2\right)r\mathrm{d}r\mathrm{d}z,$$
for our axisymmetric case with no swirl, is a conservative quantity of the 3D Euler equations. In Figure \ref{fig: stats_curve_2}, we can see that there is little change of the kinetic energy $E$ as a function of time $t$. In fact, the major reason for the change of $E$ in our computation is due to the update of adaptive mesh, where we need to interpolate $\omega_1$ and $\psi_1$ from an old mesh to a new mesh. Since the new adaptive mesh will be more focusing on the near field around the origin, the far field velocity field might lose some accuracy, leading to a change in the kinetic energy $E$. However, such an update of adaptive mesh occurs only 35 times out of the total 65000 iterations in time, and the change in the kinetic energy $E$ in each update is negligible. By the end of the computation, the change in the kinetic energy $E$ is at most $1.4\times10^{-4}$ of the magnitude of $E$.

\subsection{Evidence for a potential self-similar blow-up}
\label{sec: evidence_self_similar}

\begin{figure}[hbt!]
\centering
\includegraphics[width=.33\textwidth]{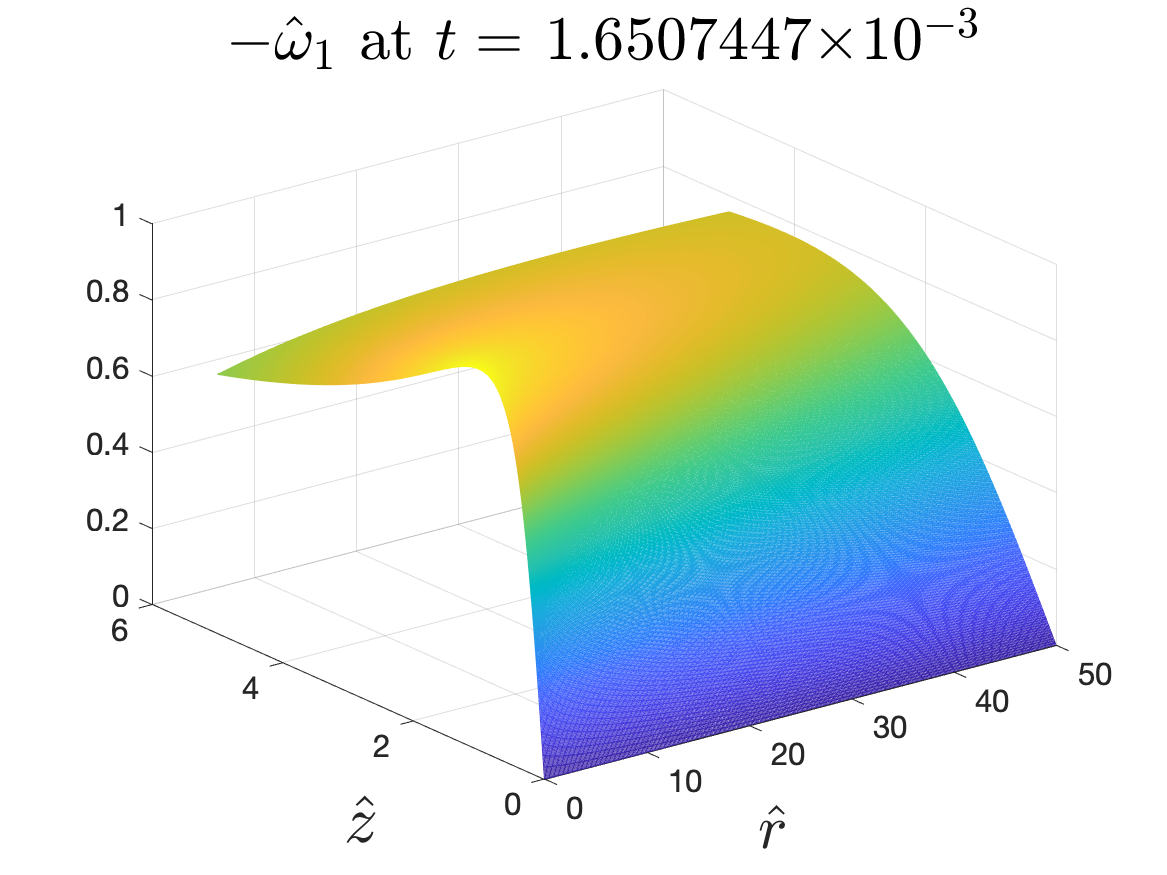}\includegraphics[width=.33\textwidth]{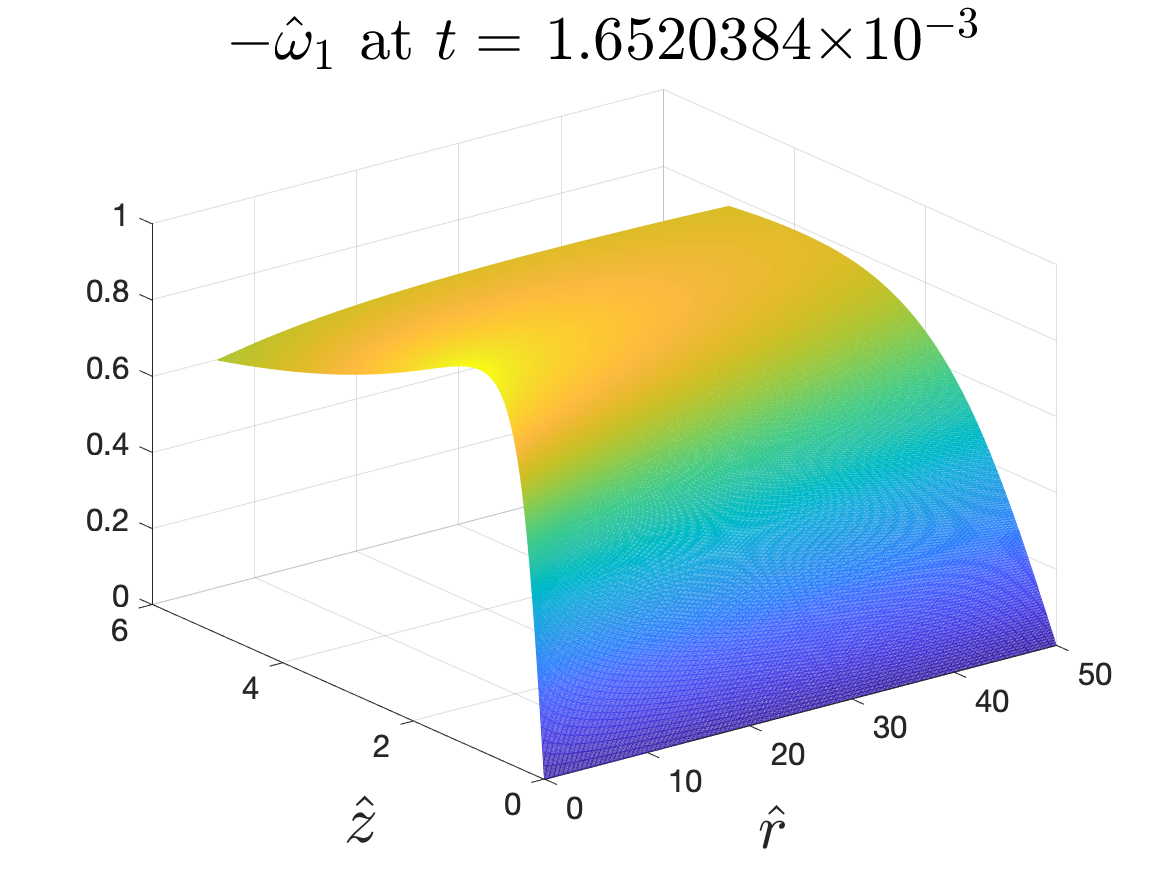}\includegraphics[width=.33\textwidth]{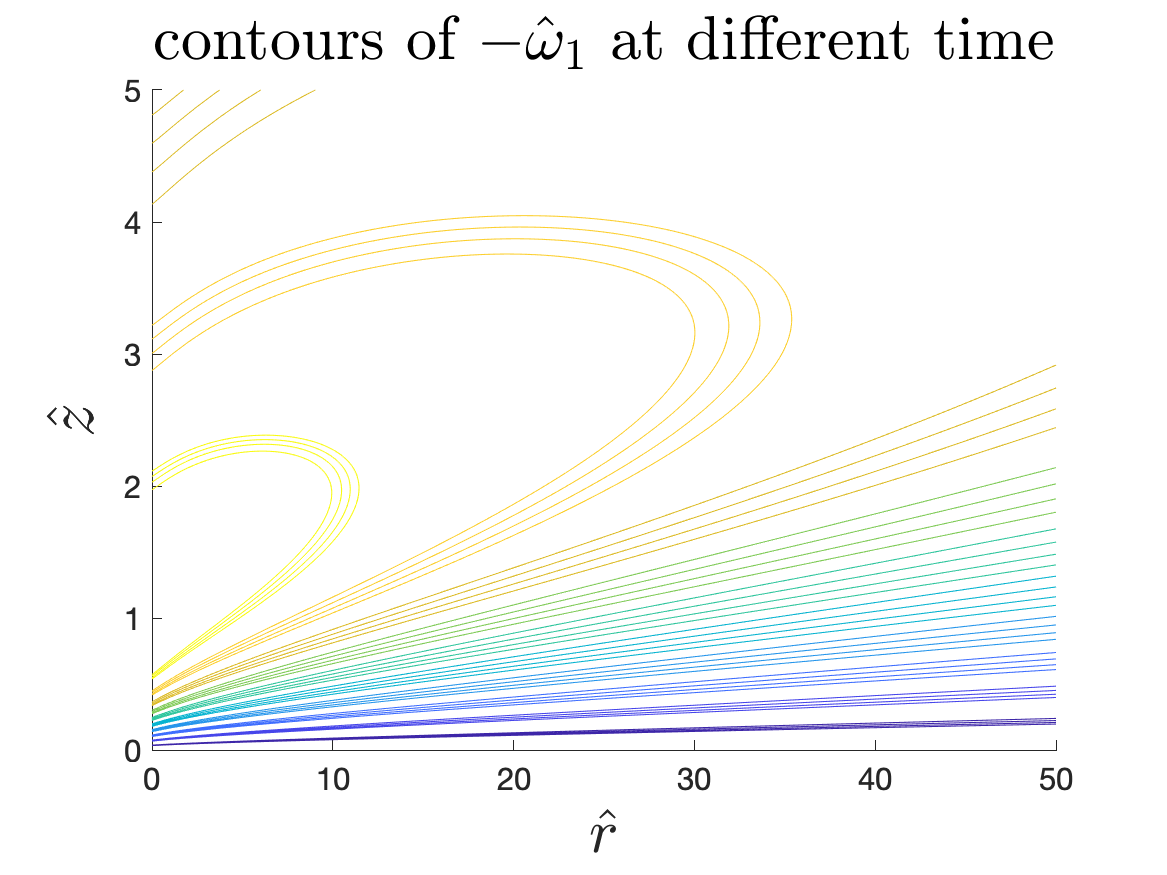}\\
\caption{Left and middle: Local profiles of $-\hat{\omega}_1$ at $t=\left\{1.6507447\right.$, $\left.1.6520384\right\}\times10^{-3}$. Right: Local contours of $-\hat{\omega}_1$ at $t=\left\{1.6507447\right.$, $1.6512953$, $1.6517173$, $\left.1.6520384\right\}\times10^{-3}$.}\label{fig: self_similar}
\end{figure}

We observe a potential self-similar blow-up in our numerical solution. To check the self-similar property, we visualize the local profile of the rescaled $\omega_1$ near the origin. Recall that $(0, Z_1)$ is the maximum location of $-\omega_1$, we define $$\hat{\omega}_1(\hat{r}, \hat{z},t)=\omega_1\left(Z_1(t)\hat{r}, Z_1(t)\hat{z}, t\right)/\|\omega_1(t)\|_{L^\infty},$$ as the rescaled version of $\omega_1$. The above definition rescales the magnitude of $|\hat{\omega}_1|$ to 1, and rescales the maximum location of $|\hat{\omega}_1|$ to $(\hat{r}, \hat{z})=(0, 1)$. We plot the profiles of $-\hat{\omega}_1$ near the origin at different time instants and the contours of $-\hat{\omega}_1(\hat{r}, \hat{z})$ at different times in Figure \ref{fig: self_similar}. The profile of $-\hat{\omega}_1$ seems to change slowly in the late time, indicating a potential self-similar structure of the blow-up profile near the origin. In other words, $x_0=0$ in the self-similar ansatz \eqref{eq: self_similar}.

\begin{figure}[hbt!]
\centering
\includegraphics[width=.4\textwidth]{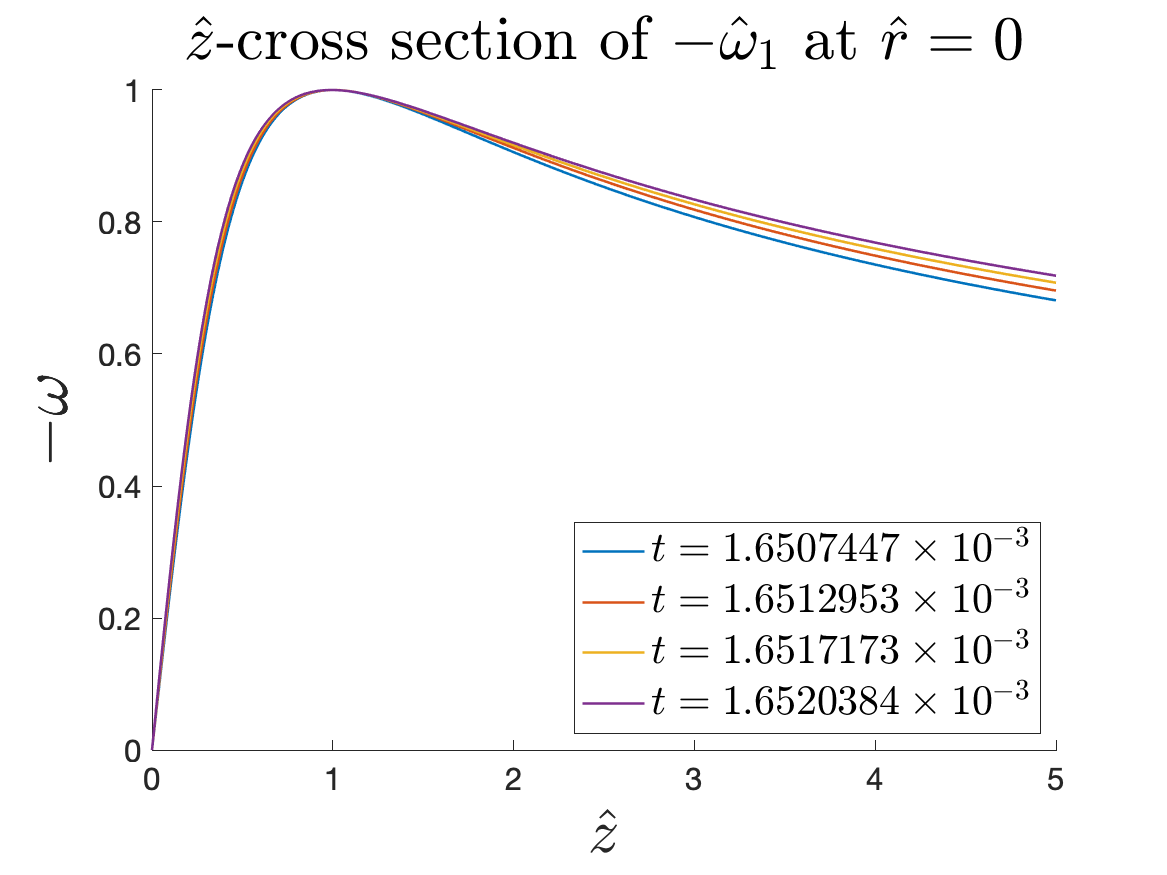}\includegraphics[width=.4\textwidth]{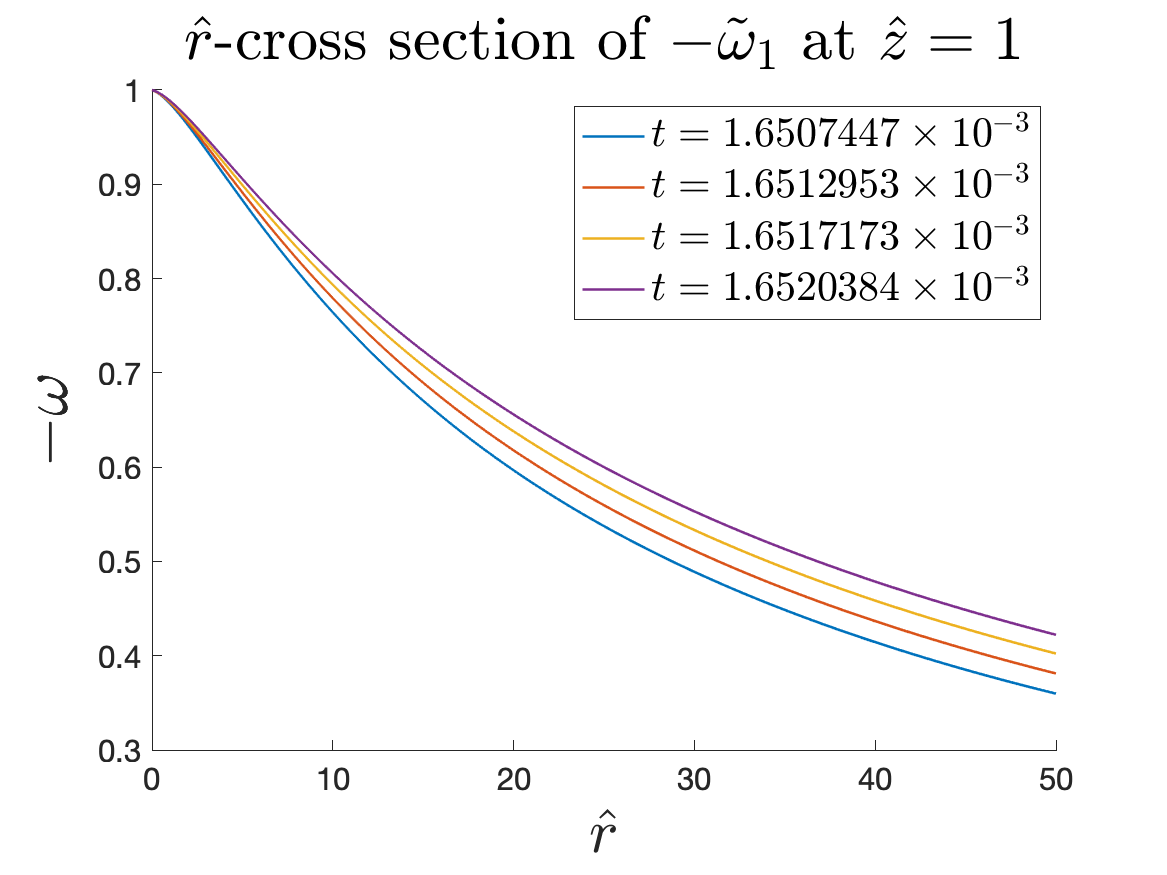}\\
\caption{Cross sections of $-\hat{\omega}_1$ at different times.}\label{fig: self_similar_cross_section}
\end{figure}

In Figure \ref{fig: self_similar_cross_section}, we plot the cross sections of $-\hat{\omega}_1$ at $\hat{r}=0$ and $\hat{z}=1$. The cross section at $\hat{r}=0$ shows a good potential for a self-similar blow-up, while the cross section at $\hat{z}=1$ shows that the blow-up profile has not converged to a self-similar profile yet. This is reasonable because although we are very close to the potential blow-up time, the strong collapsing along the $z$-direction and the effect of round-off errors prevent us from continuing the computation. We refer to Section \ref{sec: dynamic_rescaling} where we use the dynamic rescaling method and indeed observe numerically the convergence to the potential self-similar profile.

\subsection{Resolution study}
\label{sec: resolution_study}

\begin{figure}[hbt!]
\centering
\includegraphics[width=.33\textwidth]{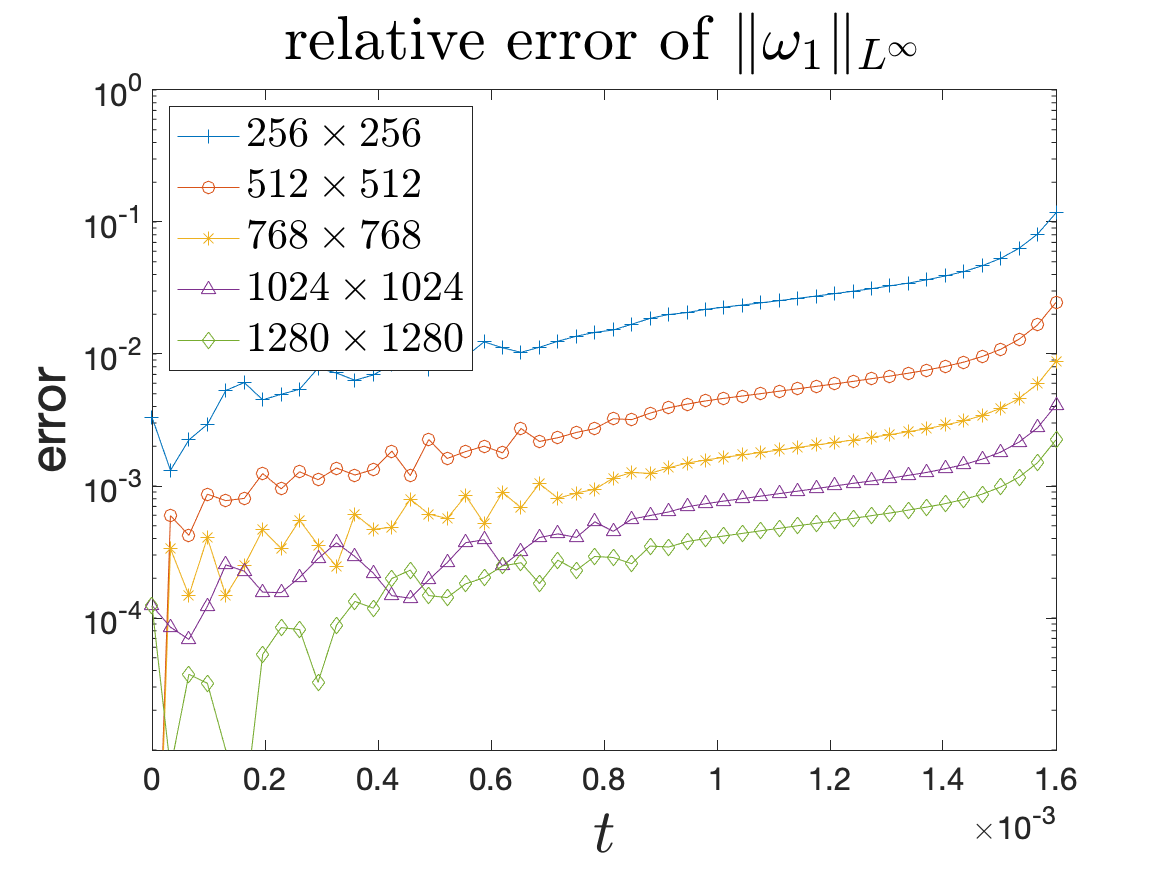}\includegraphics[width=.33\textwidth]{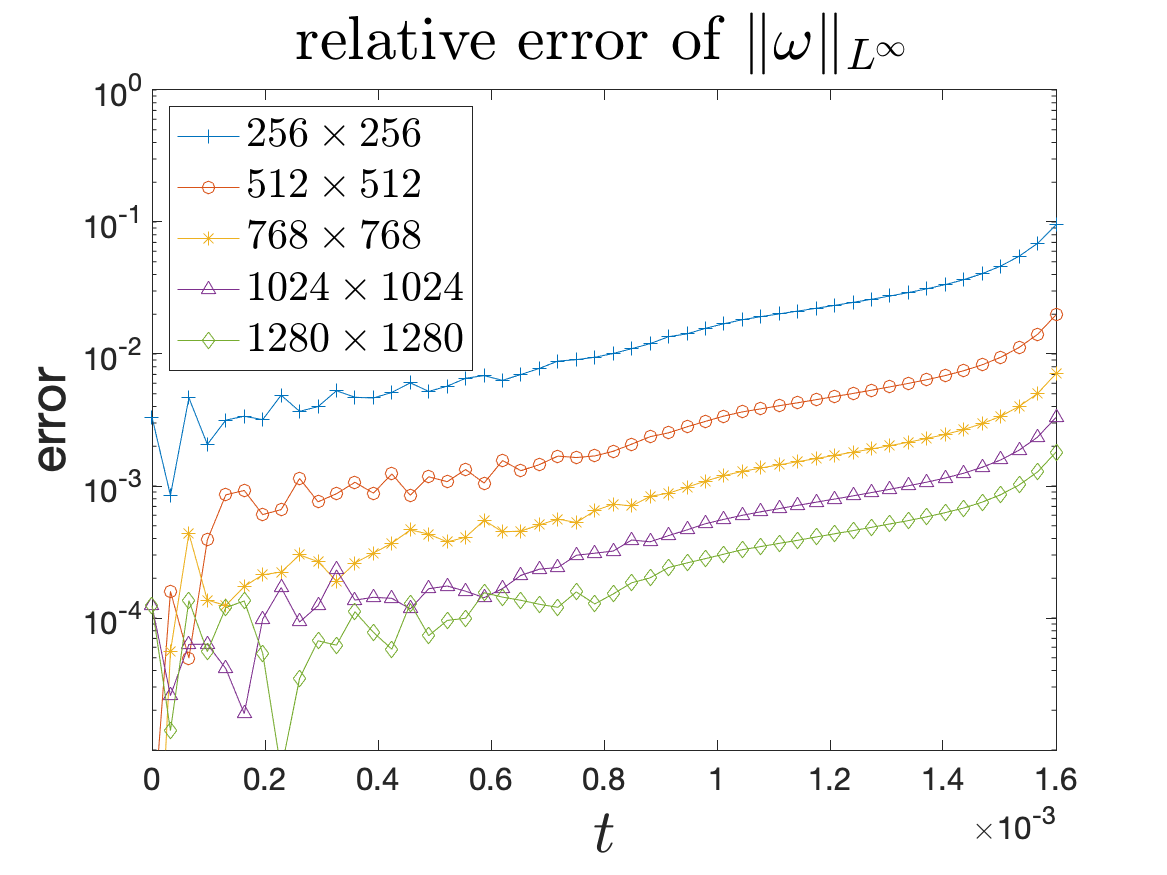}\includegraphics[width=.33\textwidth]{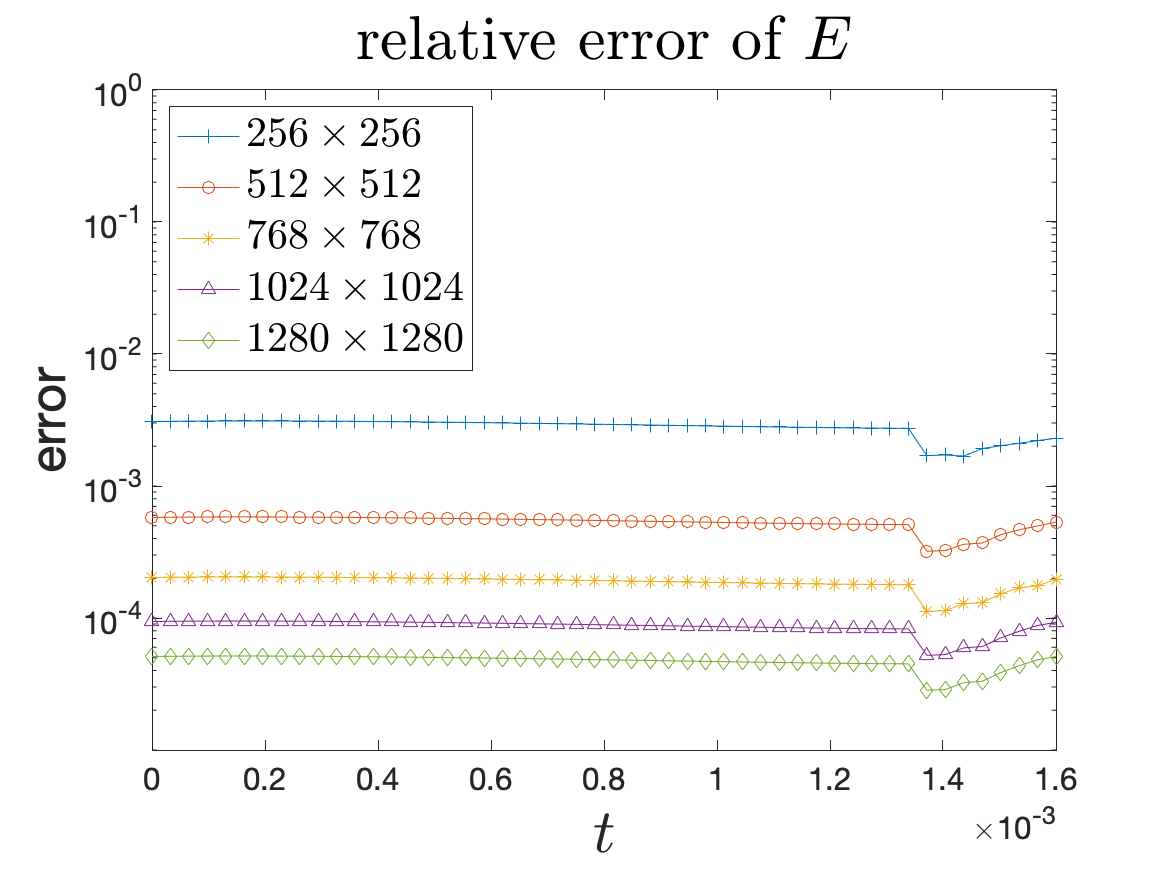}\\
\includegraphics[width=.33\textwidth]{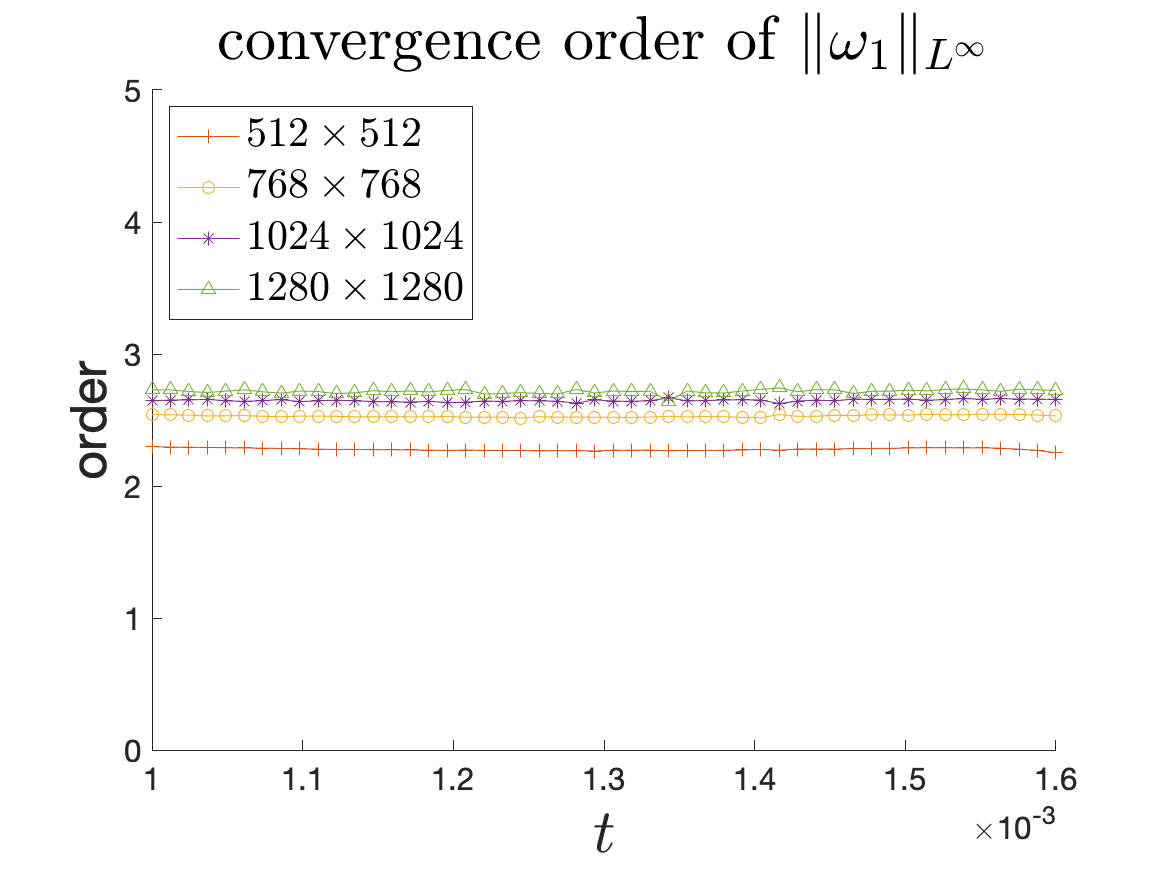}\includegraphics[width=.33\textwidth]{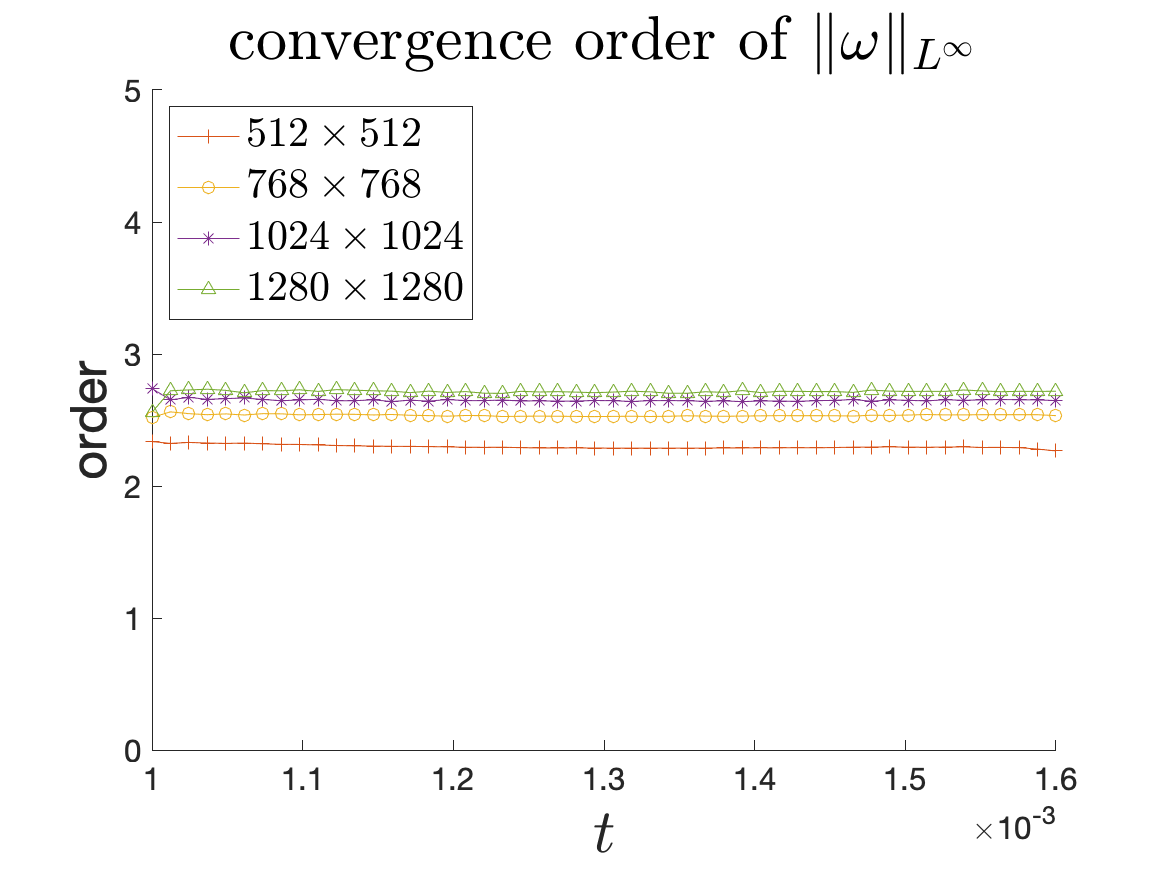}\includegraphics[width=.33\textwidth]{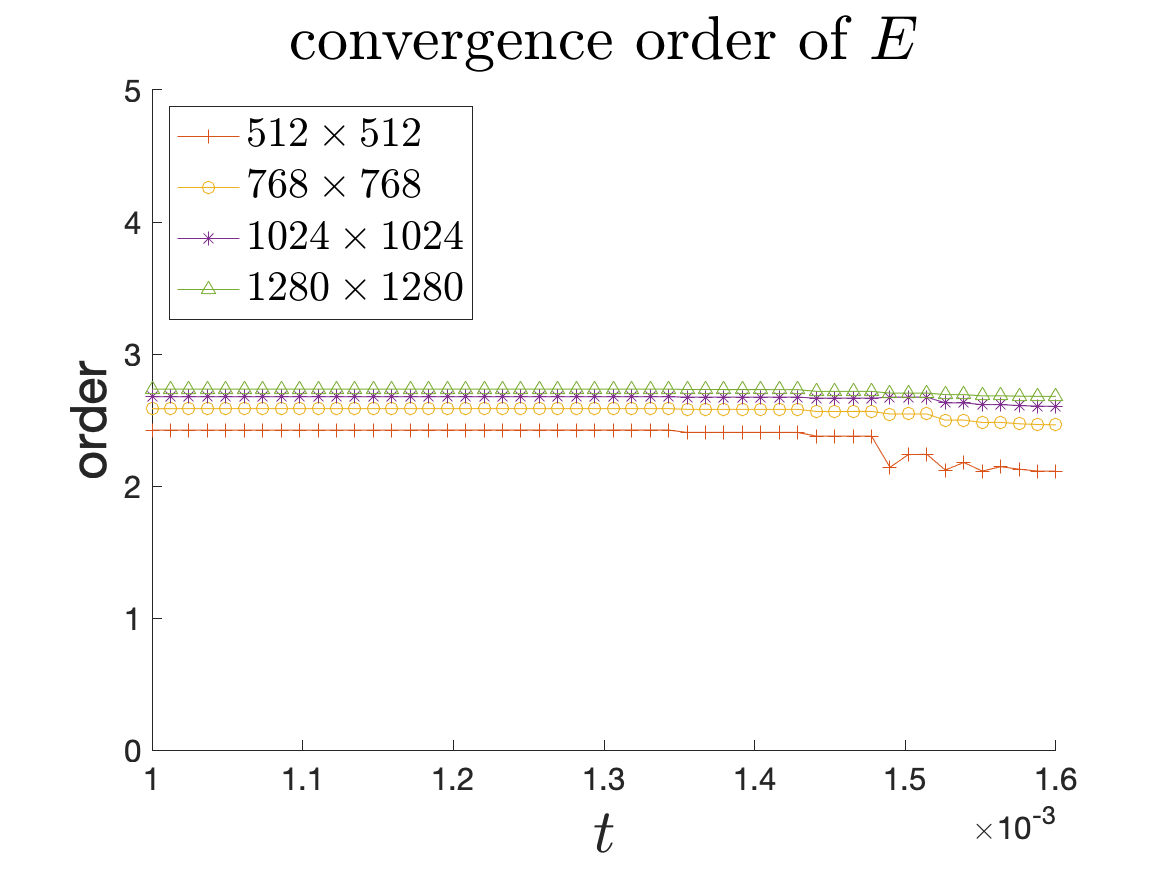}\\
\caption{Relative errors and convergence orders of $\|\omega_1\|_{L^\infty}$, $\|\omega\|_{L^\infty}$, and $E$ in sup-norm.}\label{fig: converg_order}
\end{figure}

We perform resolution study on the numerical solutions of \eqref{eq: vort_stream_1_3d_noswirl} to confirm the accuracy of our numerical solutions. We first simulate the equations on spatial resolutions of $256k\times256k$ with $k=1,2,\ldots,6$. The highest resolution we used is $1536\times1536$. Next, for the numerical solution at resolution $256k\times256k$, we compute its sup-norm relative error in several chosen quantities at selected time instants using the numerical solution at resolution $256(k+1)\times256(k+1)$ as the reference, for $k=1,2,\ldots,5$. Finally, we use the relative error obtained above to estimate the convergence order of the numerical method.

We consider two types of quantities. The first type is the function of the solutions. Here we consider the magnitude of $\omega_1$,
$\|\omega_1\|_{L^\infty}$, the maximum norm of vorticity, $\|\omega\|_{L^\infty}$, and the kinetic energy, $E$. We remark that $\|\omega_1\|_{L^\infty}$ and $\|\omega\|_{L^\infty}$ only depend on the local field near the origin, and $E$ should be considered as a global quantity. The second type is the vector fields of $\omega_1$, $\psi_1$, $u^r$, and $u^z$ that are actively participating in the simulated system \eqref{eq: vort_stream_1_3d_noswirl}. 

For each quantity, we use $q_k$ to represent the estimate we get at resolution $256k\times256k$. Then the sup-norm relative error $e_k$ is defined as
$$e_k=\|q_k-q_{k+1}\|_{L^\infty}/\|q_{k+1}\|_{L^\infty}.$$
If $q_k$ is a vector field, we first interpolate it to the reference resolution $256(k+1)\times256(k+1)$, and then compute the relative error as above. The convergence order of the error $\beta_k$ at this resolution can be estimated via 
$$\beta_k=\log\left(\frac{e_{k-1}}{e_k}\right)\Big/\log\left(\frac{k}{k-1}\right).$$

In Figure \ref{fig: converg_order}, we plot the relative error of the quantities $\|\omega_1\|_{L^\infty}$, $\|\omega\|_{L^\infty}$ and $E$ for $t\in\left[0,1.6\times10^{-3}\right]$, and the convergence order of the error in the late time $t\in\left[1\times10^{-3},1.6\times10^{-3}\right]$. We observe a numerical convergence with order slightly higher than $2$. The convergence order is quite stable in the time interval of our computation.

\begin{table}[hbt!]
\centering
\caption{Relative errors and convergence orders of $\omega_1$, $\psi_1$, $u^r$ and $u^z$ in sup-norm.}\label{tab: conv_rate}
\begin{tabular}{|c|c|c|c|c|}
\hline
\multirow{2}{*}{mesh size} & \multicolumn{4}{c|}{Sup-norm relative error at $t=1.6\times10^{-3}$} \\
\cline{2-5}
 & $\omega_1$ & order & $\psi_1$ & order \\
\hline
\hspace{0.3cm}$256\times256$\hspace{0.3cm} & \hspace{0.25cm}$2.545\times10^{-1}$\hspace{0.25cm} & \hspace{0.55cm}-\hspace{0.55cm} & \hspace{0.25cm}$5.912\times10^{-3}$\hspace{0.25cm} & \hspace{0.55cm}-\hspace{0.55cm} \\
\hline
$512\times512$ & $5.478\times10^{-2}$ & 2.216 & $1.168\times10^{-3}$ & 2.340 \\
\hline
$768\times768$ & $1.969\times10^{-2}$ & 2.524 & $4.136\times10^{-4}$ & 2.560 \\
\hline
$1024\times1024$ & $9.189\times10^{-3}$ & 2.655 & $1.926\times10^{-4}$ & 2.656 \\
\hline
$1280\times1280$ & $5.008\times10^{-3}$ & 2.720 & $1.050\times10^{-4}$ & 2.719 \\
\hline
\end{tabular}
\\
\vspace{0.3cm}
\begin{tabular}{|c|c|c|c|c|}
\hline
\multirow{2}{*}{mesh size} & \multicolumn{4}{c|}{Sup-norm relative error at $t=1.6\times10^{-3}$} \\
\cline{2-5}
 & $u^r$ & order & $u^z$ & order \\
\hline
\hspace{0.3cm}$256\times256$\hspace{0.3cm} & \hspace{0.25cm}$2.035\times10^{-2}$\hspace{0.25cm} & \hspace{0.55cm}-\hspace{0.55cm} & \hspace{0.25cm}$8.095\times10^{-3}$\hspace{0.25cm} & \hspace{0.55cm}-\hspace{0.55cm} \\
\hline
$512\times512$ & $3.954\times10^{-3}$ & 2.364 & $1.533\times10^{-3}$ & 2.310 \\
\hline
$768\times768$ & $1.405\times10^{-3}$ & 2.552 & $5.793\times10^{-4}$ & 2.556 \\
\hline
$1024\times1024$ & $6.540\times10^{-4}$ & 2.658 & $2.699\times10^{-4}$ & 2.655 \\
\hline
$1280\times1280$ & $3.594\times10^{-4}$ & 2.682 & $1.472\times10^{-4}$ & 2.719 \\
\hline
\end{tabular}
\end{table}

In Table \ref{tab: conv_rate}, we list the relative error and convergence order of the vector fields at $t=1.6\times10^{-3}$. The convergence order stays well above $2$, suggesting at least a second-order convergence for our numerical solver of the 3D axisymmetric Euler equations.

\subsection{Scaling analysis}
\label{sec: scaling-analysis}

In this section, we quantify the scaling property of the potential blow-up observed in our computation. This scaling analysis will give more supporting evidence that the potential blow-up satisfies the Beale-Kato-Majda blow-up criterion. It also uncovers more properties of the potential blow-up.

\begin{figure}[hbt!]
\centering
\includegraphics[width=.4\textwidth]{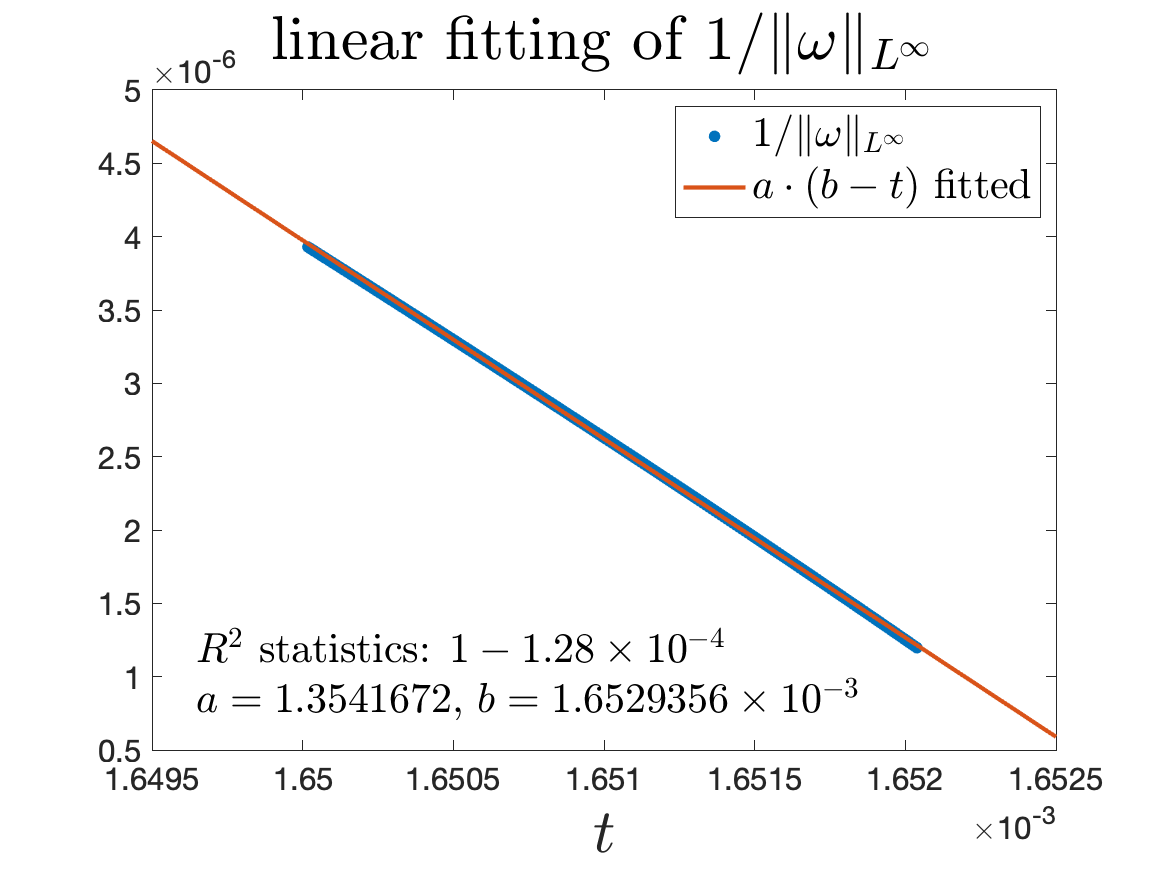}\includegraphics[width=.4\textwidth]{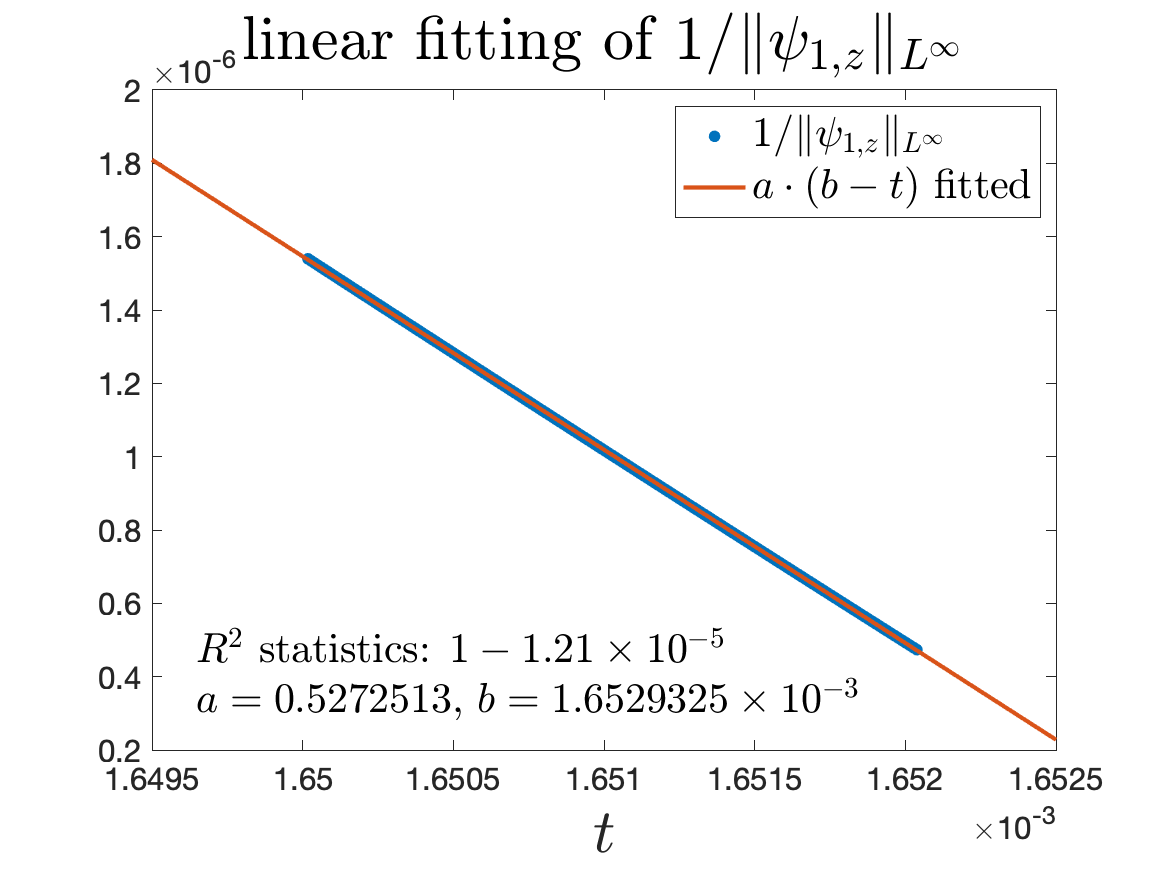}
\caption{Linear fitting of $1/\|\omega\|_{L^\infty}$ and $1/\|\psi_{1,z}\|_{L^\infty}$ with time.}\label{fig: fit_w_dpsi1dz}
\end{figure}

As discussed in \eqref{eq: self_similar} of Section \ref{sec: self_similar_solution}, if there is a self-similar blow-up, the scaling invariant property of the 3D Euler equations will ensure that $\|\omega\|_{L^\infty}\sim1/(T-t)$ and $\|\psi_{1,z}\|_{L^\infty}\sim1/(T-t)$. Therefore, we examine this property by regressing $\|\omega\|_{L^\infty}^{-1}$ and $\|\psi_{1,z}\|_{L^\infty}^{-1}$ against $t$, respectively. More specifically, for a quantity $v$, which is either $\|\omega\|_{L^\infty}^{-1}$ or $\|\psi_{1,z}\|_{L^\infty}^{-1}$, we perform the least square fitting of the model
$$v\sim a\cdot(b-t),$$
in searching for constants $a$ and $b$, where $a$ is the negated slope of the fitted line, and $b$ can be considered as the estimate time of the blow-up. In Figure \ref{fig: fit_w_dpsi1dz}, we visualize the data points and the fitted line using data between $t=1.6500174\times10^{-3}$ and $t=1.6520384\times10^{-3}$. The $R^2$ of the fitting between $\|\omega\|_{L^\infty}^{-1}$ and $t$ is $1-1.28\times10^{-4}$, and the $R^2$ of the fitting between $\|\psi_{1,z}\|_{L^\infty}^{-1}$ and $t$ is $1-1.21\times10^{-5}$. Such high $R^2$ values show strong linear relation between $\|\omega\|_{L^\infty}^{-1}$, $\|\psi_{1,z}\|_{L^\infty}^{-1}$ and $t$. Moreover, the fittings of the two quantities estimate the blow-up time to be $b=1.6529356\times10^{-3}$ and $b=1.6529325\times10^{-3}$ respectively. These two blow-up times agree with each other up to $6$ digits. Therefore, Figure \ref{fig: fit_w_dpsi1dz} provides further evidence that the 3D Euler equations develop a potential finite-time singularity.

We next move to fit the scaling factors $c_l$ and $c_\omega$ used in the self-similar ansatz \eqref{eq: self_similar} of the solutions. Since the functions $\Omega$ and $\Psi$ are time-independent in \eqref{eq: self_similar}, we should have that
$$Z_1\sim (T-t)^{c_l},\quad \|\omega_1\|_{L^\infty}^{-1}\sim(T-t)^{c_\omega},$$
where we recall that $Z_1=Z_1(t)$ is the $z$-coordinate of the maximum location of $-\omega_1$. Due to the unknown powers $c_l$ and $c_\omega$, the direct fitting of the above model is nonlinear. Therefore, we turn to a searching algorithm for the power variable. Specifically, for a quantity $v$, that is either $Z_1$ or $\|\omega_1\|_{L^\infty}^{-1}$, we search for a power $c$ such that the linear regression of
$$v^{1/c}\sim a\cdot(b-t),$$
has the largest $R^2$ value up to some error tolerance. We will start with a guessed window of the power $c$, and then exhaust the value of $c$ within the window up to some error tolerance, and choose $c$ with the largest $R^2$ value. If the optimal $c$ we searched falls on the boundary of the current window, we then adaptively adjust the window size and location, and repeat the above procedure. When the optimal searched $c$ falls within the interior of the window, we stop the searching.

\begin{figure}[hbt!]
\centering
\includegraphics[width=.4\textwidth]{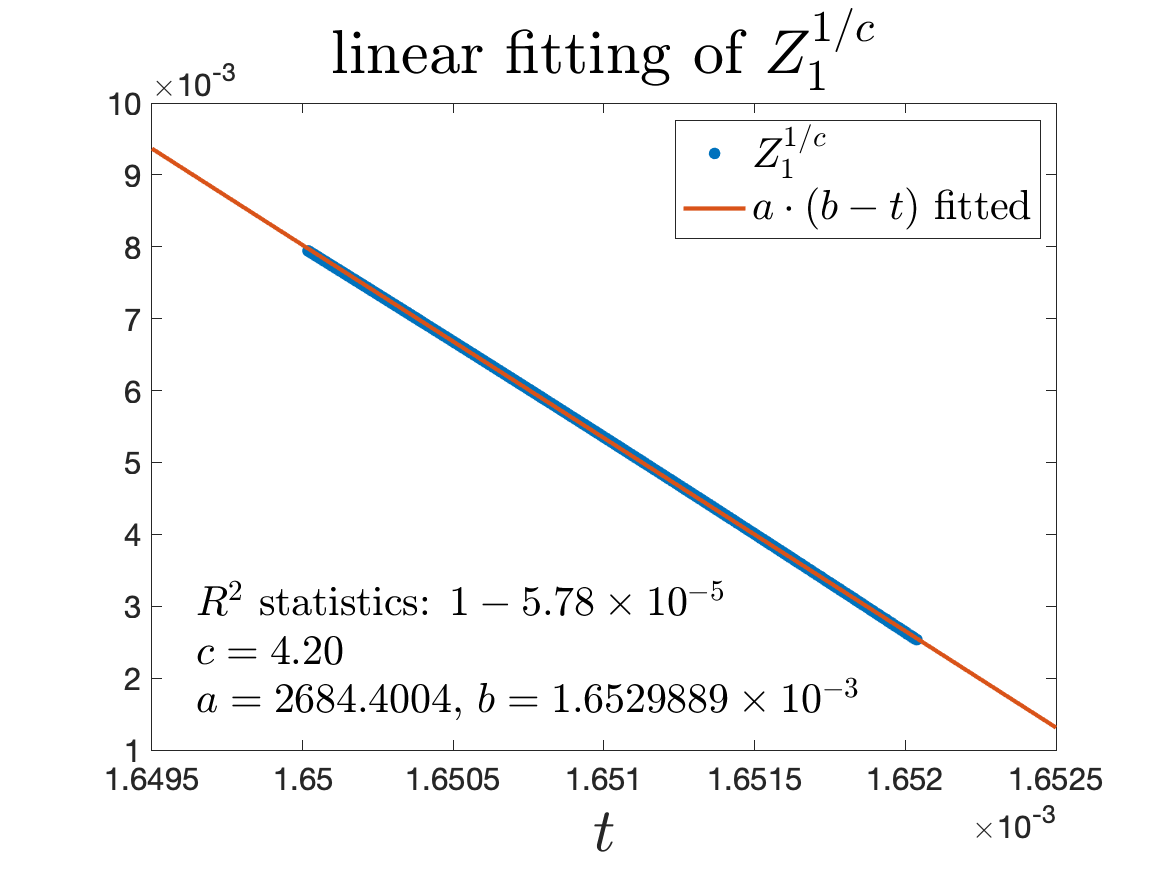}\includegraphics[width=.4\textwidth]{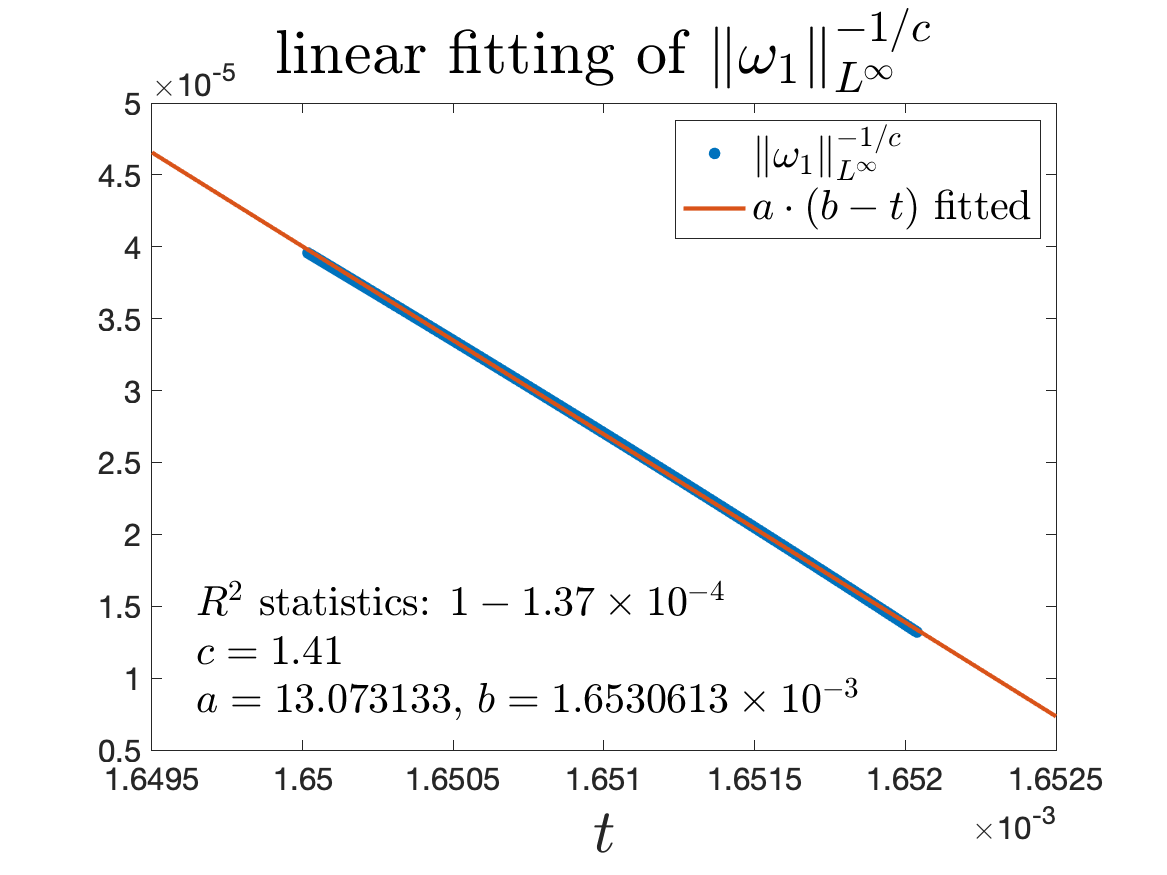}
\caption{Linear fitting of $Z_1^{1/c}$ and $\|\omega_1\|_{L^\infty}^{-1/c}$ with time.}\label{fig: fit_Z_maxw1}
\end{figure}

In Figure \ref{fig: fit_Z_maxw1}, we demonstrate the result of the searching. We can see that with the chosen $c$, the linear regression achieves a very high $R^2$ value, suggesting a strong linear relation. The relative error between the estimated blow-up time and the previous estimate smaller than  $7.8\times10^{-5}$. Moreover, the searching suggests that $c_l\approx 4.20$ and $c_\omega\approx1.41$, and these estimated values of $c_l$ and $c_\omega$ satisfy the scaling relation $c_\omega=1+\alpha c_l$  in \eqref{eq: scaling_relation} approximately.

It is worth emphasizing that the estimated $c_l$ is well above $1$, and this explains the convex curve of $Z_1(t)$ as observed in Figure \ref{fig: stats_curve_2} in Section \ref{sec: evidence_singular}. 

We remark that we did not perform the searching algorithm with $\|\psi_1\|_{L^\infty}$ to find out the scaling factor $c_\psi$, so that we could check the other scaling relation $c_\psi=1-c_l$ in \eqref{eq: scaling_relation}. This is because $\|\psi_1\|_{L^\infty}$ is mainly affected by the far field behavior of $\psi_1$, as shown in Figure \ref{fig: end_data}. However, the self-similar ansatz \eqref{eq: self_similar} is only valid in the near field, so such fitting is meaningless. In fact, the good fitting between $\|\psi_{1,z}\|_{L^\infty}^{-1}$ and $t$ already implies that $c_\psi=1-c_l$, because the self-similar ansatz suggests that $\|\psi_{1,z}\|_{L^\infty}^{-1}\sim(T-t)^{c_\psi+c_l}$.

\begin{table}[hbt!]
\centering
\caption{Fitting results of $\|\omega\|_{L^\infty}^{-1}$, $\|\psi_{1,z}\|_{L^\infty}^{-1}$, $Z_1$ and $\|\omega_1\|_{L^\infty}^{-1}$ at different mesh sizes.}\label{tab: fit_mesh}
\begin{tabular}{|c|c|c|c|c|}
\hline
\multirow{2}{*}{mesh size} & \multicolumn{2}{c|}{$1/\|\omega\|_{L^\infty}$} & \multicolumn{2}{c|}{$1/\|\psi_{1,z}\|_{L^\infty}$} \\
\cline{2-5}
 & $10^3\times b$ & $R^2$ & $10^3\times b$ & $R^2$ \\
\hline
\hspace{0.05cm}$1024\times1024$\hspace{0.05cm} &  \hspace{0.05cm}1.6529356\hspace{0.05cm} & \hspace{0.05cm}0.99987\hspace{0.05cm} & \hspace{0.05cm}1.6529325\hspace{0.05cm} & \hspace{0.05cm}0.99999\hspace{0.05cm} \\
\hline
$1280\times1280$ & 1.6527953 & 1.00000 & 1.6528189 & 1.00000 \\
\hline
$1536\times1536$ & 1.6525824 & 1.00000 & 1.6527396 & 1.00000 \\
\hline
\end{tabular}
\\
\vspace{0.3cm}
\begin{tabular}{|c|c|c|c|c|c|c|}
\hline
\multirow{2}{*}{mesh size} & \multicolumn{3}{c|}{$Z_1$} & \multicolumn{3}{c|}{$1/\|\omega_1\|_{L^\infty}$} \\
\cline{2-7}
 & $c$ & $10^3\times b$ & $R^2$ & $c$ & $10^3\times b$ & $R^2$ \\
\hline
\hspace{0.05cm}$1024\times1024$\hspace{0.05cm} & \hspace{0.05cm}4.20\hspace{0.05cm} & \hspace{0.05cm}1.6529889\hspace{0.05cm} & \hspace{0.05cm}0.99994\hspace{0.05cm} & \hspace{0.05cm}1.41\hspace{0.05cm} & \hspace{0.05cm}1.6530613\hspace{0.05cm} & \hspace{0.05cm}0.99986\hspace{0.05cm} \\
\hline
$1280\times1280$ & 4.21 & 1.6527877 & 0.99999 & 1.42 & 1.6527894 & 1.00000 \\
\hline
$1536\times1536$ & 4.25 & 1.6526864 & 1.00000 & 1.41 & 1.6526953 & 1.00000 \\
\hline
\end{tabular}
\end{table}

Finally, we perform the above fitting of different quantities using different spatial resolutions, and summarize the results in Table \ref{tab: fit_mesh}. We can see that the fitting has excellent quality at all spatial resolutions, and the fitted parameters are consistent across different spatial resolutions.

\section{The dynamic rescaling formulation}
\label{sec: dynamic_rescaling}

In order to better study the potential self-similar singularity as we have observed in Section \ref{sec: evidence_self_similar}, we add extra scaling terms to \eqref{eq: vort_stream_1_3d_noswirl} and write
\begin{subequations}
\label{eq: vort_stream_1_3d_noswirl_dr}
\begin{align}
    \tilde{\omega}_{1,\tau}+\left(\tilde{c}_l\xi+\tilde{u}^\xi\right)\tilde{\omega}_{1,\xi}+\left(\tilde{c}_l\zeta+\tilde{u}^\zeta\right)\tilde{\omega}_{1,\zeta}&=\left(c_\omega-(1-\alpha)\tilde{\psi}_{1,\zeta}\right)\tilde{\omega}_1, \label{eq: vort_1_3d_noswirl_dr}\\
    -\tilde{\psi}_{1,\xi\xi}-\tilde{\psi}_{1,\zeta\zeta}-\frac{3}{\xi}\tilde{\psi}_{1,\xi}&=\tilde{\omega}_1\xi^{\alpha-1}, \label{eq: stream_1_3d_noswirl_dr}\\
    \tilde{u}^\xi=-\xi\tilde{\psi}_{1,\zeta},\quad\tilde{u}^\zeta&=2\tilde{\psi}_1+\xi\tilde{\psi}_{1,\xi}, \label{eq: velo_rz_1_3d_noswirl_dr}
\end{align}
\end{subequations}
where $\tilde{c}_l=\tilde{c}_l(\tau)$, $\tilde{c}_\omega=\tilde{c}_\omega(\tau)$ are scalar functions of $\tau$. In \eqref{eq: vort_1_3d_noswirl_dr}, the terms $\tilde{c}_l\xi\partial_\xi$ and $\tilde{c}_l\zeta\partial_\zeta$ stretch the solutions in space to maintain a finite support of the self-similar blow-up solution. The term $\tilde{c}_\omega\tilde{\omega}_1$ acts as a damping term to ensure that the magnitude of $\tilde{\omega}_1$ remains finite. The combined effect of these terms dynamically rescales the solution to capture the potential self-similar profile. Such dynamic rescaling strategy has widely been used in the study of singularity formation of nonlinear Schr\"{o}dinger equations as in \cite{mclaughlin1986focusing, landman1988rate, lemesurier1988focusing, landman1992stability, papanicolaou1994focusing}. And recently it has been used to study singularity formation of the 3D Euler equations as in \cite{hou2018potential, chen2021finite, chen2021finite2}.

If we define
\begin{align}
\label{eq: dr_relation}
    \tilde{c}_\psi(\tau)=\tilde{c}_\omega(\tau)+(1+\alpha)\tilde{c}_l(\tau),
\end{align}
we can check that \eqref{eq: vort_stream_1_3d_noswirl_dr} admits the following solution
\begin{equation}
\begin{aligned}
\label{eq: solu_dr}
    \tilde{\omega}_1(\xi,\zeta,\tau)&=\tilde{C}_\omega(\tau)\omega_1\left(\tilde{C}_l(\tau)\xi,\tilde{C}_l(\tau)\zeta,t(\tau)\right),\\
    \tilde{\psi}_1(\xi,\zeta,\tau)&=\tilde{C}_\psi(\tau)\psi_1\left(\tilde{C}_l(\tau)\xi,\tilde{C}_l(\tau)\zeta,t(\tau)\right),
\end{aligned}
\end{equation}
where $(\omega_1, \psi_1)$ is the solution to \eqref{eq: vort_stream_1_3d_noswirl}, and
\begin{align*}
    \tilde{C}_\omega(\tau)&=\exp\left(\int_0^\tau \tilde{c}_\omega(s)\mathrm{d}s\right),\quad \tilde{C}_\psi(\tau)=\exp\left(\int_0^\tau \tilde{c}_\psi(s)\mathrm{d}s\right),\\
    \tilde{C}_l(\tau)&=\exp\left(-\int_0^\tau \tilde{c}_l(s)\mathrm{d}s\right),\quad t^\prime(\tau)=\tilde{C}_\psi(\tau)\tilde{C}_l(\tau)=\tilde{C}_\omega(\tau)\tilde{C}_l^{-\alpha}(\tau).
\end{align*}

The new equations \eqref{eq: vort_stream_1_3d_noswirl_dr} leave us with two degrees of freedom: we are free to choose $\left\{\tilde{c}_l(\tau), \tilde{c}_\omega(\tau)\right\}$. This allows us to impose the following normalization conditions
\begin{align}
\label{eq: norm_cond}
    \tilde{\omega}_{1}(0,1,\tau)=-1,\quad\tilde{\omega}_{1,\zeta}(0,1,\tau)=0,\quad\text{for }\tau\geq0.
\end{align}
These two conditions will help fix the maximum value of $-\tilde{\omega}_1$ at $1$ and the maximum location at $(\xi, \zeta)=(0, 1)$.

One way to enforce the normalization conditions, as used in many literatures like \cite{hou2018potential, liu2017spatial}, is to first enforce them at $\tau=0$ using the scaling invariant relation \eqref{eq: vort_stream_invariant}, and then enforce their time derivatives to be zero
\begin{align}
\label{eq: norm_cond_dr}
    \frac{\partial}{\partial\tau}\tilde{\omega}_{1}(0,1,\tau)=0,\quad\frac{\partial}{\partial\tau}\tilde{\omega}_{1,\zeta}(0,1,\tau)=0,\quad\text{for }\tau\geq0.
\end{align}
Using \eqref{eq: vort_1_3d_noswirl_dr}, the above conditions are equivalent to
\begin{equation}
\label{eq: norm_cond_close}
\begin{aligned}
\tilde{c}_l(\tau)=&-2\tilde{\psi}_1(0,1,\tau)-(1-\alpha)\tilde{\psi}_{1,\zeta\zeta}(0,1,\tau)\frac{\tilde{\omega}_1(0,1,\tau)}{\tilde{\omega}_{1,\zeta\zeta}(0,1,\tau)},\\
\tilde{c}_\omega(\tau)=&(1-\alpha)\tilde{\psi}_{1,\zeta}(0,1,\tau).
\end{aligned}
\end{equation}
However, it is hard to evaluate \eqref{eq: norm_cond_close} accurately, because it requires calculating second-order derivatives. More importantly, due to the complicated nonlinear nature of \eqref{eq: vort_1_3d_noswirl_dr}, even if \eqref{eq: norm_cond_close} can be accurately evaluated, the temporal discretization (Runge-Kutta method) makes it difficult to enforce \eqref{eq: norm_cond} exactly for the next time step. As a result, imposing \eqref{eq: norm_cond_dr} is not as helpful to preserve the normalization conditions \eqref{eq: norm_cond} in the following time steps. The maximum magnitude and location will gradually change in time, which makes it difficult to compute the self-similar profile numerically.

\subsection{The operator splitting strategy}

To enforce the normalization conditions \eqref{eq: norm_cond} accurately at every time step, we utilize the operator splitting method and rewrite \eqref{eq: vort_1_3d_noswirl_dr} as
\begin{align}
\label{eq: os_ode}
    \tilde{\omega}_{1,\tau}=F(\tilde{\omega}_{1})+G(\tilde{\omega}_{1}),
\end{align}
where $F(\tilde{\omega}_{1})=-\tilde{u}^\xi\tilde{\omega}_{1,\xi}-\tilde{u}^\zeta\tilde{\omega}_{1,\zeta}-(1-\alpha)\tilde{\psi}_{1,\zeta}\tilde{\omega}_1$ contains the original terms in \eqref{eq: vort_1_3d_noswirl}, and
$G(\tilde{\omega}_{1})=-\tilde{c}_l\xi\tilde{\omega}_{1,\xi}-\tilde{c}_l\zeta\tilde{\omega}_{1,\zeta}+\tilde{c}_\omega\tilde{\omega}_1$ is the linear part that controls the rescaling. Here we view $\tilde{\psi}_{1}$ as a function of $\tilde{\omega}_1$ through the Poisson equation \eqref{eq: stream_1_3d_noswirl_dr}. The operator splitting method allows us to solve \eqref{eq: vort_1_3d_noswirl_dr} by solving $\tilde{\omega}_{1,\tau}=F(\tilde{\omega}_{1})$ and $\tilde{\omega}_{1,\tau}=G(\tilde{\omega}_{1})$ alternatively.

We can use the standard Runge-Kutta method to solve $\tilde{\omega}_{1,\tau}=F(\tilde{\omega}_{1})$. As for $\tilde{\omega}_{1,\tau}=G(\tilde{\omega}_{1})$, we notice that there is a closed form solution for the initial value problem
\begin{align}
\label{eq: op_spl_close}
    \tilde{\omega}_{1}(\xi,\zeta,\tau)=\tilde{C}_\omega(\tau)\tilde{\omega}_1\left(\tilde{C}_l(\tau)\xi,\tilde{C}_l(\tau)\zeta, 0\right),
\end{align}
where $\tilde{C}_\omega(\tau)=\exp\left(\int_0^\tau \tilde{c}_\omega(s)\mathrm{d}s\right)$ and $\tilde{C}_l(\tau)=\exp\left(-\int_0^\tau \tilde{c}_l(s)\mathrm{d}s\right)$.

In the first step, solving $\tilde{\omega}_{1,\tau}=F(\tilde{\omega}_{1})$ will violate the normalization conditions \eqref{eq: norm_cond}. But we will correct this error in the second step by solving $\tilde{\omega}_{1,\tau}=G(\tilde{\omega}_{1})$ with a smart choice of $\tilde{C}_l$ and $\tilde{C}_\omega$ in \eqref{eq: op_spl_close}. In other words, at every time step when we solve $\tilde{\omega}_{1,\tau}=G(\tilde{\omega}_{1})$, we can exactly enforce \eqref{eq: norm_cond} by properly choosing $\tilde{C}_l$ and $\tilde{C}_\omega$ in \eqref{eq: op_spl_close}. We could also adopt Strang's splitting \cite{strang1968construction} for better temporal accuracy.

\subsection{Numerical settings}
\label{sec: dynamic_rescaling_numeric}

\begin{figure}[hbt!]
\centering
\includegraphics[width=.4\textwidth]{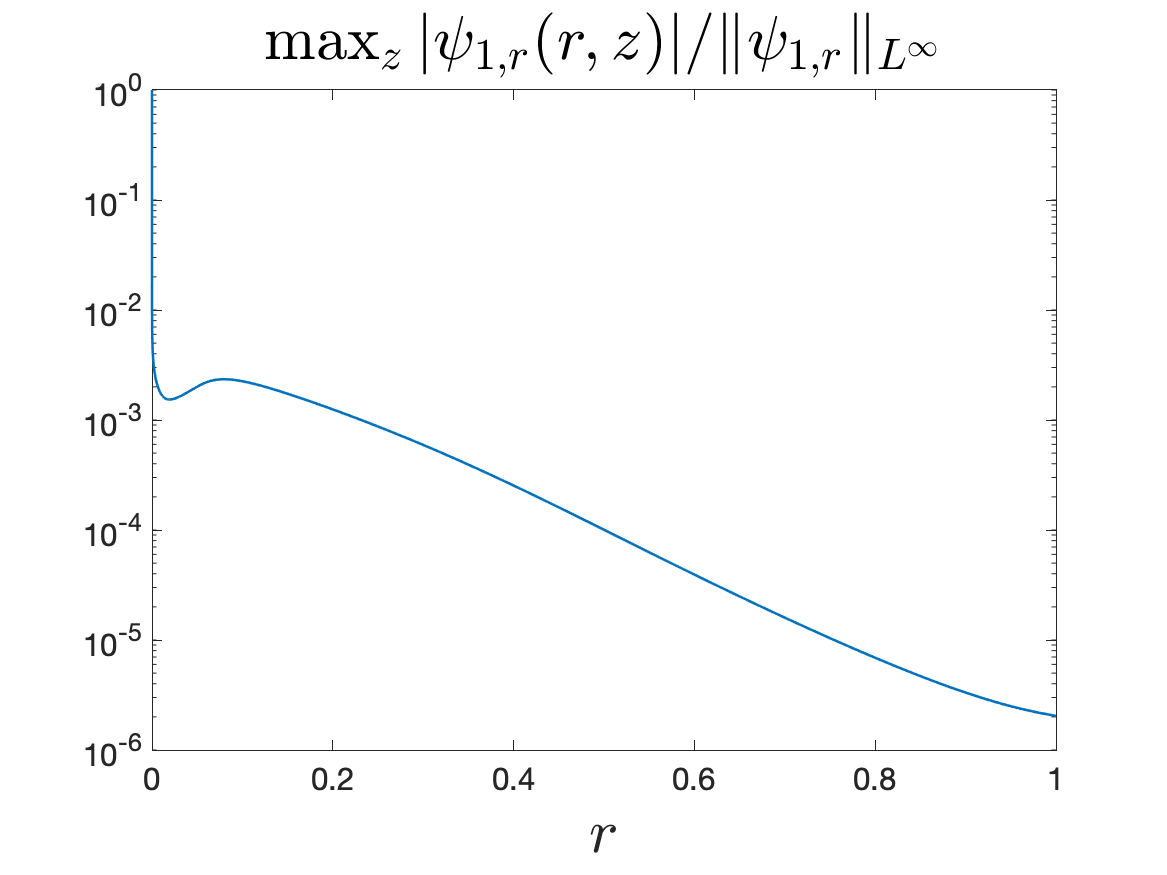}\includegraphics[width=.4\textwidth]{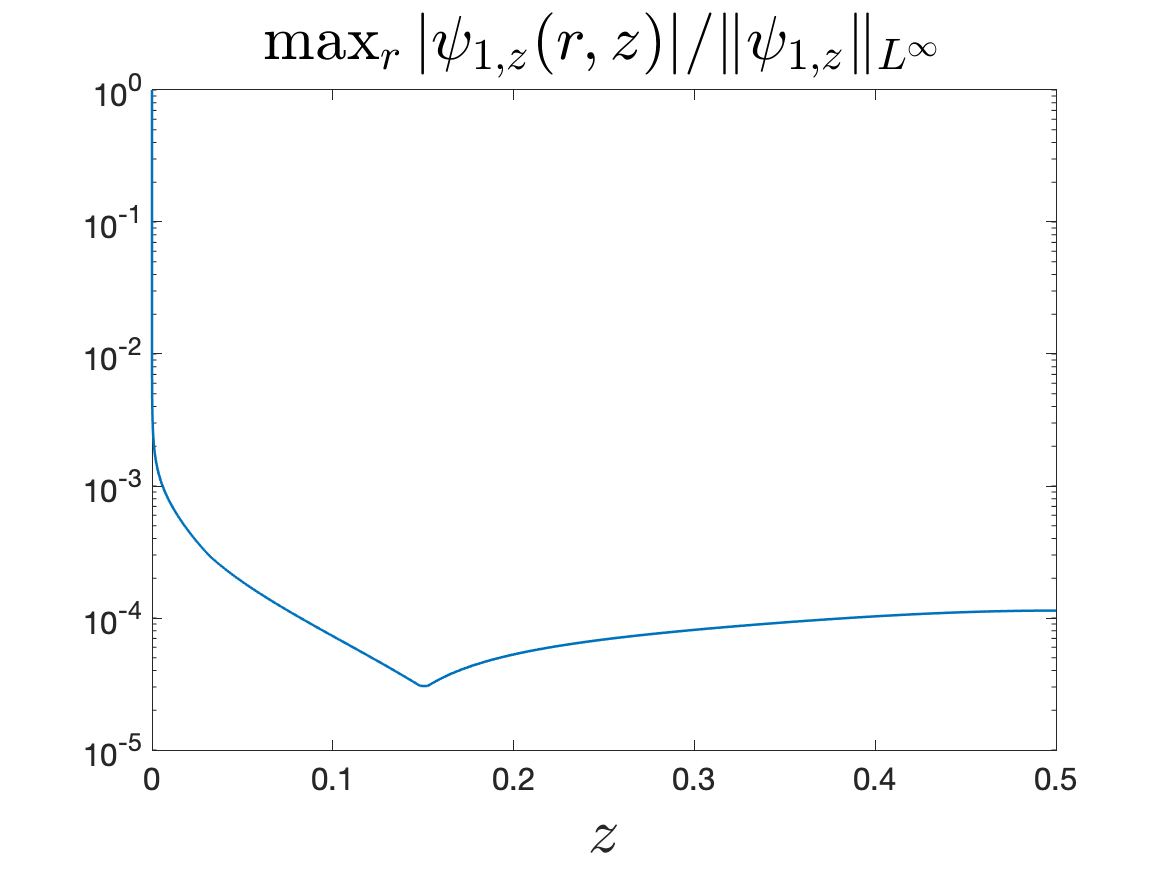}
\caption{Decay of the derivatives of $\psi_1$.}\label{fig: dr_decay}
\end{figure}

Now we numerically solve the dynamic rescaling formulation \eqref{eq: vort_stream_1_3d_noswirl_dr}. For the initial condition, we use the solution obtained from the final iteration of the adaptive mesh method in Section \ref{sec: evidence_singular}, and use the relation \eqref{eq: vort_stream_invariant} to enforce the normalization conditions \eqref{eq: norm_cond_dr}. Now that the maximum location of $\tilde{\omega}_1$ is pinned at $(\xi,\zeta)=(0,1)$, we focus on a large computational domain
$$\mathcal{D}^\prime=\left\{(\xi,\zeta): 0\leq\xi\leq D, 0\leq\zeta\leq D/2\right\}\,,$$
with domain size $D=1\times10^5$. The reason for such large domain size is because under the normalization of the dynamic rescaling formulation, the original computational domain $(r, z)\in\left[0, 1\right]\times\left[0, 1/2\right]$ is now equivalent to $(\xi, \zeta)\in\left[0, 1/Z_1\right]\times\left[0, 1/2/Z_1\right]$. This domain grows large quickly as $Z_1$ tends to zero. While we adopt the boundary conditions at $\xi=0$ and $\zeta=0$ of \eqref{eq: vort_stream_1_3d_noswirl} in Section \ref{sec: boundary_symmetry}, we need to find a good far field boundary conditions for $\tilde{\psi}_1$. Due to extra stretching terms, the far field boundary for $\tilde{\psi}_1$ will no longer correspond to the far field boundary for $\psi_1$, namely $r=1$ and $z=1/2$. However, we notice that $\psi_{1,r}$ decays rapidly with respect to $r$, and $\psi_{1,z}$ decays rapidly with respect to $z$. For example, Figure \ref{fig: dr_decay} shows the decay of $\psi_{1,r}$ as $r\rightarrow1$ and the decay of $\psi_{1,z}$ as $z\rightarrow1/2$ for the solution to \eqref{eq: vort_stream_1_3d_noswirl} at $t=1.6524635\times10^{-3}$. Therefore, it is reasonable to impose the zero Neumann boundary condition at the far field boundaries of $\mathcal{D}^\prime$: $\xi=D$ and $\zeta=D/2$. Due to the size of the computation domain $\mathcal{D}^\prime$ and the presence of the vortex stretching terms, the error introduced by this boundary condition will have little influence on the near field around $(\xi,\zeta)=(0,1)$. We will discuss in details about the influence of the domain size $D$ in Section \ref{sec: domain-size} and in the end of Section \ref{sec: holder_exponent}.

We remark that we still need the adaptive mesh in the $r$- and $z$-directions, because we not only need to cover a very large field, but also need to focus around $(\xi,\zeta)=(0,1)$. The adaptive mesh that we use to solve the dynamic rescaling formulation will not change during the computation, since the dynamically rescaled vorticity has its maximum location fixed at $(\xi,\zeta)=(0,1)$ for all times instead of traveling toward the origin.

\subsection{Convergence to the steady state}

\begin{figure}[hbt!]
\centering
\includegraphics[width=.4\textwidth]{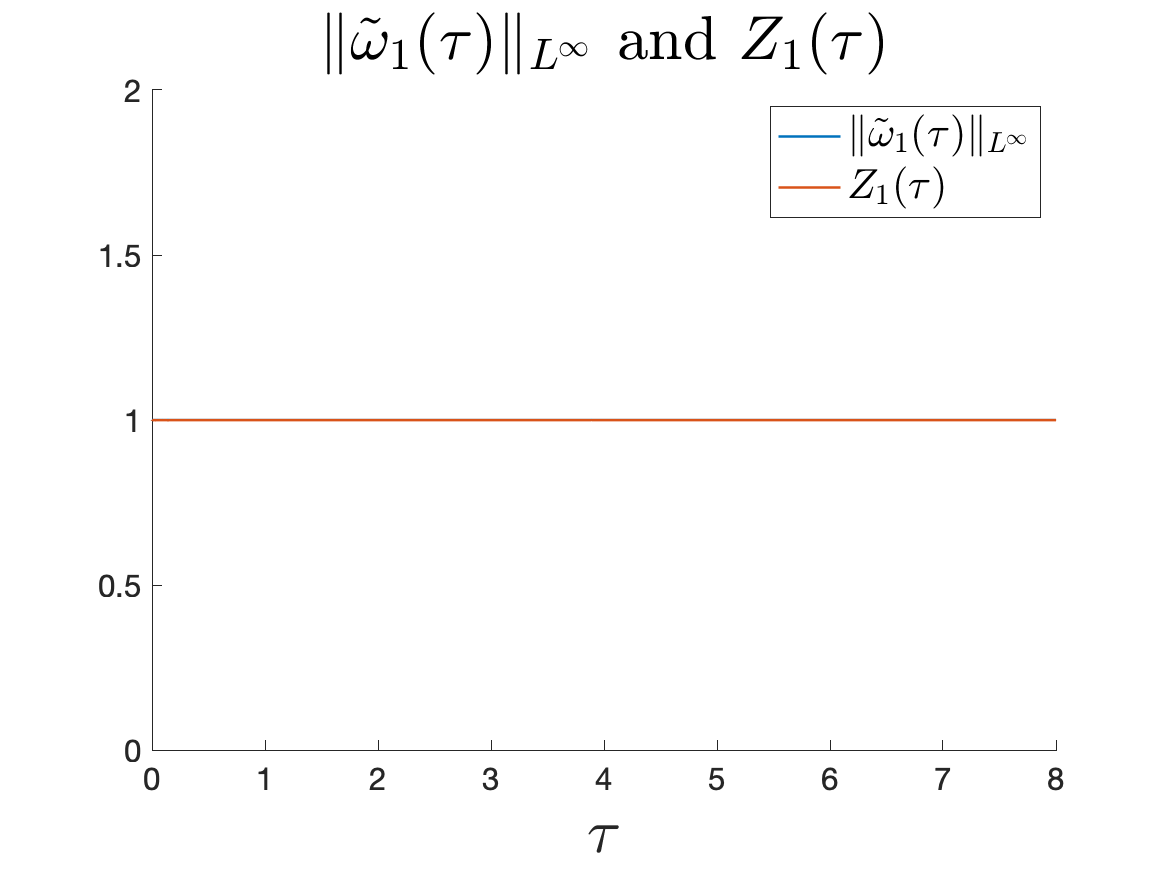}\includegraphics[width=.4\textwidth]{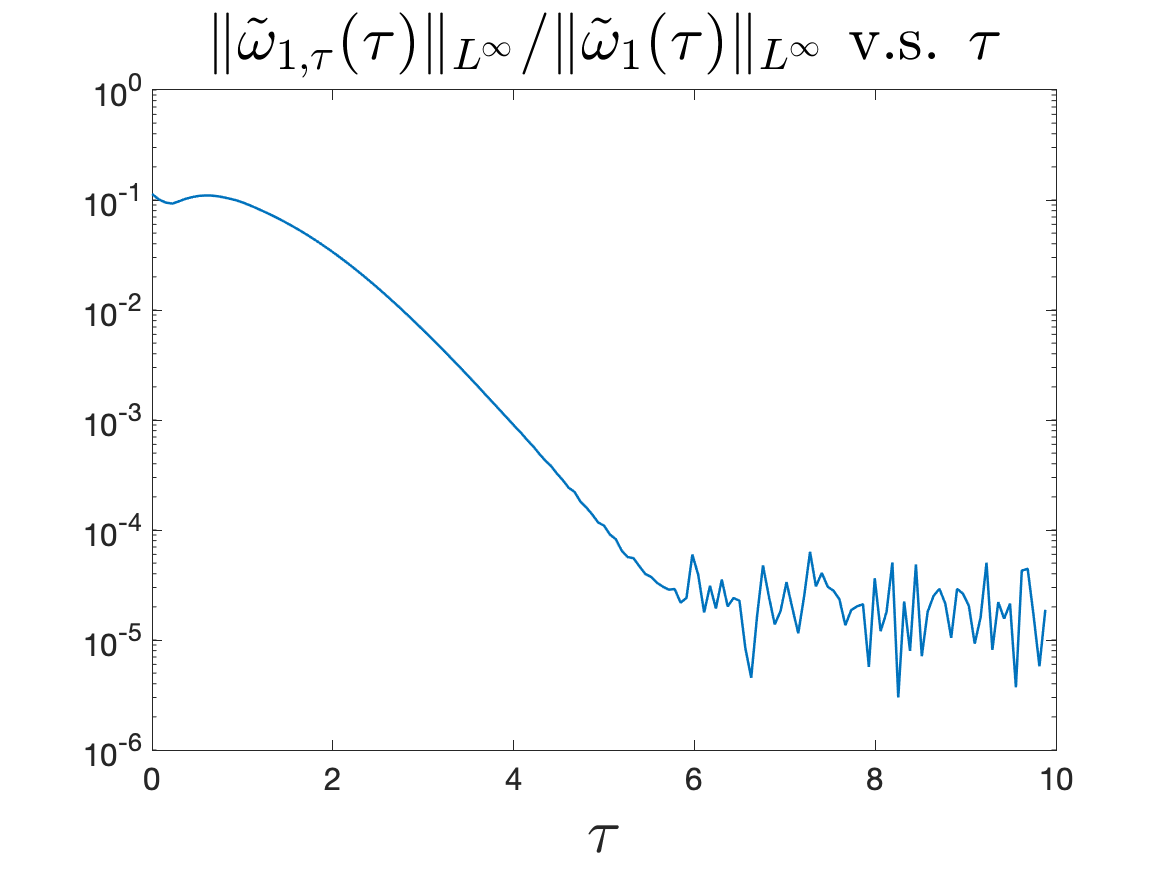}\\
\caption{Left: curves of the normalized quantities $\|\tilde{\omega}_1(\tau)\|_{L^\infty}$ and $Z_1(\tau)$. Right: Curve of the relative time derivative strength $\|\tilde{\omega}_{1,\tau}(\tau)\|_{L^\infty}/\|\tilde{\omega}_1(\tau)\|_{L^\infty}$.}\label{fig: dr_norm_tau_and_dr_time_deri}
\end{figure}

We solve \eqref{eq: vort_stream_1_3d_noswirl_dr} until it converges to a steady state. In the left subplot of Figure \ref{fig: dr_norm_tau_and_dr_time_deri}, we monitor how the normalization conditions \eqref{eq: norm_cond_dr} are enforced. The two normalized quantities, $\|\tilde{\omega}_1(\tau)\|_{L^\infty}$ and $Z_1(\tau)$,  are essentially fixed at $1$, and in fact, they deviate from $1$ by less than $5.14\times10^{-4}$. In the right subplot of Figure \ref{fig: dr_norm_tau_and_dr_time_deri}, We view the system \eqref{eq: vort_stream_1_3d_noswirl_dr} as an ODE of $\tilde{\omega}_1$ as in \eqref{eq: os_ode}, and plot the relative strength of the time derivative
$$\|\tilde{\omega}_{1,\tau}\|_{L^\infty}/\|\tilde{\omega}_1\|_{L^\infty}=\|F(\tilde{\omega}_1)+G(\tilde{\omega}_1)\|_{L^\infty}/\|\tilde{\omega}_1\|_{L^\infty},$$
as a function of time $\tau$. This relative strength of the time derivative has a decreasing trend and drops below $8.18\times10^{-6}$ near the end of the computation, which implies that we are very close to the steady state.

When the solution of \eqref{eq: vort_stream_1_3d_noswirl_dr} converges to a steady state, $\tilde{\omega}_1$ and $\tilde{\psi}_1$ are independent of the time $\tau$. Therefore, we should have the following relation from \eqref{eq: solu_dr}
\begin{align*}
    \omega_1(r, z, t)&\sim\tilde{C}_\omega^{-1}(\tau(t))\tilde{\omega}_1\left(\tilde{C}_l^{-1}(\tau(t))r, \tilde{C}_l^{-1}(\tau(t))z\right),\\
    \psi_1(r, z, t)&\sim\tilde{C}_\psi^{-1}(\tau(t))\tilde{\psi}_1\left(\tilde{C}_l^{-1}(\tau(t))r, \tilde{C}_l^{-1}(\tau(t))z\right),
\end{align*}
where $\tau=\tau(t)$ is the rescaled time variable. Comparing the above relation with the ansatz stated in \eqref{eq: self_similar}, we conclude that
\begin{align}
\label{eq: scale_conversion}
    c_l=-\frac{\tilde{c}_l}{\tilde{c}_\omega+\alpha\tilde{c}_l},\quad c_\omega=\frac{\tilde{c}_\omega}{\tilde{c}_\omega+\alpha\tilde{c}_l},\quad c_\psi=\frac{\tilde{c}_\psi}{\tilde{c}_\omega+\alpha\tilde{c}_l}.
\end{align}
We remark that assuming \eqref{eq: dr_relation}, the above relation naturally guarantees that the scaling relation \eqref{eq: scaling_relation} holds true.

\begin{figure}[hbt!]
\centering
\includegraphics[width=.4\textwidth]{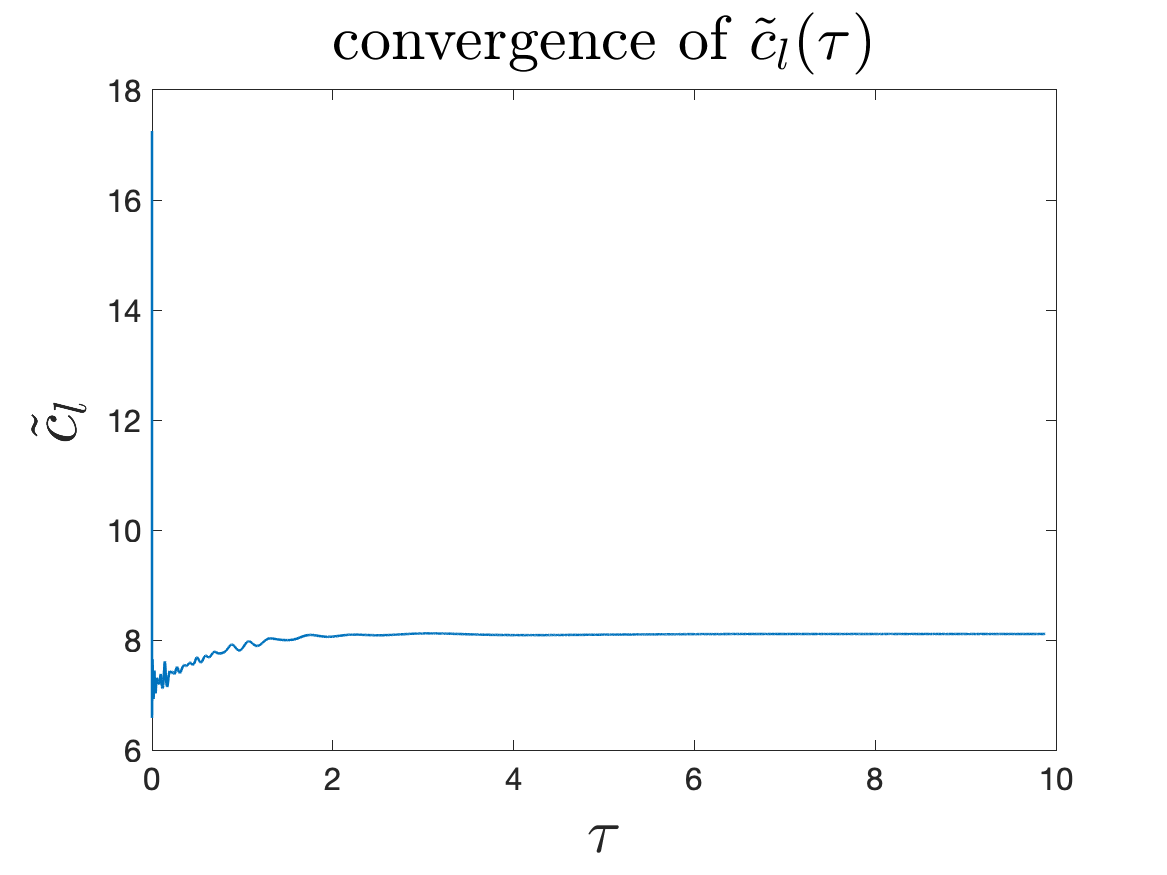}\includegraphics[width=.4\textwidth]{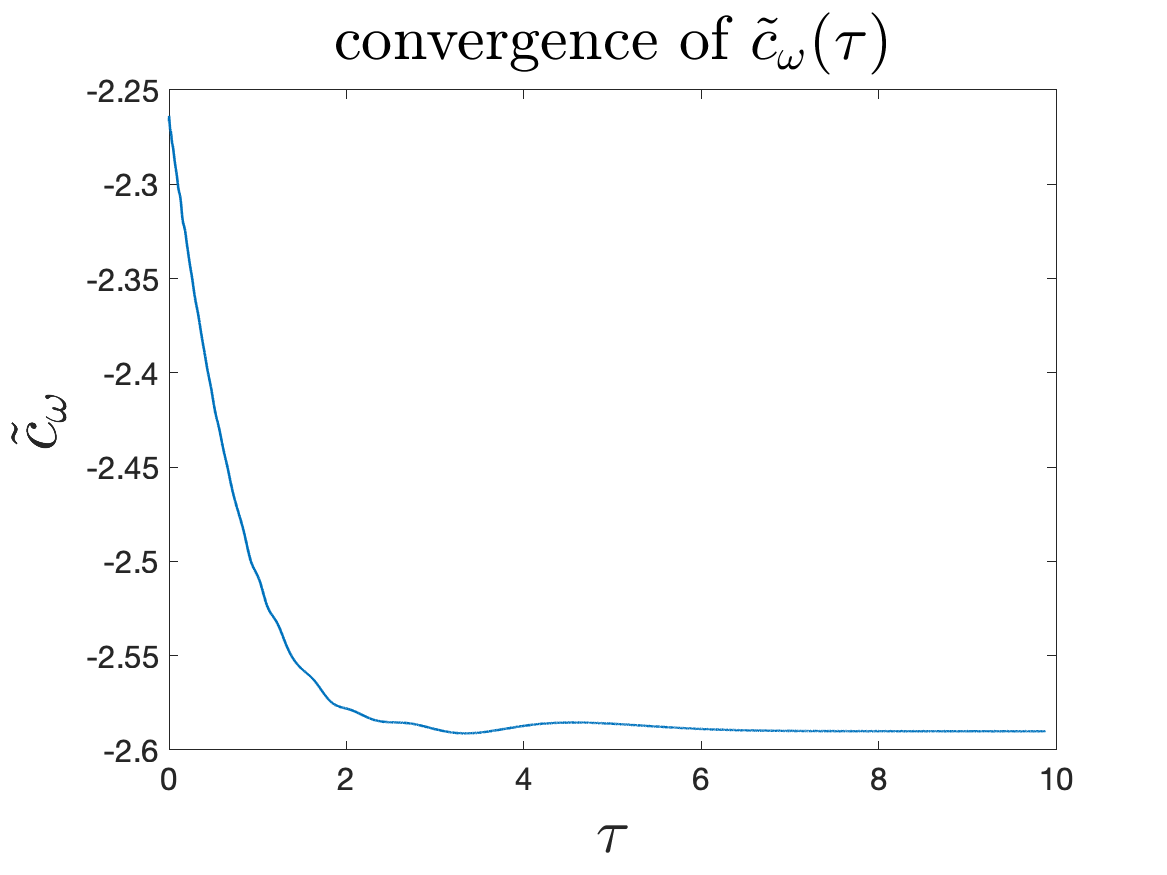}\\
\includegraphics[width=.4\textwidth]{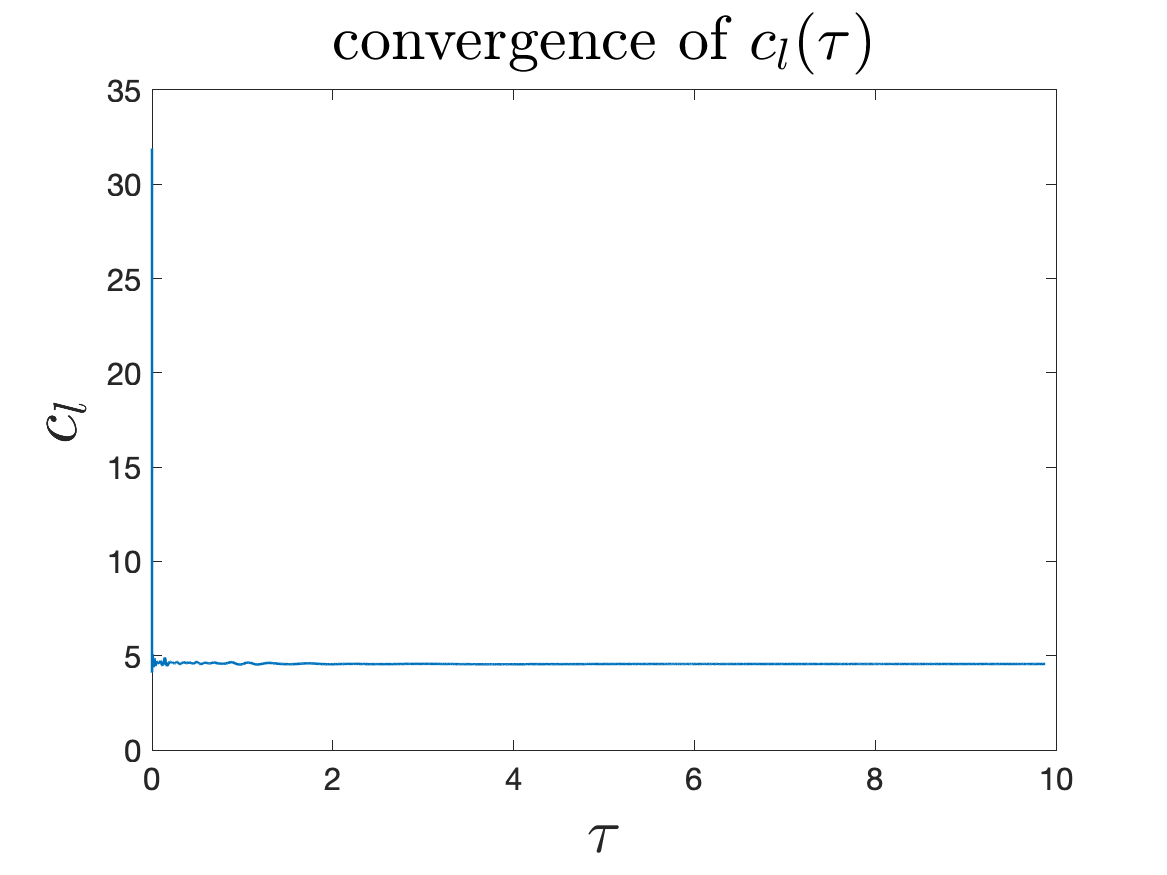}\includegraphics[width=.4\textwidth]{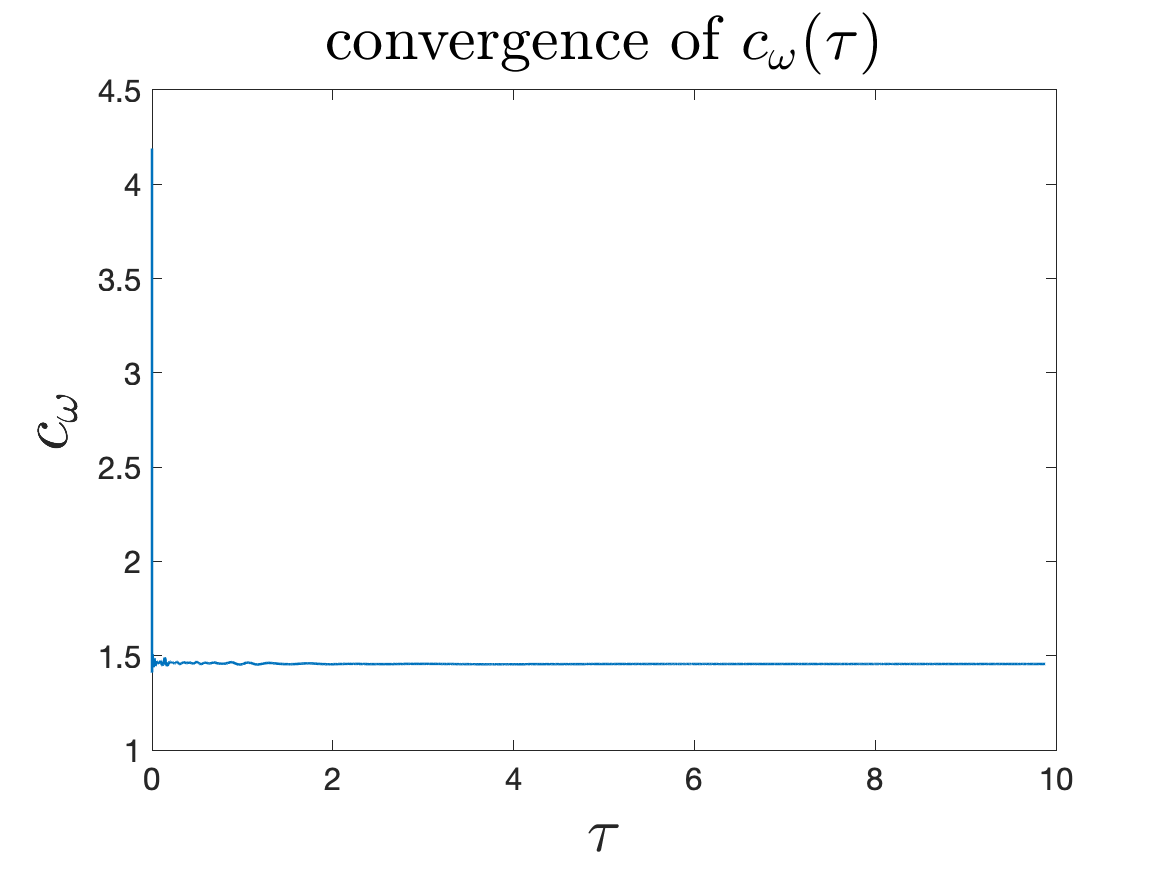}
\caption{Convergence curves of the scaling factors using dynamic rescaling method. Top row: $\tilde{c}_l$ and $\tilde{c}_\omega$. Bottom row: $c_l$ and $c_\omega$.}\label{fig: dr_scaling_factor}
\end{figure}

In Figure \ref{fig: dr_scaling_factor}, we show the curves of scaling factors $\tilde{c}_l$, $\tilde{c}_\omega$ for the dynamic rescaling formulation \eqref{eq: vort_stream_1_3d_noswirl_dr} and $c_l$, $c_\omega$ for the self-similar ansatz \eqref{eq: self_similar}. We observe a relatively fast convergence to the steady state as time increases. The converged values $c_l=4.549$ and $c_\omega=1.455$ are close to the approximate values obtained in Section \ref{sec: scaling-analysis}. Moreover, they also satisfy the relation \eqref{eq: scaling_relation}.

\begin{figure}[hbt!]
\centering
\includegraphics[width=.4\textwidth]{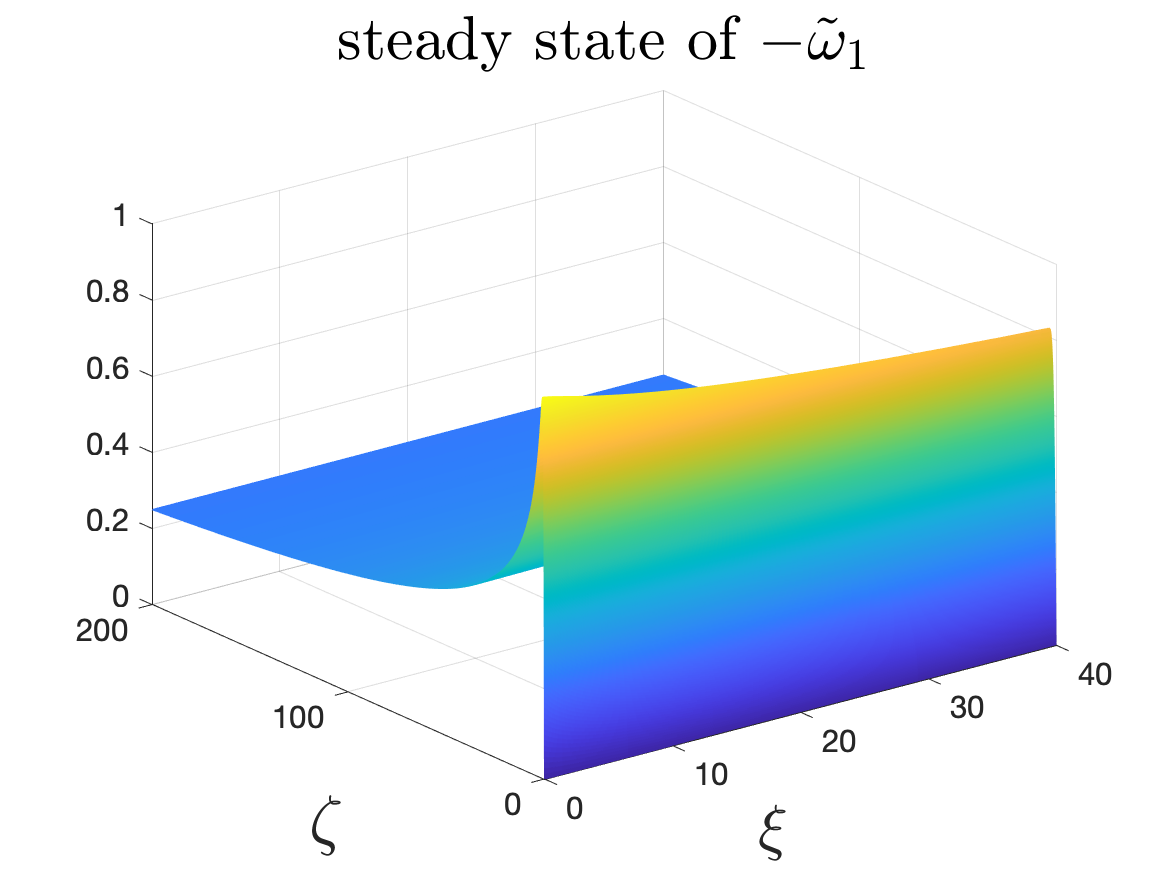}\includegraphics[width=.4\textwidth]{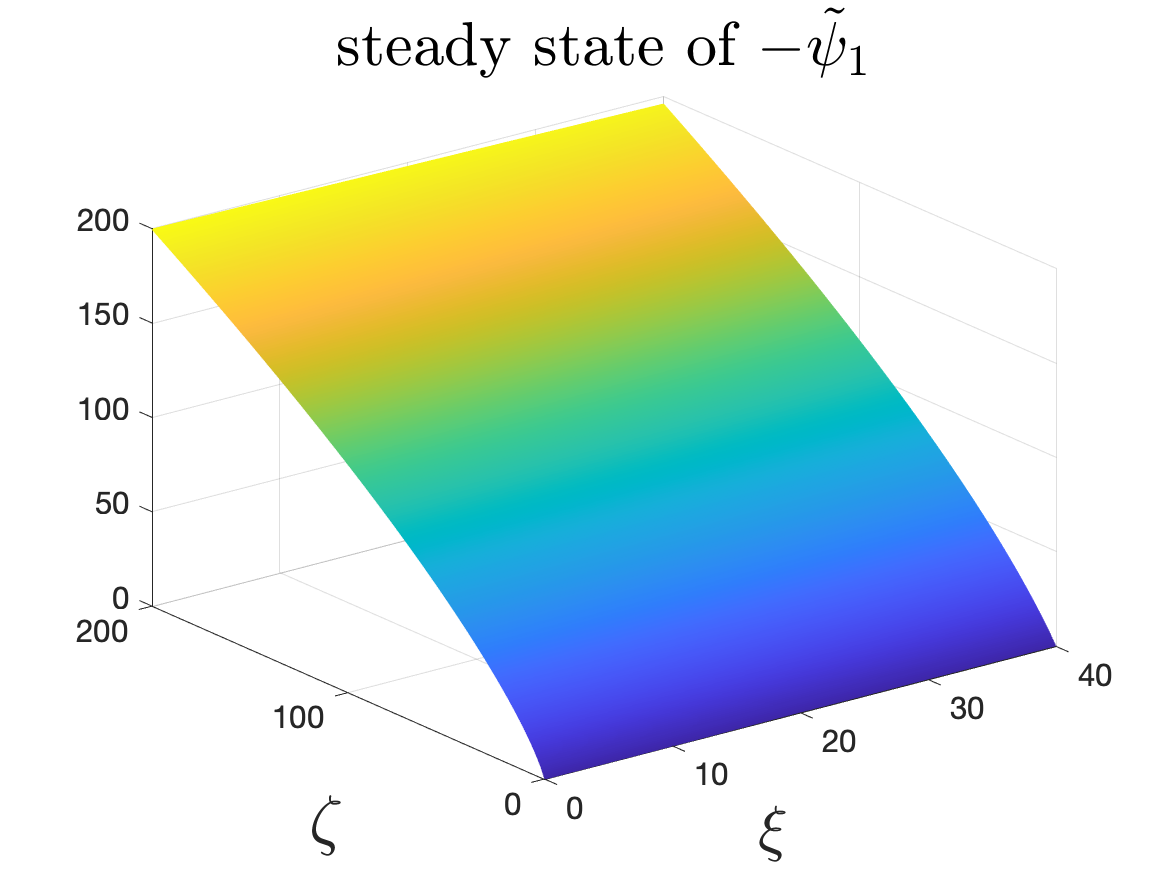}
\caption{Steady states of $-\tilde{\omega}_1$ and $-\tilde{\psi}_1$.}\label{fig: dr_steadystate}
\end{figure}

The approximate steady states of $\tilde{\omega}_1$ and $\tilde{\psi}_1$ are plotted in Figure \ref{fig: dr_steadystate}. We see that both $\tilde{\omega}_1$ and $\tilde{\psi}_1$ are relatively flat in $\xi$, suggesting a possible 1D structure of their profiles. While both functions have weak dependence on $\xi$, $-\tilde{\omega}_1$ seems to tilt up around $\xi=0$ a little bit. The shape of the steady states looks similar to the shape of the profiles we obtained via the adaptive mesh at the stopping time in Figure \ref{fig: end_data_zoom}.

\begin{figure}[hbt!]
\centering
\includegraphics[width=.4\textwidth]{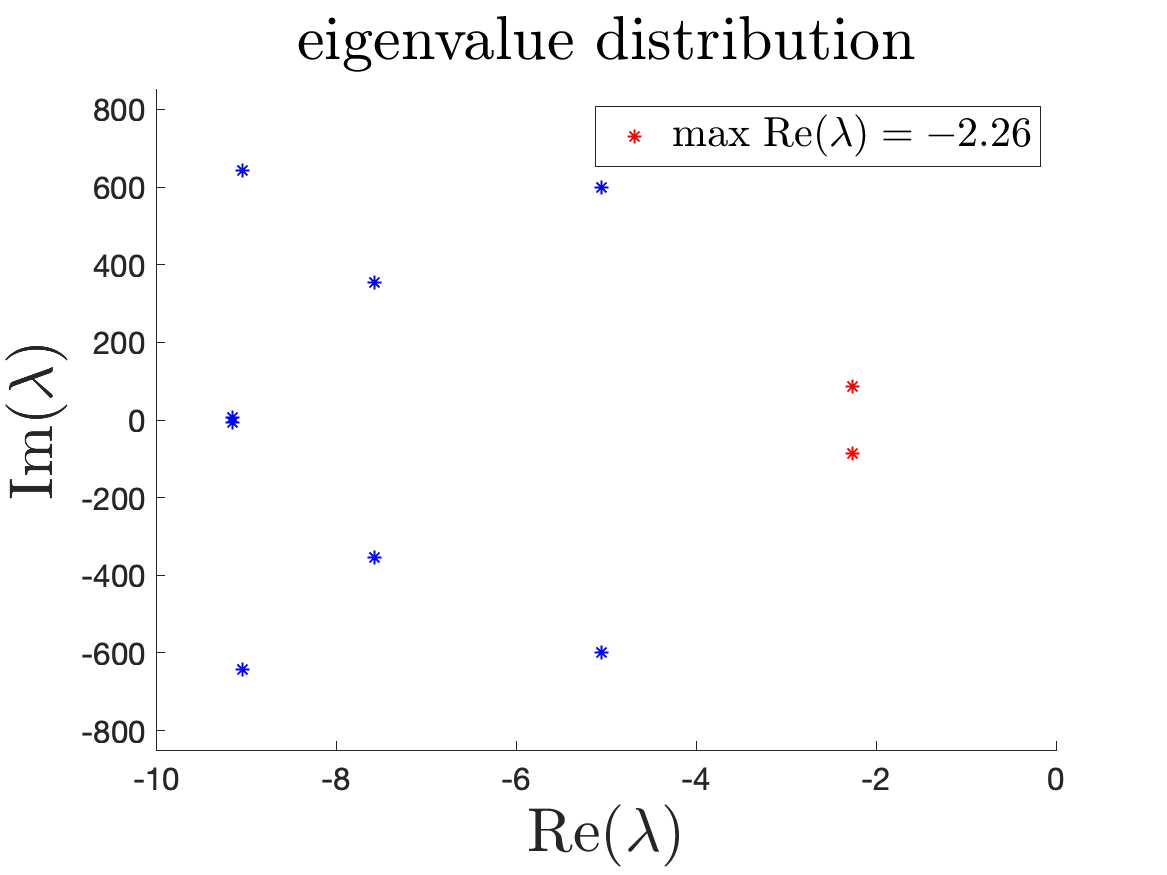}
\caption{Top 10 eigenvalues of the Jacobian with the largest real part.}\label{fig: jacobian_eigen}
\end{figure}

To study the stability of the steady states, we linearize the time-evolution equation \eqref{eq: vort_stream_1_3d_noswirl_dr} of $\tilde{\omega}_1$ near the steady states, and compute the eigenvalues of the Jacobian matrix. The top few eigenvalues with largest real part are plotted in Figure \ref{fig: jacobian_eigen}. We see that the eigenvalues have negative real part, which demonstrates the stability of our potential self-similar finite-time blow-up.

\subsection{The domain size study}
\label{sec: domain-size}

Ideally, the dynamic rescaling formulation \eqref{eq: vort_stream_1_3d_noswirl_dr} should be solved in the first quadrant $\mathcal{D}^\prime_\infty=\left\{(\xi,\zeta): \xi\geq0, \zeta\geq0\right\}$. In Section \ref{sec: setting_numeric}, we use a large rectangular region $\mathcal{D}^\prime=\left[0, D\right]\times\left[0, D/2\right]$, with domain size $D=10^5$, to approximate the unbounded domain $\mathcal{D}^\prime_\infty$ and propose to use the Neumann boundary condition at the far field.

\begin{figure}[hbt!]
\centering
\includegraphics[width=.4\textwidth]{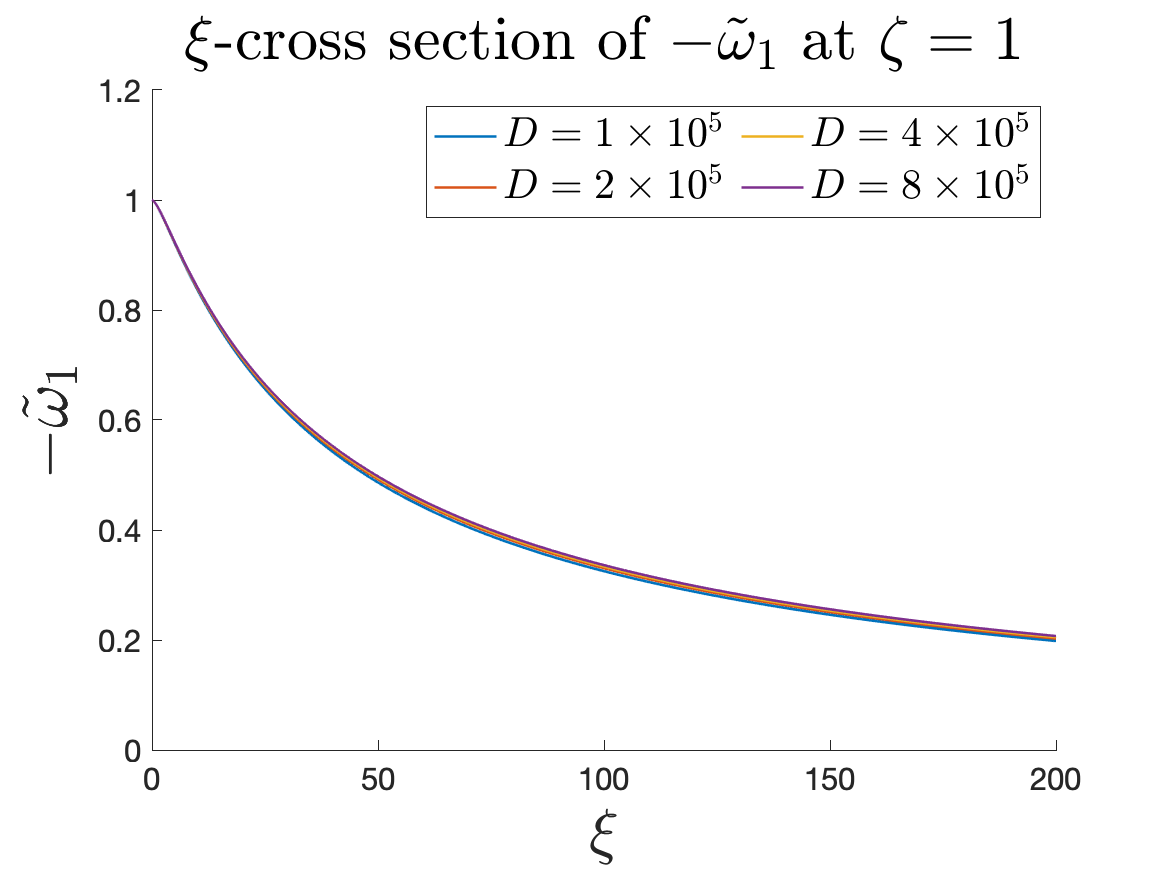}\includegraphics[width=.4\textwidth]{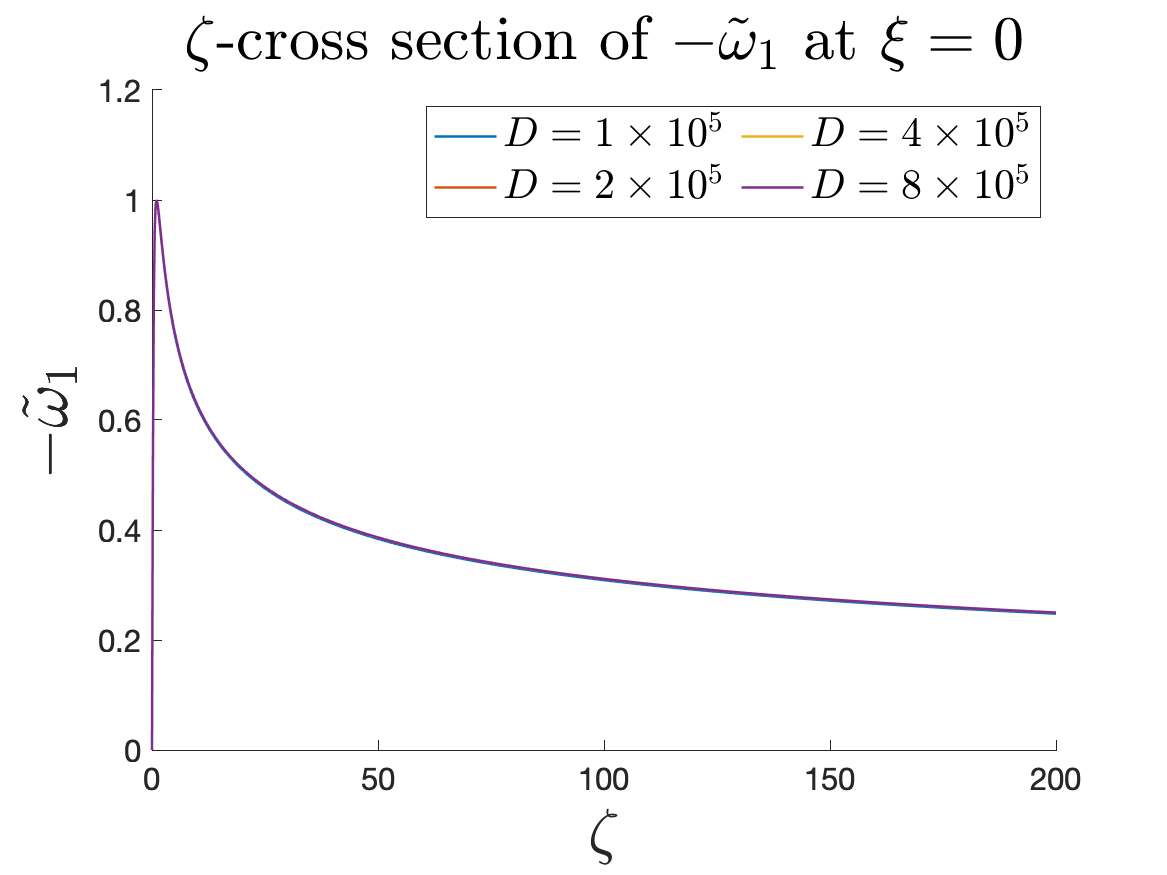}
\caption{Cross sections of steady states of $-\tilde{\omega}_1$ with different domain size $D$.}\label{fig: steadystate_cross_section_size}
\end{figure}

We study how the domain size would influence our steady state solution by extending the domain size and considering $D=1\times10^5, 2\times10^5, 4\times10^5, 8\times10^5$. Using these four different domain sizes, we solve the dynamic rescaling formulation \eqref{eq: vort_stream_1_3d_noswirl_dr} to its steady state. Figure \ref{fig: steadystate_cross_section_size} is the cross section comparison of the steady state of $\tilde{\omega}_1$. We also list the scaling factor $c_l$ with different domain size $D$ in Table \ref{tab: cl_D}. We can see that, with significant larger domain size, the steady state is nearly the same, and the scaling factor is almost the same. This shows our choice of domain size $D=10^5$ is large enough to approximate the steady state well.

\begin{table}[hbt!]
\centering
\caption{The scaling factor $c_l$ with different domain size $D$.}\label{tab: cl_D}
\begin{tabular}{|c|c|c|c|c|}
\hline
\hspace{0.5cm}$D$\hspace{0.5cm} & \hspace{0.1cm}$1\times10^5$\hspace{0.1cm} & \hspace{0.1cm}$2\times10^5$\hspace{0.1cm} & \hspace{0.1cm}$4\times10^5$\hspace{0.1cm} & \hspace{0.1cm}$8\times10^5$\hspace{0.1cm} \\
\hline
$c_l$ & $4.549$ & $4.546$ & $4.545$ & $4.545$ \\
\hline
\end{tabular}
\index{tables}
\end{table}

\section{The H\"{o}lder exponent and the dimension in the potential blow-up}
\label{sec: holder_dimension}

Starting this section, we will no longer fix the H\"{o}lder exponent $\alpha=0.1$.

\subsection{The H\"{o}lder exponent $\alpha$}
\label{sec: holder_exponent}

In his study of the finite-time blow-up of the axisymmetric Euler equations with no swirl and with H\"{o}lder continuous initial data \cite{elgindi2021finite}, Elgindi assumes that $\alpha$ is very close to zero, smaller than $10^{-14}$. Such small value of $\alpha$ is used to control the higher order terms of $\alpha$ in Elgindi's proof. However, as stated in the Conjecture 8 of \cite{drivas2022singularity} by Drivas and Elgindi, such a blow-up may still hold for a range of $\alpha\in(0,1/3)$ for the 3D Euler equations. For $\alpha>1/3$, it has been shown by \cite{ukhovskii1968axially, serfati1994regularite, shirota1994note, saint1994remarks, danchin2007axisymmetric, abidi2010global} that the solution will be globally regular. 

Therefore, we try different H\"{o}lder exponent $\alpha$ and explore the window of $\alpha$ that admits potential finite-time blow-up. For each $\alpha$, we first use the adaptive mesh method to solve the equations \eqref{eq: vort_stream_1_3d_noswirl} close enough to its potential blow-up time, and then use the dynamic rescaling method \eqref{eq: vort_stream_1_3d_noswirl_dr} to continue the computation and capture the self-similar profile.

\begin{figure}[hbt!]
\centering
\includegraphics[width=.4\textwidth]{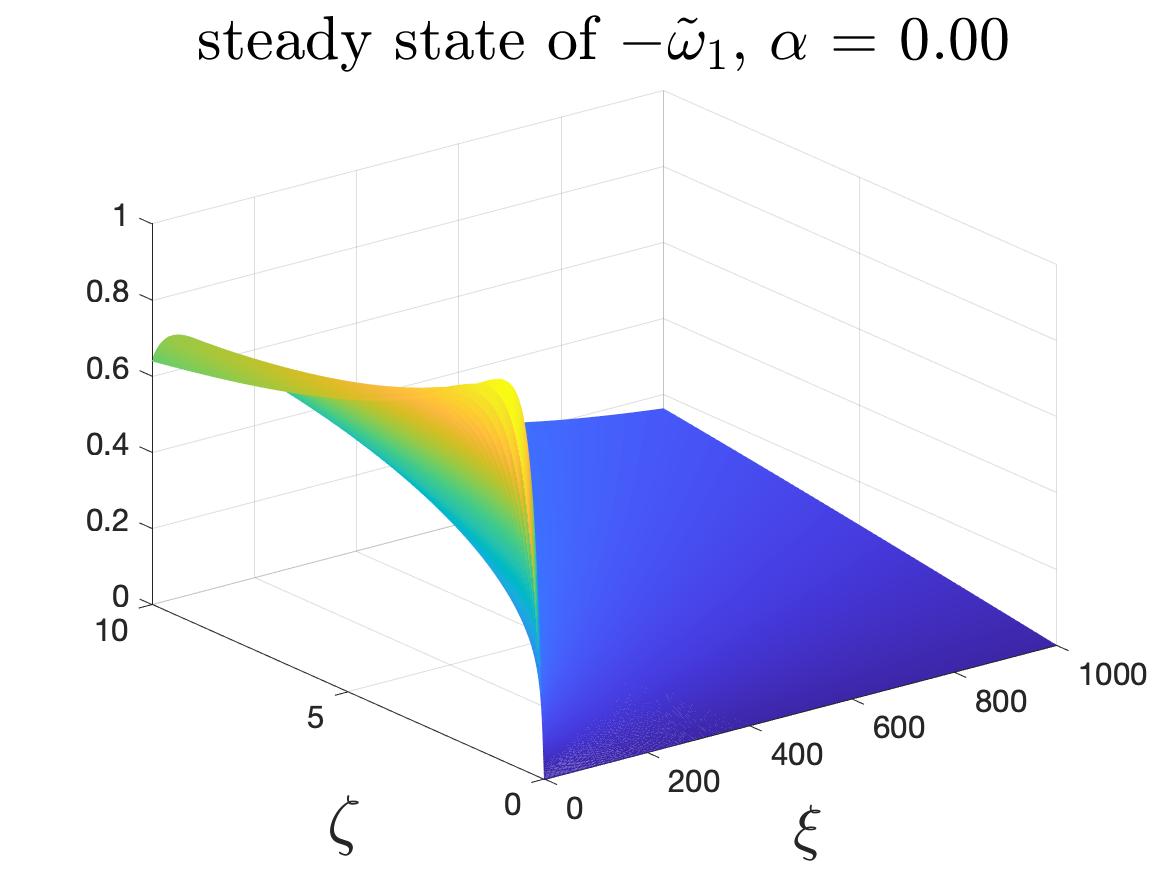}\includegraphics[width=.4\textwidth]{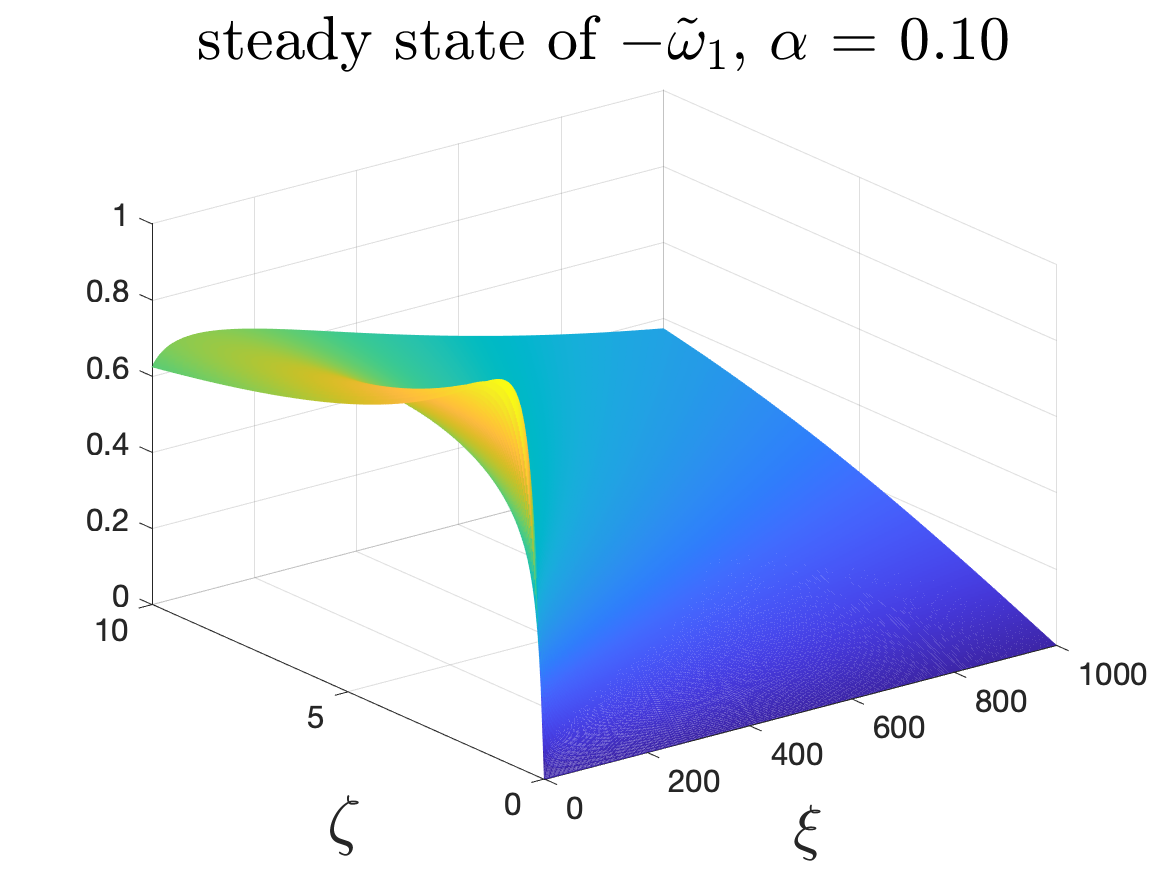}\\
\includegraphics[width=.4\textwidth]{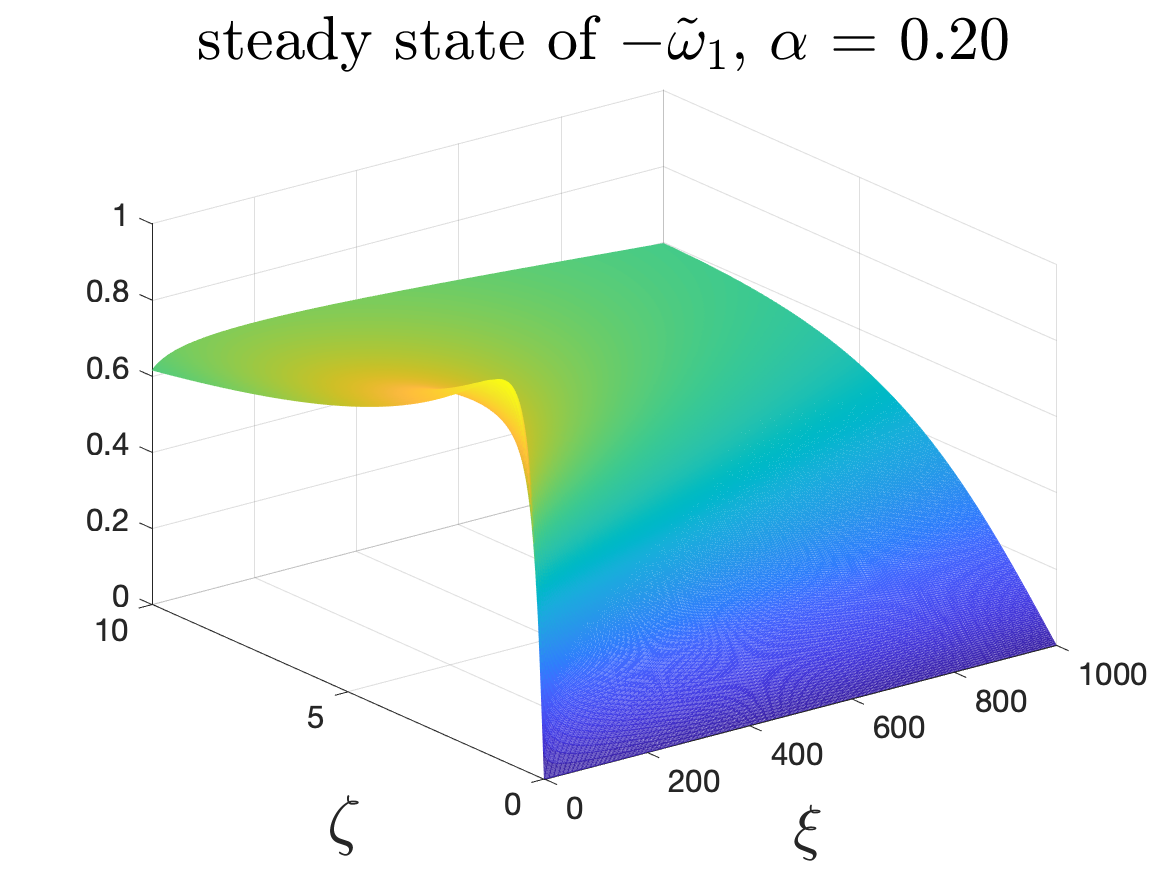}\includegraphics[width=.4\textwidth]{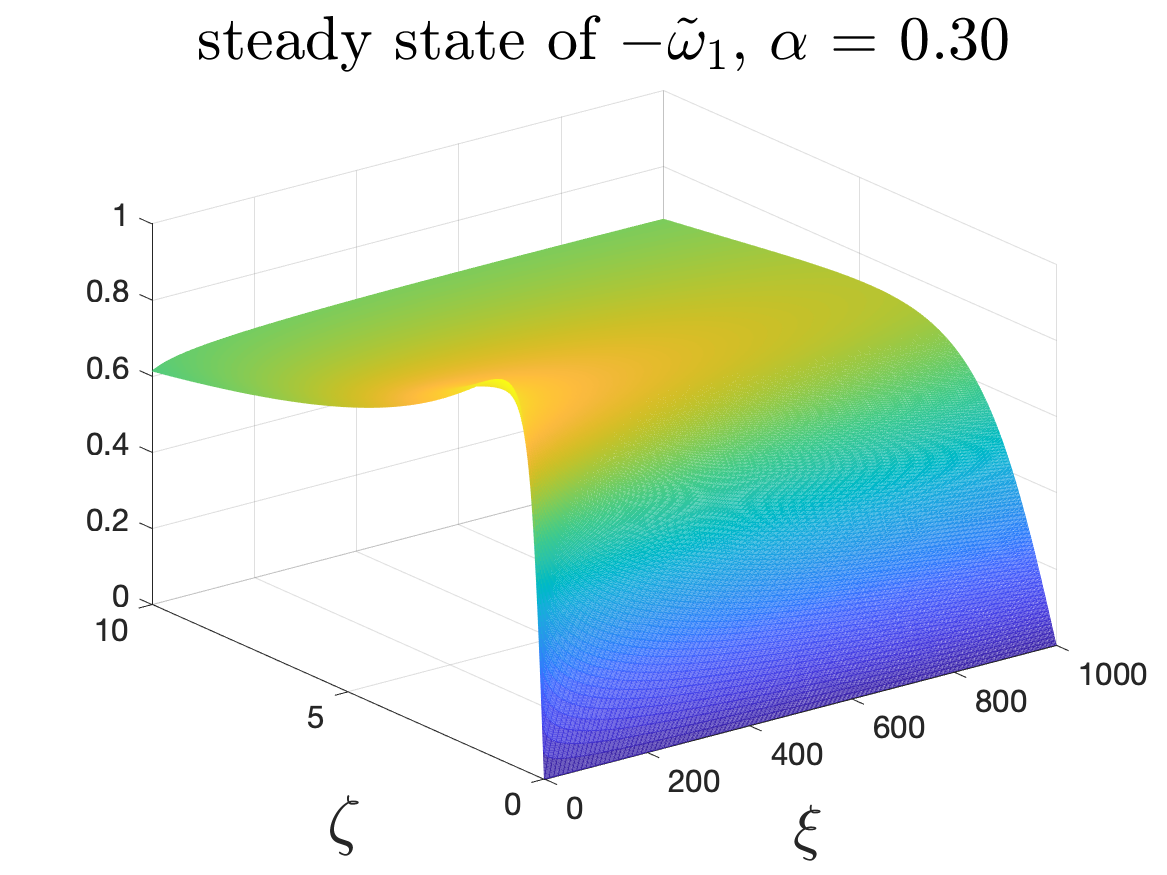}
\caption{Steady states of $-\tilde{\omega}_1$ with different $\alpha$ in $\mathbb{R}^3$.}\label{fig: steadystate_alpha}
\end{figure}

\begin{figure}[hbt!]
\centering
\includegraphics[width=.4\textwidth]{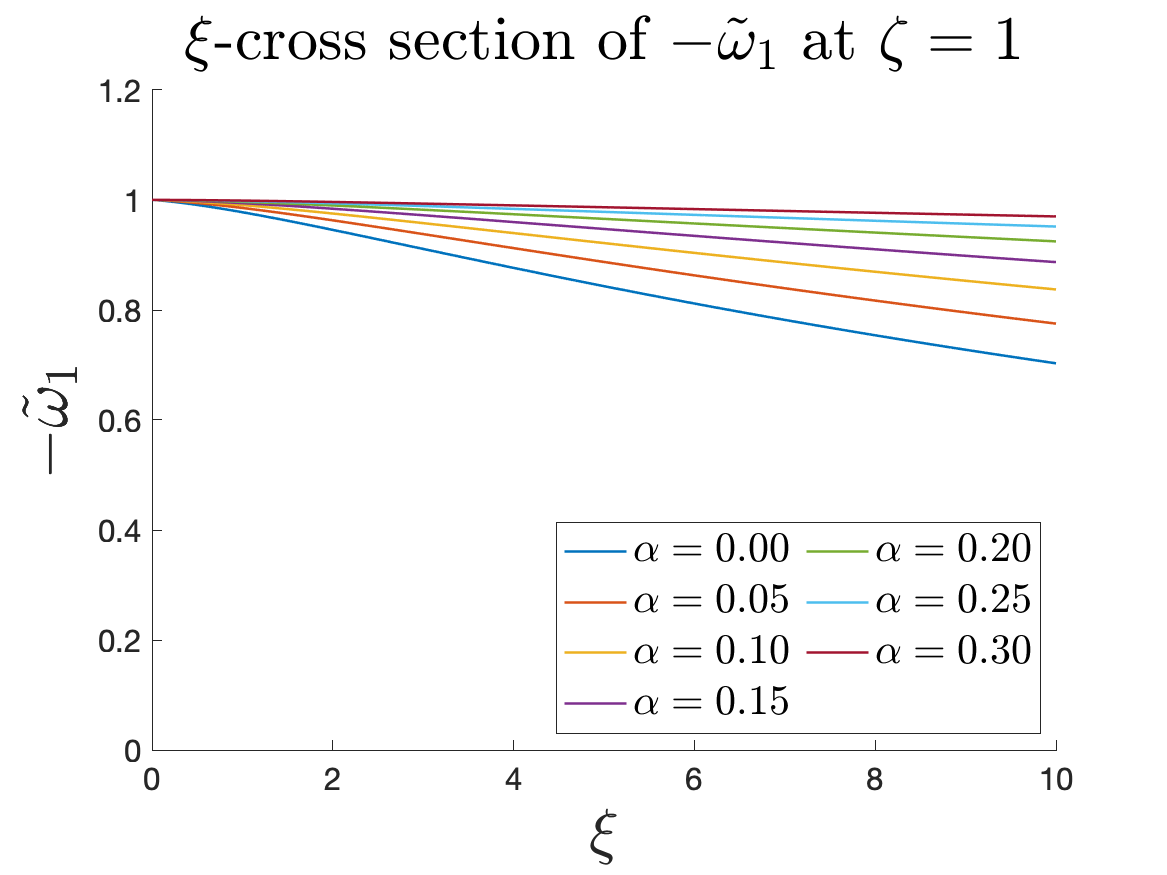}\includegraphics[width=.4\textwidth]{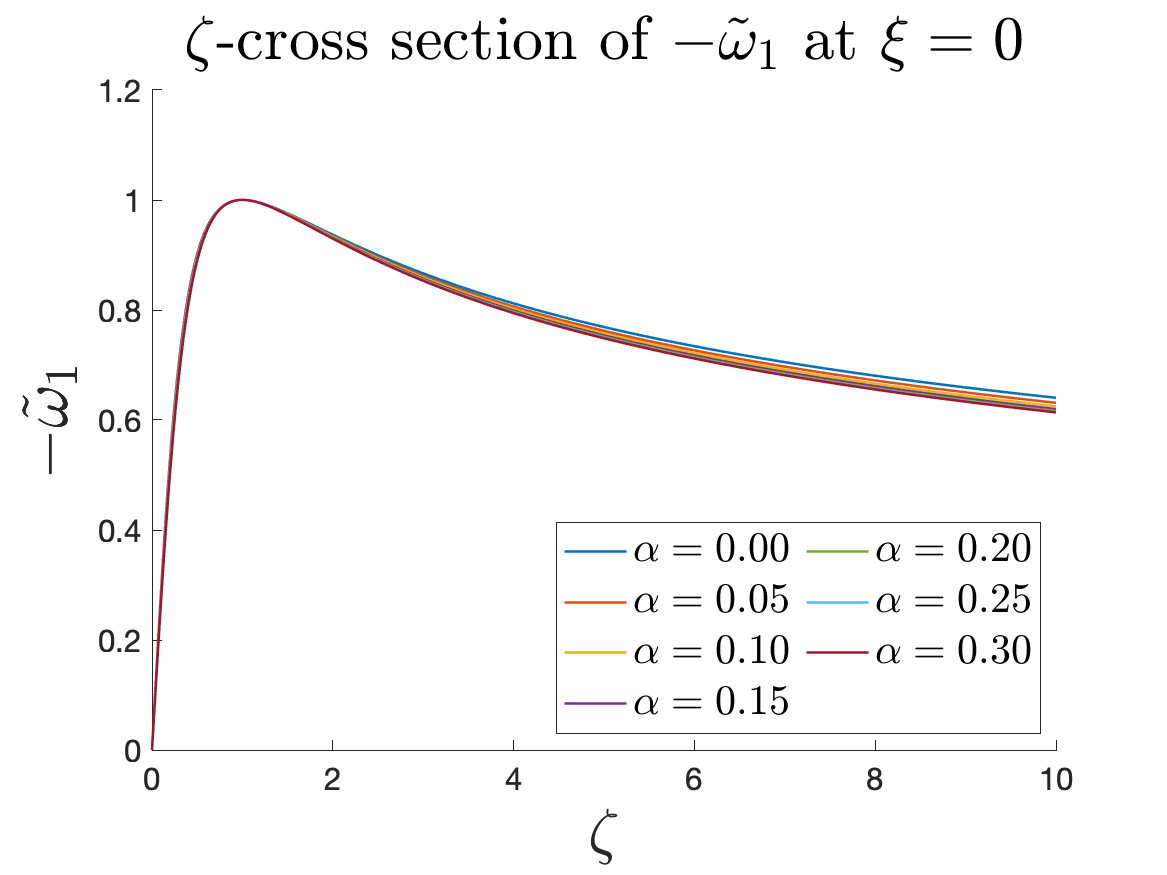}\\
\includegraphics[width=.4\textwidth]{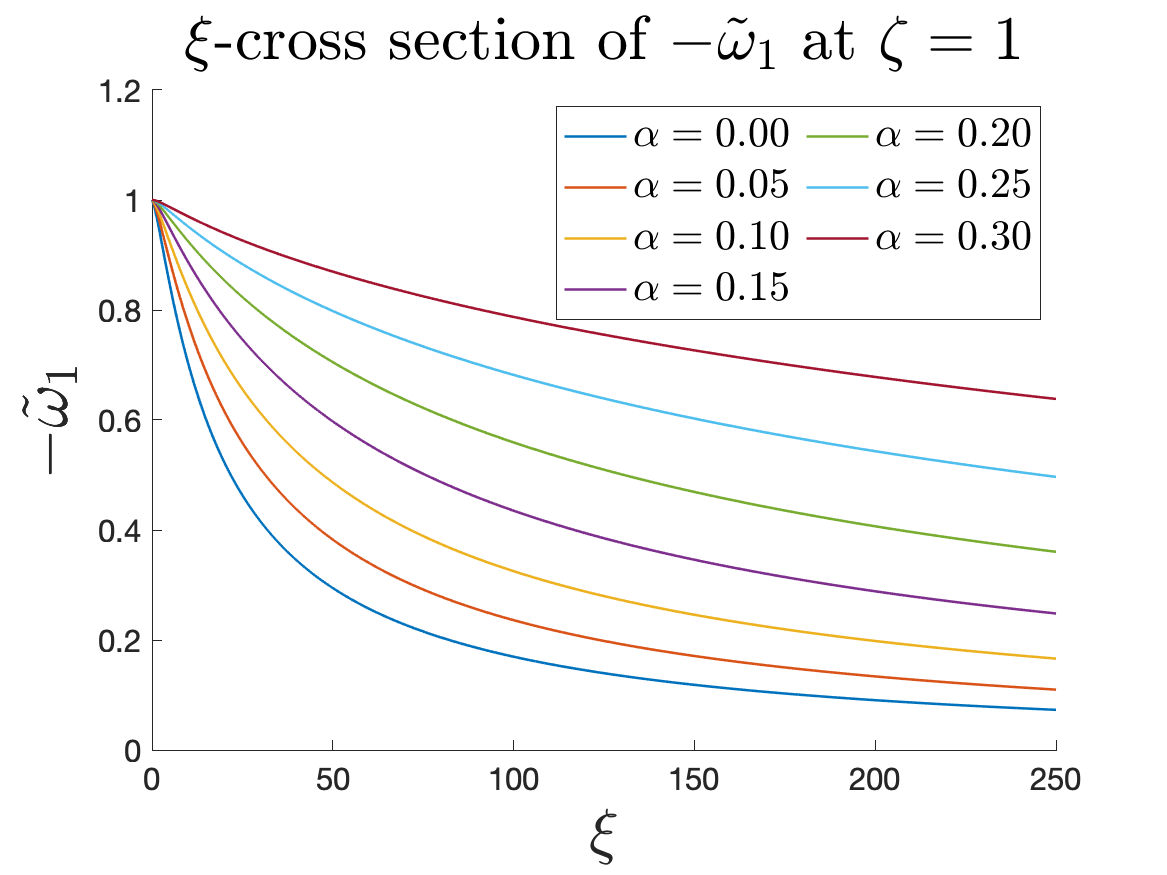}\includegraphics[width=.4\textwidth]{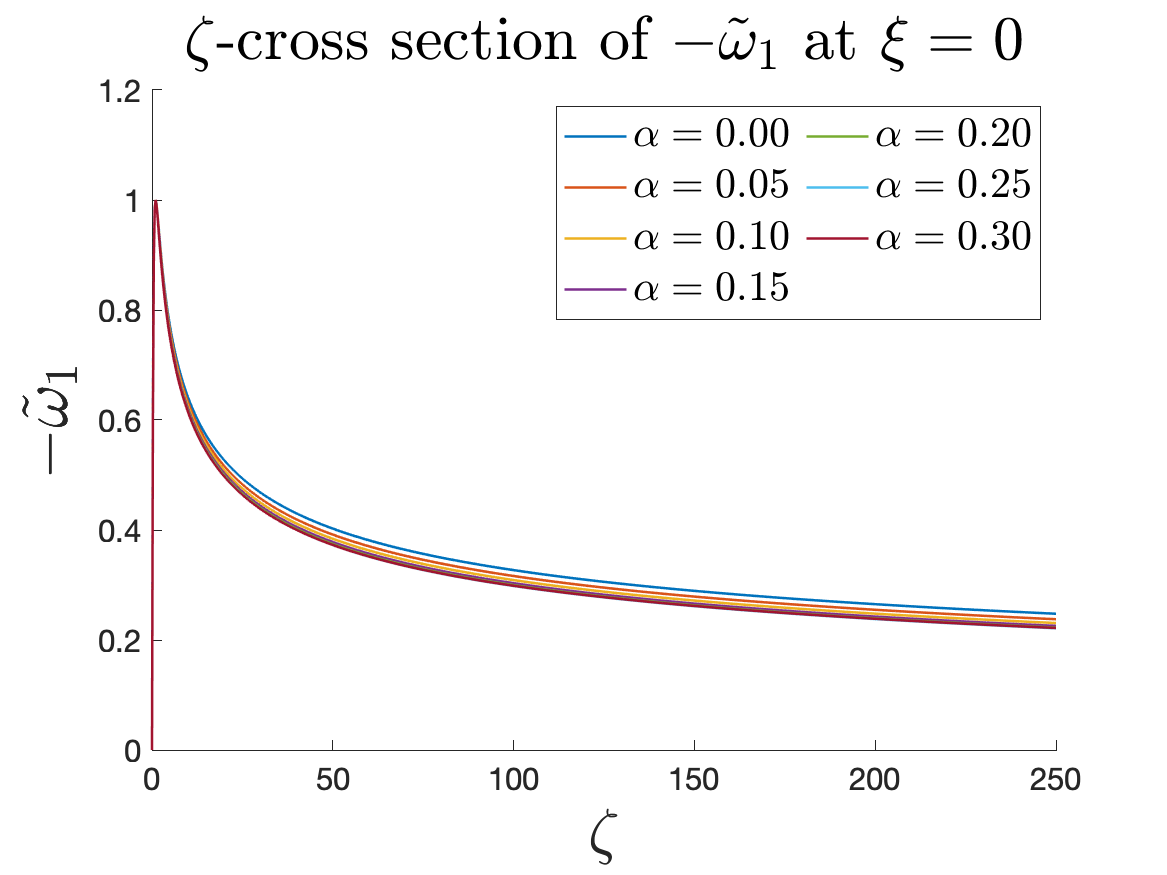}
\caption{Cross sections of steady states of $-\tilde{\omega}_1$ with different $\alpha$. Top row: on a local window. Bottom row: on a larger window.}\label{fig: steadystate_cross_section_alpha}
\end{figure}

\begin{figure}[hbt!]
\centering
\includegraphics[width=.4\textwidth]{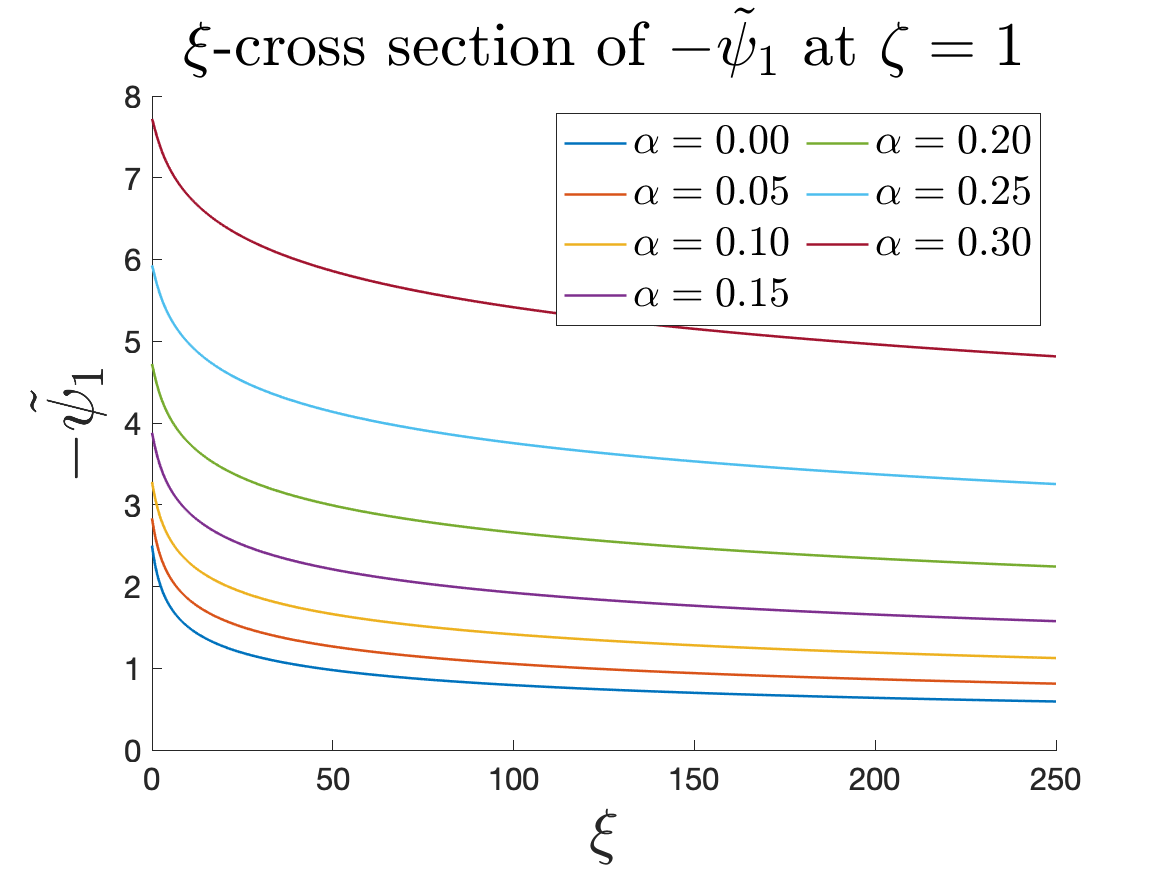}\includegraphics[width=.4\textwidth]{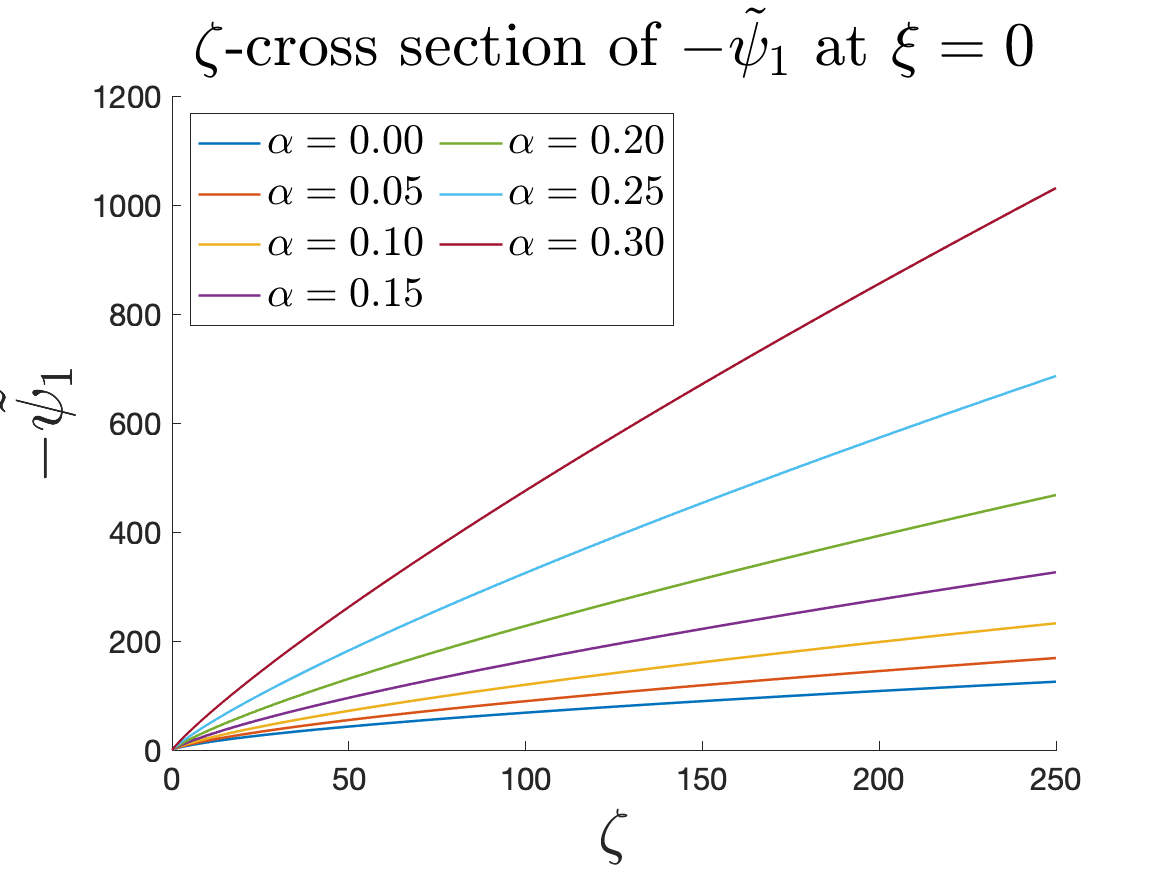}
\caption{Cross sections of steady states of $-\tilde{\psi}_1$ with different $\alpha$.}\label{fig: steadystate_cross_section_alpha_psi1}
\index{figures}
\end{figure}

For our 3D axisymmetric Euler equations with initial data \eqref{eq: inital_data} with $\alpha=0.0$, $0.1$, $0.2$, $0.3$, we obtain strong evidence for the formation of self-similar singularity. The steady states of the solutions are plotted in Figure \ref{fig: steadystate_alpha}. We can see that as $\alpha$ increases, $\tilde{\omega}_1$ will have weaker dependence on $\xi$, and the self-similar profile becomes more and more one-dimensional. We plot the cross sections of the steady states of $\tilde{\omega}_1$ in Figure \ref{fig: steadystate_cross_section_alpha}. As $\alpha$ increases, $-\tilde{\omega}_1(\xi, 1)$ becomes more and more flat, especially in the local window around $\xi=0$. Moreover, $-\tilde{\omega}_1(0, \zeta)$ seems to be insensitive to the value of $\alpha$. The cross sections of the stead states of $\tilde{\psi}_1$ in Figure \ref{fig: steadystate_cross_section_alpha_psi1} shows that as $\alpha$ increases, $-\tilde{\psi}_1(0, \zeta)$ becomes more and more like a linear function.

As $\alpha$ increases, $c_l$ increases rapidly. We can see from Table \ref{tab: cl_n3} that, $c_l$ is more than $100$ when $\alpha=0.3$. Such large $c_l$ can cause a lot of troubles for our adaptive mesh method, as the collapsing speed of the solution is extremely fast. Fortunately, the dynamic rescaling method is stable with large $c_l$, as the extra stretching term can control the rate of collapse. Based on observation of the rapid increase of $c_l$ as $\alpha$ approaches $\alpha^*$, we conjecture that $c_l$ will tend to infinity as $\alpha$ tends to $\alpha^*$. It is interesting to notice that when $\alpha=0.0$, the scaling factor $c_l$ is approximately 3.248, which doesn't seem special nor implies any degeneracy. In fact, when $\alpha=0$, our initial data \eqref{eq: inital_data} of $\omega_1$ is smooth in the axisymmetric variable $(r, z)$, but when we lift it to $\mathbb{R}^3$, the vorticity $\omega=\omega^\theta e_\theta$ will have singularity due to the coordinate singularity of $e_\theta$. In this case, the symmetry axis can be viewed as a boundary in the $(r, z)$-plane.

\begin{table}[hbt!]
\centering
\caption{The scaling factor $c_l$ with different $\alpha$ in the 3D case.}\label{tab: cl_n3}
\begin{tabular}{|c|c|c|c|c|c|c|c|}
\hline
\hspace{0.5cm}$\alpha$\hspace{0.5cm} & \hspace{0.2cm}$0.00$\hspace{0.2cm} & \hspace{0.2cm}$0.05$\hspace{0.2cm} & \hspace{0.2cm}$0.10$\hspace{0.2cm} & \hspace{0.2cm}$0.15$\hspace{0.2cm} & \hspace{0.2cm}$0.20$\hspace{0.2cm} & \hspace{0.2cm}$0.25$\hspace{0.2cm} & \hspace{0.2cm}$0.30$\hspace{0.2cm} \\
\hline
$c_l$ & $3.248$ & $3.771$ & $4.549$ & $5.818$ & $8.270$ & $15.00$ & $112.8$ \\
\hline
\end{tabular}
\index{tables}
\end{table}

\begin{figure}[hbt!]
\centering
\includegraphics[width=.4\textwidth]{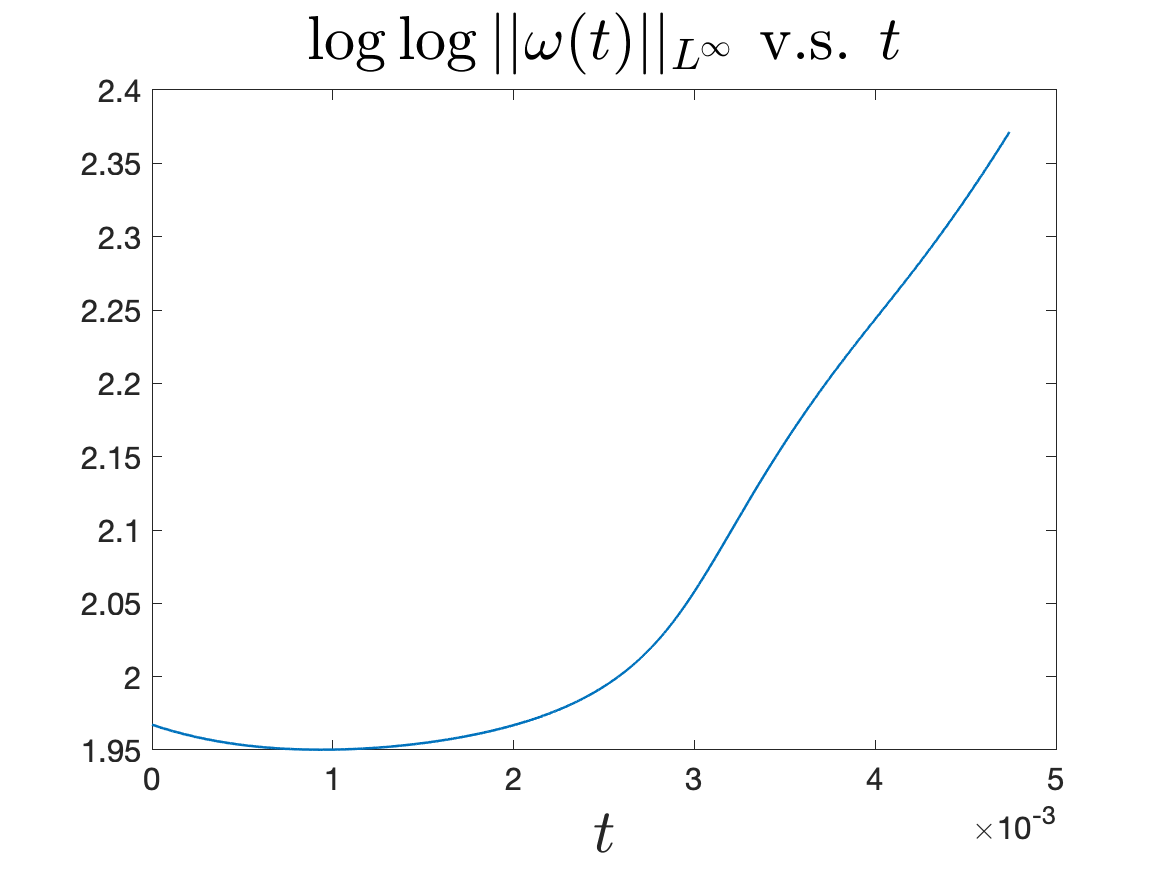}\includegraphics[width=.4\textwidth]{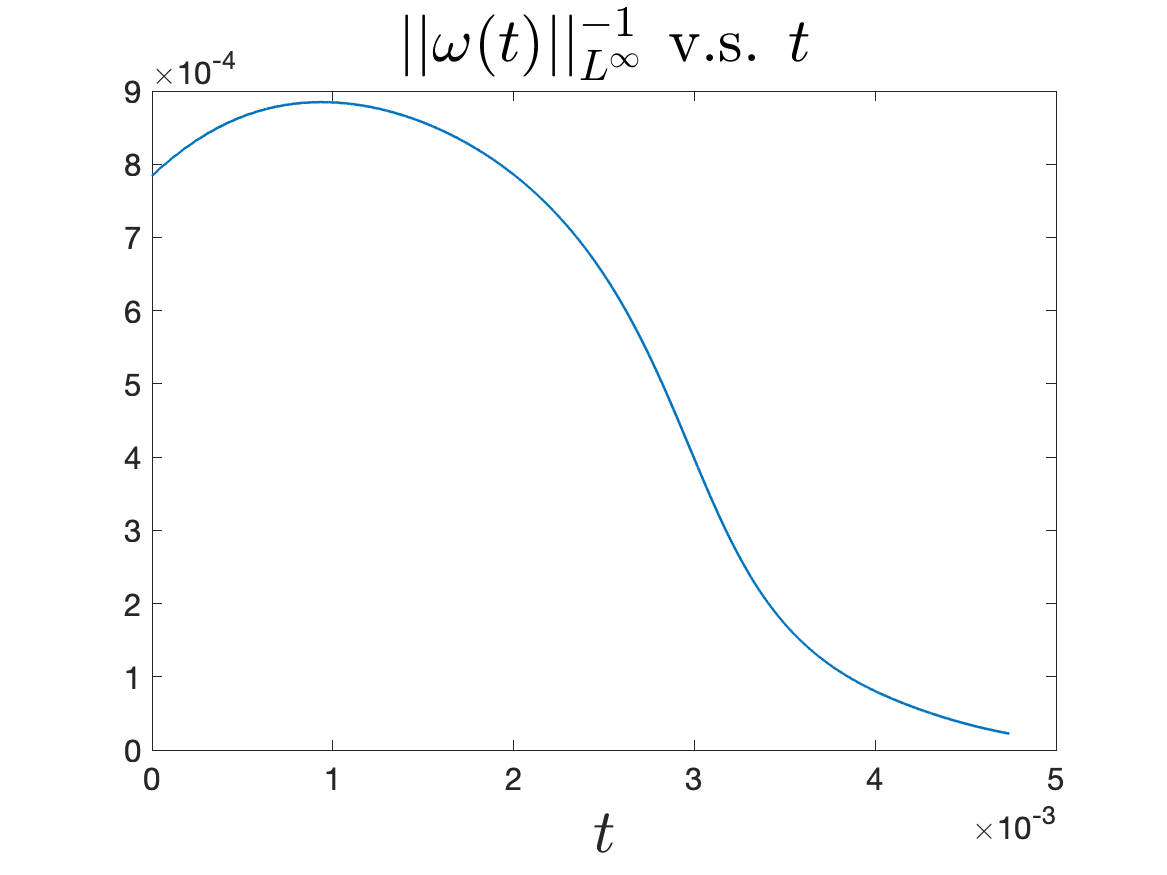}
\caption{Evidence of non-blow-up for $\alpha=0.31$ in $\mathbb{R}^3$.}\label{fig: no_blow_up}
\end{figure}

\begin{figure}[hbt!]
\centering
\includegraphics[width=.4\textwidth]{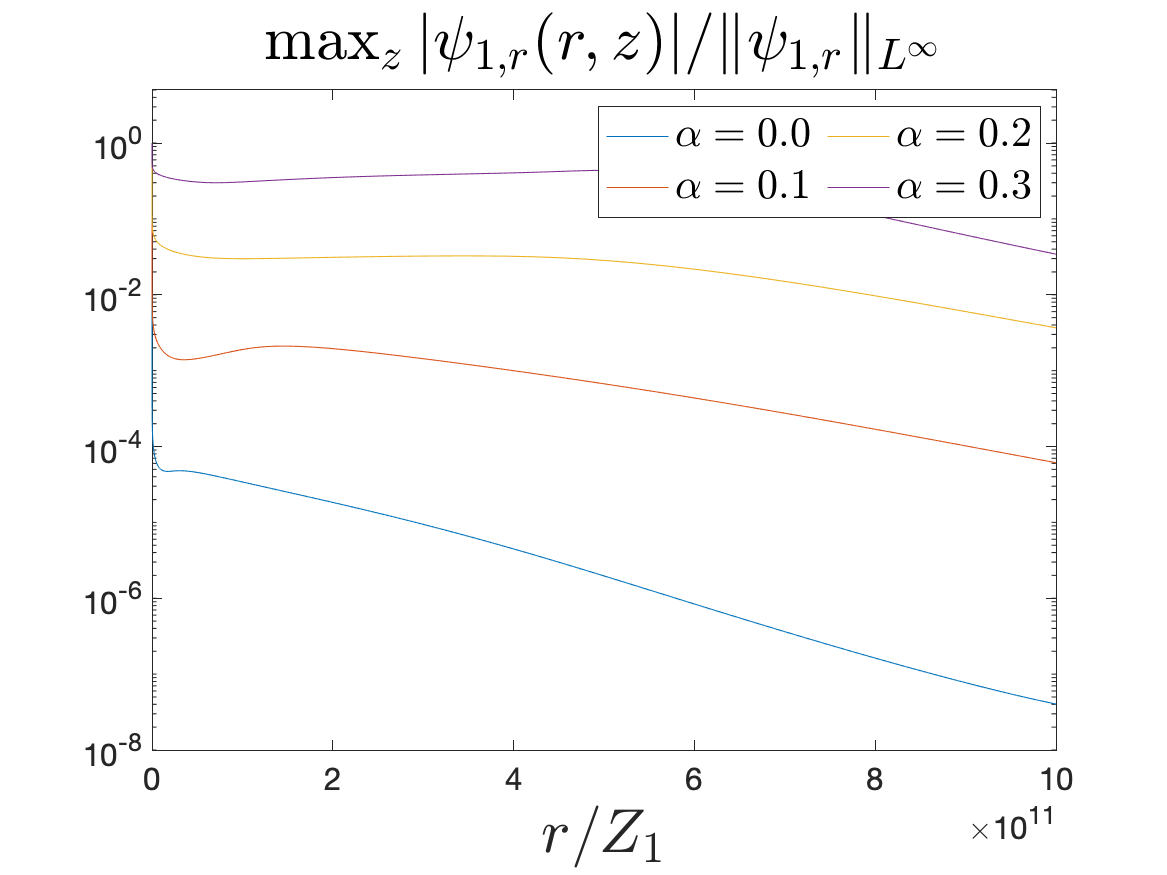}\includegraphics[width=.4\textwidth]{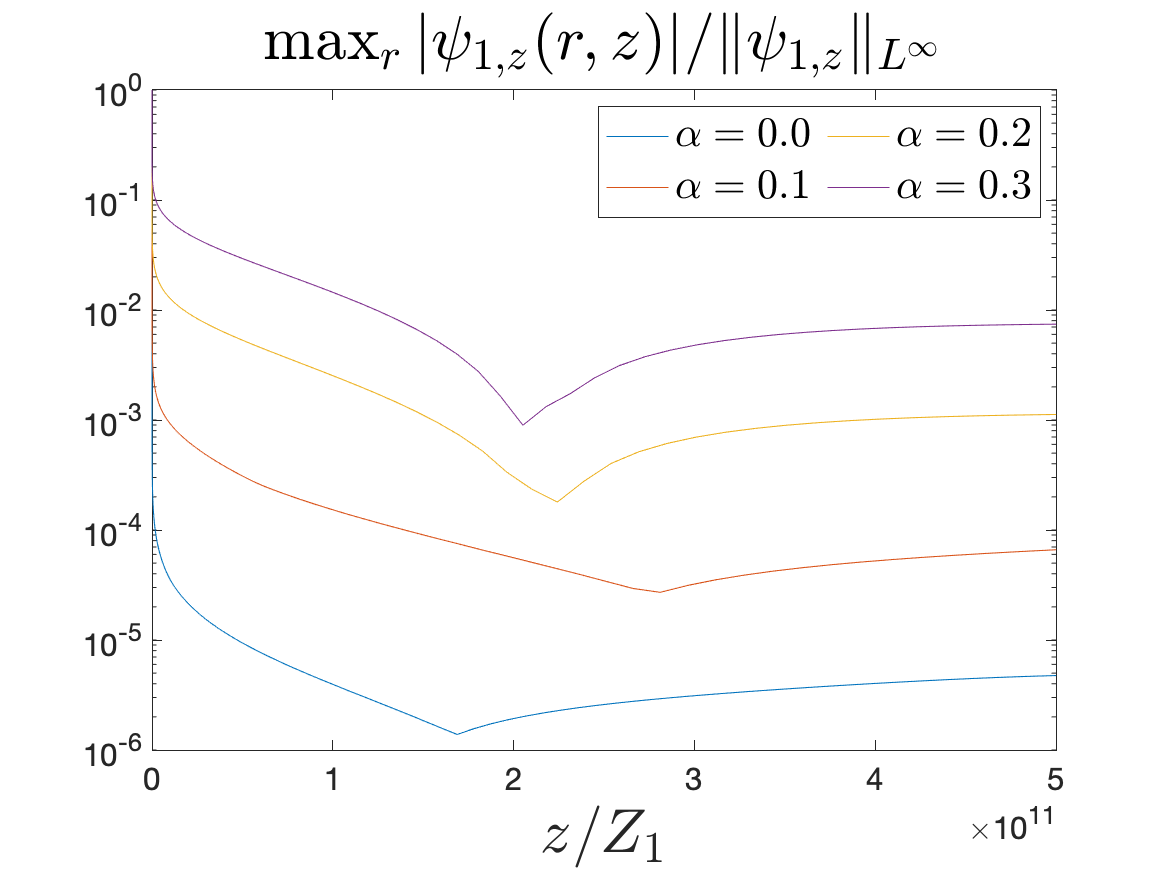}
\caption{Decay of the derivatives of $\psi_1$ with different $\alpha$.}\label{fig: dr_decay_alpha}
\index{figures}
\end{figure}

For $\alpha>0.30$, like $\alpha=0.31, 0.40, 0.50$, we observe that although $\|\tilde{\omega}\|_{L^\infty}$ grows rapidly in the initial stage, it eventually slows down and starts to decrease, and the dynamic rescaling formulation fails to converge to a steady state. For example, in the case of $\alpha=0.31$ shown in Figure \ref{fig: no_blow_up}, the double logarithm of $\|\omega\|_{L^\infty}$ becomes sublinear in the late stage, and $\|\omega\|_{L^\infty}^{-1}$ seems to decay slowly to zero, which would violate the Beale-Kato-Majda blow-up criterion. While the value $\alpha=0.31$ is still far from the critical case of $\alpha=1/3$, we remark that this could be due to the fact that the stability of the steady states becomes weaker as $\alpha$ tends to the critical value $\alpha^*$. The largest real part of the eigenvalues of the Jacobian matrix, like we plotted in Figure \ref{fig: jacobian_eigen}, is $-2.26$ when $\alpha=0.1$, $-1.89$ when $\alpha=0.2$, and $-0.53$ when $\alpha=0.3$. This shows that the steady states are less stable when $\alpha$ approaches $\alpha^*$. Another reason is that $\tilde{c}_\omega + \alpha \tilde{c}_l$ is very close to zero as $\alpha$ tends to the critical value $\alpha^*$. Thus, numerical errors may cause $\tilde{c}_\omega + \alpha \tilde{c}_l$ to change sign dynamically when $\alpha$ approaches $\alpha^*$, which leads to non-blow-up. In order to capture the blow-up behavior as $\alpha$ approaches $\alpha^*$, we need much higher resolution to prevent $\tilde{c}_\omega + \alpha \tilde{c}_l$ to change sign, which poses great numerical challenges in accuracy and computational time. We also note that the domain size $D$ of the computational domain $\mathcal{D}^\prime$ needs to be enlarged as $\alpha$ approaches the critical value $\alpha^*$. In Figure \ref{fig: dr_decay_alpha}, we plot the decay of derivatives of $\psi_1$ at the end of adaptive mesh method computation, which is also the starting point of the dynamic rescaling formulation computation. We can see that $\psi_{1,r}$ decays slower with a larger $\alpha$, which implies that we need a larger domain size $D$ to ensure that the zero Neumann boundary condition is appropriate. This further poses challenges in accuracy and computational time. Despite the numerical difficulty, the consistent potential blow-up behavior for $\alpha\leq 0.3$ makes us believe that  the 3D Euler equations could also develop potential self-similar blow-up for all $\alpha < 1/3$.

\subsection{The dimension $n$}
\label{sec: dimension}

We extend the numerical study from the three-dimensional case to the higher dimensional case, and show that the critical value $\alpha^*$ is close to $1-\frac{2}{n}$ in the $n$-dimensional case. It is expected that the higher the dimension, the easier it is for a perfect fluid to form a singularity in finite time, see Section 6 of \cite{drivas2022singularity}. We provide some numerical evidence to support such point of view.

In the $n$-dimensional case, we still use
\begin{align*}
    u(x,t):&\mathbb{R}^n\times\left[0, T\right)\rightarrow\mathbb{R}^n,\qquad\text{ and }\qquad p(x,t):\mathbb{R}^n\times\left[0, T\right)\rightarrow\mathbb{R},
\end{align*}
to denote the $n$-D vector field of the velocity and the $n$-D scalar field of the pressure respectively, where $x=(x_1, x_2, \ldots, x_n)\in\mathbb{R}^n$. Then the $n$-dimensional Euler equations can be written as
\begin{subequations}
\label{eq: euler_nd}
\begin{align}
    u_t + u\cdot\nabla u &= -\nabla p, \label{eq: velo_nd}\\
    \nabla\cdot u&=0. \label{eq: incomp_nd}
\end{align}
\end{subequations}

Next, we introduce $r=\sqrt{\sum_{k=1}^{n-1}x_k^2}$, $z=x_n$, and the unit vectors
\begin{align*}
    e_r=(x_1/r, x_2/r, \ldots, x_{n-1}/r, 0),\qquad e_z=(0, 0, \ldots, 0, 1)\,.
\end{align*}
Similar to the 3D case, we call an $n$-D vector field $v: \mathbb{R}^n\rightarrow\mathbb{R}^n$ to be axisymmetric and no swirl if the following ansatz applies
\begin{align}
    v=v^r(r,z)e_r+v^z(r,z)e_z.
    \label{eq: axisymmetric_no_swirl}
\end{align}
The axisymmetric $n$-D Euler equations with no swirl can be written in the vorticity-stream function form as
\begin{subequations}
\label{eq: vort_stream_nd}
\begin{align}
    \omega^\theta_t + u^r\omega^\theta_r + u^z\omega^\theta_z &= \frac{n-2}{r}u^r\omega^\theta, \label{eq: vort_theta_nd}\\
    -\psi^\theta_{rr}-\psi^\theta_{zz}-\frac{n-2}{r}\psi^\theta_{r}+\frac{n-2}{r^2}\psi^\theta&=\omega^\theta, \label{eq: stream_theta_nd}\\
    u^r=-\psi^\theta_z,\quad u^z&=\frac{n-2}{r}\psi^\theta+\psi^\theta_r , \label{eq: velo_rz_nd}
\end{align}
\end{subequations}
where we introduce the angular vorticity $\omega^\theta$ and angular stream function $\psi^\theta$ as $\omega^\theta=u^r_z-u^z_r$ and $-\Delta\psi^\theta=\omega^\theta$, similar to the 3D axisymmetric Euler equations.

We would like to note that here we consider axisymmetric and no swirl condition together in \eqref{eq: axisymmetric_no_swirl} because, if the velocity has component on a direction orthogonal to $e_r$ and $e_z$, then the incompressibility condition $\nabla\cdot u=0$ will inevitably introduce dependence on variable other than $r$ and $z$ when the dimension $n$ is greater than $3$, even if this component only depends on $r$ and $z$ at time $t=0$.

Since we focus on $C^\alpha$ continuous initial data for the angular vorticity $\omega^\theta$, we introduce $(\omega_1, \psi_1)$ similarly as in \eqref{eq: new_variable_3d_noswirl}, and re-write the $n$-D axisymmetric Euler equations with no swirl in the below
\begin{subequations}
\label{eq: vort_stream_nd_1_noswirl}
\begin{align}
    \omega_{1,t} + u^r\omega_{1,r} + u^z\omega_{1,z} &= -(n-2-\alpha)\psi_{1,z}\omega_1, \label{eq: vort_theta_nd_1_noswirl}\\
    -\psi_{1,rr}-\psi_{1,zz}-\frac{n}{r}\psi_{1,r}&=\omega_1r^{\alpha-1}, \label{eq: stream_theta_nd_1_noswirl}\\
    u^r=-r\psi_{1,z},\quad u^z&=(n-1)\psi_1+r\psi_{1,r}. \label{eq: velo_rz_1_nd_noswirl}
\end{align}
\end{subequations}

Roughly speaking, the dimension $n$ controls the strength of the vortex stretching term $-(n-2-\alpha)\psi_{1,z}\omega_1$ and the velocity in the $z$-direction $u^z=(n-1)\psi_1+r\psi_{1,r}$. It also modifies the Poisson equation for $\psi_1$. It seems natural to conjecture that the singularity formation will be more likely in the high-dimensional case because of the stronger vortex stretching term $-(n-2-\alpha)\psi_{1,z}\omega_1$.

In the following, for each combination of $\alpha$ and $n$, we first use the adaptive mesh method to solve the equations \eqref{eq: vort_stream_nd_1_noswirl} close enough to its possible blow-up time, and then use the dynamic rescaling formulation to continue the computation to capture the potential self-similar structure.

\begin{figure}[hbt!]
\centering
\includegraphics[width=.4\textwidth]{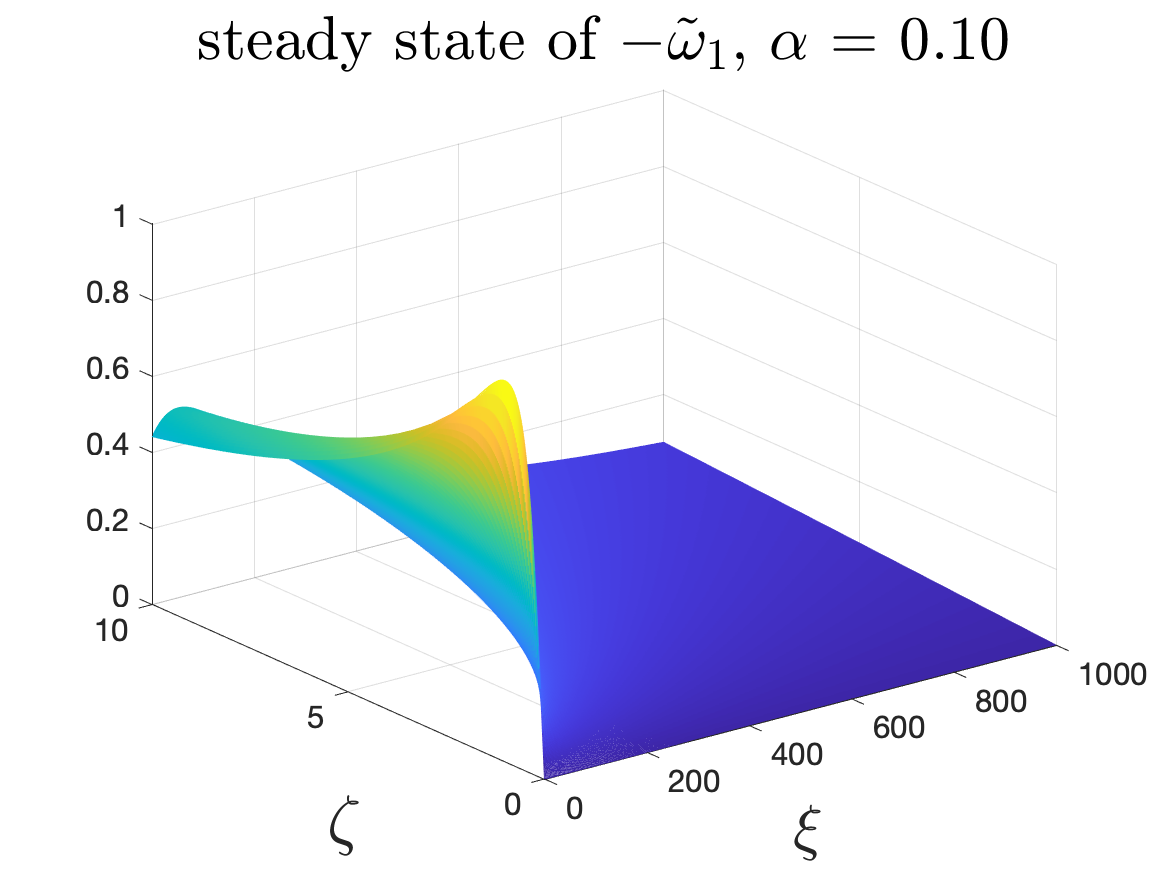}\includegraphics[width=.4\textwidth]{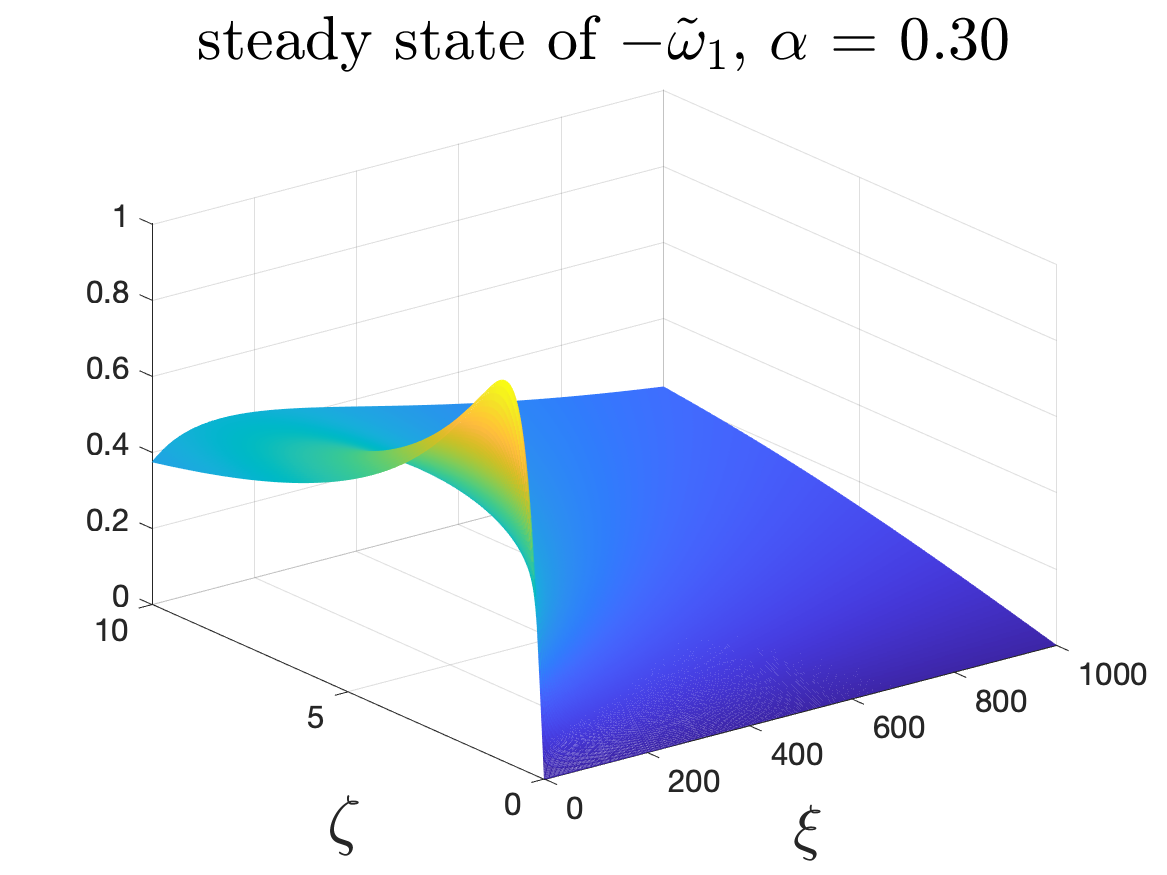}\\
\includegraphics[width=.4\textwidth]{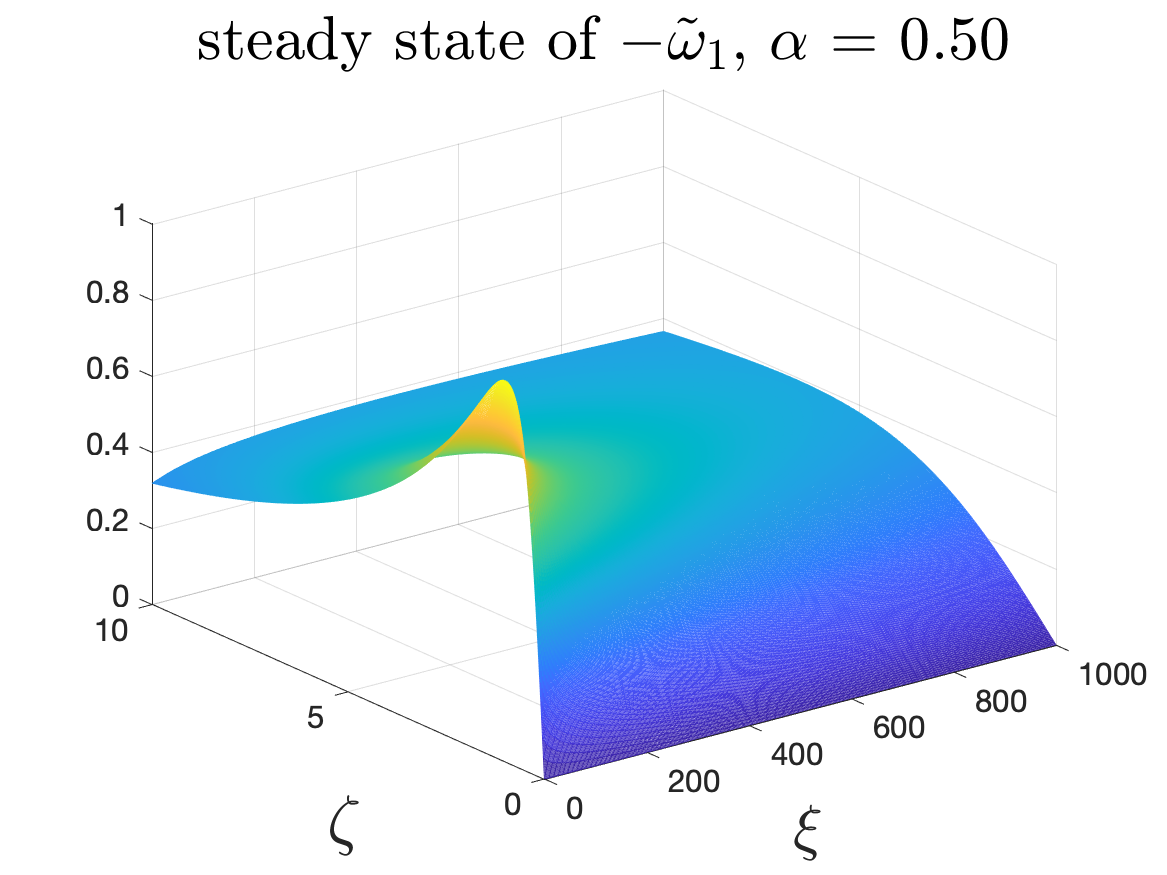}\includegraphics[width=.4\textwidth]{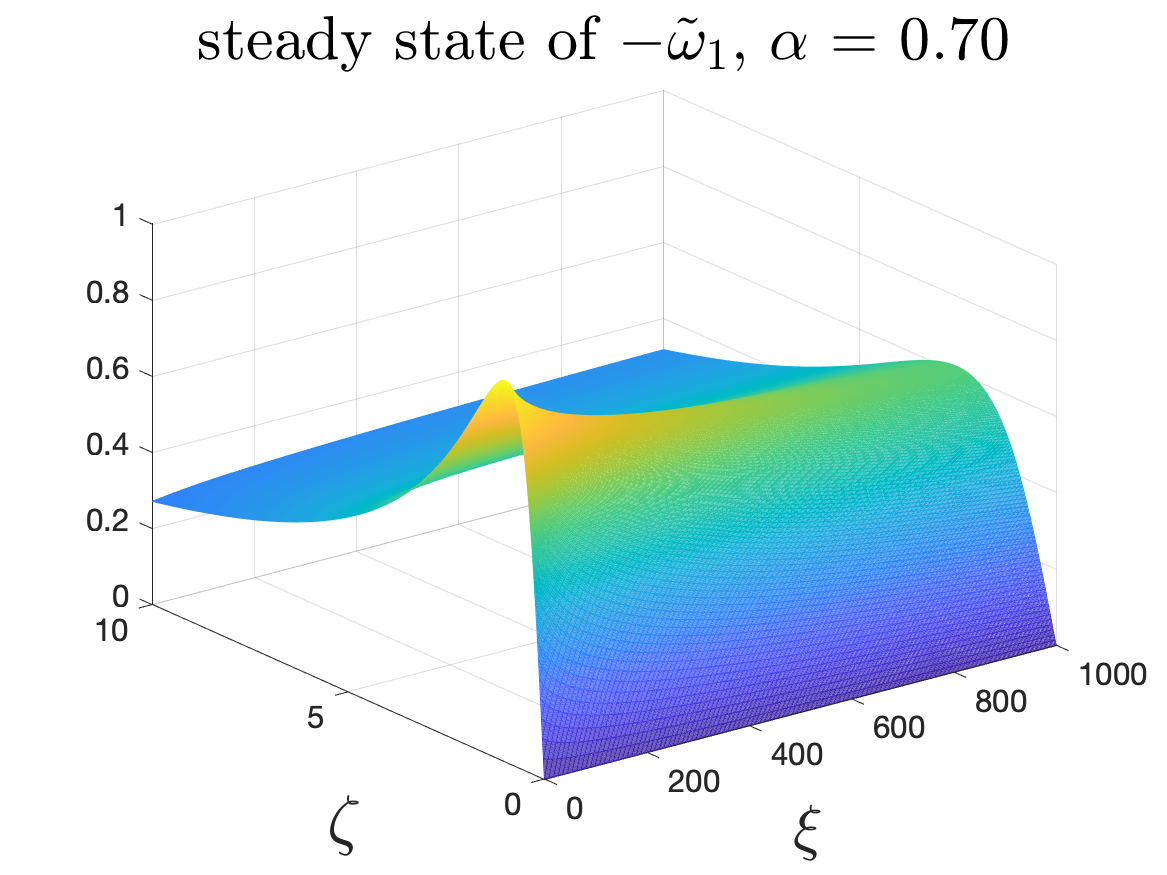}
\caption{Steady states of $-\tilde{\omega}_1$ with different $\alpha$ in $\mathbb{R}^{10}$.}\label{fig: steadystate_alpha_n10}
\end{figure}

\begin{figure}[hbt!]
\centering
\includegraphics[width=.4\textwidth]{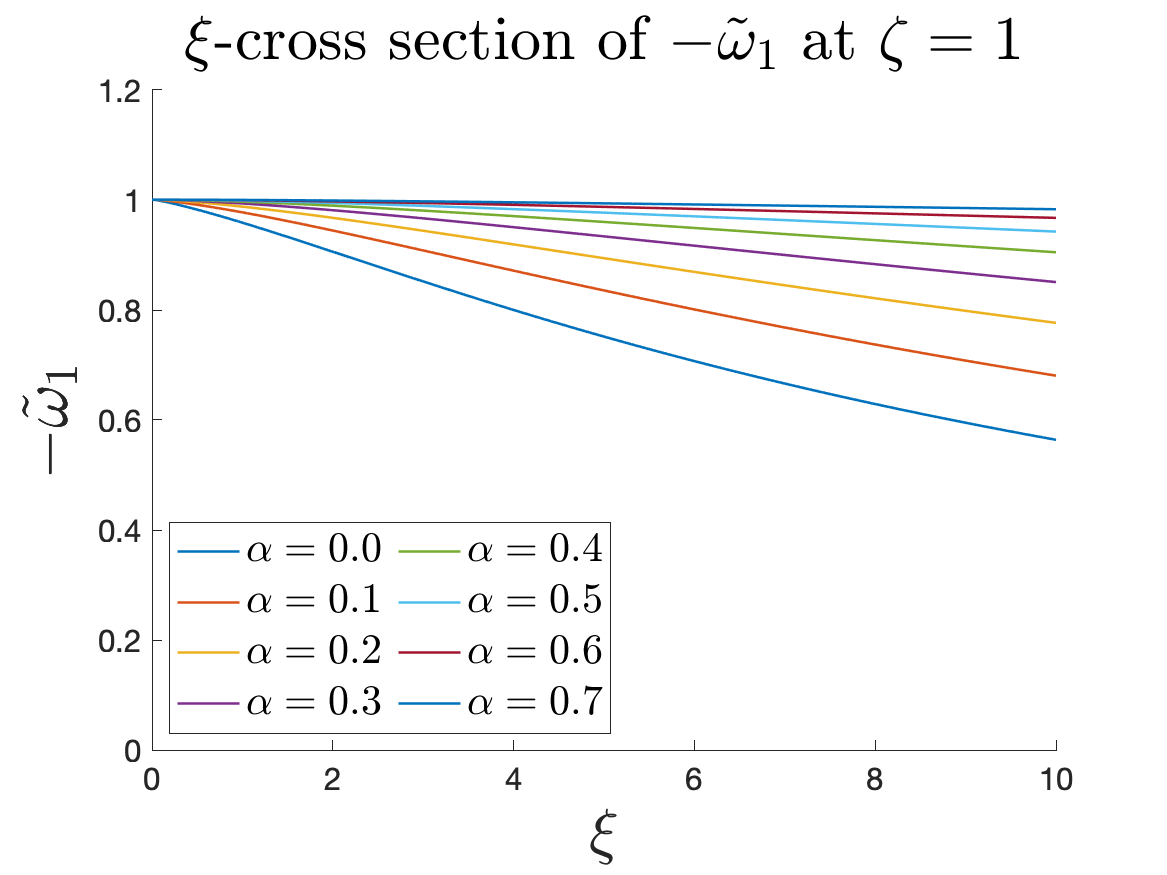}\includegraphics[width=.4\textwidth]{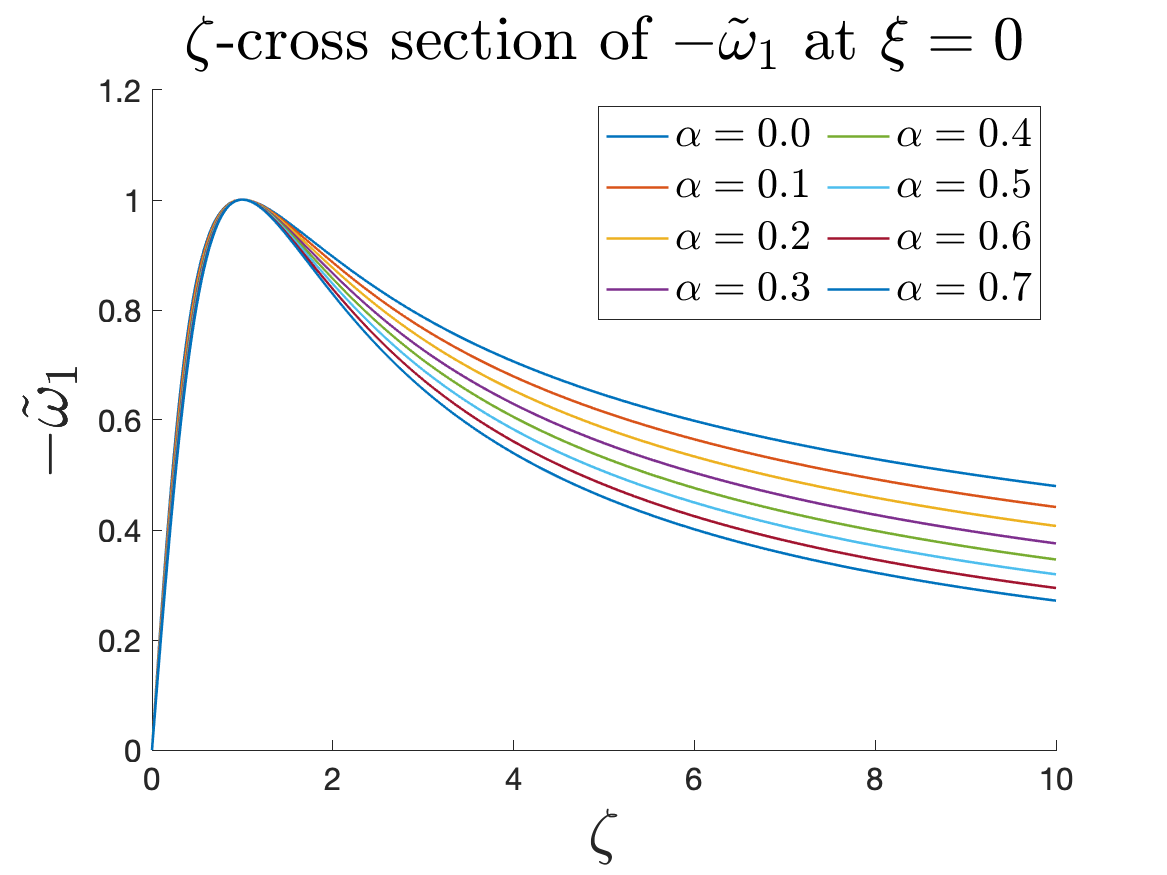}\\
\includegraphics[width=.4\textwidth]{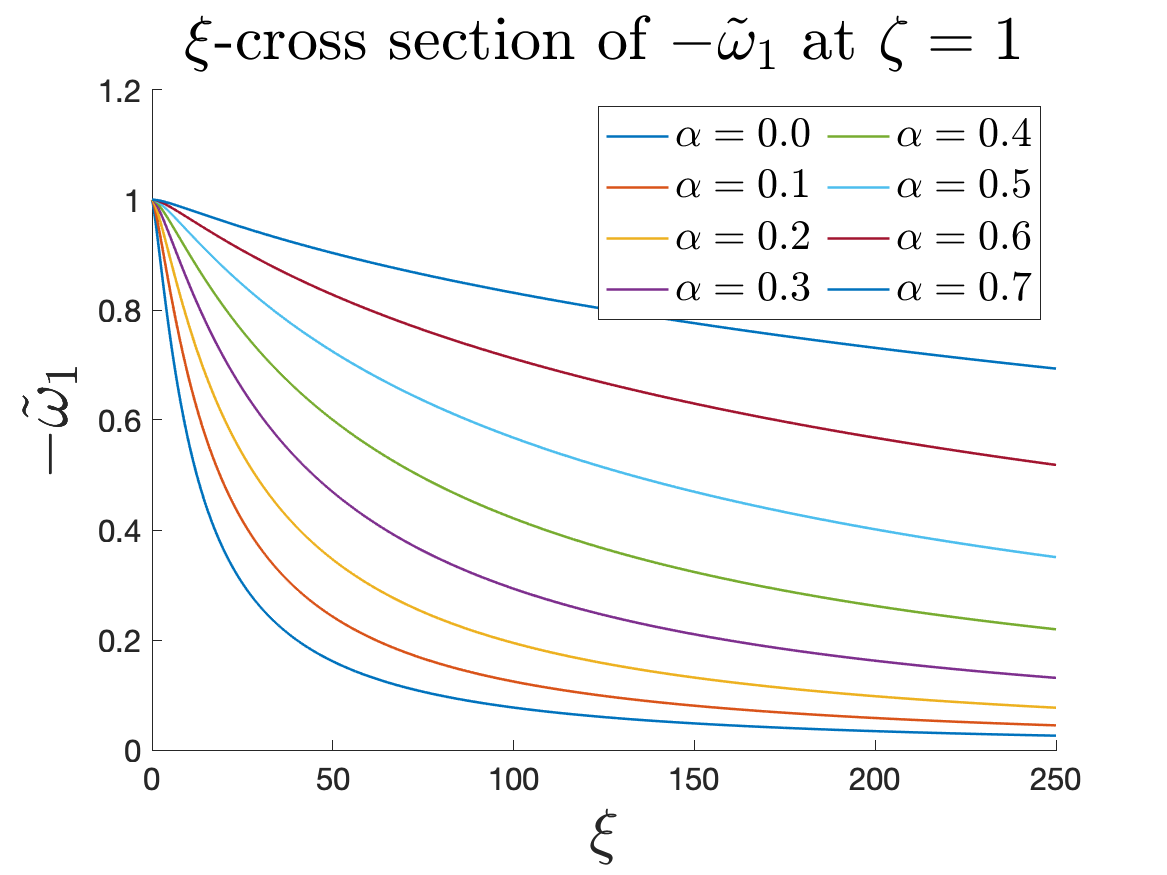}\includegraphics[width=.4\textwidth]{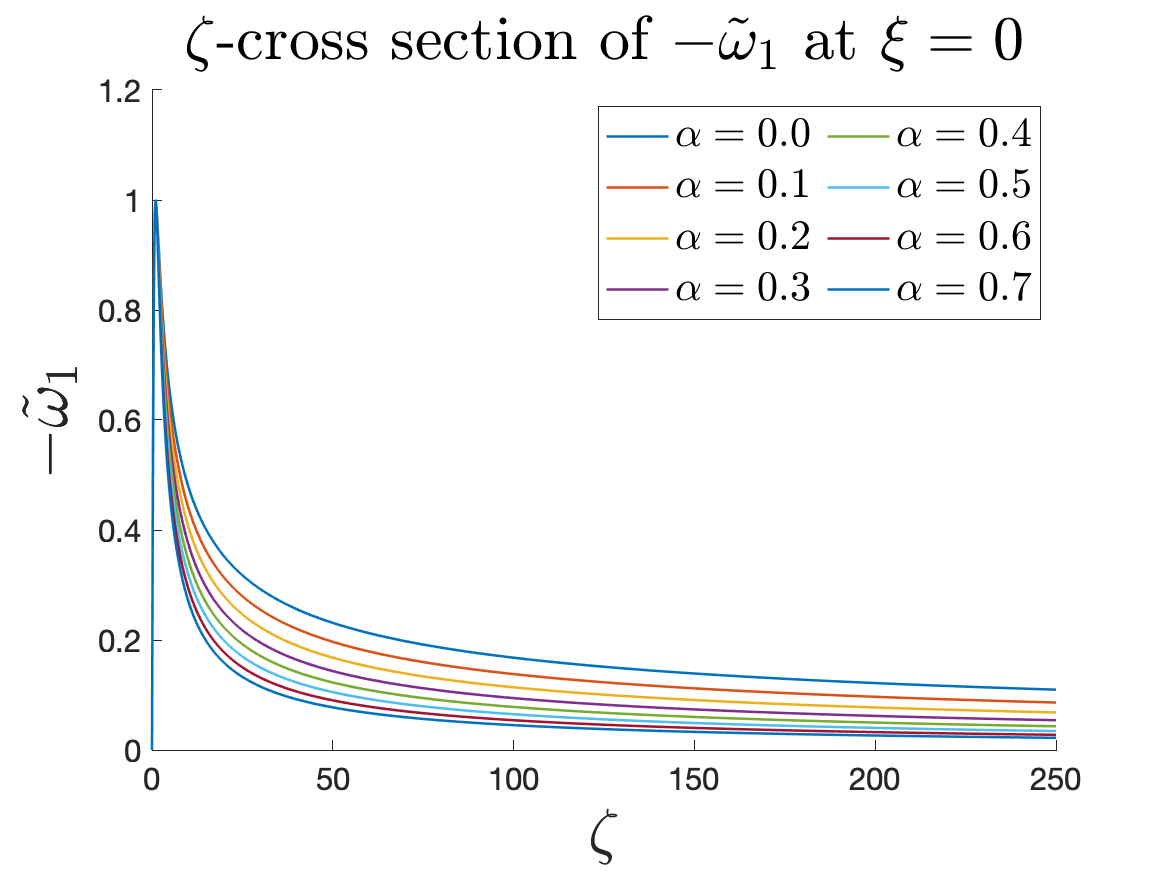}
\caption{Cross sections of steady states of $-\tilde{\omega}_1$ with different $\alpha$ in $\mathbb{R}^{10}$. Top row: on a local window. Bottom row: on a larger window.}\label{fig: steadystate_cross_section_alpha_n10}
\index{figures}
\end{figure}

\begin{figure}[hbt!]
\centering
\includegraphics[width=.4\textwidth]{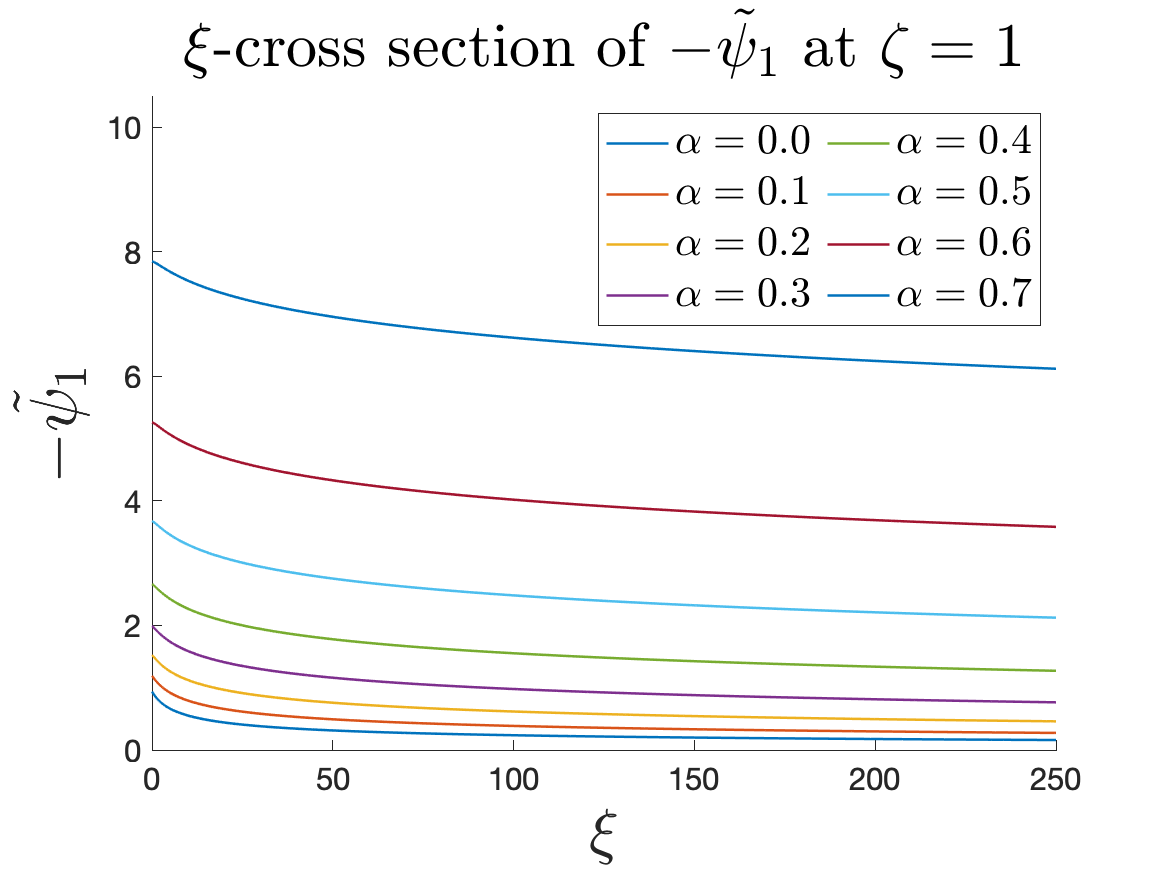}\includegraphics[width=.4\textwidth]{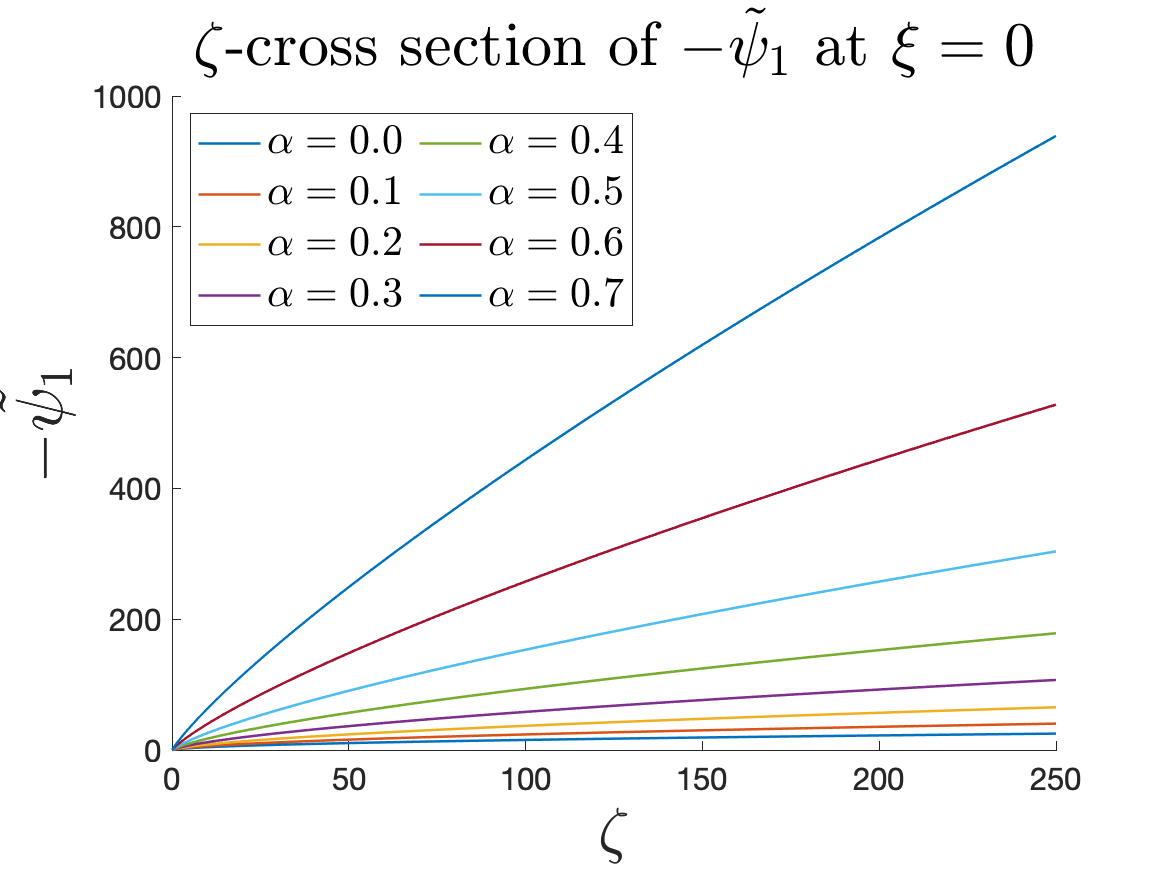}
\caption{Cross sections of steady states of $-\tilde{\psi}_1$ with different $\alpha$ in $\mathbb{R}^{10}$.}\label{fig: steadystate_cross_section_alpha_n10_psi1}
\index{figures}
\end{figure}

The steady states of $-\tilde{\omega}_1$ and $-\tilde{\psi}_1$ for $\alpha=0.1, 0.3, 0.5, 0.7$ when the dimension $n=10$ are plotted in Figure \ref{fig: steadystate_alpha_n10}. In Figures \ref{fig: steadystate_cross_section_alpha_n10} and \ref{fig: steadystate_cross_section_alpha_n10_psi1}, we provide the cross sections of the steady states $-\tilde{\omega}_1$ and $-\tilde{\psi}_1$ for different $\alpha$ when $n=10$. Similar to the $n=3$ case, we see that as $\alpha$ increases, the steady state becomes flatter in $\xi$. The cross section in $\zeta$ shows the different decay rates for different values of $\alpha$. It is also very interesting to see that the $\zeta$-cross section of $-\tilde{\psi}_1$ seems to be well approximated by a linear function of $\zeta$ in the near field.

In Table \ref{tab: cl_n10}, we showed how $c_l$ grows with $\alpha$ when $n=10$. Similar to the observation in Section \ref{sec: holder_exponent}, we see $c_l$ quickly increases with $\alpha$, and has the trend to go to infinity as $\alpha$ approaches $\alpha^*$. We remark that in our computation, we found $\alpha^*>0.72$. As in the 3D case in Section \ref{sec: holder_exponent}, we need to solve the dynamic rescaling formulation on a larger computational domain with higher accuracy and longer physical time to reach the steady states for larger $\alpha$. Based on our observation in the $n=3$ case and $n=10$ case, we conjecture that for our example, the critical value $\alpha^*=1-\frac{2}{n}$, which means our example would potentially support Conjecture 8 of \cite{drivas2022singularity} and extend it to the $n$-dimensional case.

\begin{table}[hbt!]
\centering
\caption{The scaling factor $c_l$ with different $\alpha$ in the $n$-D case, where $n=10$.}\label{tab: cl_n10}
\begin{tabular}{|c|c|c|c|c|c|c|c|c|}
\hline
\hspace{0.5cm}$\alpha$\hspace{0.5cm} & \hspace{0.25cm}$0.0$\hspace{0.25cm} & \hspace{0.25cm}$0.1$\hspace{0.25cm} & \hspace{0.25cm}$0.2$\hspace{0.25cm} & \hspace{0.25cm}$0.3$\hspace{0.25cm} & \hspace{0.25cm}$0.4$\hspace{0.25cm} & \hspace{0.25cm}$0.5$\hspace{0.25cm} & \hspace{0.25cm}$0.6$\hspace{0.25cm} & \hspace{0.25cm}$0.7$\hspace{0.25cm} \\
\hline
$c_l$ & $2.155$ & $2.432$ & $2.811$ & $3.363$ & $4.244$ & $5.897$ & $10.16$ & $47.14$ \\
\hline
\end{tabular}
\index{tables}
\end{table}

We next study how the dimension $n$ influences the finite-time blow-up. We fix $\alpha=0.1$ and try different choices of dimensions $n=3, 4, 5, 6, 8, 10$ using the same initial data
\begin{align*}
    \omega^\circ_1=\frac{-12000\left(1-r^2\right)^{18}\sin(2\pi z)}{1+12.5\sin^2(\pi z)}\,.
\end{align*}
We note that this initial data will lead to the same steady state phenomenon as the initial data in \eqref{eq: inital_data}, as we will study in detail in Section \ref{sec: sensitivity}. We choose to report results in the $n$-dimensional case on this initial data because it takes less time to develop a potential blow-up, since the initial vorticity concentrates more near the origin, which helps reduce our computational cost.

In Table \ref{tab: stats_dim}, we report the estimated blow-up times $T$ and scaling factors $c_l$ for different dimensions. It is not surprising that the blow-up time of the same initial data is shorter for the higher dimensional case, because the vortex stretching term has a larger amplification coefficient. However, the scaling factor $c_l$ is smaller for larger $n$. Intuitively, the velocity component $u^z$ seems to be stronger with larger $n$. This phenomenon suggests that in the high dimensional case, the dimension-related term $-\frac{n}{r}\partial_r$ in the Poisson equation controls $\psi_1$ and therefore weakens the collapsing speed of the solution. We also remark that the decay of $c_l$ with $n$ significantly slows down in Table \ref{tab: stats_dim}. It is tempting to speculate if there is a limit of $c_l$ as $n$ approaches infinity.

\begin{figure}[hbt!]
\centering
\includegraphics[width=.4\textwidth]{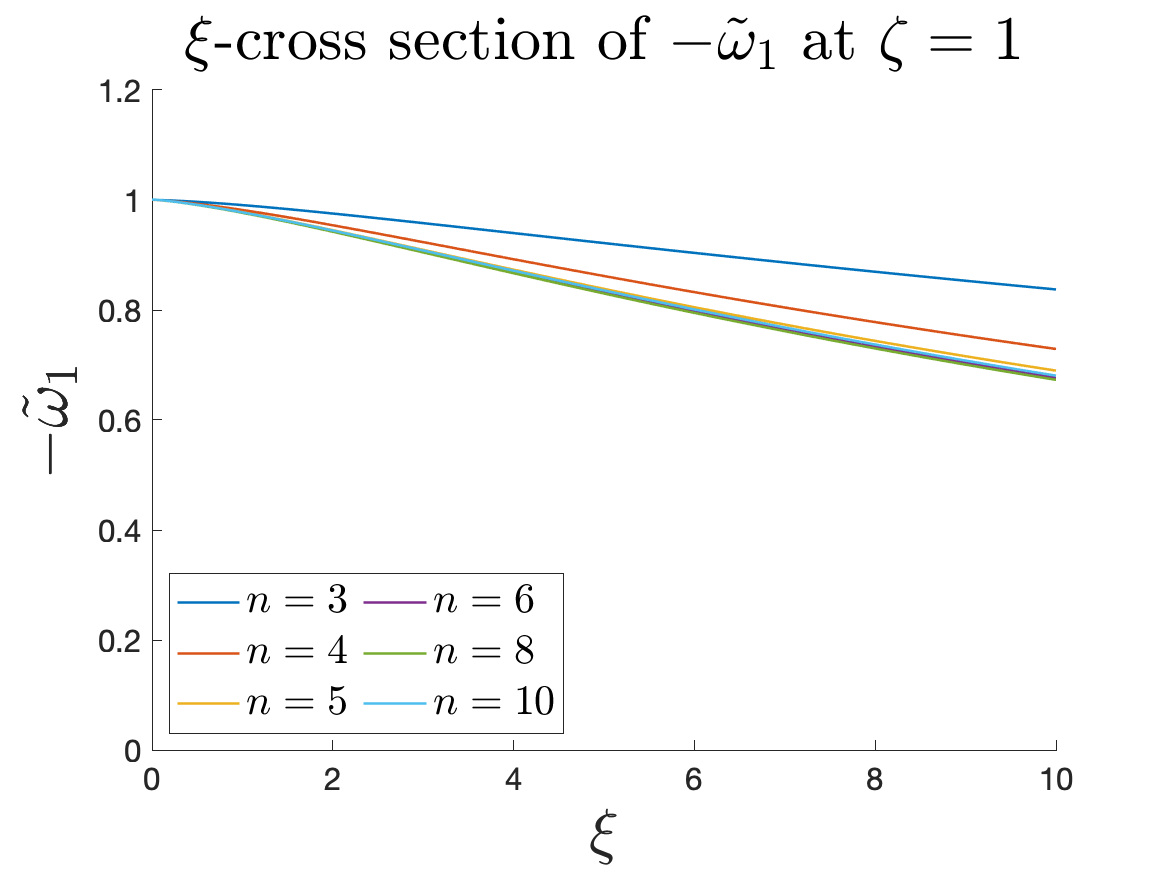}\includegraphics[width=.4\textwidth]{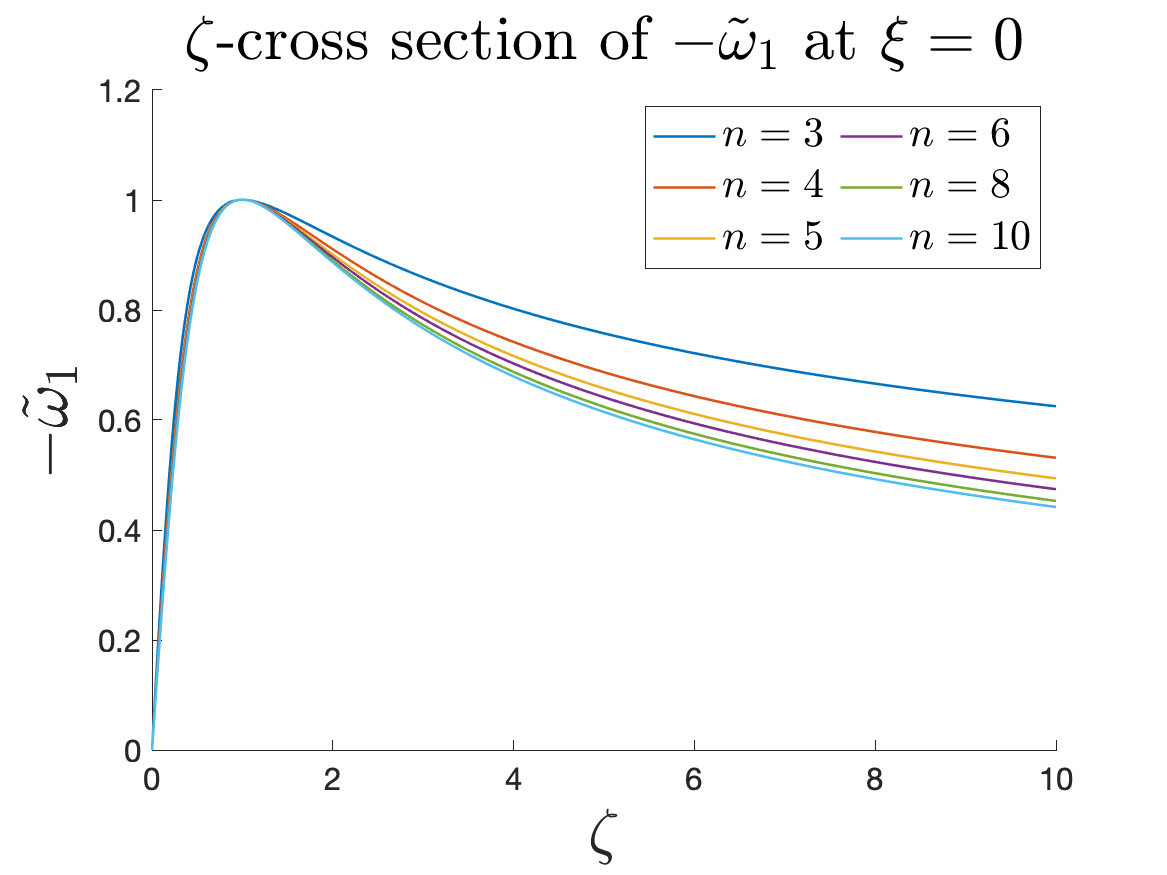}\\
\includegraphics[width=.4\textwidth]{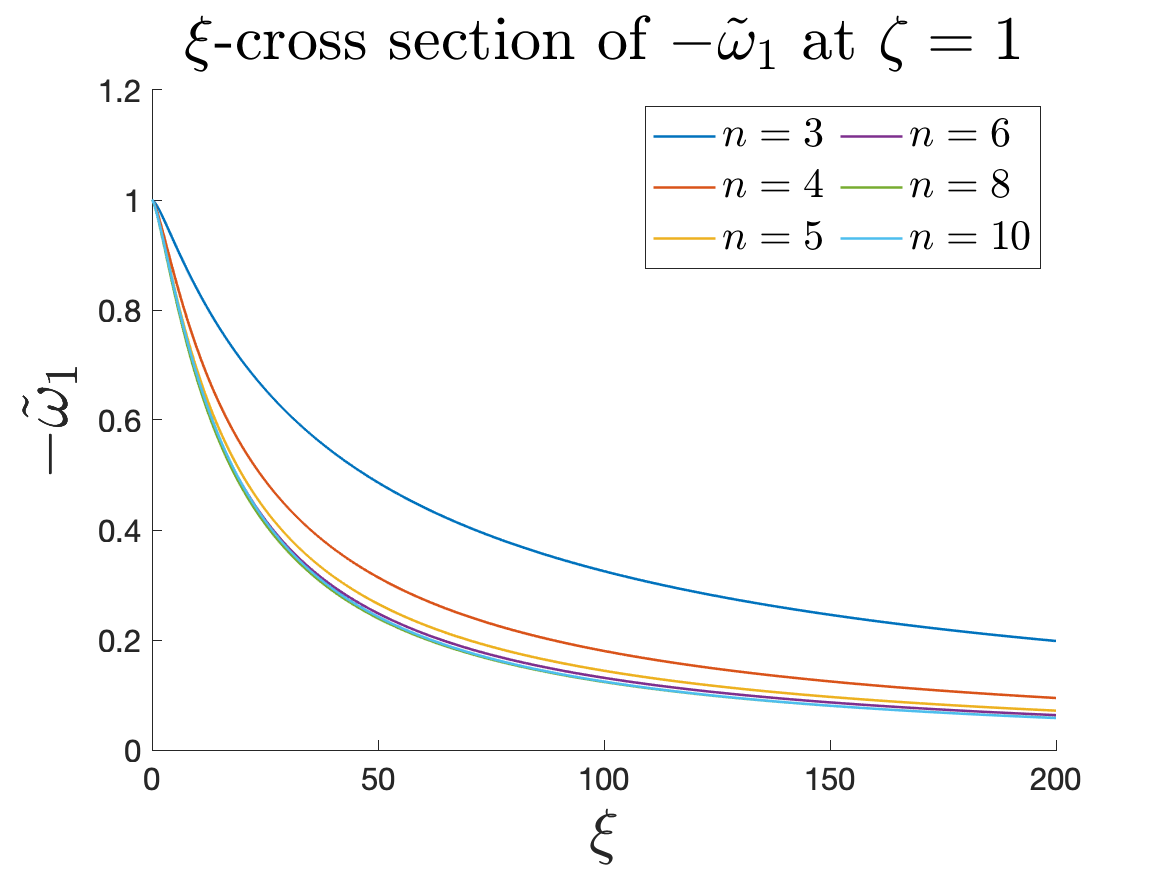}\includegraphics[width=.4\textwidth]{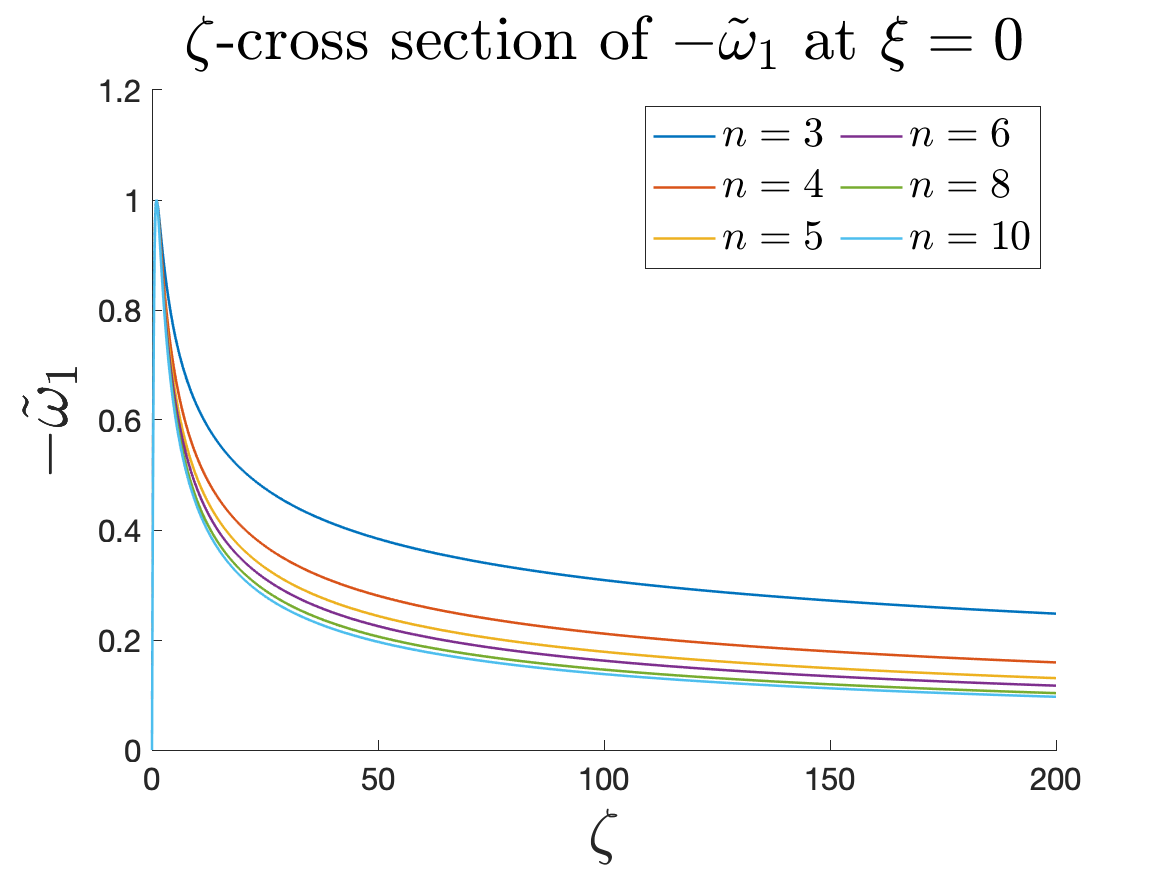}
\caption{Cross sections of steady states of $-\tilde{\omega}_1$ with different $n$ with $\alpha=0.1$. Top row: on a local window. Bottom row: on a larger window.}\label{fig: steadystate_cross_section_dim}
\index{figures}
\end{figure}

\begin{figure}[hbt!]
\centering
\includegraphics[width=.4\textwidth]{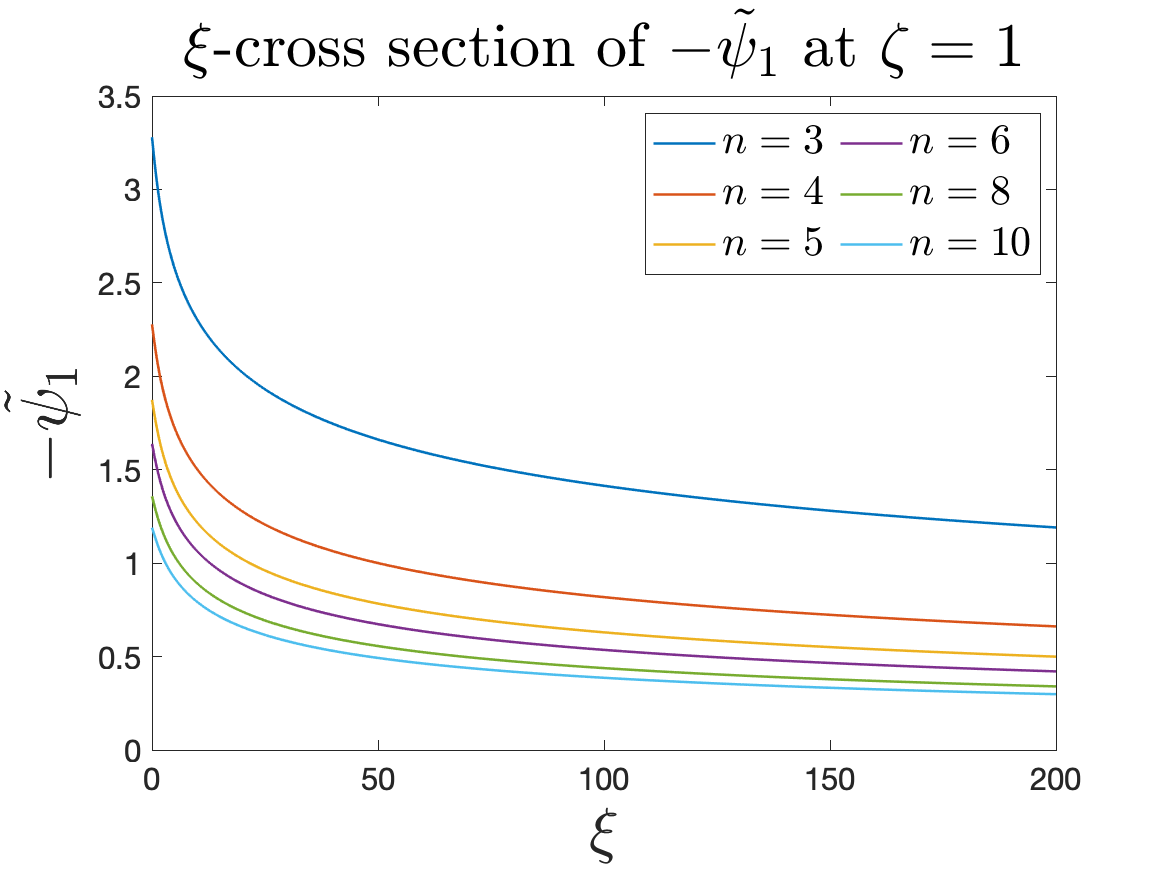}\includegraphics[width=.4\textwidth]{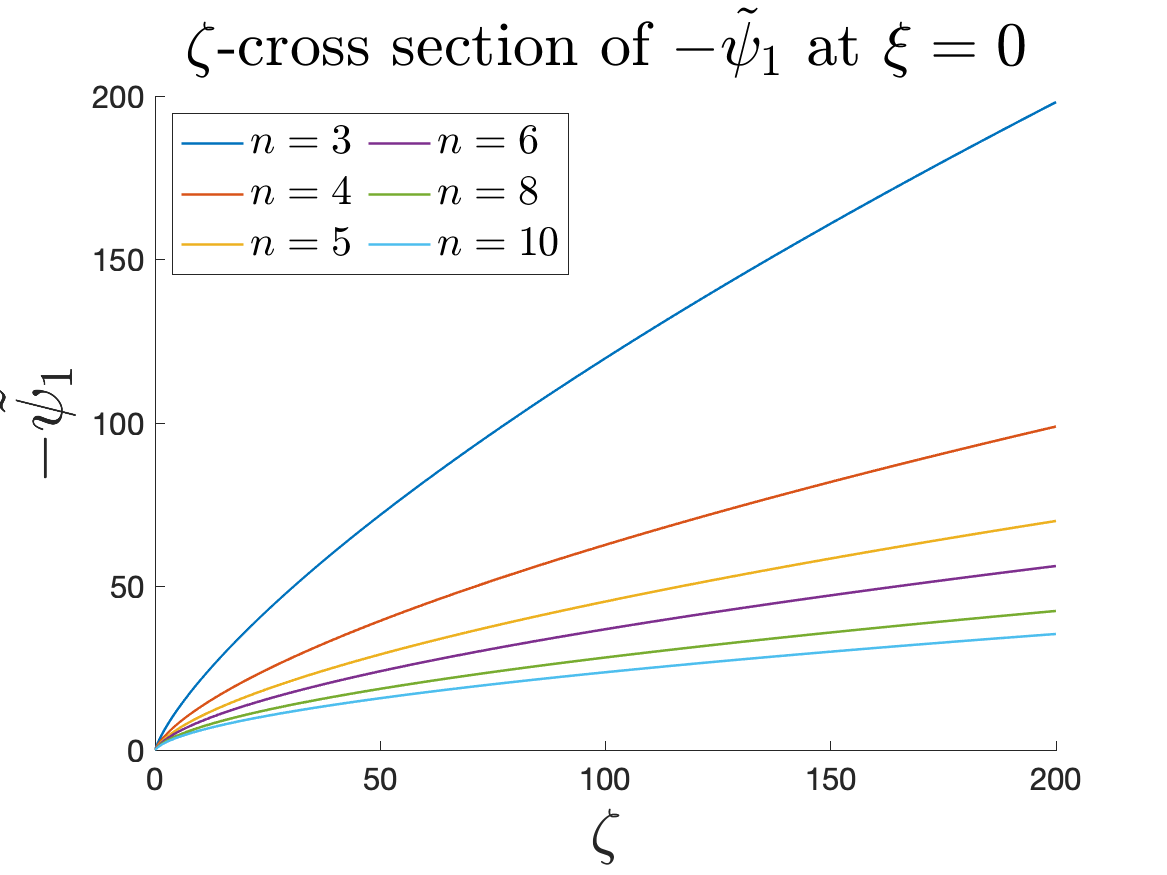}
\caption{Cross sections of steady states of $-\tilde{\psi}_1$ with different $n$ with $\alpha=0.1$.}\label{fig: steadystate_cross_section_dim_psi1}
\index{figures}
\end{figure}

In Figure \ref{fig: steadystate_cross_section_dim} and \ref{fig: steadystate_cross_section_dim_psi1}, we provide the cross sections of the steady states of $-\tilde{\omega}_1$ and $-\tilde{\psi}_1$ from the dynamic rescaling formulation. We observe that the cross sections change with $n$. But as $n$ becomes larger than $5$, the difference quickly narrows down. This would give more evidence that there is some non-trivial limit as $n$ goes to infinity. It would be very interesting to further explore this infinite dimension limit in the future.

\begin{table}[hbt!]
\centering
\caption{Estimated blow-up times $T$ and scaling factors $c_l$ in different dimensions with $\alpha=0.1$.}\label{tab: stats_dim}
\begin{tabular}{|c|c|c|c|c|c|c|}
\hline
\hspace{0.5cm}$n$\hspace{0.5cm} & \hspace{0.5cm}$3$\hspace{0.5cm} & \hspace{0.5cm}$4$\hspace{0.5cm} & \hspace{0.5cm}$5$\hspace{0.5cm} & \hspace{0.5cm}$6$\hspace{0.5cm} & \hspace{0.5cm}$8$\hspace{0.5cm} & \hspace{0.5cm}$10$\hspace{0.5cm} \\
\hline
$10^4\times T$ & $4.183$ & $3.589$ & $2.877$ & $2.498$ & $2.089$ & $1.866$ \\
\hline
$c_l$ & $4.549$ & $3.218$ & $2.851$ & $2.680$ & $2.514$ & $2.432$ \\
\hline
\end{tabular}
\index{tables}
\end{table}

\section{Sensitivity of the potential blow-up to initial data}
\label{sec: sensitivity}

We study the sensitivity of the potential self-similar blow-up to initial data. In addition to the initial data \eqref{eq: inital_data}, we consider the following cases,
\begin{equation}
\label{eq: inital_data_2}
    \begin{aligned}
    \omega^{\circ,1}_1&=\frac{-12000\left(1-r^2\right)^{18}\sin(2\pi z)}{1+12.5\sin^2(\pi z)},\\
    \omega^{\circ,2}_1&=-6000\cos\left(\frac{\pi r}{2}\right)\sin(2\pi z)\left(2+\exp\left(-r^2\sin^2(\pi z)\right)\right),\\
    \omega^{\circ,3}_1&=\frac{-12000\left(1-r^2\right)^{18}\sin(2\pi z)^3}{1+12.5\sin^2(\pi z)}.
    \end{aligned}
\end{equation}

\begin{figure}[hbt!]
\centering
\includegraphics[width=.32\textwidth]{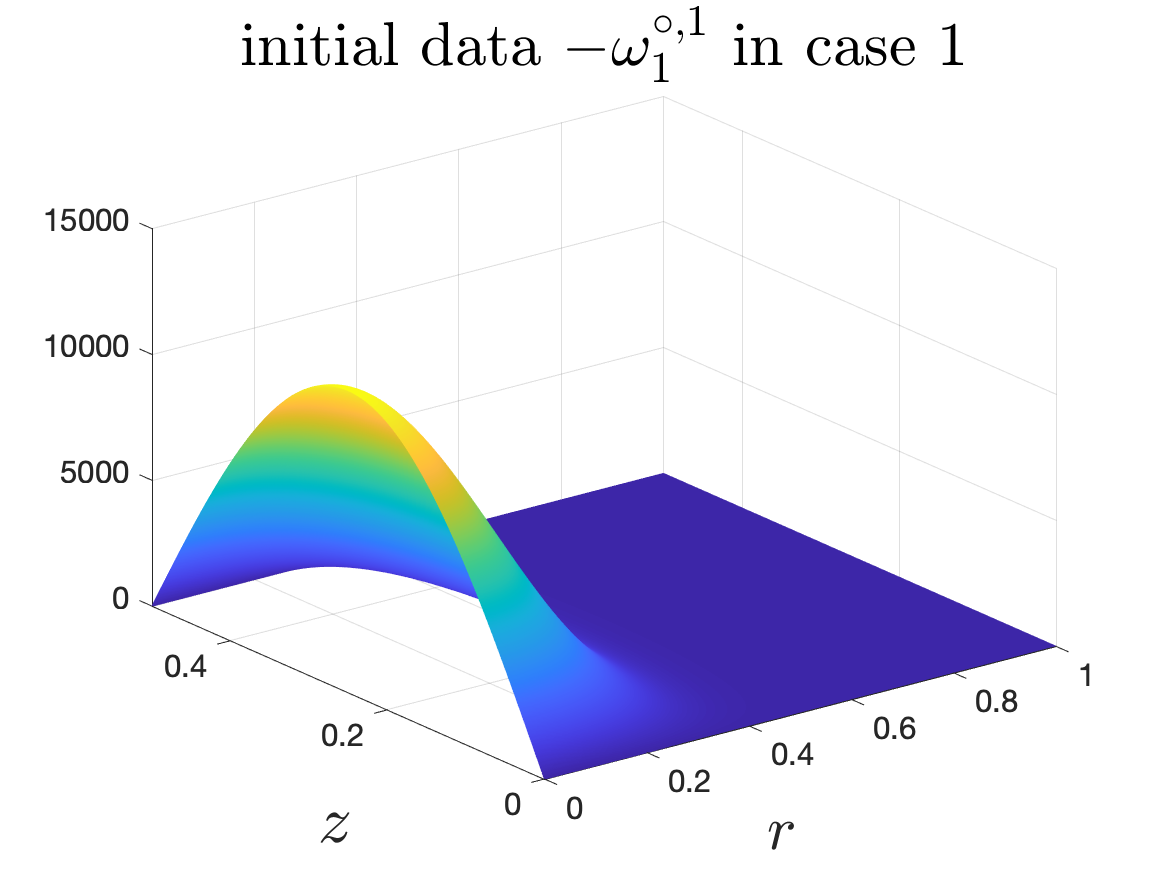}\includegraphics[width=.32\textwidth]{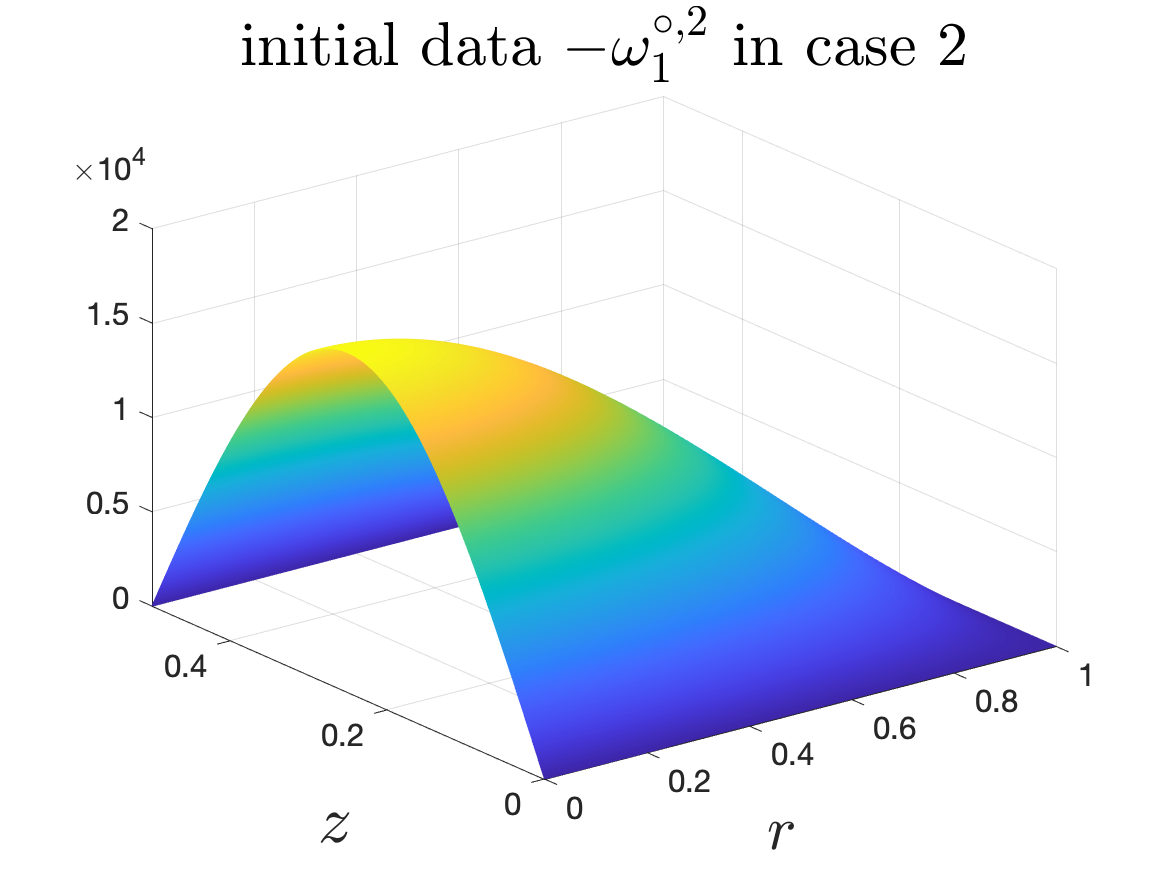}\includegraphics[width=.32\textwidth]{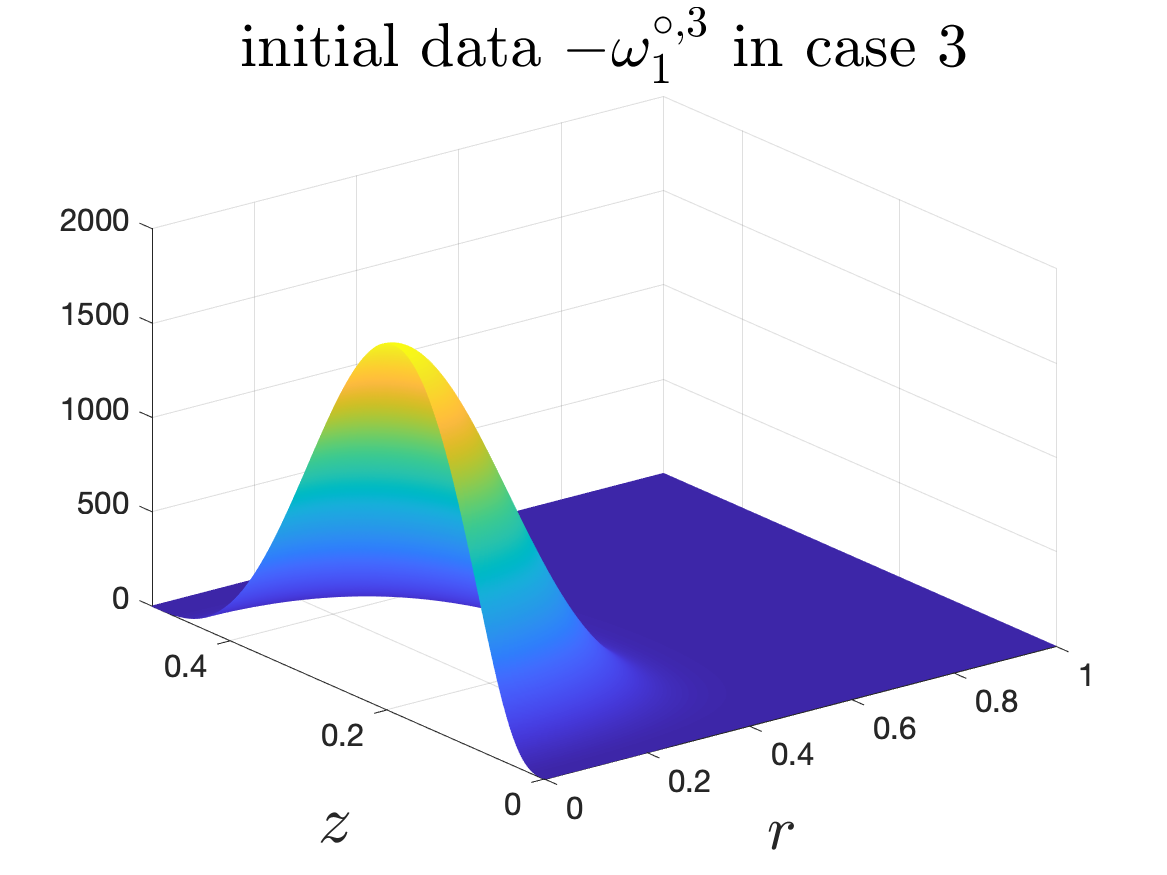}
\caption{Profiles of the initial data in all three cases.}\label{fig: init_data_cases}
\end{figure}

We show the profiles of these three initial data in Figure \ref{fig: init_data_cases}. In Case 1, $\omega^{\circ,1}_1$ is a perturbation of $\omega^\circ_1$ by setting cosine in the denominator to be sine, which is adopted as the initial data for the $n$-D case in Section \ref{sec: dimension}. In Case 2, $\omega^{\circ,2}_1$ has a decay rate in $r$ slower than $\left(1-r^2\right)^{18}$, and is no longer a simple tensor product of $r$ and $z$. In Case 3, $\omega^{\circ,3}_1$ has an improved regularity in $\rho=\sqrt{r^2+z^2}$ near the origin. Indeed, we have, with $\omega_1(r,z,0)=\omega^{\circ,3}_1(r,z)$,
$$\omega^\theta(r,z,0)=r^\alpha\omega^{\circ,3}_1(r,z)\sim r^\alpha z^3=\rho^{3+\alpha}\cos^\alpha\theta\sin^3\theta.$$
While for the original choice of the initial data \eqref{eq: inital_data}, $\omega^\theta(r,z,0)\sim\rho^{1+\alpha}\cos^\alpha\theta\sin\theta$.

For all three cases, we only solve the 3D axisymmetric Euler equations with $\alpha=0.3$, due to the limited computational resources. As shown in Table \ref{tab: cl_n3}, for our original initial data, $c_l=112.8$ is very large, which suggests that our choice of $\alpha$ is very close to the borderline between the blow-up and non-blow-up. If the blow-up profile of the above initial data agrees with our original initial data well, we then have good confidence that they should have the same behavior for other settings of $\alpha$.

We solve the 3D axisymmetric Euler equations with the above initial data by first using the adaptive mesh method to get close enough to the potential blow-up time, and then using the dynamic rescaling method to capture the potential self-similar solution.

\begin{figure}[hbt!]
\centering
\includegraphics[width=.4\textwidth]{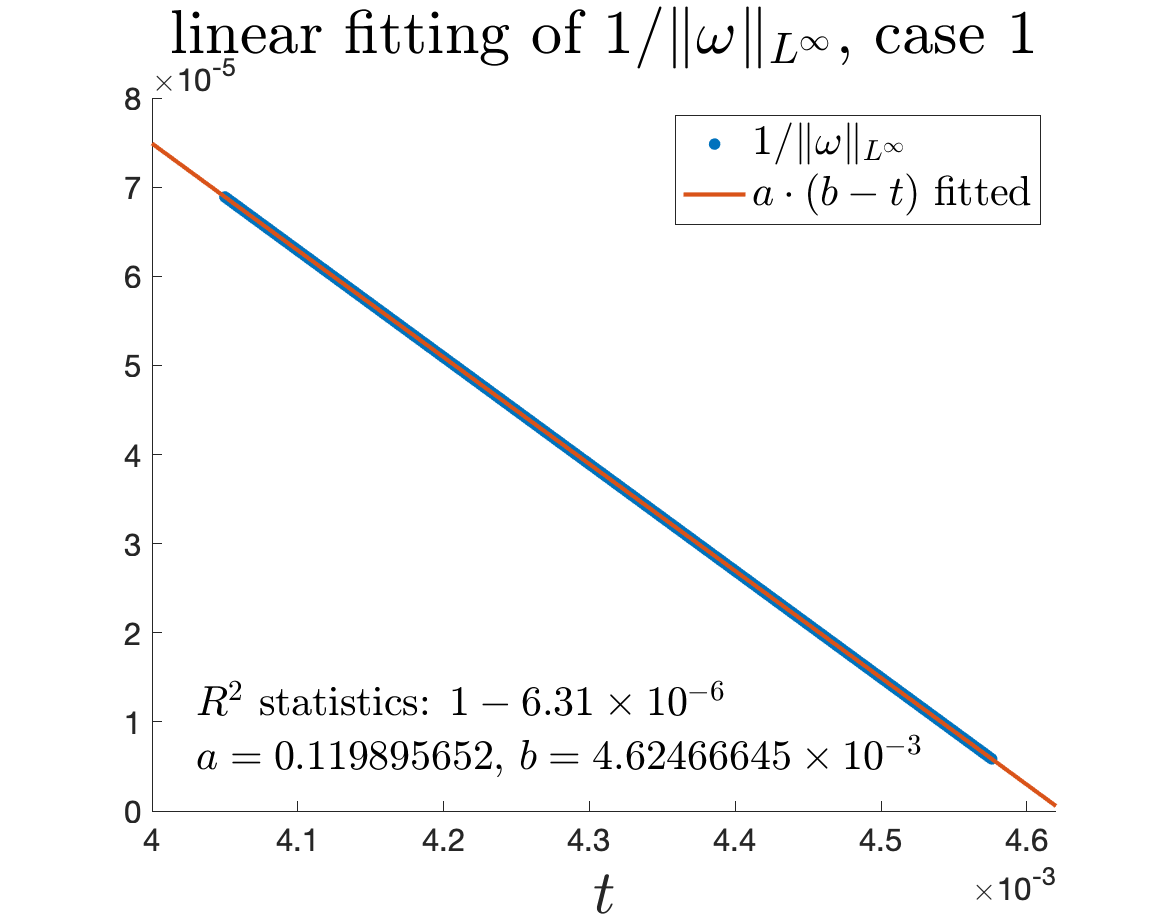}\includegraphics[width=.4\textwidth]{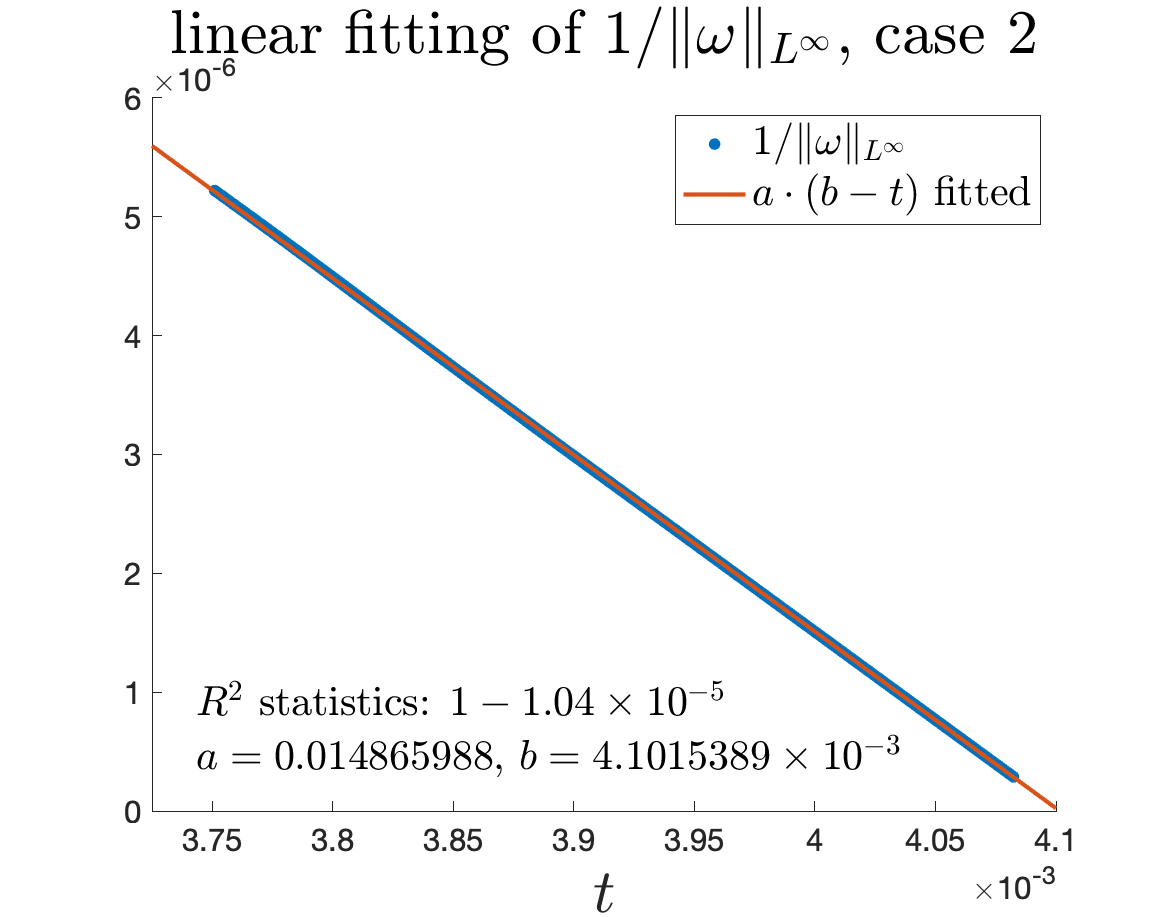}
\caption{Fitting of $1/\|\omega\|_{L^\infty}$ with time $t$ in the first and second cases.}\label{fig: fit_12}
\end{figure}

\begin{figure}[hbt!]
\centering
\includegraphics[width=.4\textwidth]{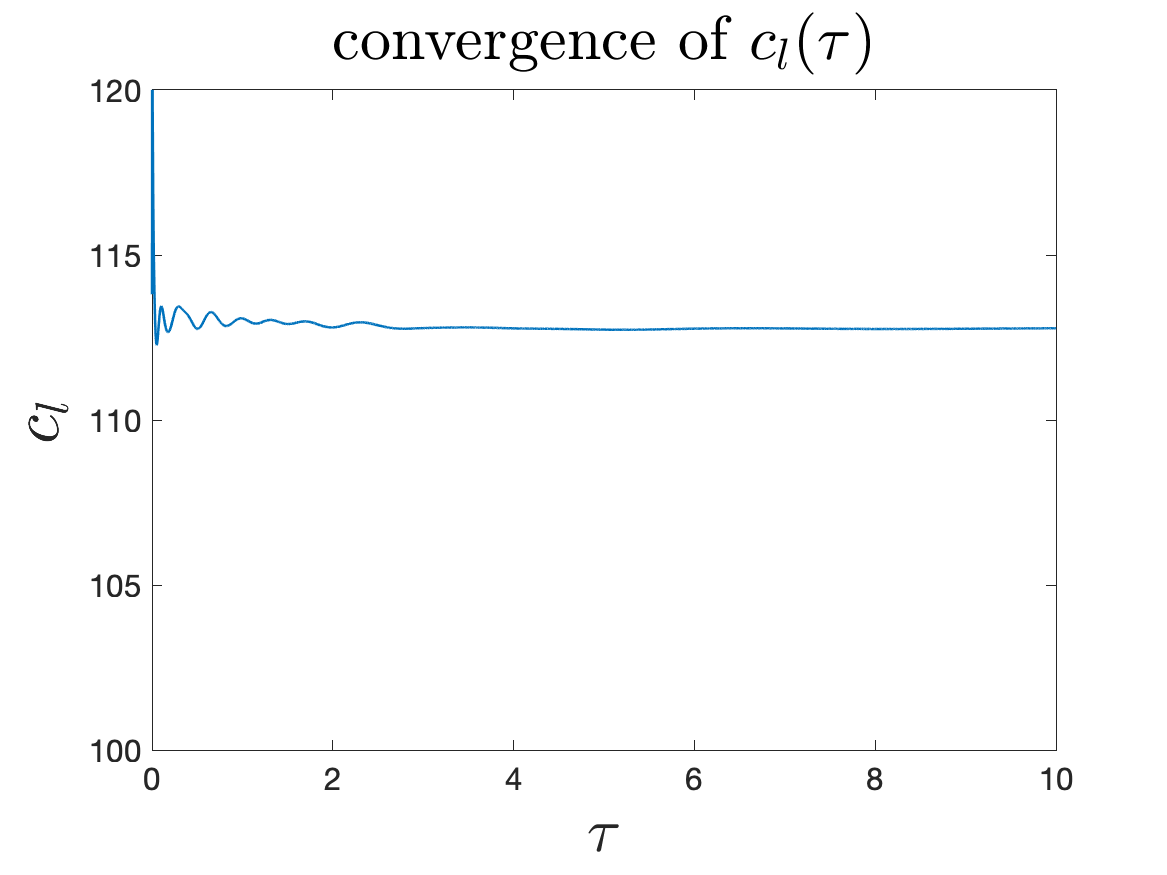}\includegraphics[width=.4\textwidth]{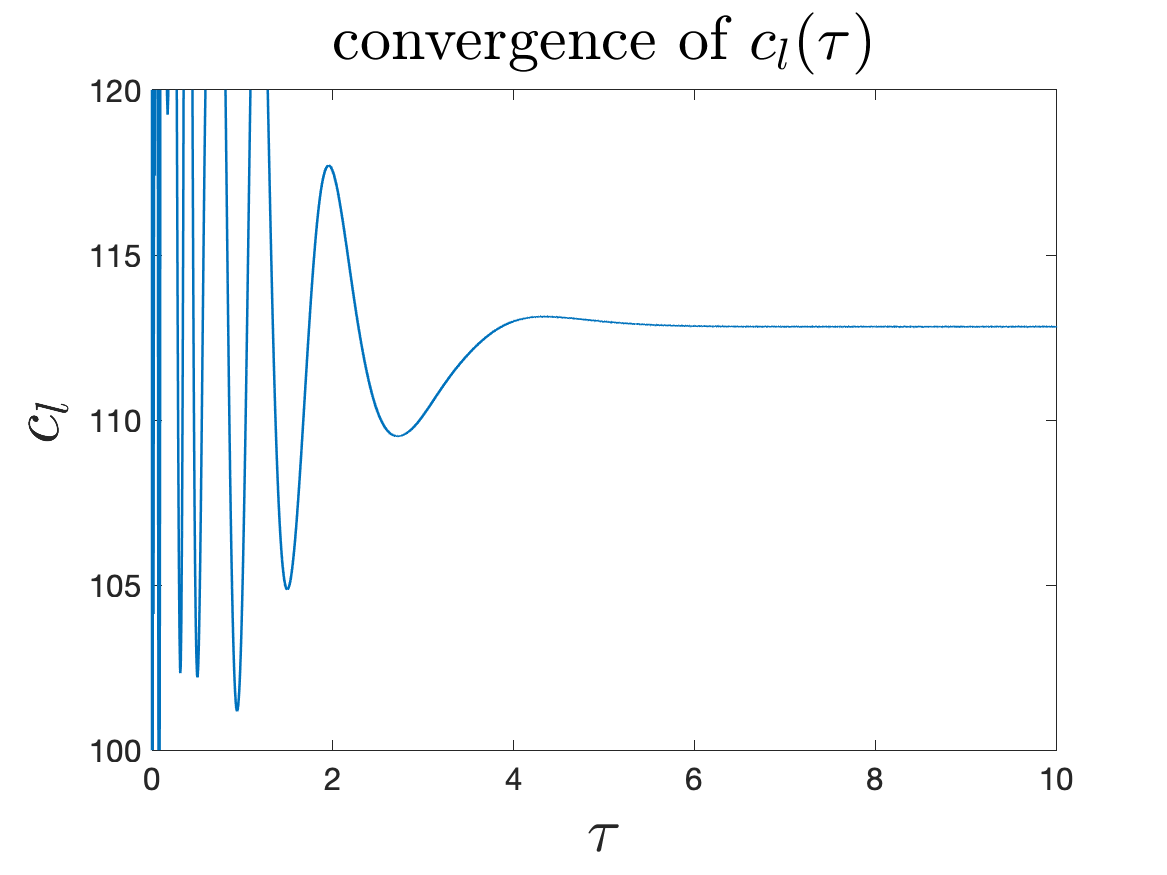}
\caption{Curves of the scaling factor $c_l$ in the first and second cases.}\label{fig: cl_12}
\end{figure}

\begin{figure}[hbt!]
\centering
\includegraphics[width=.4\textwidth]{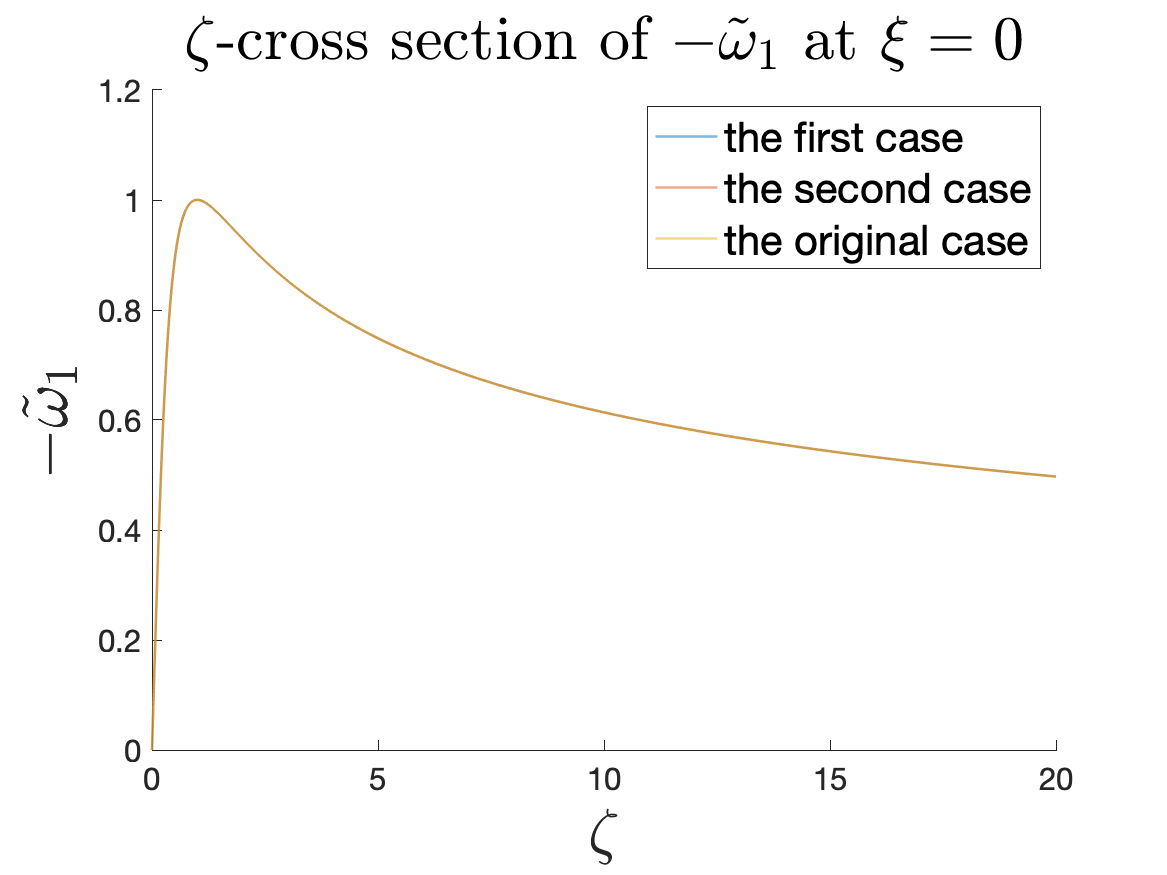}\includegraphics[width=.4\textwidth]{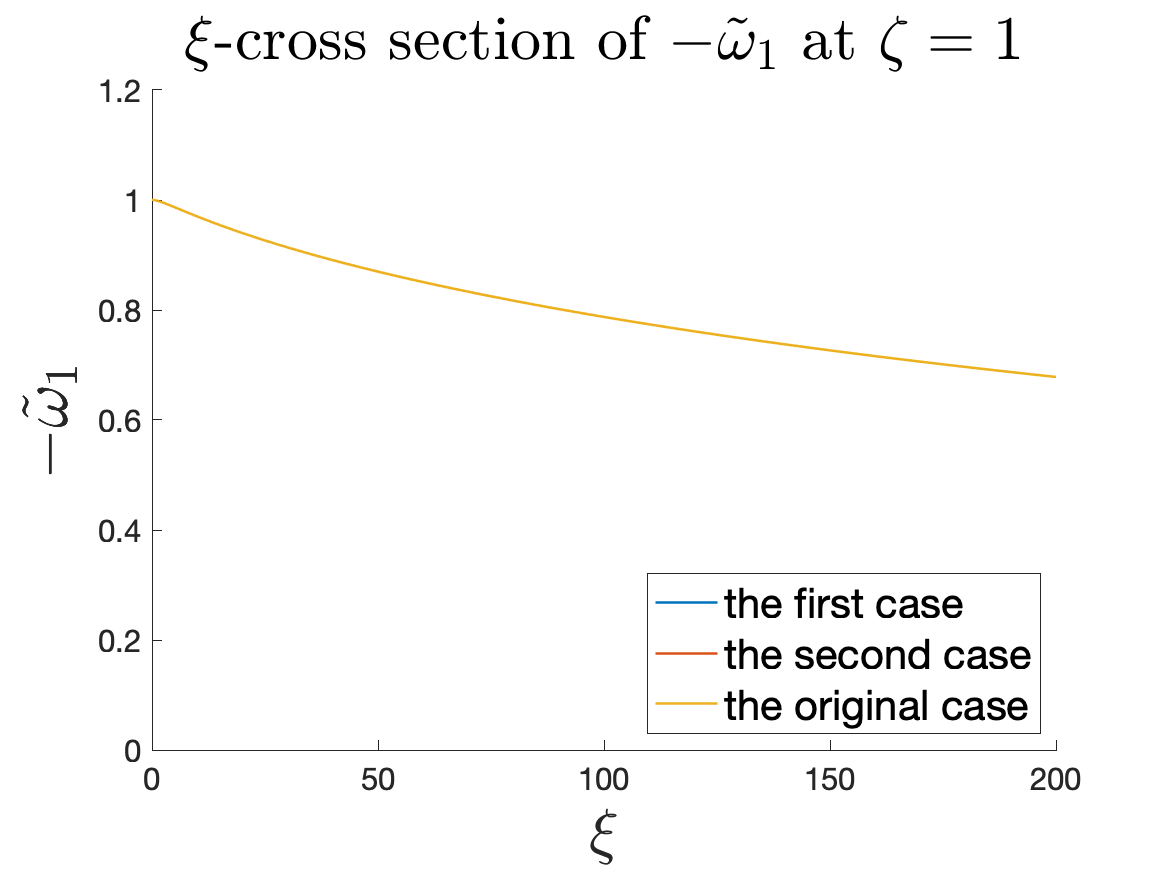}
\caption{Cross sections of the steady states of $-\tilde{\omega}_1$ in the first and second cases.}\label{fig: cross_section_12}
\end{figure}

\subsection{The first two cases: same regularity near the origin}

For the first and second cases, we show the fitting of $1/\|\omega\|_{L^\infty}$ with time $t$ in Figure \ref{fig: fit_12}, and the curve of the scaling factor $c_l$ in Figure \ref{fig: cl_12}. We can see that in both cases, $\|\omega\|_{L^\infty}$ scales like $1/(T-t)$, which implies a finite-time blow-up. Moreover, $c_l$ converges to $112.8$, matching the value of $c_l$ we obtained using the original initial data well. In Figure \ref{fig: cross_section_12}, we show the cross sections of the steady state of $-\tilde{\omega}_1$ in comparison with the result obtained using the original initial data. There is no visible difference between the three steady states presented. In fact, even on the whole computational domain $\mathcal{D}^\prime=\left\{(\xi,\zeta): 0\leq\xi\leq 1\times10^5, 0\leq\zeta\leq 5\times10^4\right\}$ in the dynamic rescaling computation, the steady states in the first and second cases only differ by $3.13\times10^{-10}$ and $5.29\times10^{-10}$ respectively from the steady state using our original initial data $\omega^{\circ}_1$ in the relative sup-norm.

\begin{figure}[hbt!]
\centering
\includegraphics[width=.4\textwidth]{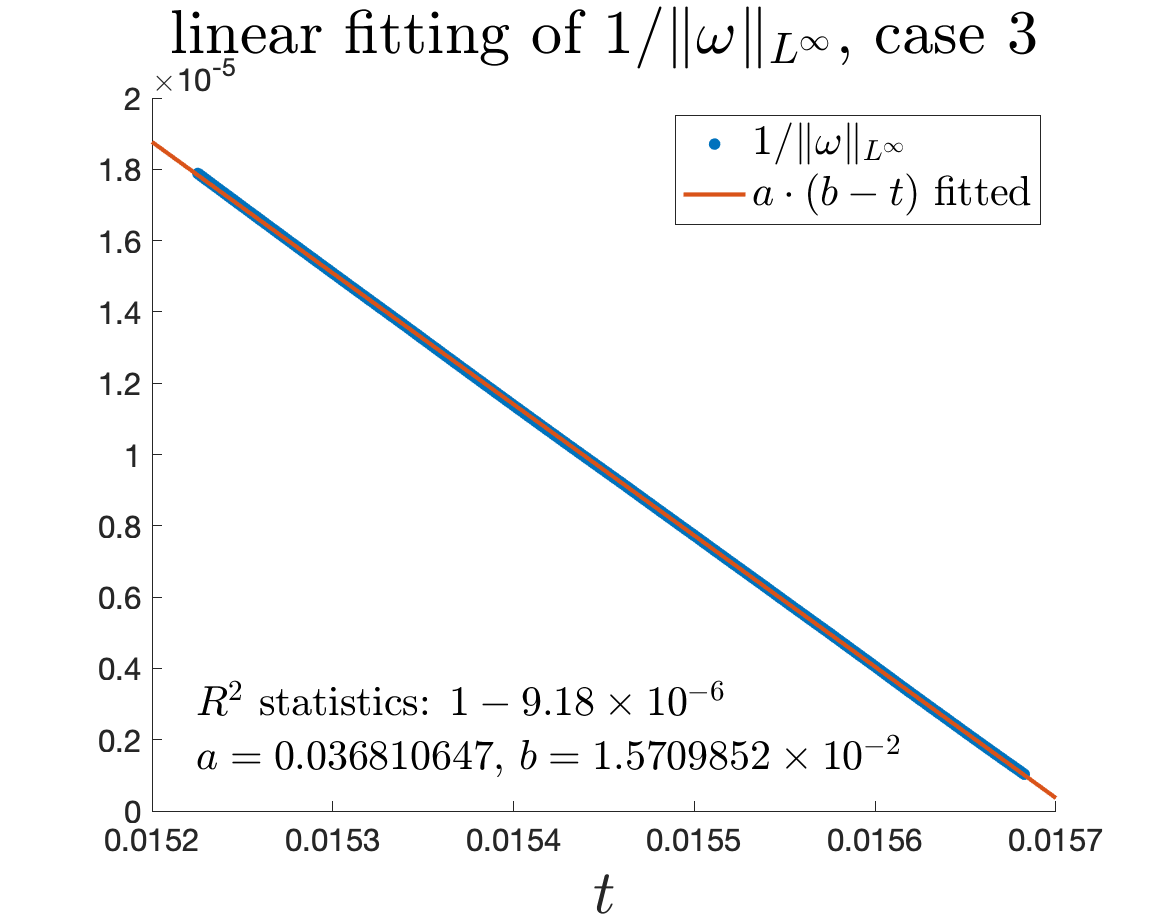}\includegraphics[width=.4\textwidth]{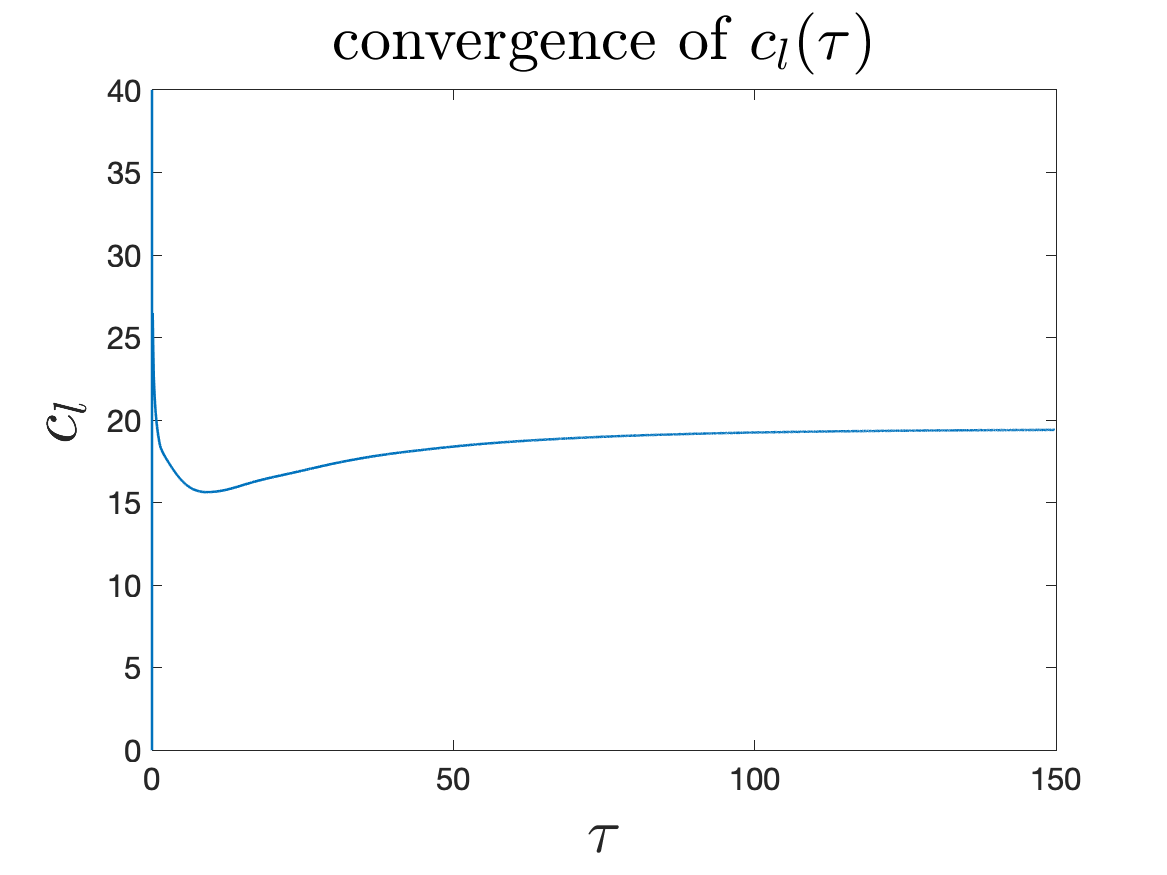}
\caption{Fitting of $1/\|\omega\|_{L^\infty}$ and curve of the scaling factor $c_l$ in the third case.}\label{fig: cl_3}
\end{figure}

\begin{figure}[hbt!]
\centering
\includegraphics[width=.4\textwidth]{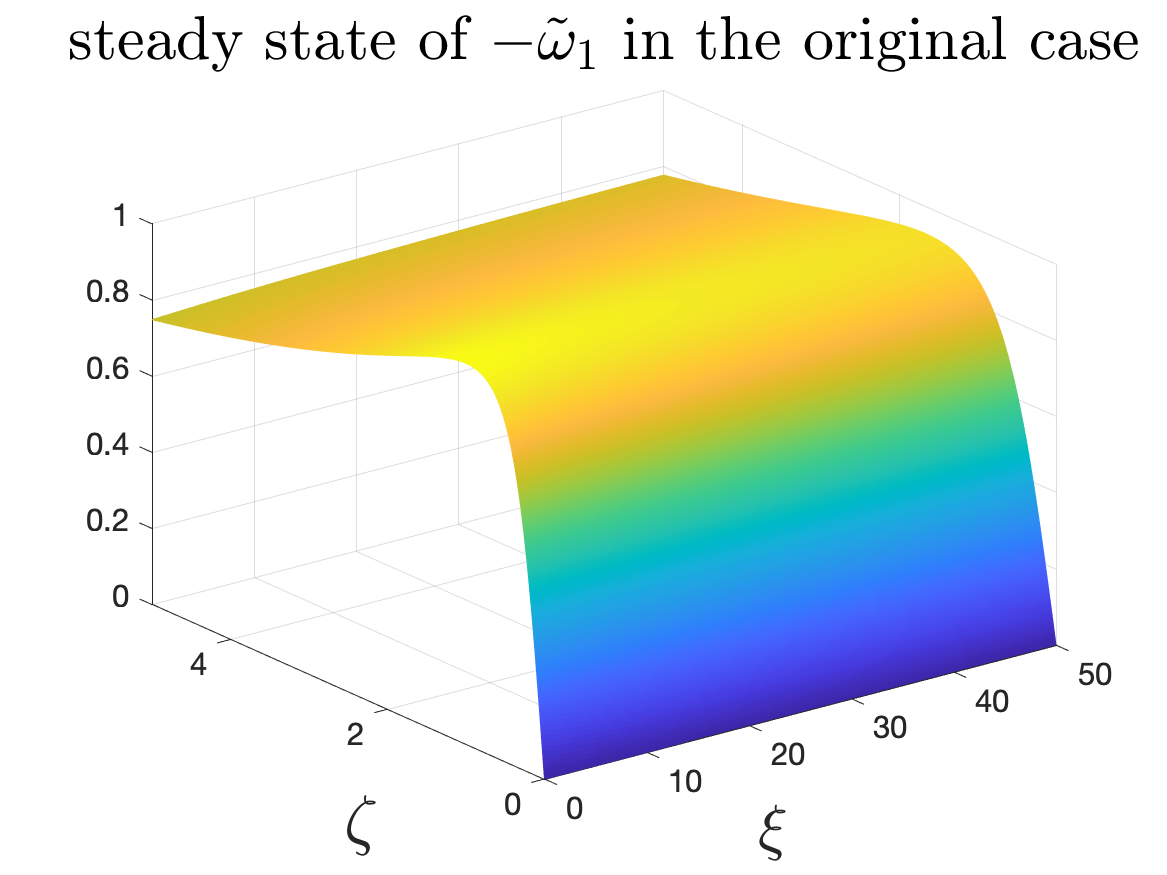}\includegraphics[width=.4\textwidth]{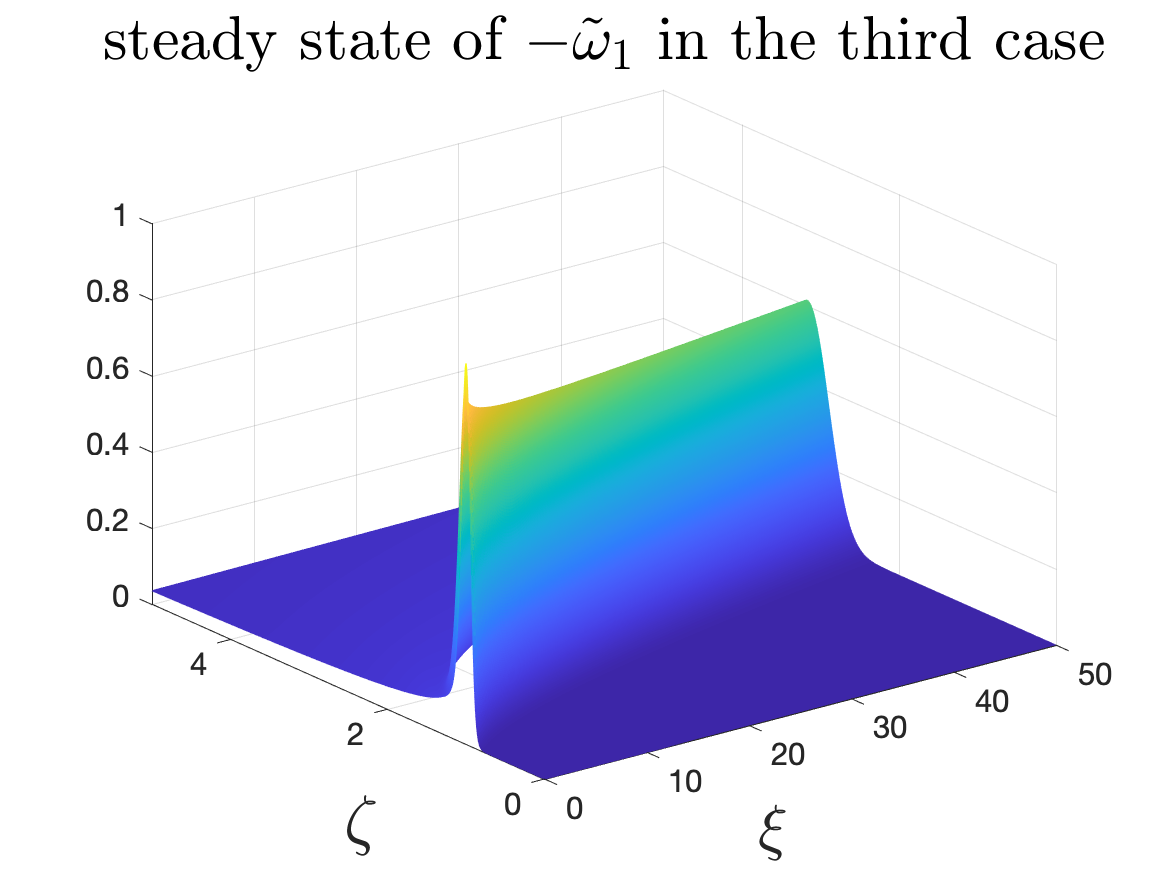}
\includegraphics[width=.4\textwidth]{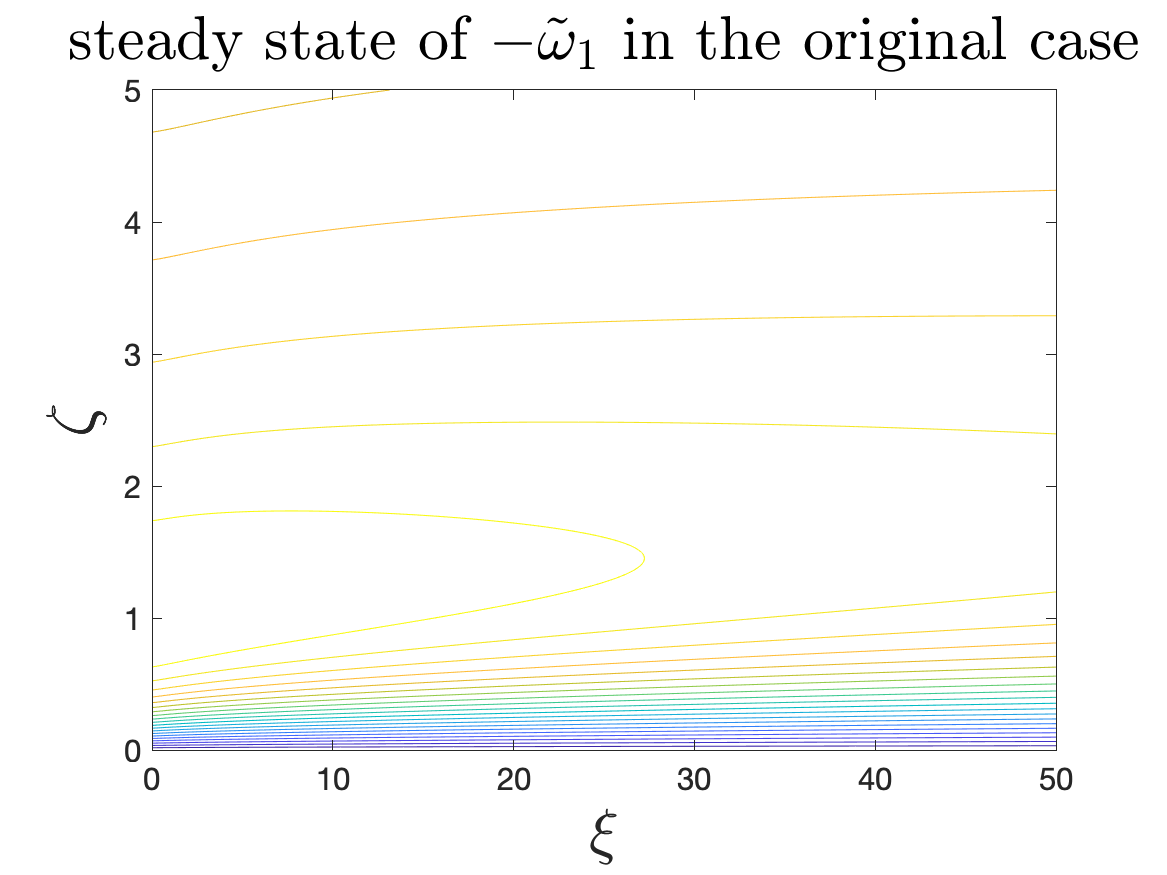}\includegraphics[width=.4\textwidth]{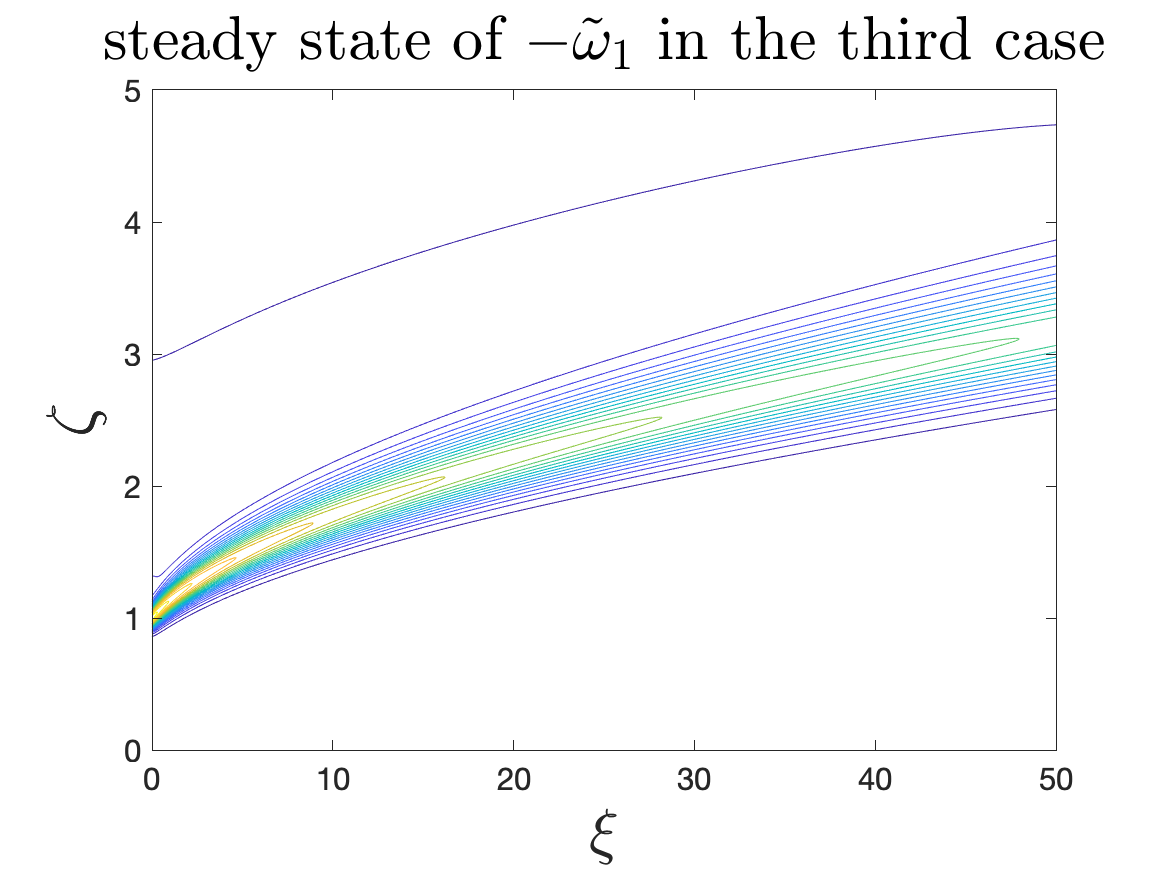}
\caption{Profiles and contours of the steady states of $-\tilde{\omega}_1$ in the original and third cases.}\label{fig: profile_3}
\end{figure}

\subsection{The third case: improved regularity near the origin}
\label{sec: sensitivity_third}

For the third case, the fitting of $1/\|\omega\|_{L^\infty}$ and the curve of the scaling factor $c_l$ is shown in Figure \ref{fig: cl_3}. We observe that $1/\|\omega\|_{L^\infty}$ has a good linear fitting with time, suggesting a finite-time blow-up. However, $c_l$ converges to $19.44$ which is clearly different from $112.8$, suggesting that there might be a new blow-up mechanism. In Figure \ref{fig: profile_3}, we compare the steady states of $\omega^{\circ,3}_1$ and $\omega^{\circ}_1$ in the 3D profiles and the 2D contours. The steady state of $\omega^{\circ,3}_1$ has a slower change near $z=0$. This might be caused by the smoothness of the initial data near $z=0$, because we have $\omega^{\circ,3}_1\sim r^\alpha z^3$, in contrast to $\omega^{\circ}_1\sim r^\alpha z$ near $(r,z)=(0,0)$. The steady state of the third case develops a channel-like structure that is not parallel to either axis.

The new blow-up scenario in the third case provides some support of Conjecture 9 of \cite{drivas2022singularity}, in which the authors conjectured that the 3D Euler equations could still develop a finite-time blow-up for initial data that are $C^\infty$ in $\rho$. In our future study, we plan to investigate the potential blow-up using a class of initial data of the form
\begin{align*}
    \omega^{\circ,4}_1&=-12000\left(1-r^2\right)^{18}\sin(2\pi z)^{2k+1},
\end{align*}
with a positive integer $k$, so that $\omega^{\circ,4}_1\sim r^\alpha z^{2k+1}=\rho^{2k+1+\alpha}\cos^\alpha\theta\sin^{2k+1}\theta$ is $C^{2k+1}$ in $\rho$.

\section{Comparison with Eligindi's singularity}
\label{sec: compare_elgindi}

In this section, we compare our blow-up scenario with the scenario in \cite{elgindi2021finite} studied by Elgindi.

Elgindi introduced a polar coordinate system on the $(r, z)$-plane to construct his blow-up solution. More specifically, he introduced
$$\rho=\sqrt{r^2+z^2}, \quad\theta=\arctan\left(\frac{z}{r}\right).$$
Then for a H\"{o}lder exponent $\alpha$, he introduced a change of variable $R=\rho^\alpha$ and defined the variables
$$\Omega(R,\theta)=\omega^\theta(r,z),\quad\Psi(R,\theta)=\frac{1}{\rho^2}\psi^
\theta(r,z).$$
In this setting, \eqref{eq: vort_stream_3d_noswirl} can be rewritten as
\begin{subequations}
\label{eq: elgindi}
\begin{align}
    &\Omega_{t} + \left(3\Psi+\alpha R\Psi_R\right)\Omega_\theta - \left(\Psi_\theta-\Psi\tan\theta\right)\Omega_R = \left(2\Psi\tan\theta+\alpha R\Psi_R\tan\theta+\Psi_\theta\right)\Omega, \label{eq: vort_elgindi}\\
    &-\alpha^2R^2\Psi_{RR}-\alpha(5+\alpha)R\Psi_{R}-\Psi_{\theta\theta}+\left(\Psi\tan\theta\right)_\theta-6\Psi=\Omega. \label{eq: stream_elgindi}
\end{align}
\end{subequations}

Elgindi's analysis of \eqref{eq: stream_elgindi} establishes the following leading order approximation for small $\alpha$
\begin{align}
\label{eq: elgindi_leading}
    \Psi(R,\theta)=\frac{1}{4\alpha}\sin(2\theta)L_{12}(\Omega)(R)+\text{lower order terms},
\end{align}
where
$$L_{12}(\Omega)(R)=\int_R^\infty\int_0^{\frac{\pi}{2}}\Omega(s,\theta)\frac{K(\theta)}{s}\mathrm{d}s\mathrm{d}\theta,$$
with $K(\theta)=3\sin\theta\cos^2\theta$. If we plug in the approximation \eqref{eq: elgindi_leading} to \eqref{eq: vort_elgindi}, neglecting lower order terms of $\alpha$, and (time) scaling out some constant factor, we arrive at Elgindi's fundamental model
\begin{align}
\label{eq: elgindi_fund}
    \Omega_{t} = \frac{1}{\alpha}L_{12}(\Omega)\Omega,
\end{align}
which admits self-similar finite-time blow-up. In his analysis, Elgindi chose the following self-similar solution of the fundamental model \eqref{eq: elgindi_fund}
\begin{align}
\label{eq: elgindi_solution}
    \Omega(R,\theta,t)=\frac{c}{1-t}F\left(\frac{R}{1-t}\right)\left(\sin\theta\cos^2\theta\right)^{\alpha/3},
\end{align}
where $c>0$ is some fixed constant, and $F(z)=2z/(1+z)^2$.

One difference between our blow-up scenario and Engindi's blow-up scenario is how the scaling factor $c_l$ depends on $\alpha$. We rewrite \eqref{eq: elgindi_solution} as
\begin{align*}
    \Omega&=\frac{c}{1-t}F\left(\frac{\rho^\alpha}{1-t}\right)\left(\frac{r^2z}{\rho^3}\right)^{\alpha/3}=\frac{c}{1-t}F\left(\left(\frac{\rho}{\left(1-t\right)^{1/\alpha}}\right)^\alpha\right)\left(\frac{r^{2/3}z^{1/3}}{\rho}\right)^{\alpha}.
\end{align*}
If we let $G(z)=F(z^\alpha)$, we see
$$\Omega=\frac{c}{1-t}G\left(\frac{\rho}{\left(1-t\right)^{1/\alpha}}\right)\left(\frac{r^{2/3}z^{1/3}}{\rho}\right)^{\alpha}.$$
Since $r^{2/3}z^{1/3}/\rho$ is homogeneous, we may conclude that the scaling factors for the self-similar blow-up solution \eqref{eq: elgindi_solution} are
$$c_l=1/\alpha,\quad c_\omega=2.$$
Note that this also satisfies the relation $c_\omega=1+\alpha c_l$ in \eqref{eq: scaling_relation}. This implies that $c_l$ decreases as $\alpha$ increases, and $c_l$ will tend to infinity as $\alpha\rightarrow0$. However, as shown in Table \ref{tab: cl_n3} and \ref{tab: cl_n10}, our $c_l$ increases as $\alpha$ increases. In the limit of $\alpha\rightarrow0$, our solution is well behaved and we observed an finite value of the scaling factor $c_l$ when $\alpha=0$, as reported in Table \ref{tab: cl_n3}, while Elgindi's case needs some renormalization and $c_l$ tends to infinity. The scaling factor $c_l$ in our case goes to infinity when $\alpha$ approaches to the critical value $\alpha^*$ which is close to $1/3$.

Furthermore, the regularity of our initial data as a function of $\rho$ is different from that of Elgindi's initial data. Around $(r,z)=(0,0)$, Elgindi's initial condition has the following leading order behavior
\begin{align*}
    \Omega\sim \rho^\alpha\left(\sin\theta\cos^2\theta\right)^{\alpha/3}=r^{2\alpha/3}z^{\alpha/3}.
\end{align*}
However, our initial condition gives
\begin{align*}
    \omega^\theta=r^\alpha\omega^\circ_1\sim r^\alpha z=\rho^{1+\alpha}\cos^{\alpha}\theta\sin\theta.
\end{align*}
These two leading order scaling properties differ from each other in that
\begin{itemize}
    \item Elgindi's initial condition of $\omega^\theta$ has a $C^\alpha$ H\"{o}lder continuity in $\rho$, whereas ours is $C^{1,\alpha}$ in $\rho$,
    \item Elgindi's initial condition of $\omega^\theta$ is H\"{o}lder continuous in both $z=0$ and $r=0$, whereas our initial condition is H\"{o}lder continuous in $r=0$ but smooth in $z$.
\end{itemize}

In Conjecture 8 of \cite{drivas2022singularity}, the authors conjectured that the initial data could be $C^\infty$ in $\rho$ for finite-time blow-up of the 3D axisymmetric Euler equations with no swirl. Our initial data slightly improves the regularity of the initial data in $\rho$.  In Section \ref{sec: sensitivity_third}, we also briefly explored the initial data with higher regularity in $\rho$. 

In Lemma 4.33 of \cite{drivas2022singularity}, the authors stated that the limiting equations at $\alpha=0$ of \eqref{eq: elgindi} can blow up in finite time for initial data of $\Omega$ that only has a $C^\alpha$-H\"{o}lder continuity near $r=0$ for $\alpha<1/3$. Our study shows that the blow-up of the axisymmetric Euler equations does not require to have H\"{o}lder continuity of the initial vorticity along the $z$-direction. The essential driving force for the finite-time blow-up comes from the H\"{o}lder continuity of the initial vorticity along the $r$-direction.

\section{A one-dimensional model of the potential self-similar blow-up}
\label{sec: 1d_model}

From Figure \ref{fig: steadystate_cross_section_alpha} and \ref{fig: steadystate_cross_section_alpha_n10} in Section \ref{sec: holder_dimension}, we observe that as $\alpha$ approaches the critical value $\alpha^*$, $-\tilde{\omega}_1$ will become very flat in $\xi$. This inspires us to conjecture that in the $\alpha\rightarrow\alpha^*$ limit, $-\tilde{\omega}_1$ will eventually become a function of $\zeta$ only in a relatively large domain. Based on this observation, we assume that
\begin{align}
    \omega_1(r,z)=\omega_1(0, z),
    \label{eq: 1d_assumption}
\end{align}
and derive a one-dimensional model for the $n$-D Euler equations \eqref{eq: vort_stream_nd_1_noswirl}.

At $r=0$, the velocity fields \eqref{eq: velo_rz_1_nd_noswirl} become $u^r=0$, $u^z=(n-1)\psi_1$. Therefore, the vorticity equation \eqref{eq: vort_theta_nd_1_noswirl} becomes
$$\omega_{1,t}(0,z)+(n-1)\psi_1(0,z)\omega_{1,z}(0,z)=-(n-2-\alpha)\psi_{1,z}(0,z)\omega_1(0,z).$$
As for the Poisson equation \eqref{eq: stream_theta_nd_1_noswirl}, we use the Green's function $G_{n,\alpha}(r,r^\prime, z, z^\prime)$ for the operator $L_{n,\alpha}=r^{1-\alpha}\left(-\partial_{rr}-\frac{n}{r}\partial_r-\partial_{zz}\right)$. We have
\begin{align*}
    \psi_1(r,z)&=\int_{(r^\prime,z^\prime)\in\mathcal{D}}G_{n,\alpha}(r,r^\prime, z, z^\prime)\omega_1(r^\prime,z^\prime)\mathrm{d}r^\prime\mathrm{d}z^\prime,\\
    &=\int_{(r^\prime,z^\prime)\in\mathcal{D}}G_{n,\alpha}(r,r^\prime, z, z^\prime)\omega_1(0,z^\prime)\mathrm{d}r^\prime\mathrm{d}z^\prime,
\end{align*}
and therefore
\begin{align*}
    \psi_1(0,z)&=\int_{0}^{1/2}H_{n,\alpha}(z, z^\prime)\omega_1(0,z^\prime)\mathrm{d}z^\prime,
\end{align*}
where
\begin{align*}
    H_{n,\alpha}(z, z^\prime)=\int_0^1G_{n,\alpha}(0,r^\prime,z,z^\prime)\mathrm{d}r^\prime.
\end{align*}

Putting these equations together, and omitting the $r$-coordinate when there is no ambiguity, we have the following closed system in 1D: for $z\in\left[0,1/2\right]$, $\omega_1$ and $\psi_1$ are functions of $z$ whose evolution in time is governed by the equations
\begin{subequations}
\label{eq: 1d_model}
\begin{align}
    \omega_{1,t} + (n-1)\psi_1\omega_{1,z} &= -(n-2-\alpha)\psi_{1,z}\omega_1, \label{eq: 1d_model_vort}\\
    \psi_1&=T_{n,\alpha}\omega_1, \label{eq: 1d_model_steam}
\end{align}
\end{subequations}
where $T_{n,\alpha}$ is an integral transform with kernel function $H_{n,\alpha}$:
\begin{align*}
    T_{n,\alpha}\omega_1&=\int_{0}^{1/2}H_{n,\alpha}(z, z^\prime)\omega_1(z^\prime)\mathrm{d}z^\prime.
\end{align*}

\subsection{The Kernel Function $H_{n,\alpha}$}

We look for a more explicit expression for the kernel $H_{n,\alpha}$.

Following the idea in \cite{hou2018potential}, we view $-\partial_{rr}-\frac{n}{r}\partial_r-\partial_{zz}$ as the Laplacian operator in the $(n+2)$-dimensional space for axisymmetric functions. The fundamental solution for the $(n+2)$-dimensional Laplace equation is
$$\Phi_0(x)=\frac{\Gamma(n/2)}{4\pi^{n/2+1}}\frac{1}{|x|^n},$$ for $x\in\mathbb{R}^{n+2}$. Now, since we have zero Dirichlet boundary conditions at $r=1$, $z=0$, $z=1/2$, we can obtain the Green's function for the above equation by properly symmetrizing the fundamental solution of the Laplace equation, which gives us
\begin{align*}
    G_{n,\alpha}(r,r^\prime,z,z^\prime)=&\sum_{m\in\mathbb{Z}}\left(G^\circ_{n,\alpha}\left(r,r^\prime,z+m,z^\prime\right)-G^\circ_{n,\alpha}\left(r,r^\prime,-z+m,z^\prime\right)\right.\\
    &-\left.G^\circ_{n,\alpha}\left(1,rr^\prime,z+m,z^\prime\right)+G^\circ_{n,\alpha}\left(1,rr^\prime,-z+m,z^\prime\right)\right),
\end{align*}
where
$$G^\circ_{n,\alpha}(r,r^\prime,z,z^\prime)=C(n)\frac{{r^\prime}^{n+\alpha-1}}{A^{n/2}}{}_2F_1\left(n/2,n/2,n,B\right),$$
with $\Gamma$ being the Gamma function, ${}_2F_1$ being the Gauss hypergeometric function, and
$$A=\left(r+r^\prime\right)^2+\left(z-z^\prime\right)^2,\quad B=4rr^\prime/A,\quad C(n)=\frac{2^{n-2}}{\pi}\frac{\Gamma(n/2)^2}{\Gamma(n)}.$$
In fact, it is easy to check that $G_{n,\alpha}$ satisfies the boundary conditions:
$$G_{n,\alpha}(1,r^\prime,z,z^\prime)=G_{n,\alpha}(r,r^\prime,0,z^\prime)=G_{n,\alpha}(r,r^\prime,1/2,z^\prime)=0.$$

We notice that the Gaussian hypergeometric function ${}_2F_1$ has the property that ${}_2F_1(n/2,n/2,n,0)=1$. Therefore, we know that
\begin{align*}
    G^\circ_{n,\alpha}(0,r^\prime,z,z^\prime)=\frac{C(n){r^\prime}^{n+\alpha-1}}{\left({r^\prime}^2+\left(z-z^\prime\right)^2\right)^{n/2}},\quad G^\circ_{n,\alpha}(1,0,z,z^\prime)=0.
\end{align*}
Therefore, we arrive at an expression for $H_{n,\alpha}$:
\begin{align*}
    H_{n,\alpha}(z,z^\prime)&=\int_0^1G_{n,\alpha}(0,r^\prime,z,z^\prime)\mathrm{d}r^\prime\\
    =C(n)&\int_0^{1}\sum_{m\in\mathbb{Z}}\left(\frac{r^{n+\alpha-1}}{\left(r^2+\left(z+m-z^\prime\right)^2\right)^{n/2}}-\frac{r^{n+\alpha-1}}{\left(r^2+\left(z-m+z^\prime\right)^2\right)^{n/2}}\right)\mathrm{d}r.
\end{align*}
From $H_{n,\alpha}$, we can also find the integral transform for $\psi_{1,z}$:
$$\psi_{1,z}(z)=\int_0^{1/2}\partial_zH_{n,\alpha}(z,z^\prime)\omega_1(z^\prime)\mathrm{d}z^\prime,$$
with
\begin{align*}
    &\partial_zH_{n,\alpha}(z,z^\prime)=\\
    nC(n)&\int_0^{1}\sum_{m\in\mathbb{Z}}\left(\frac{r^{n+\alpha-1}\left(z+m-z^\prime\right)}{\left(r^2+\left(z+m-z^\prime\right)^2\right)^{n/2+1}}-\frac{r^{n+\alpha-1}\left(z-m+z^\prime\right)}{\left(r^2+\left(z-m+z^\prime\right)^2\right)^{n/2+1}}\right)\mathrm{d}r.
\end{align*}

\begin{figure}[hbt!]
\centering
\includegraphics[width=.4\textwidth]{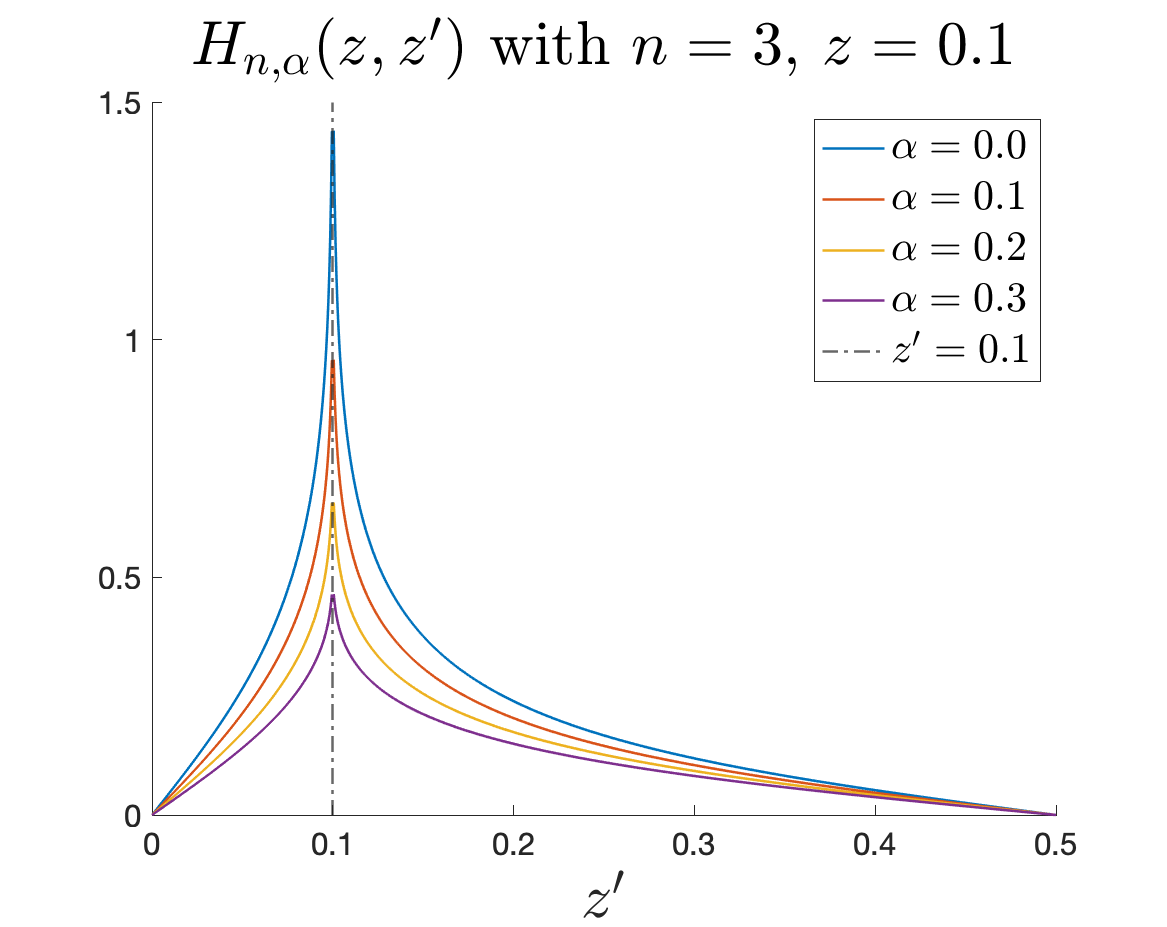}\includegraphics[width=.4\textwidth]{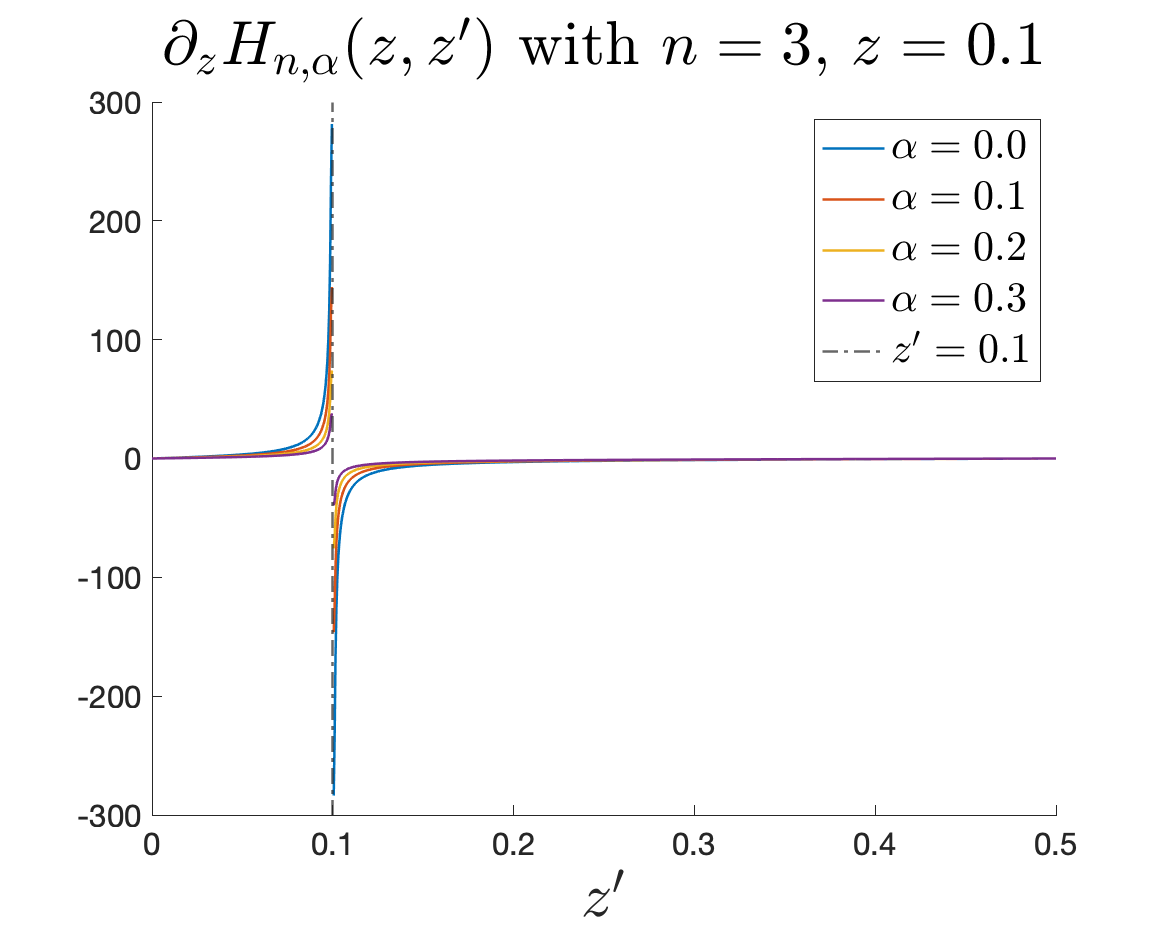}
\caption{Profiles of $H_{n,\alpha}$ and $\partial_zH_{n,\alpha}$ as functions of $z^\prime$ for different $\alpha$ with $n=3$ and $z=0.1$.}\label{fig: H_a}
\index{figures}
\end{figure}

\begin{figure}[hbt!]
\centering
\includegraphics[width=.4\textwidth]{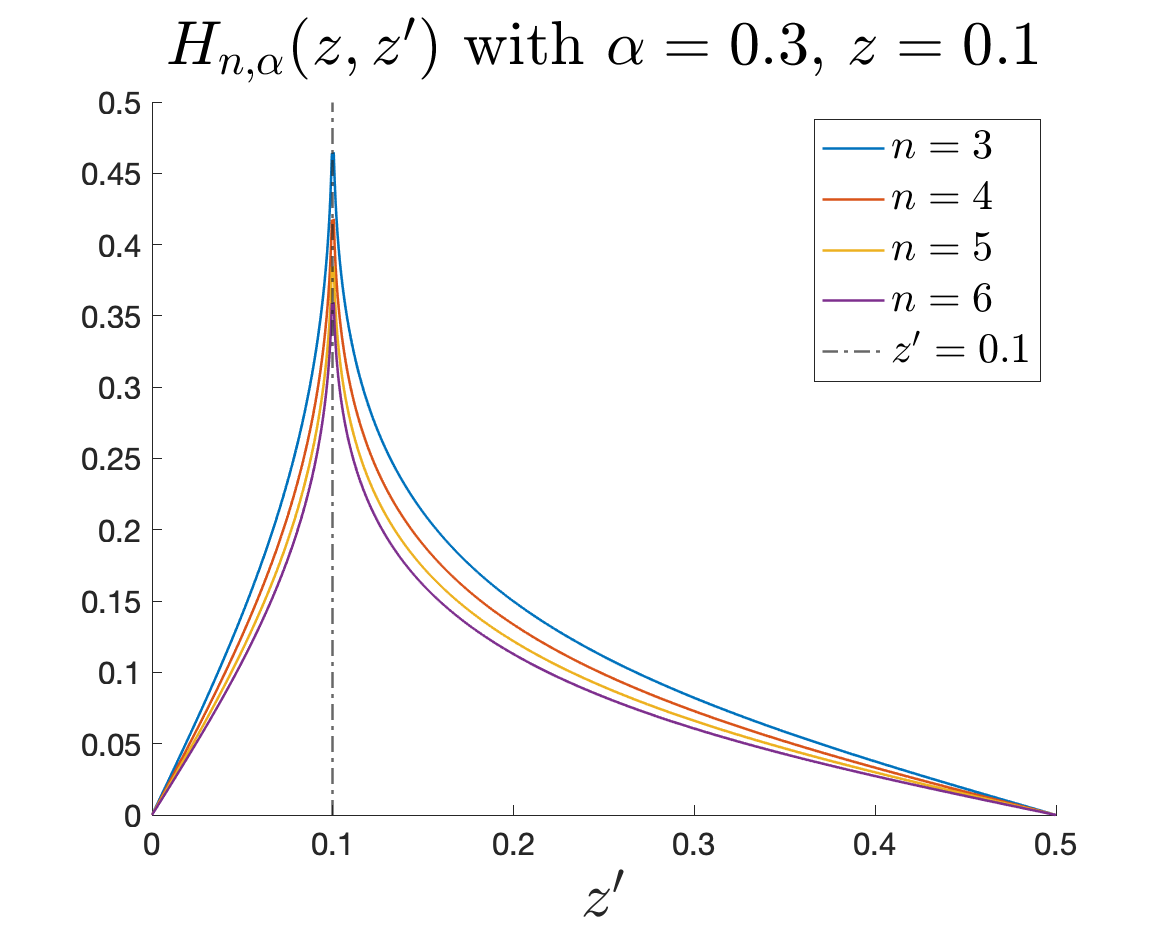}\includegraphics[width=.4\textwidth]{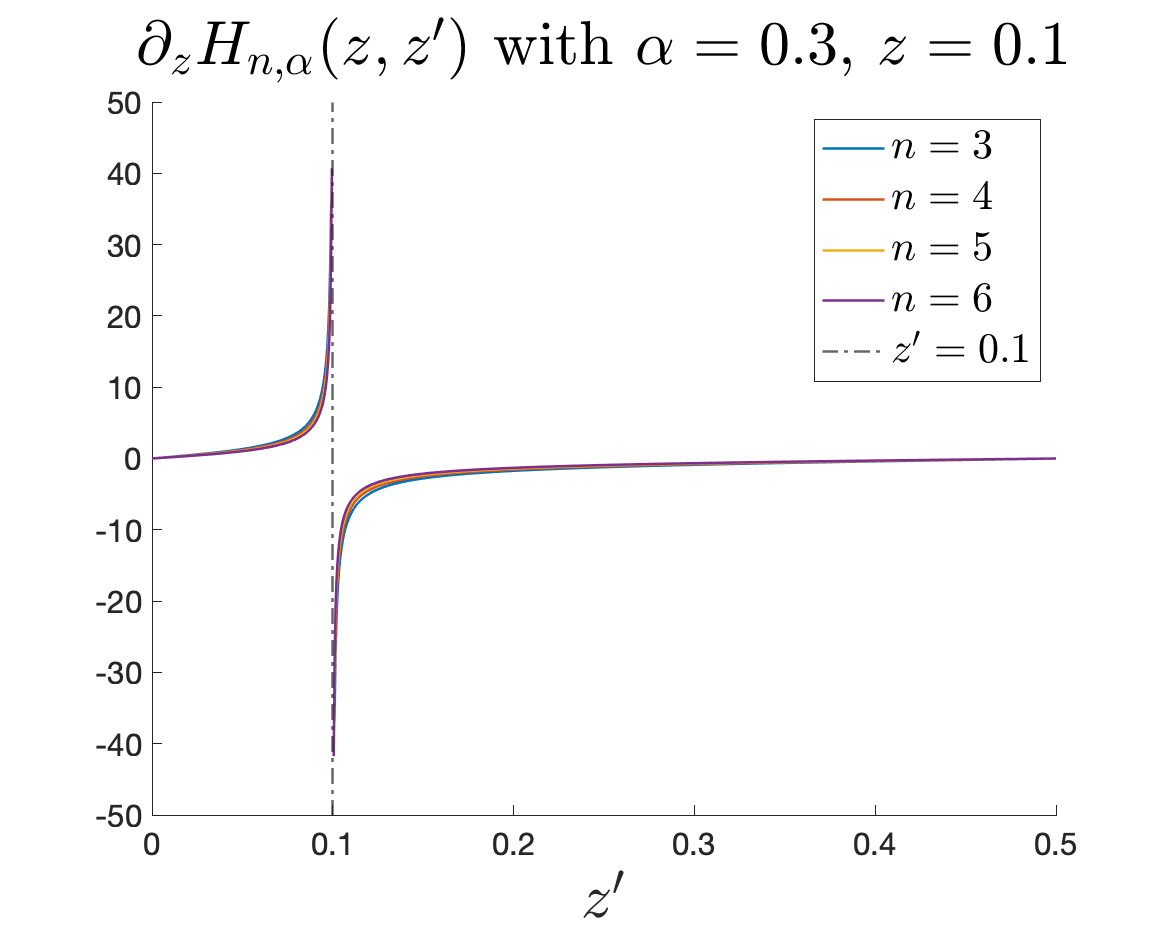}
\caption{Profiles of $H_{n,\alpha}$ and $\partial_zH_{n,\alpha}$ as functions of $z^\prime$ for different $n$ with $\alpha=0.3$ and $z=0.1$.}\label{fig: H_n}
\index{figures}
\end{figure}

In Figure \ref{fig: H_a} and \ref{fig: H_n}, we show the profiles of the kernel functions $H_{n,\alpha}$ and $\partial_zH_{n,\alpha}$ for various combinations of $(n,\alpha)$. Figure \ref{fig: H_a} shows that the smaller $\alpha$ will make $H_{n,\alpha}$ and $\partial_zH_{n,\alpha}$ larger in scale. And this corresponds to the fact that smaller $\alpha$ is easier to develop the blow-up. In Figure \ref{fig: H_n}, we see that as the dimension $n$ increases, the profile of $H_{n,\alpha}$ seems to become shorter and thinner, which makes the velocity component $u^z=(n-1)\psi_1$ smaller. This phenomenon is consistent with Table \ref{tab: stats_dim} where larger $n$ tends to have a slower collapsing rate.

\subsection{Numerical Simulation}

Although the assumption \eqref{eq: 1d_assumption} of our 1D model is based on our observation in the large $\alpha$ scenario, we simulate the 1D system \eqref{eq: 1d_model} numerically with different choices of $\alpha$ to show its approximation to the $n$-D axisymmetric Euler equations \eqref{eq: vort_stream_nd_1_noswirl}.

Similar to our study in Section \ref{sec: holder_exponent}, we use the adaptive mesh method to solve \eqref{eq: 1d_model}, and then compute the potential self-similar blow-up profile by solving the steady state of the following dynamic rescaling formulation:
\begin{subequations}
\label{eq: 1d_model_dr}
\begin{align}
    \tilde{\omega}_{1,\tau} + \left(\tilde{c}_l\zeta+(n-1)\tilde{\psi}_1\right)\tilde{\omega}_{1,\zeta} &=\left(\tilde{c}_\omega-(n-2-\alpha)\tilde{\psi}_{1,\zeta}\right)\tilde{\omega}_1, \label{eq: 1d_model_vort_dr}\\
    \tilde{\psi}_1&=T^\prime_{n,\alpha}\tilde{\omega}_1. \label{eq: 1d_model_steam_dr}
\end{align}
\end{subequations}
Here we use $T^\prime_{n,\alpha}$ to distinguish from the integral transform $T_{n,\alpha}$ in \eqref{eq: 1d_model}, because $T^\prime_{n,\alpha}$ is on the large domain $\mathcal{D}^\prime$ with zero Neumann boundary condition in the far field, the same as the setup in Section \ref{sec: dynamic_rescaling_numeric}. In our computation, we use the late stage solution of the equations \eqref{eq: vort_stream_nd_1_noswirl} as the initial data to \eqref{eq: 1d_model}. To calculate the integral transform $T_{n,\alpha}$ numerically, we first temporarily lift up $\omega_1$ from the 1D axis of $z$ to the 2D plane of $(r, z)$ by relation \eqref{eq: 1d_assumption}, then solve the Poisson equation \eqref{eq: stream_theta_nd_1_noswirl} on the $(r, z)$-plane for $\psi_1$, and finally restrict $\psi_1$ to the symmetry axis to obtain its value on the $z$-axis. The same method goes for the integral transform $T^\prime_{n,\alpha}$, but on different computational domain ($\mathcal{D}^\prime$) and with a different boundary condition in the far field.

\begin{figure}[hbt!]
\centering
\includegraphics[width=.4\textwidth]{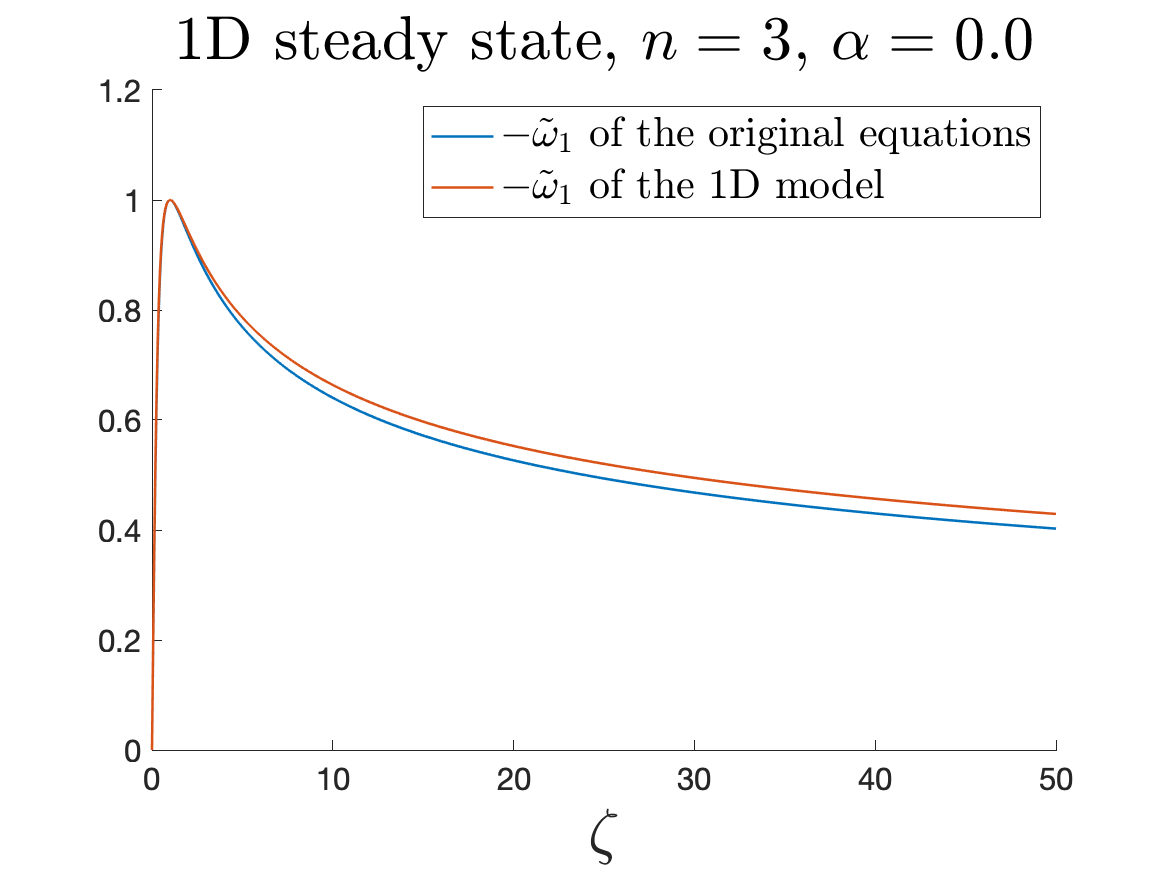}\includegraphics[width=.4\textwidth]{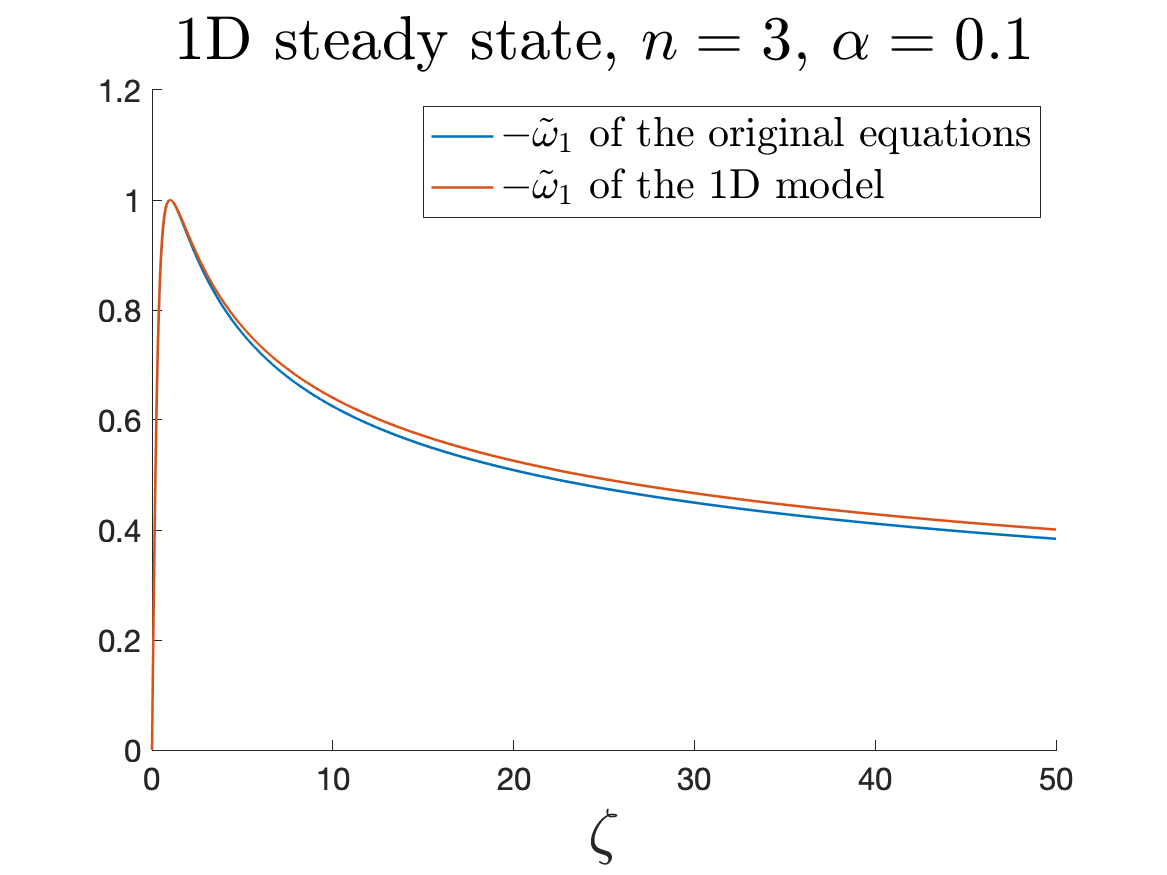}\\
\includegraphics[width=.4\textwidth]{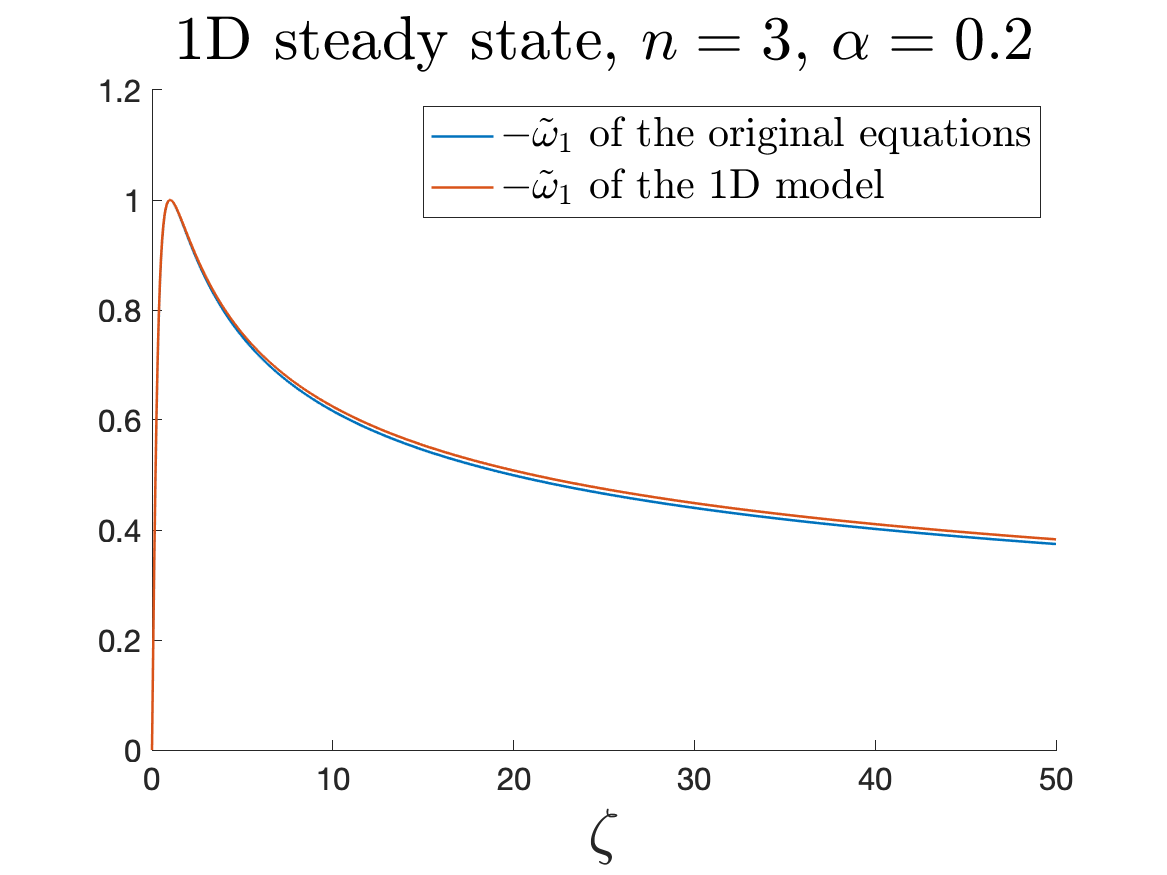}\includegraphics[width=.4\textwidth]{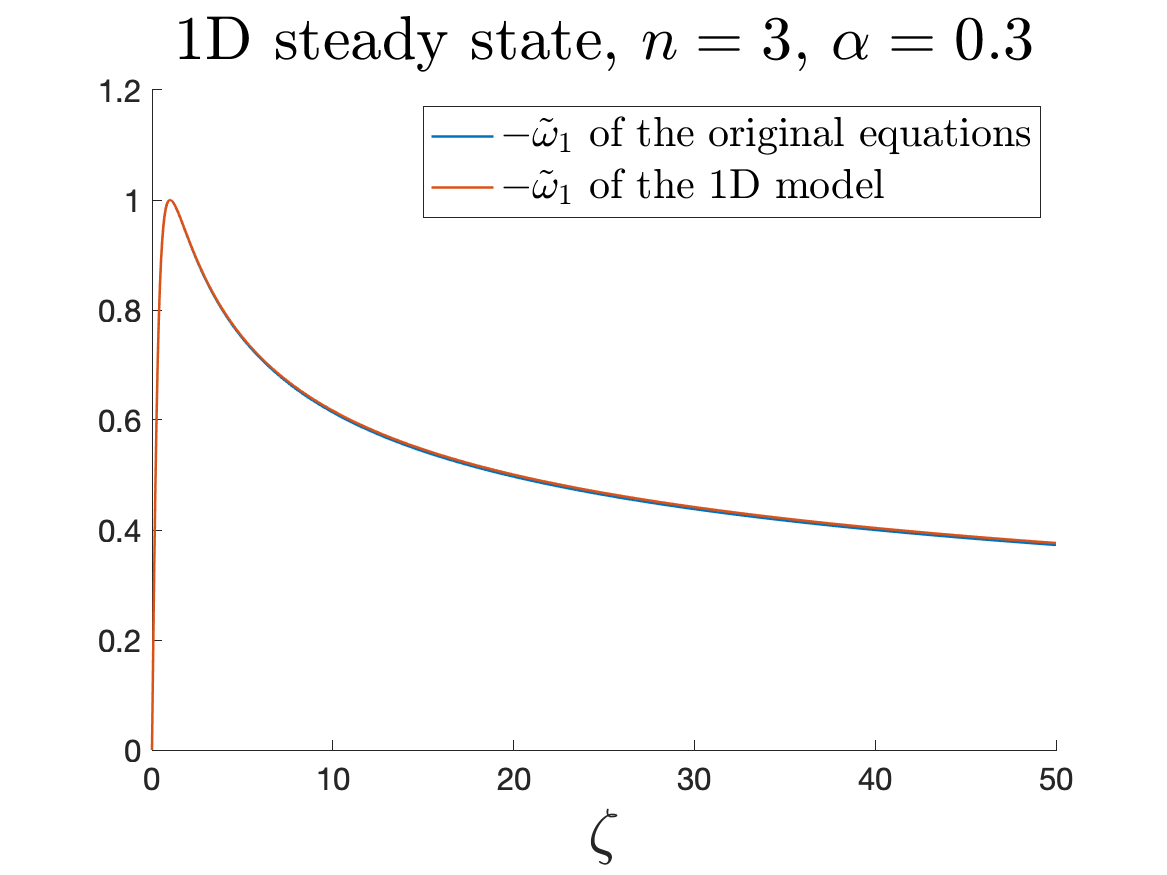}
\caption{Comparison of the steady state of $-\tilde{\omega}_1(0, \zeta)$ with different $\alpha$.}\label{fig: 1d_steadystate_compare}
\index{figures}
\end{figure}

Since the dimension $n$ would not be critical in the comparison, we only consider the 3D case ($n=3$) here, with H\"{o}lder exponent $\alpha=0.0, 0.1, 0.2, 0.3$. In Figure \ref{fig: 1d_steadystate_compare}, we compare the steady state of of $-\tilde{\omega}_1(0, \zeta)$ from the 3D axisymmetric Euler equations and the 1D model. When $\alpha=0.3$, the two steady states match with each other quite well, with a relative sup-norm error of $1.5\times10^{-3}$ on the $\zeta$-axis. It is interesting to see that at small $\alpha$, the two steady states are also close to each other. The relative sup-norm error at $\alpha=0.0$ is smaller than $2.7\times10^{-2}$ on the $\zeta$-axis.

In Table \ref{tab: cl_1d_compare}, we report the comparison of the scaling factor $c_l^{\text{1D}}$ of our 1D model and the scaling factor $c_l$ for the 3D axisymmetric Euler equations. As $\alpha$ increases, the 1D model \eqref{eq: 1d_model} can better approximate the 3D axisymmetric Euler equation in the relative error sense. This shows that our 1D model \eqref{eq: 1d_model} can serve as a good model to understand the blow-up mechanism for the potential self-similar blow-up of the $n$-dimensional Euler equations. Although the 1D model is derived based on the flatness of $\omega_1$ as a function of $\xi$ as $\alpha$ is close to the critical value $\alpha^*$, the agreement of $c_l$ between the 1D model and the 3D Euler equations is still quite good for smaller values of $\alpha$, which is quite surprising.

\begin{table}[hbt!]
\centering
\caption{Scaling factor comparison between the 1D model and the 3D axisymmetric Euler equations.}\label{tab: cl_1d_compare}
\begin{tabular}{|c|c|c|c|c|c|c|}
\hline
\hspace{0.5cm}$\alpha$\hspace{0.5cm} & \hspace{0.5cm}$0.0$\hspace{0.5cm} & \hspace{0.5cm}$0.1$\hspace{0.5cm} & \hspace{0.5cm}$0.2$\hspace{0.5cm} & \hspace{0.5cm}$0.3$\hspace{0.5cm} \\
\hline
$c_l$ & $3.248$ & $4.549$ & $8.270$ & $112.8$ \\
\hline
$c_l^{\text{1D}}$ & $3.374$ & $4.682$ & $8.464$ & $114.8$ \\
\hline
\end{tabular}
\index{tables}
\end{table}

\section{Concluding remarks}
\label{sec: conclusions}

In this paper, we have numerically studied the singularity formation in the axisymmetric Euler equations with no swirl when the initial condition for the angular vorticity is $C^\alpha$ H\"{o}lder continuous. With carefully-chosen initial data and specially-designed adaptive mesh, we have solved the solution very close to the potential blow-up time, and obtained strong convincing numerical evidence for the singularity formation by numerically examining the Beale-Kato-Majda blow-up criterion. Scaling analysis and dynamic rescaling method have further suggested the potential self-similar blow-up. We observed the potential self-similar blow-up in finite time when the H\"{o}lder exponent $\alpha$ is greater or equal to $0$, and is smaller than a critical value $\alpha^*$, and this upper bound $\alpha^*$ is larger than $0.3$, and have the potential to be $1/3$. This result supports Conjecture 8 of \cite{drivas2022singularity}. Since when $\alpha>1/3$, the axisymmetric Euler equations with no swirl admit global regularity \cite{ukhovskii1968axially, serfati1994regularite, shirota1994note, saint1994remarks, danchin2007axisymmetric, abidi2010global}, this would potentially close gap between blow-up and non-blow-up, leaving only the critical case of $\alpha=1/3$. We also extend this result to the high-dimensional case, and find that in general the critical value $\alpha^*$ for the $n$-D axisymmetric Euler equation is close to $1-\frac{2}{n}$.

The potential blow-up observed in this paper is insensitive to the perturbation of initial data. And our initial study suggested that the regularity of the initial data around the origin would determine its scaling properties and the shape of the self-similar blow-up profile. Compared with Elgindi's blow-up result reported in \cite{elgindi2021finite}, our potential blow-up scenario has very different scaling properties. The  regularity properties of the initial condition of the two initial data are also quite different. 

Inspired by our numerical observations, we proposed a simple one-dimensional model to capture the leading order behavior of the $n$-dimensional Euler equations. Our numerical experiments showed that the one-dimensional model can develop approximately the same potential finite-time blow-up as the original $n$-dimensional Euler equations. This one-dimensional model could play a role similar to the leading order system derived by Elgindi in \cite{elgindi2021finite} in the analysis of the finite-time singularity of the 3D Euler equations.

\section*{Acknowledgments}
The research was in part supported by DMS-2205590. We would like to acknowledge the generous support from Mr. K. C. Choi through the Choi Family Gift Fund and the Choi Family Postdoc Gift Fund. We would also like to thank Mr. Xiang Qin and Mr. Xiuyuan Wang from Peking University for helpful discussion.

\bibliographystyle{siamplain}
\bibliography{references}
\end{document}